\documentclass[preprint,plainstyle]{imsart} 

\RequirePackage[colorlinks,citecolor=blue,urlcolor=blue]{hyperref}
\usepackage{url}

\usepackage[bottom]{footmisc}

\usepackage{amsmath}
\allowdisplaybreaks
\usepackage{graphicx}
\usepackage{enumerate}

\RequirePackage{natbib}
\usepackage{url} % not crucial - just used below for the URL  

\startlocaldefs

\usepackage{amsmath,amssymb,amsthm,graphicx}
\usepackage{comment}
\usepackage{color} 
\usepackage{mathtools}
\usepackage{enumitem}
\usepackage{subcaption} 
\usepackage[english]{babel}

\usepackage{macros}

\definecolor{green}{rgb}{.1,.4,.1}

\newcommand{\sw}{\mathrm{SW}}
\newcommand{\sj}{\mathrm{SJ}}  
 
\makeatletter
\let\save@mathaccent\mathaccent
\newcommand*\if@single[3]{%
  \setbox0\hbox{${\mathaccent"0362{#1}}^H$}%
  \setbox2\hbox{${\mathaccent"0362{\kern0pt#1}}^H$}%
  \ifdim\ht0=\ht2 #3\else #2\fi
  }
\newcommand*\rel@kern[1]{\kern#1\dimexpr\macc@kerna}
\newcommand*\widebar[1]{\@ifnextchar^{{\wide@bar{#1}{0}}}{\wide@bar{#1}{1}}}
\newcommand*\wide@bar[2]{\if@single{#1}{\wide@bar@{#1}{#2}{1}}{\wide@bar@{#1}{#2}{2}}}
\newcommand*\wide@bar@[3]{%
  \begingroup
  \def\mathaccent##1##2{%
    \let\mathaccent\save@mathaccent
    \if#32 \let\macc@nucleus\first@char \fi
    \setbox\z@\hbox{$\macc@style{\macc@nucleus}_{}$}%
    \setbox\tw@\hbox{$\macc@style{\macc@nucleus}{}_{}$}%
    \dimen@\wd\tw@
    \advance\dimen@-\wd\z@
    \divide\dimen@ 3
    \@tempdima\wd\tw@
    \advance\@tempdima-\scriptspace
    \divide\@tempdima 10
    \advance\dimen@-\@tempdima
    \ifdim\dimen@>\z@ \dimen@0pt\fi
    \rel@kern{0.6}\kern-\dimen@
    \if#31
      \overline{\rel@kern{-0.6}\kern\dimen@\macc@nucleus\rel@kern{0.4}\kern\dimen@}%
      \advance\dimen@0.4\dimexpr\macc@kerna
      \let\final@kern#2%
      \ifdim\dimen@<\z@ \let\final@kern1\fi
      \if\final@kern1 \kern-\dimen@\fi
    \else
      \overline{\rel@kern{-0.6}\kern\dimen@#1}%
    \fi
  }%
  \macc@depth\@ne
  \let\math@bgroup\@empty \let\math@egroup\macc@set@skewchar
  \mathsurround\z@ \frozen@everymath{\mathgroup\macc@group\relax}%
  \macc@set@skewchar\relax
  \let\mathaccentV\macc@nested@a
  \if#31
    \macc@nested@a\relax111{#1}%
  \else
    \def\gobble@till@marker##1\endmarker{}%
    \futurelet\first@char\gobble@till@marker#1\endmarker
    \ifcat\noexpand\first@char A\else
      \def\first@char{}%
    \fi
    \macc@nested@a\relax111{\first@char}%
  \fi
  \endgroup
}
\makeatother

\theoremstyle{plain}

\newtheorem{theorem}{Theorem}
\newtheorem{lemma}{Lemma}
\newtheorem{example}{Example}
\newtheorem{corollary}{Corollary}
\newtheorem{proposition}{Proposition}

\endlocaldefs

\begin{document}

\begin{frontmatter}

% "Title of the Paper"
\title{Minimax Confidence Intervals for the Sliced Wasserstein Distance}%thanksref{t1}}

\runtitle{Minimax Confidence Intervals for the Sliced Wasserstein Distance}

% indicate corresponding author with \corref{}
% \author{\fnms{John} \snm{Smith}\thanksref{t2}\corref{}\ead[label=e1]{smith@foo.com}\ead[label=e2,url]{www.foo.com}}
% \thankstext{t2}{Thanks to somebody} 
% \address{line 1\\ line 2\\ \printead{e1}\\ \printead{e2}}

\author{\fnms{Tudor} \snm{Manole}\corref{}\ead[label=e1]{tmanole@andrew.cmu.edu}},   %\thanksref{t1,t2,t3}
\author{\fnms{Sivaraman} \snm{Balakrishnan}\ead[label=e2]{siva@stat.cmu.edu}},    %\thanksref{t1,t3,t4}
\author{\fnms{Larry} \snm{Wasserman}%\thanksref{t1}
\ead[label=e3]{larry@stat.cmu.edu}}
\address{Department of Statistics and Data Science \\Carnegie Mellon University \\
\printead{e1,e2,e3}}

\runauthor{Manole et al.}

\begin{abstract}
Motivated by the growing popularity of variants of  
the Wasserstein distance in statistics and machine learning,  
we study statistical inference for the Sliced Wasserstein distance---an easily
computable variant of the Wasserstein distance. Specifically, we construct
confidence intervals for the Sliced Wasserstein distance which
have finite-sample validity under no assumptions or under mild moment assumptions.
These intervals 
are adaptive in length to the regularity of the
underlying distributions. 
We also bound the minimax risk of estimating the Sliced Wasserstein distance, 
and as a consequence establish 
that the lengths of our proposed confidence intervals are minimax optimal over appropriate 
distribution classes. 
To motivate the choice of these classes, we also study minimax 
rates of estimating a distribution under the Sliced Wasserstein distance. 
These theoretical findings are complemented with a simulation study demonstrating the deficiencies
of the classical bootstrap, and the advantages of our proposed methods. We also show 
strong correspondences between our theoretical
predictions and the adaptivity of our confidence interval lengths in simulations.
We conclude 
by demonstrating the use of our confidence intervals in the setting of 
simulator-based likelihood-free inference. 
In this 
%likelihood-free 
setting, 
%in contrast to 
contrasting popular approximate Bayesian computation methods, 
we develop uncertainty quantification methods with rigorous frequentist coverage guarantees.
\end{abstract}

%Primary: ; Secondary: 
 
\end{frontmatter}

% Main text entry area
\begingroup
%\hypersetup{linkcolor=black}
%\tableofcontents
\endgroup

\section{Introduction} 
The Wasserstein distance is a metric between probability distributions
which has received a surge of interest in statistics and machine learning
\citep{panaretos2019, kolouri2017}. This distance arises from the 
optimal transport problem \citep{villani2003}, and measures the work required 
to couple one distribution with another. Specifically, given probability distributions $P$ and $Q$
admitting at least $r \geq 1$ moments, with support in $\bbR^d$ for $d \geq 1$,
the $r$-th order Wasserstein distance between $P$ and $Q$ is defined by
\begin{equation}
\label{eq::wasserstein}
W_r(P,Q) = \left(\inf_{\gamma\in \Pi(P,Q)}\int\norm{x-y}^r d\gamma(x,y)\right)^{1/r},
\end{equation}
where
$\Pi(P,Q)$ denotes the set of joint probability distributions
 with marginals $P$ and $Q$,   known as couplings. Any minimizer $\gamma$ is called an optimal
 coupling between $P$ and $Q$. The norm $\norm\cdot$ is taken to be Euclidean in this
 paper, but may more generally be replaced by any metric on~$\bbR^d$.
 
The Wasserstein distance has broadly
served two uses in the statistics literature (see the review article of~\cite{panaretos2019a} and references therein).
On the one hand, it has been used
as a theoretical tool for asymptotic theory
\citep{shorack2009, shao2012}, since
convergence in the $r$-Wasserstein distance is equivalent 
to weak convergence of probability measures and their $r$-th moments \citep{villani2003}.
Wasserstein distances also play a prominent role in the analysis of 
mixture models~\citep{nguyen2013a,ho2019}. 
On the other hand, increasingly many statistical applications 
employ the Wasserstein distance as a methodological
tool in its own right.
Unlike many common metrics between probability distributions, the Wasserstein
distance does not presume distributions which are absolutely continuous
with respect to a common dominating measure, and is sensitive to the underlying geometry
of their support, due to the $\ell_2$-norm embedded in its definition. These considerations
make it a natural and powerful data analytic tool---see for 
instance~\cite{delbarrio1999,  delbarrio2005,
courty2016a,ramdas2017,arjovsky2017,ho2017c,bernton2019a,bernton2019,
verdinelli2019} and references therein.

Despite the popularity of the Wasserstein distance,  
its high computational complexity often limits its applicability to large-scale problems.
Developing
efficient numerical approximations  of the distance
remains an active research area---see 
 \cite{peyre2019} for a recent review. 
A key exception to the high computational cost is the univariate case, 
in which the Wasserstein distance admits a closed form
as the $L^r$ norm between the quantile functions of~$P$ and~$Q$,
which can be easily computed. This fact 
has led to the study of an alternate metric, known as the Sliced Wasserstein distance 
\citep{rabin2011, bonneel2015}
obtained by averaging      
the Wasserstein distance between random one-dimensional projections of the distributions
$P$ and $Q$. 
\begin{comment}

Specifically, let $\bbS^{d-1} = \{x \in \bbR^d: \norm x = 1\}$, 
and let $\mu$ denote the uniform probability measure on $\bbS^{d-1}$. 
The $r$-th order Sliced Wasserstein distance
between $P$ and $Q$ is given by
$$\sw_r(P,Q) = \left(\int_{\bbS^{d-1}} W_r^r(P_\theta, Q_\theta) d\mu(\theta)\right)^{\frac 1 r}=\left(\int_{\bbS^{d-1}} \int_0^1 \big|F_\theta^{-1}(u) - G_\theta^{-1}(u)\big|^r dud\mu(\theta)\right)^{\frac 1 r},$$
where for any $\theta \in \bbS^{d-1}$, $P_\theta$ and $Q_\theta$ denote the respective
probability distributions of $X^\top\theta$ and $Y^\top\theta$ where $X \sim P$ and $Y \sim Q$.
Furthermore, $F_\theta:x \in \bbR \mapsto \bbP(X^\top\theta \leq x)$ is the cumulative distribution function (CDF) of $P_\theta$, 
and $F_\theta^{-1} : u \in [0,1] \mapsto \inf \big\{ x \in \bbR: F_\theta(x) \geq u\big\}$
denotes its quantile function (and similarly for $G_\theta$ and $G_\theta^{-1}$).

\end{comment}
The Sliced Wasserstein distance is a generally weaker metric than the Wasserstein distance \citep{bonnotte2013}, 
%as shown by the inequality $\sw_r(P,Q) \leq \frac 1 d W_r(P,Q)$ for all $r \geq 2$ .
but nevertheless preserves many qualitatively similar properties which % to the Wasserstein distance, 
make it an attractive and easily computable 
alternative in many applications.  

Motivated by the fact that the Wasserstein distance and its sliced analogue
are sensitive to outliers and heavy tails, we introduce a trimmed
version of the Sliced Wasserstein distance, denoted by $\sw_{r,\delta}(P,Q)$ for some trimming constant $\delta \in [0,1/2)$ and defined formally in
equation~\eqref{eq::trimmed_SW}.  
This robustification of the Sliced Wasserstein distance compares 
distributions up to a $2\delta$ fraction of their probability mass, thereby generalizing 
the one-dimensional trimmed Wasserstein distance introduced by \cite{munk1998}. % (see also \cite{alvarez-esteban2008,alvarez-esteban2012}). 
One of the aims of our 
paper is to derive
confidence intervals for the trimmed Sliced Wasserstein distance which make either no
assumptions or mild moment assumptions  
on the unknown distributions $P$ and $Q$. Specifically, 
given a level $\alpha\in(0,1)$ and
i.i.d. samples $X_1, \dots, X_n \sim P$ and $Y_1, \dots, Y_m \sim Q$, we derive confidence sets
$C_{nm} \subseteq \bbR$ such that
\begin{equation}
\label{eq::two-sample}
\inf_{P,Q} \bbP\big( \sw_{r,\delta}(P,Q) \in C_{nm}\big) \geq 1-\alpha,
\end{equation}
where the infimum is over a suitable family of distributions $P,Q$.

One of the main reasons that the Wasserstein distance has found many applications is the
fact that it is a useful notion of distance under weak %moment 
assumptions. Unlike the Total Variation, Hellinger, Kullback-Leibler and other 
divergences, the Wasserstein distance between a pair of distributions can be estimated from samples (optimally) 
under mild assumptions without requiring any smoothing. However, existing results on \emph{inference} for the Wasserstein distance 
\citep{munk1998,freitag2003,freitag2007,freitag2005} typically require strong 
smoothness assumptions and suggest different inferential procedures
when $P=Q$ as compared to when $P \neq Q$. In contrast, we construct various assumption-light confidence intervals $C_{nm}$ which 
have finite-sample validity under weak moment assumptions.  

The confidence intervals we construct
are adaptive to the regularity of the distributions $P$ and $Q$, as measured by %the magnitude of 
a functional $\sj_{r,\delta}$ introduced formally in Section \ref{sec::minimax_distribution} (equation \eqref{eq::sj}).
The magnitude of $\sj_{r,\delta}(P)$ is largely controlled by the tails of $P$ and by whether its one-dimensional projections
have connected support. %In rough terms, the $\sj$ functional has a small value when the distribution has light tails,
%and has connected support. 
The one-dimensional counterpart of this functional was identified in the work of~\cite{bobkov2019} who showed that when this
functional is finite, the empirical measure of $X_1, \dots, X_n$ 
converges to $P$ under the Wasserstein distance
at the fast rate of $O(1/\sqrt{n})$, assuming $d=1$. On the other hand, when this functional is infinite~\cite{bobkov2019} showed that 
this convergence happens at a slower rate of $O((1/n)^{1/2r})$. 
Our work shows that the role of the $\sj_{r,\delta}$ functional in inference is
more nuanced. When  
$\sj_{r,\delta}(P)$ and $\sj_{r,\delta}(Q)$ 
are finite, our confidence intervals have length scaling at the fast rate of $O(1/\sqrt{n \wedge m})$, mirroring the
rates of convergence in the work of~\citet{bobkov2019}.  
On the other hand, when these values are infinite, a dichotomy arises: in full generality,
 when $\sw_{r,\delta}(P,Q)$ is allowed to 
take arbitrary (small) values uncertainty quantification is difficult and our intervals can have lengths scaling as $O((1/n \wedge m)^{1/2r})$ in the worst case. However,
we find, somewhat surprisingly, even when the $\sj_{r,\delta}$ functional is infinite, accurate $O(1/\sqrt{n \wedge m})$-inference is possible so long as 
$\sw_{r,\delta}(P,Q)$ is bounded away from 0. 
We emphasize that the intervals we construct are \emph{adaptive,}
i.e. they have small lengths under appropriate conditions on the $\sj_{r,\delta}$ functional and $\sw_{r,\delta}(P,Q)$, without needing the statistician to specify or have knowledge of these quantities. 
We also show that our confidence intervals have minimax optimal length over classes of distributions
with varying magnitudes of $\sj_{r,\delta}(P)$. 
 
To complement our results on confidence intervals for the Sliced Wasserstein distance we also consider the problem of estimating the Sliced Wasserstein 
distance between two distributions, given samples from each of them. 
We provide minimax upper and lower bounds for this problem as well. Indeed, our 
minimax lower bounds for confidence interval length are derived directly from 
minimax lower bounds for estimating the Sliced Wasserstein distance 
by noting that the minimax length of a confidence interval
is bounded from below by the corresponding minimax estimation rate. 

We illustrate the practical significance of our methodology via %several applications. First, we consider
an application to likelihood-free inference \citep{sisson2018}, 
in which a parametrized stochastic simulator for the data-generating process is available, but its underlying distribution is intractable. 
Here, our goal is to construct confidence intervals for unknown parameters of the simulator, on the basis of minimizing its Sliced Wasserstein
distance from an observed sample. 
Distributional assumptions such as those made in past work   on inference for 
the one-dimensional Wasserstein distance
\citep{munk1998,freitag2003,freitag2007,freitag2005}
are typically unverifiable in such applications.

\textbf{Our Contributions.} 
We summarize the contributions of this paper as follows.
\begin{itemize}
\item We define the $\delta$-trimmed Sliced Wasserstein distance $\sw_{r,\delta}$, and 
the functional $\sj_{r,\delta}$, generalizing the functional $J_r$ of \cite{bobkov2019}. 
We show that the finiteness of $\sj_{r,\delta}(P)$
is a sufficient condition for the empirical measure to estimate $P$ at the parametric
rate under the trimmed Sliced Wasserstein distance, and we prove corresponding minimax lower bounds.
We also derive minimax rates of estimating the Sliced Wasserstein distance between two distributions, 
both in the trimmed and untrimmed settings. These  rates
are sensitive to the magnitude of the $\sj_{r,\delta}$ functional. 
\item We propose two-sample
confidence intervals for $\sw_{r,\delta}(P,Q)$ which have finite-sample coverage under minimal moment assumptions. 
 We bound the length of our confidence intervals, showing that they are adaptive both
 to the magnitude of $\sj_{r,\delta}(P), \sj_{r,\delta}(Q)$ and to whether or not $P = Q$.
 These lengths achieve the minimax rate of estimating the Sliced Wasserstein distance,
 up to polylogarithmic factors. 
\item We further contrast our finite-sample confidence intervals with asymptotic methods.
In particular, under certain regularity conditions, we derive limit laws and show that 
the bootstrap is consistent in estimating the distribution of the empirical
$r$-Sliced Wasserstein distance for all $r > 1$, whenever $P\neq Q$.  
We then show how this last assumption may be removed by
combining the strengths of our finite-sample intervals and the bootstrap.

\item We illustrate our theoretical findings with a simulation study and an application
to likelihood-free inference. 
\end{itemize}

% 
%\textbf{Paper Outline.} The rest of this paper is organized as follows. Section \ref{sec::background}
%contains background on Wasserstein distances, and defines the trimmed Sliced Wasserstein distance. In
% Section \ref{sec::minimax},  
%we define the $\sj_{r,\delta}$ functional, and we 
%establish minimax rates of estimating a distribution under the Sliced Wasserstein distance. We also
%bound the minimax risk of estimating the Sliced Wasserstein distance. 
%In Section \ref{sec::cis}, we derive confidence intervals for the one-dimensional and Sliced Wasserstein distances
%which have finite-sample coverage under minimal assumptions. We discuss asymptotically-valid approaches
%in Section \ref{sec::asymptotic}.  
%We illustrate the performance of our confidence intervals via a
%simulation study in Section \ref{sec::simulations}, and in Section \ref{sec::applications},
%we describe applications of our methodology to likelihood-free inference.
%We close with discussions in Section \ref{sec:conclusion}.

\textbf{Notation.} In what follows, given a vector $x = (x_1, \dots, x_d) \in \bbR^d$, $\norm{x} = (\sum_{i=1}^d x_i^2)^{1/2}$
denotes the $\ell_2$ norm of $x$. For any $a, b\in \bbR$, $a \vee b$ denotes the maximum of $a$ and $b$, and $a \wedge b$ denotes the
minimum of $a$ and $b$. For any function $f$ mapping a set $A$ to $\bbR$, its supremum norm is denoted by
$\norm f_\infty = \sup_{x \in A} |f(x)|$. For any sequences of real numbers $\{a_n\}_{n=1}^{\infty}$ and $\{b_n\}_{n=1}^{\infty}$, 
we write $a_n \lesssim b_n$ if there exists a universal constant $C > 0$ such that $a_n \leq C b_n$, and we write
$a_n \asymp b_n$ if $a_n \lesssim b_n \lesssim a_n$. When the constant
$C$ depends on another numerical constant $x \in \bbR$, 
we either state this dependence explicitly, or write $a_n\lesssim_x b_n$
and $a_n \asymp_x b_n$.
 $\delta_x$
denotes the Dirac delta measure at a point $x \in \bbR^d$.
The Lebesgue measure on $\bbR^k$, for an integer $k\geq 1$ to be understood
from context, is denoted $\lambda$.
Given a map $T : \bbR^d \to \bbR^d$ and a Borel probability measure
$P$ supported in $\bbR^d$, $T_{\#} P$ denotes the pushforward
of $P$ under $T$, defined by $T_{\#} P(B) = P(T^{-1}(B))$
for all Borel sets $B \subseteq \bbR^d$. We also denote
by $P^{\otimes n}$ the $n$-fold product measure of $P$. 
For any set $A \subseteq \bbR^d$, its diameter
is denoted
$\diam(A) = \sup\{\norm{x-y}: x,y \in A\}$. 
For any real numbers $a, b \in \bbR$, $I(a \leq b)$
is the indicator function equal to $1$ if $a\leq b$ and 0 otherwise.

\section{Background and Related Work}
\label{sec::background}
In this section we first provide some background on the Wasserstein distance and its sliced counterpart before turning our attention to a detailed discussion of related work.

\subsection{The Wasserstein Distance}
\label{sec::background_wasserstein}
Let $\calP(\calX)$ denote the set of Borel probability
measures whose support is contained in  $\calX \subseteq\bbR^d$. For all $r \geq 1$, let
$\calP_r(\calX)$ denote the subset of measures in $\calP(\calX)$
admitting finite $r$-th moment. 

\paragraph{The One-Dimensional Wasserstein Distance.} 
The infimum in the definition
of the Wasserstein distance \eqref{eq::wasserstein} is always achieved in the setting of this paper (cf. Theorem 4.1, \cite{villani2003}).
Closed form expressions for the minimizer are, however, unavailable in general. The one-dimensional case is a key exception.

Let $\calX \subseteq \bbR$ and $P,Q \in \calP_r(\calX)$. 
Let $F, G$ denote the cumulative distribution functions (CDFs) of $P$ and $Q$, and denote
their respective quantile functions by $F^{-1}$ and $G^{-1}$, where
$F^{-1}(u)= \inf\{x \in \bbR: F(x) \geq u\}$ for all $u \in [0,1]$.
We extend $F^{-1}$ to be defined over the entire real line under the convention
$F^{-1}(u) = \inf(\calX)$ for all $u < 0$ and $F^{-1}(u) = \sup(\calX)$ for all $u > 1$, and similarly for $G^{-1}$.  
The one-dimensional Wasserstein distance admits the closed form \citep{bobkov2019} 
%One-dimensional probability distributions form a key special case in which the %infimum in the definition of
%Wasserstein distance (equation \eqref{eq::wasserstein}) admits a closed form \citep{bobkov2019}, given by
\begin{equation}
\label{eq::wass-onedim}
W_r(P,Q) = \left(\int_0^1 \big|F^{-1}(u) - G^{-1}(u)\big|^rdu \right)^{1/r}.
\end{equation} 

\paragraph{$\infty$-Wasserstein Distance.}
Let $\calX \subseteq \bbR^d$ be a bounded set, and $P,Q \in \calP(\calX)$. 
In this case, the limit,
\begin{equation}
\label{eq::W_infty}
W_\infty(P,Q) := \lim_{r \to \infty} W_r(P,Q) = \sup_{r \geq 1} W_r(P,Q)
\end{equation}
exists, and defines a new metric $W_\infty$ on $\calP(\calX)$.
In the special case $\calX \subseteq \bbR$, we have,
$W_\infty(P,Q) = \sup_{0\leq u \leq 1} \big|F^{-1}(u) - G^{-1}(u)\big|.$
The relationship $W_r(P,Q) \leq W_\infty(P,Q)$ shows that $W_\infty$ is a stronger metric than $W_r$ for any $r \geq 1$. 
In fact, it is strictly stronger: 
for instance, $W_r(\delta_0, (1-\epsilon) \delta_0 + \epsilon\delta_1) \to 0$ 
as $\epsilon \to 0$ if and only if $r$ is finite.
In contrast, the metrics $W_r$ induce the same (weak) topology for all finite $r$, when $\diam(\calX) < \infty$ \citep{villani2003}.

\paragraph{The One-Dimensional Trimmed Wasserstein Distance.} 
Given distributions $P,Q \in \calP(\bbR)$ and a trimming constant $\delta \in [0,1/2)$, \cite{munk1998} define the $\delta$-trimmed Wasserstein
distance (up to rescaling) by
\begin{equation}
\label{eq::trimmed_W}
W_{r,\delta}(P,Q) = \left(\frac 1 {1-2\delta} \int_\delta^{1-\delta} \big|F^{-1}(u) - G^{-1}(u)\big|^r du \right)^{\frac 1 r}.
\end{equation}
When $\delta =0$, $W_{r,\delta}$ reduces to the original Wasserstein distance $W_r$, 
and when $\delta > 0$, $W_{r,\delta}$ compares the distributions $P$ and $Q$
up to a $2\delta$ fraction of their tail mass. %, thereby providing a robustification of the Wasserstein distance. 
Specifically, let $P^\delta$ denote the distribution with CDF
$F^\delta(x) = (F(x)-\delta)/(1-2\delta)$, for all $ F^{-1}(\delta) \leq x \leq F^{-1}(1-\delta)$, 
and similarly for $Q^\delta$. Then, \cite{alvarez-esteban2008} note that $
W_{r,\delta}(P,Q) = W_r(P^\delta, Q^\delta).
$

In addition, we define the trimmed $\infty$-Wasserstein distance by
$$W_{\infty,\delta}(P,Q) := \lim_{r \to \infty} W_{r,\delta}(P,Q) =  \sup_{\delta \leq u \leq 1-\delta} \big|F^{-1}(u) - G^{-1}(u)\big|.$$
 
\subsection{The Sliced Wasserstein Distance}
The Sliced Wasserstein distance 
\citep{rabin2011} %,bonneel2015} 
is defined as the mean of Wasserstein
distances between one-dimensional projections of the distributions $P$ and $Q$.
Specifically, let $\bbS^{d-1} = \{x \in \bbR^d: \norm x = 1\}$, 
and let $\pi_\theta : x \in \bbR^d \mapsto x^\top\theta$, for all $ \theta \in \bbS^{d-1}$.
Let $P_\theta = {\pi_\theta}_\# P$ and $Q_\theta = {\pi_\theta}_\# Q$, that is,
$P_\theta$ and $Q_\theta$ are the respective probability distributions of $X^\top\theta$ and $Y^\top\theta$,
for $X \sim P$ and $Y \sim Q$. 
Let $\mu$ denote the uniform probability measure on $\bbS^{d-1}$. %, that is, 
The $r$-th order Sliced Wasserstein distance between two distributions $P,Q \in \calP_r(\bbR^d)$
is given by
\begin{equation}
\label{eq::SW}
\sw_r(P,Q) = \left(\int_{\bbS^{d-1}} W_r^r(P_\theta, Q_\theta) d\mu(\theta)\right)^{\frac 1 r}.
\end{equation}
Since 
$P_\theta$ and $Q_\theta$ are one-dimensional distributions, 
equation \eqref{eq::SW} admits the closed form
\begin{equation}
\sw_r(P,Q) =   \left( 
\int_{\bbS^{d-1}} \int_0^1 \big| F_\theta^{-1}(u)- G_\theta^{-1}(u)\big|^r du d\mu(\theta)\right)^{\frac 1 r},
\end{equation} 
where $F_\theta^{-1}$ and $G_\theta^{-1}$ are the respective quantile functions of $P_\theta$ and $Q_\theta$.
Both integrals of the above expression can be approximated via Monte Carlo 
sampling from $\bbS^{d-1}$ and from the unit interval $[0,1]$.
This fact makes the computation of the Sliced
Wasserstein distance significantly simpler than that of the Wasserstein distance. 
Moreover, the Sliced Wasserstein distance  retains some of the qualitative behaviour of the Wasserstein distance, 
at least for compactly-supported distributions. Indeed, 
\cite{bonnotte2013} showed that for any distributions  $P,Q \in \calP(\{x \in \bbR^d: \norm x \leq M\})$,
where $M > 0$, we have
\begin{equation}
\label{eq::sw_w_sw}
\sw_r^r(P, Q) \leq c_{d,r} W_r^r(P,Q) \leq C_{d,r} M^{r-1/(d+1)} \sw_r^{1/(d+1)}(P,Q),
\end{equation}
where
$C_{d,r} > 0$ is a constant depending on $d$ and $r$, but not depending on $M$, and
$c_{d,r} = \frac 1 d \int_{\bbS^{d-1}} \norm \theta_r^r d\mu(\theta)$, which
is bounded above by $1/d$ whenever $r \geq 2$.
In particular, it follows that the metrics $W_r$ and $\sw_r$ are topologically equivalent
over $\calP(\calX)$ when $\diam(\calX) < \infty$. As we shall see, however, the statistical
behaviour of the Wasserstein and Sliced Wasserstein distances can differ dramatically
for large dimensions $d$. 

Though the original Sliced Wasserstein distance of~\cite{rabin2011}
was defined in terms of the uniform distribution $\mu$ over $\bbS^{d-1}$, 
a straightforward
adaptation of Proposition 5.12 of~\cite{bonnotte2013} shows that 
$\sw_{r}$ remains a metric over $\calP(\bbR^d)$ when $\mu$ is replaced by
any probability measure which is absolutely continuous with respect to the Hausdorff measure on $\bbS^{d-1}$. 
We allow $\mu$ to be any such measure throughout the sequel.

\paragraph{The Trimmed Sliced Wasserstein Distance.} In analogy to the trimmed Wasserstein distance in equation \eqref{eq::trimmed_W}, 
we further define %the trimmed Sliced Wasserstein distance by
\begin{equation}
\label{eq::trimmed_SW}
\sw_{r,\delta}(P,Q) = \left(\frac 1 {1-2\delta} \int_{\bbS^{d-1}} \int_\delta^{1-\delta} \big|F_\theta^{-1}(u) - G_\theta^{-1}(u) \big|^rdud\mu(\theta)\right)^{\frac 1 r},
\end{equation}
for some $\delta \in [0,1/2)$. We also define the trimmed $\infty$-Sliced Wasserstein distance by
$\sw_{\infty,\delta}(P,Q) %:= \sup_{r \geq 1}~\sw_{r,\delta}(P,Q) 
= \int_{\bbS^{d-1}} W_{\infty,\delta}(P_\theta, Q_\theta) d\mu(\theta)$,
and more generally,  we write 
$$\sw_{\infty,\delta}^{(r)}(P,Q) %= \norm{W_{\infty,\delta}(P_\theta,Q_\theta)}_{L^r(\mu)}^r = 
=\int_{\bbS^{d-1}} W_{\infty,\delta}^r(P_\theta, Q_\theta) d\mu(\theta).$$
 When $\delta>0$, the trimmed-Sliced Wasserstein  distance 
is well-defined and finite for all distributions $P,Q \in \calP(\bbR^d)$, including those
admitting fewer than $r$ moments. Nevertheless, it can be easily seen that 
$\sup_{P,Q \in \calP(\bbR^d)} \sw_{r,\delta}(P,Q)=\infty$. It will 
be fruitful in our development to impose moment conditions which ensure that
the quantity $\sw_{r,\delta}(P,Q)$ is uniformly bounded---one such condition is
given in terms of the class 
\begin{equation}
\label{eq::calK}
\calK_{r,\rho}(b) = \left\{P \in \calP(\bbR^d): 
  \int_{\bbS^{d-1}}  \bbE_{X \sim P}\big[ |X^\top\theta|^\rho\big]^{\frac  r \rho}d\mu(\theta) \leq b\right\}, \quad b,\rho,r \geq 1.
\end{equation} 
We shall use the special case $\rho=2$ most often, in which case we drop the subscript $\rho$ and simply write
$\calK_{r}(b): = \calK_{r,2}(b)$. 
It follows from Lemma~\ref{lem::calK} in Appendix~\ref{app::preliminary}
that $\sw_{r,\delta}(P,Q)$ is uniformly bounded over distributions $P,Q \in \calK_{r}(b)$,
by a constant depending only on $b,r$ and $\delta$. Notice  that if $\bar b = b^{\rho/r}$, then $\calK_{r,\rho}(b)$ 
contains the class %$\widebar \calK_\rho(\bar b)$ of distributions with  second moment bounded above by~$\bar b$,
\begin{equation}
\label{eq::bar_calK} 
\widebar \calK_\rho(\bar b) = \left\{P \in \calP_\rho(\bbR^d): 
   \bbE_{X \sim P}\big[\norm{X}^\rho\big] \leq \bar b\right\}.
\end{equation} 
Finally, we write  $\calK_{r,\rho} = \bigcup_{b \geq 1} \calK_{r,\rho}(b)$,
$\calK_r = \bigcup_{b \geq 1} \calK_r(b)$ and $\widebar\calK_\rho = \bigcup_{\bar b \geq 1} \widebar\calK_\rho(\bar b)$. 

\subsection{Related Work}
\label{sec::related_work}
We are unaware of any other work regarding statistical
inference for the Sliced Wasserstein distance, except in the special
case $d=1$ when it coincides with the one-dimensional Wasserstein distance. 
In this case, \cite{munk1998} study limiting distributions of the 
empirical (plug-in) Wasserstein distance estimator, and \cite{freitag2003,freitag2007,freitag2005} establish sufficient conditions
for the validity of the bootstrap in estimating the distribution of the empirical  second-order trimmed Wasserstein distance. 
While these results are very useful, they assume that (i) $P$ and $Q$ are absolutely continuous, (ii) with densities supported on connected sets,
and (iii) require different inferential procedures at the classical null $(P=Q)$ and away from the null $(P\neq Q)$.
In contrast, the confidence intervals derived in the present paper are valid under either 
no assumptions or mild moment assumptions on $P$ and $Q$,
and are applied more generally to the Sliced Wasserstein distance in arbitrary dimension.
Though our methodology is assumption-light, our confidence intervals are adaptive to (iii), and
assumptions (i) and (ii) are closely related to the finiteness
of $\sj_{r,\delta}(P),\sj_{r,\delta}(Q)$, 
to which our confidence intervals are also adaptive. 
 
The Sliced Wasserstein distance is one of many modifications of the Wasserstein
distance based on low-dimensional projections. We mention here 
the Generalized Sliced \citep{kolouri2019}, Tree-Sliced \citep{le2019},
max-Sliced \citep{deshpande2019}, Subspace Robust \citep{paty2019a,niles-weed2019}, and Distributional Sliced \citep{nguyen2020} Wasserstein distances. 
It is also possible to define various other interesting distances by slicing (averaging along univariate projections---see \cite{kim2020}).

Beyond the aforementioned inferential results %of \cite{munk1998, freitag2003, freitag2007,freitag2005} 
for the one-dimensional Wasserstein distance, 
statistical inference for Wasserstein distances over finite or countable spaces has been studied by 
\cite{sommerfeld2018a, tameling2019, klatt2020a, klatt2020b}.
For distributions with multidimensional support, \cite{rippl2016} consider the situation where $P$ and $Q$ only differ by a location-scale transformation.
\cite{imaizumi2019} study the validity of the multiplier bootstrap for estimating the distribution of 
the plug-in estimator of an approximation of the Wasserstein distance.
Central Limit Theorems for empirical Wasserstein distances in general dimension
have been established by \cite{delbarrio2019a,delbarrio2021}, but with unknown centering constants which are a barrier to using these results for statistical inference.

Rates of convergence for the problem of estimating a distribution under the Wasserstein distance (\cite{dudley1969, %bolley2006, 
boissard2014, 
fournier2015, bobkov2019, weed2019,singh2019, lei2020}, and references therein) 
have received significantly more attention than the problem
of estimating the Wasserstein distance, the latter being more closely related to our work. 
Minimax rates of estimating the Wasserstein distance between two distributions
have been established 
by \cite{niles-weed2019}, as well as by \cite{liang2019} when $r=1$. 
In the special case $d=1$, where
the Sliced Wasserstein distance coincides with the Wasserstein distance, our
results refine those of \cite{niles-weed2019} by showing that faster rates can be achieved depending on the 
finiteness of the $\sj_{r,\delta}$ functional, and on the magnitude of $\sw_{r,\delta}(P,Q)$.

Likelihood-free inference methodology with respect to the Wasserstein and Sliced Wasserstein distances has recently been developed
by \cite{bernton2019} and \cite{nadjahi2020}. %, respectively. 
In contrast to these methods, both of which employ approximate Bayesian computation, 
our work provides frequentist coverage guarantees under minimal assumptions.

\section{Estimating the Sliced Wasserstein Distance}
\label{sec::minimax}

The goal of this section is to bound the minimax risk of estimating the Sliced
Wasserstein distance between two distributions, that is 
\begin{equation}
\label{eq::minimax_distance}
\calR_{nm} \equiv \calR_{nm}(\calO; r) = \inf_{\hat S_{nm}} \sup_{(P,Q) \in \calO} \bbE_{P^{\otimes n} \otimes Q^{\otimes m}} \big| \hat S_{nm} - \sw_{r,\delta}(P,Q)\big|,
\end{equation}
where the infimum is over all estimators $\hat S_{nm}$ of the Sliced Wasserstein
distance based on a sample of size $n$ from $P$ and a sample of size $m$ from $Q$, 
and $\calO\subseteq \calP(\bbR^d) \times \calP(\bbR^d)$ is a %n appropriately chosen 
collection of pairs of distributions. 
Our motivation for studying this quantity is the observation that 
$\calR_{nm}$ lower bounds the minimax length of 
a confidence interval for the Sliced Wasserstein distance.  
We construct confidence
intervals with matching length in Section \ref{sec::cis}.

The estimation problem in equation \eqref{eq::minimax_distance} is related to, but distinct from, the
problem of estimating a distribution under the Sliced Wasserstein distance. The minimax
risk associated with this problem is given by
\begin{equation}
\label{eq::minimax_distribution}
\calM_n \equiv \calM_n(\calJ;r) = \inf_{\hat P_n} \sup_{P\in\calJ} \bbE_{P^{\otimes n}} \Big\{ \sw_{r,\delta}(\hat P_n,P)\Big\},
\end{equation}
where the infimum is over all estimators $\hat P_n$ of Borel probability distributions~$P$,
based on a sample of size $n$ from $P$, and $\calJ \subseteq \calP(\bbR^d)$.
Problems~\eqref{eq::minimax_distance} and~\eqref{eq::minimax_distribution} 
are related as follows: Given estimators $\hat P_n$ and $\hat Q_m$ for
two distributions $P,Q \in \calP(\bbR^d)$, which are minimax optimal in the sense of equation
\eqref{eq::minimax_distribution},
we have, by the triangle inequality,
\begin{align}
\nonumber \calR_{nm}(\calO; r) 
           &\lesssim \bbE \big| \sw_{r,\delta}(\hat P_n, \hat Q_m) - \sw_{r,\delta}(P,Q) \big| \\
  &\leq \bbE  \sw_{r,\delta}(\hat P_n, P)  + \bbE \sw_{r,\delta}(\hat Q_m, Q)   % \label{eq::minimax_loose} 
% &\asymp \calM_n(\calJ; r) + \calM_m(\calJ; r) \\ 
 \lesssim \calM_{n\wedge m}(\calJ;r),
 \label{eq::minimax_relation}
\end{align}
for suitable families $\calJ$ and $\calO$ (typically $\calO \subseteq \calJ \times \calJ$). 
Inequality \eqref{eq::minimax_relation}
implies that estimating a distribution under $\sw_{r,\delta}$
is a more challenging problem, statistically, than that of estimating  
the Sliced Wasserstein distance between two distributions. 
It is unclear, however, whether the rate $\calM_{n\wedge m}$
is a tight upper bound on $\calR_{nm}$, or whether the latter can be further reduced. 
For the Wasserstein distance $W_r$ in general dimension,  
\cite{liang2019} and \cite{niles-weed2019} showed that there is no gap between these minimax risks (ignoring polylogarithmic
factors)
for compactly supported distributions.   

Let us now briefly summarize the main results of this section. 
We bound $\calM_n$ and $\calR_{nm}$, and show that 
there can be a large gap between these minimax risks when the pairs of distributions
in $\calO$ are appropriately separated.
In the special case $d=1$, $\sw_{r,\delta}$ reduces to the 
 (trimmed) Wasserstein distance, and our results imply faster rates 
 than those of \cite{liang2019} and  \cite{niles-weed2019}, 
 for estimating
 the Wasserstein distance between distributions bounded away from each other. 
 Furthermore,  in   contrast to the minimax risk for estimating the Wasserstein distance and estimating
under  the Wasserstein distance, the minimax risks 
 we obtain for the Sliced Wasserstein distance when $d > 1$ are dimension-free.
  
 Though our primary
 interest is in   $\calR_{nm}$ (due to its direct connection to confidence
 intervals) we begin by studying $\calM_n$ to motivate our choices
 of families~$\calO$. Inspired by \cite{bobkov2019}, in Section \ref{sec::minimax_distribution} 
 we define a functional $\sj_{r,\delta}$, whose magnitude
 is related to the regularity of the supports of $P$ and $Q$, and whose finiteness implies
 improved rates of decay for $\calM_n$. We then study the minimax risk $\calR_{nm}$
 over various families~$\calO$ 
 in Section \ref{sec::minimax_distance}.

\subsection{Minimax  Estimation under the Sliced Wasserstein Distance} 
\label{sec::minimax_distribution} 

Let $\delta \in [0,1/2)$, $P \in \calP(\bbR^d)$, and let $X_1, \dots, X_n \sim P$ be an i.i.d. sample. 
Let $P_n = \frac 1 n \sum_{i=1}^n \delta_{X_i}$
denote the corresponding empirical measure. 
The goal of this section is to characterize the rates of convergence of $P_n$ to the distribution $P$
under the (trimmed) Sliced Wasserstein distance, extending the comprehensive treatment by \cite{bobkov2019}.
We then provide corresponding minimax lower bounds on $\calM_n$.

For any $\theta \in \bbS^{d-1}$, let $p_\theta$ denote the density of the absolutely continuous 
component in the Lebesgue decomposition of the measure $P_\theta = {\pi_\theta}_{\#} P$. 
Define the functional
\begin{equation}
\label{eq::sj}
\sj_{r,\delta}(P) = \frac 1 {1-2\delta}\int_{\bbS^{d-1}} 
\int_{\delta}^{1-\delta} \left(\frac{\sqrt{u(1-u)}}{p_\theta(F_\theta^{-1}(u))}\right)^r du d\mu(\theta),
\end{equation}
with the convention that $0/0=0$. When $d = 1$, we write $\text{J}_{r,\delta}$ instead of $\sj_{r,\delta}$, 
and in the untrimmed case $\delta = 0$, 
we omit the subscript $\delta$ and write $\text{SJ}_r$ or $\text{J}_r$. 
When $d=1$ and $\delta = 0$, \cite{bobkov2019} prove that 
the finiteness of $\text{J}_r(P)$ is a necessary and sufficient condition for $\bbE\big[W_r(P_n,P)\big]$
to decay at the parametric rate $n^{-1/2}$. 
The magnitude of $\text{J}_r$ is thus closely related
to the convergence behaviour of empirical measures under 
one-dimensional Wasserstein distances, and we show below that the same is true for the $\sj_{r,\delta}$ functional
with respect to trimmed Sliced Wasserstein distances, using distinct proof techniques.

It can be seen that a necessary condition for the finiteness of $\sj_{r,\delta}(P)$ is that
for $\mu$-almost every $\theta \in \bbS^{d-1}$, the density $p_\theta$ is supported
on a (possibly infinite) interval.  
When $\delta$ vanishes, 
the value of $\sj_{r,\delta}(P)$ also depends on the tail behaviour of $P$ 
and the value of $r$. For example, 
if $P = N(0, I_d)$ is the standard Gaussian distribution, 
it can be shown that $\sj_{r,\delta}(P) < \infty$ whenever $\delta > 0$, 
whereas for $\delta=0$, $\text{SJ}_r(P) < \infty$ if and only if $1 \leq r < 2$ by a similar argument as \citeauthor{bobkov2019}
(\citeyear{bobkov2019}, p. 46).
On the other hand, 
if $P = \frac 1 2 U(0,\Delta_1) + \frac 1 2 U(\Delta_2, 1)$, for some $0 < \Delta_1 \leq \Delta_2 < 1$, 
where $U(a, b)$ denotes the uniform distribution on the interval $(a,b) \subseteq \bbR$, 
%then a direct calculation reveals that
one has $\sj_{r,\delta}(P) < \infty$
if and only if $\Delta_1 = \Delta_2$, for every $\delta \in [0,1/2)$. 

We now provide two upper bounds on $\bbE\big[\sw_{r,\delta}(P_n,P)\big]$, 
which are effective when $\sj_{r,\delta}(P) < \infty$ and $\sj_{r,\delta}(P) = \infty$ respectively.
Recall the class $\calK_{r}(b)$ from equation~\eqref{eq::calK}.
\begin{proposition} 
\label{prop::ub_empirical}
Let $b,r\geq 1$ and $\delta \in (0,1/2)$. Assume $P \in \calK_r(b)$,
and that $\delta \geq 2(r+2)/n$.
\begin{itemize}
\item[(i)] There exist constants $c_r,c_r' > 0$ depending only on $r$ such that
\begin{align}
\label{eq::ub_empirical_fast}
\bbE\big[\sw_{r,\delta}(P_n,P)\big]\leq c_r \sj_{r,\frac \delta 2}^{\frac 1 r}(P)
\sqrt{\frac{\log n}{n}}+ \frac{c_r(b n  e^{-c_r' n\delta})^{\frac 1 r}}{\sqrt \delta}.
\end{align}
%\item[(ii)] A necessary and sufficient condition for the rate $\bbE \big[\sw_{r,\delta}(P,Q)\big] = O(n^{-1/2})$ 
%to hold is $\sj_{r,\delta}(P) < \infty$.
\item[(ii)]   There exists a constant $k_r > 0$ depending only on $r$ such that 
%$$\bbE\big[\sw_{r,\delta}(P_n,P)\big] \leq C n^{-1/2r}.$$
\begin{align}
\bbE\big[\sw_{r,\delta}(P_n,P)\big] \leq \frac{k_r}{\sqrt\delta} \left(\frac b {1-2\delta}\right)^{\frac 1 r}n^{-1/2r}.
\end{align} 
\end{itemize}
\end{proposition} 
Proposition
\ref{prop::ub_empirical} provides two upper bounds on the 
rate of convergence of the empirical measure under $\sw_{r,\delta}$, 
which are closely related to those of Theorem~5.3 and Theorem~7.16 
of \cite{bobkov2019} for the one-dimensional untrimmed Wasserstein distance. 
\cite{bobkov2019} established these results by hinging upon a representation
of the empirical one-dimensional Wasserstein distance  
in terms of the so-called mean square beta distribution, coupled with Poincar\'e-type inequalities
for such measures. While extensions of these techniques to the {\it untrimmed} Sliced
Wasserstein distance are straightforward,  and will be stated for completeness in Section~\ref{sec:untrimmed},
it was unclear to us whether they may be adapted to the {\it trimmed} setting $\delta > 0$. 
Our proof of Proposition~\ref{prop::ub_empirical} is instead based on uniform concentration inequalities
for the empirical quantile process---an approach which we now describe, and which foreshadows 
the construction of our confidence intervals
in Section~\ref{sec::cis}.

Proposition~\ref{prop::ub_empirical}(i) is proved using a uniform bound
for self-normalized empirical processes, known as the relative Vapnik-Chervonenkis (VC) inequality~\citep{vapnik2013,bousquet2003}, which
will be further described in Example~\ref{ex::rel_VC} below. This result
shows that the empirical measure converges at the parametric rate under $\sw_{r,\delta}$, up to a polylogarithmic factor, 
provided $\sj_{r,\delta/2}(P)$ is bounded, and provided, for instance, that the trimming constant $\delta$ does not vanish at a rate
faster than $n^{-\beta}$ for some $\beta \in (0,1)$. We emphasize that this convergence is uniform in $P$---for instance, 
one has that for all $s \geq 1$,  
$$\sup_{\substack{P \in \calK_r(b) \\ \sj_{r,\delta/2}(P) \leq s}} \bbE\big[\sw_{r,\delta}(P_n,P)\big] 
 \underset{b,r}{\lesssim} s^{\frac 1 r} \sqrt{\frac{\log n}{n}},
~\text{when } \delta \asymp n^{-\beta}, \text{for some } \beta \in (0,1).$$
In fact, the above bound continues to hold when $s \leq 1$---a regime relevant for distributions with vanishing
variances---so long as $s^{1/r}\sqrt{\log n/n}$ remains greater than the second term in equation~\eqref{eq::ub_empirical_fast}. 
Finally, we note that the polylogarithmic factor in the above display arises from the relative VC inequality. 
Example 2.7 of \cite{gine2006} suggests that this factor may be improved to $\sqrt{\log\log n}$, but not further
if $\delta\gtrsim 1/n$. We do not know whether this factor can be omitted if stronger conditions are placed on $\delta$
while allowing it to vanish.

Proposition~\ref{prop::ub_empirical}(ii) is primarily of interest for distributions such that
$\sj_{r,\delta}(P) = \infty$, 
and is proved using the Dvoretzky-Kiefer-Wolfowitz inequality \citep{dvoretzky1956, massart1990}. 
This result shows that the empirical measure converges to $P$ at the nonparametric rate $n^{-1/2r}$
under no assumptions on $P$ apart from the mild moment assumption $P \in \calK_{r}(b)$. 
In contrast to Proposition~\ref{prop::ub_empirical}(i), however, this result suffers from a markedly worse
dependence on $\delta$. Indeed, the resulting rate of convergence deteriorates as soon as $\delta = o(1)$, 
and we conjecture that this behaviour is necessary under the stated assumptions.

The rates in Proposition \ref{prop::ub_empirical} 
do not depend on the dimension $d$, contrasting
 generic rates of convergence of the empirical measure under the Wasserstein distance. 
For instance, if $P$ is supported on a bounded set in $\bbR^d$,
\cite{lei2020} (see also \cite{fournier2015}, \cite{weed2019}) shows that 
$\bbE [W_r(P_n,P)] \lesssim n^{-1/d}$ whenever $d > 2r$, and this rate is known to be tight~\citep{dudley1969,singh2019}. 
Thus, estimating a distribution in 
the Sliced Wasserstein distance does not suffer from the
 curse of dimensionality 
 despite
 metrizing the same topology on $\calP(\bbR^d)$---see equation \eqref{eq::sw_w_sw}.

For completeness, we close this subsection 
by stating a lower bound on the minimax risk $\calM_n$ in equation \eqref{eq::minimax_distribution}. 
In view of Proposition \ref{prop::ub_empirical},
it is natural to carry out our analysis over the  
class of distributions
$$\calJ(s) = \big\{P \in \calK_{r}(b): \sj_{r,\delta}(P) \leq s \big\}, \quad s \in [0,\infty].$$
\begin{proposition}
\label{thm::minimax_distr}
Let $b,r\geq 1$ and $\delta \in (0,1/2)$. 
 Then, there exist constants $C_1, C_2 > 0$, possibly 
depending on $b, r, \delta$, such that for all $s > 0$ satisfying $b \geq (2s)^{1/r}$,
$$\calM_n(\calJ(s);r) \geq C_1 s^{\frac 1 r}n^{-1/2}, \quad \text{and} \quad 
  \calM_n(\calJ(\infty);r)\geq C_2 n^{-1/2r}.$$
\end{proposition}
Proposition \ref{thm::minimax_distr} implies that the rates 
achieved by the empirical measure in Proposition \ref{prop::ub_empirical} are minimax optimal
over the classes considered above, up to polylogarithmic factors.
The proof of this result will follow as a special case of our bounds on the minimax risk $\calR_{nm}$, to which we turn our attention next.

\subsection{Minimax Estimation of the Sliced Wasserstein Distance}
\label{sec::minimax_distance}
In this section, we bound the minimax risk $\calR_{nm}$  defined in equation \eqref{eq::minimax_distance}. 
We begin by providing upper bounds on the estimation error of the empirical Sliced Wasserstein distance, $\sw_{r,\delta}(P_n, Q_m)$.
Recall that, 
$P_n$ and  $Q_m$ denote the empirical measures of i.i.d. samples $X_1, \dots, X_n \sim P$ and $Y_1, \dots, Y_m \sim Q$, respectively. 
\begin{proposition}
\label{prop::empirical_sw} 
Let $b,r\geq 1$ and $\delta \in (0,1/2)$. Assume $P,Q \in \calK_r(b)$,
and that $\delta \geq 2(r+2)/(n\wedge m)$.  
\begin{enumerate}
\item[(i)] There exists a constant $c_r > 0$, possibly depending on $r$, such that
\begin{align*}
\bbE \big| \sw_{r,\delta}&(P_n,Q_m) - \sw_{r,\delta}(P,Q)\big| \\& \underset{r}{\lesssim} 
\left(\frac{b}{1-2\delta}\right)^{\frac 1 r} \frac{n^{-\frac 1 {2r}}}{\sqrt\delta} \wedge\left( \sj_{r,\frac \delta 2}^{\frac 1 r}(P)\sqrt{\frac{\log n}{n}}+
\frac{(b n  e^{-c_r n\delta})^{\frac 1 r}}{\sqrt \delta}\right)  \\ & \ + 
\left(\frac{b}{1-2\delta}\right)^{\frac 1 r} \frac{m^{-\frac 1 {2r}}}{\sqrt\delta}\wedge\left( \sj_{r,\frac \delta 2}^{\frac 1 r}(Q)\sqrt{\frac{\log m}{m}}+
\frac{(b m  e^{-c_r m\delta})^{\frac 1 r}}{\sqrt \delta}\right).
\end{align*}

\item[(ii)] Suppose $\sw_{r,\delta}(P,Q) \geq \Gamma$, for some real number $\Gamma > 0$. Then,  
$$\bbE \big| \sw_{r,\delta}(P_n,Q_m) - \sw_{r,\delta}(P,Q)\big| \underset{\Gamma,r}{\lesssim}   \frac{b}{\delta^{r/2}(1-2\delta)}\left(n^{-\frac 1 2} + m^{-\frac 1 2}\right).$$
\end{enumerate}
\end{proposition}
 
Proposition 
\ref{prop::empirical_sw}(i) is an immediate consequence of inequality \eqref{eq::minimax_relation}, which implies that the rate of estimating the Sliced
Wasserstein distance with the plug-in estimator $\sw_{r,\delta}(P_n, Q_m)$ 
is no worse than the rate of convergence of the empirical measures under
$\sw_{r,\delta}$ given in Proposition~\ref{prop::ub_empirical}. 
In particular, these results
show that the parametric rate for estimating $\sw_{r,\delta}$ 
is achievable for distributions satisfying $\sj_{r,\delta/2}(P), \sj_{r,\delta/2}(Q) < \infty$,
while the rate $n^{-1/2r} + m^{-1/2r}$ is otherwise achievable. 
%\item 
On the other hand, Proposition \ref{prop::empirical_sw}(ii) implies that the 
parametric rate of estimating $\sw_{r,\delta}(P,Q)$ is \textit{always} achievable
when $P$ and $Q$ are bounded away from each other under $\sw_{r,\delta}$. This fast rate of
convergence is obtained irrespective of the values of $\sj_{r,\delta}(P)$ and $\sj_{r,\delta}(Q)$. 
Discrepancies between rates of convergence at the null ($P=Q$) and away from the null ($P\neq Q$) have
previously been noted by \cite{sommerfeld2018a} for Wasserstein distances over finite spaces---indeed, 
their rates match those of Proposition \ref{prop::empirical_sw} when $\sj_{r,\delta/2}(P), \sj_{r,\delta/2}(Q) = \infty$.
Finally, we note that the natural estimator $\sw_{r,\delta}(P_n,Q_m)$ is \emph{adaptive} to the typically unknown 
quantities $\sj_{r,\delta/2}(P)$ and $\sj_{r,\delta}(Q)$, and 
does not require the statistician to specify if $P = Q$ or $P \neq Q$. Instead, the estimator adapts and yields favorable rates in 
favorable situations---when either
the $\sj_{r,\delta/2}$ functionals are finite, or when $P$ and $Q$ are sufficiently well-separated.

We now provide corresponding lower bounds on the minimax risk $\calR_{nm}$. Inspired by Proposition \ref{prop::empirical_sw}, we define the following
collection of pairs of distributions,
\begin{align*}
\calO(\Gamma; s_1, s_2) {=} \Big\{(P,Q) {\in} \calK_r^2(b): \sj_{r,\delta}(P)\leq s_1, \sj_{r,\delta}(Q) \leq s_2,  \sw_{r,\delta}(P,Q) \geq \Gamma\Big\}, 
\end{align*}
where $s_1,s_2 \in [0,\infty]$ and $\Gamma \geq 0$. 
%%In constructing our lower bounds, one minor concern that we have to address is that if $\Gamma$ is chosen sufficiently large relative to $b$
%%then the class $\calO(\Gamma; s_1, s_2)$ might be empty. Consequently, 
To ensure that the class $\calO(\Gamma; s_1, s_2)$ is nonempty, we assume in what follows that $\Gamma^r \leq c_r b,$  
for some sufficiently small constant $c_r > 0$ depending only on $r$. 
With these definitions in place we now state our minimax lower bounds on the risk $\calR_{nm}$.
\begin{theorem}
\label{thm::minimax_distance}
Let $b,r\geq 1$ and $\delta \in (0,1/2)$. Fix  $s > 0$, and assume $b \geq (2s)^{1/r}$.  
\begin{enumerate}
\item[(i)] There exists a constant $C_1 > 0$, possibly depending on $\delta,r,b$, such that
for any $s_1, s_2 \in [0,\infty]$,
$$\calR_{nm}(\calO(0; s_1, s_2);r) \geq C_1 
\begin{cases}
n^{-\frac 1 {2r}} + m^{-\frac 1 {2r}}, & s_1 =  s_2=\infty\\
\frac{s_1^{\frac 1 r}}{\sqrt n} + \frac{s_2^{\frac 1 r }}{\sqrt m}, & s_1\vee s_2 \leq s.
\end{cases} $$
\item[(ii)] For any $\Gamma > 0$ such that $\Gamma^r \leq c_r b$, there exists a constant $C_2 > 0$ 
possibly depending on $\delta, r, b, \Gamma$ such that
$$\calR_{nm}(\calO(\Gamma; \infty,\infty);r) \geq C_2\left(n^{-\frac 1 2} + m^{-\frac 1 2}\right).$$
\end{enumerate}
\end{theorem} 
Theorem \ref{thm::minimax_distance} implies that the rates achieved by the empirical Sliced
Wasserstein distance $\sw_r(P_n,Q_m)$
in Proposition \ref{prop::empirical_sw}, including their dependence on the $\sj_{r,\delta}$
functional, are minimax optimal (ignoring polylogarithmic factors). 
We defer its proof to Appendix~\ref{app::main_lowerbound}. This result is proved by a standard information-theoretic technique of 
constructing pairs of distributions %(each satisfying the various hypotheses). These distributions are carefully chosen 
which are statistically indistinguishable
but have very different Sliced Wasserstein distances. We then obtain lower bounds
via an application of Le Cam's Lemma (see, for instance, Theorem 2.2 of \citet{tsybakov2008}). Beyond this careful choice of 
distributions, the bulk of our technical effort
lies in tightly bounding the various Sliced Wasserstein distances (see Lemma~\ref{lem::main_univariate} in Appendix~\ref{app::main_lowerbound}).

In Section~\ref{sec::cis}, we construct finite-sample confidence intervals for $\sw_{r,\delta}(P,Q)$
whose lengths achieve these same rates of convergence, up to polylogarithmic factors.
Before turning to these results, we discuss estimation rates in the untrimmed case when $\delta=0$.

\subsection{Minimax Estimation of the Untrimmed Sliced Wasserstein Distance}
\label{sec:untrimmed}	
Though our  results below on finite-sample and asymptotic inference will focus on the trimmed Sliced Wasserstein distance, 
as it is an estimand for which assumption-free inference is possible, 
we end this section by deriving convergence rates for estimating the untrimmed Sliced Wasserstein distance. 
In this setting, a straightforward extension of Theorem~5.3 and Theorem~7.16  of~\cite{bobkov2019} already leads to the following
untrimmed analogue of Proposition~\ref{prop::ub_empirical}, which we state for completeness.
We recall that the class $\calK_{r,\rho}(b)$ is defined in equation~\eqref{eq::calK}.  
\begin{proposition}
\label{prop:bobkov_untrimmed}
Let $r \geq 1$ and $s  > 0$. Then, 
$$\sup_{\substack{P \in \calP(\bbR^d) \\ \sj_r(P) \leq s}} \bbE \sw_r(P_n, P) \lesssim_r s^{\frac 1 r} n^{-\frac 1 2} .$$
Furthermore, for any $\rho > 2r$ and $b > 0$, 
$$\sup_{P \in \calK_{r,\rho}(b)} \bbE \sw_r(P_n, P) \lesssim_{b,\rho,r} n^{-\frac 1 {2r}}.$$
\end{proposition}
Convergence rates for estimating $\sw_{r}(P,Q)$ immediately follow
from Proposition~\ref{prop:bobkov_untrimmed}. 
For example, we obtain
$$\sup_{P,Q \in \calK_{r,\rho}(b)} \bbE \big|\sw_r(P_n,Q_m) - \sw_r(P,Q)\big| 
\lesssim_{b,\rho,r} n^{-\frac 1 {2r}} + m^{-\frac 1 {2r}}.$$
By analogy with Proposition~\ref{prop::empirical_sw}(ii), 
it is natural to ask whether the above convergence rate may be improved to the parametric rate
when $P$ and $Q$ are separated in Sliced Wasserstein distance. Such an assertion cannot be deduced from the work of~\cite{bobkov2019}, 
and is the subject of the following main result.
\begin{theorem}
\label{thm:untrimmed_rate}
For any $r,b \geq 1$, $\Gamma > 0$, and any $\rho > 2r$, it holds that
\begin{equation}
\label{eq:untrimmed_rate}
\sup_{\substack{P,Q \in \calK_{r,\rho}(b) \\ \sw_r(P,Q) \geq \Gamma}} 
\bbE \big|\sw_r(P_n,Q_m) - \sw_r(P,Q)\big| \lesssim_{\Gamma,\rho,r}  b\left(\sqrt{\frac{\log n}{n}} + 
\sqrt{\frac{\log m}{m}}\right).
\end{equation}
\end{theorem}
Theorem~\ref{thm:untrimmed_rate} proves that
the parametric rate for estimating the Sliced Wasserstein distance between well-separated distributions continues to hold in the absence of trimming, 
at the price of a polylogarithmic factor. In fact, our proof shows more generally that
the following bound holds without separation conditions on $P$ and $Q$, 
\begin{equation}
\label{eq:untrimmed_rate_power_r}
\sup_{P,Q \in \calK_{r,\rho}(b)} \bbE \big|\sw_r^r(P_n,Q_m) - \sw_r^r(P,Q)\big| \lesssim_{\rho,r} b\left(\sqrt{\frac{\log n}{n}} + 
\sqrt{\frac{\log m}{m}}\right).
\end{equation}

The only regularity condition required for these bounds is that $P,Q \in \calK_{r,\rho}(b)$,
for some $\rho > 2r$. When $d=1$, this condition is equivalent
to assuming that $P$ and $Q$ have finite moments of order $\rho$, and is otherwise weaker when $d > 1$. 
The threshold $\rho > 2r$ appears to be nearly sharp, at least for the conclusion 
of equation~\eqref{eq:untrimmed_rate_power_r} to hold. One clearly requires $\rho \geq r$, as otherwise $\sw_r(P,Q)$ may be infinite. 
In the range $r < \rho < 2r$, for the special case $d=1$ and $P=Q$, \cite{fournier2015} argue that the rate in equation~\eqref{eq:untrimmed_rate_power_r} cannot be improved beyond $n^{-\frac{\rho-r}{\rho}}$, which is polynomially slower than the parametric~rate. 
While we do not know the sharp rate in this regime when $P \neq Q$, we expect that the parametric
rate is not achievable even under this restriction, for~$\rho < 2r$. 

Theorem~\ref{thm:untrimmed_rate} is proved using a peeling argument, coupled with a uniform self-normalized
concentration inequality for the empirical quantile process, 
which was already discussed following the statement of Proposition~\ref{prop::ub_empirical}. Unlike the latter result, where this inequality allowed us to obtain rates which adapt to the $\sj_{r,\delta}$
functional, its use here is essential for obtaining a nearly sharp rate without unnecessary moment assumptions, 
as it allows us to tightly control the behaviour of extremal empirical quantiles. We defer the proof to Appendix~\ref{app:untrimmed}.

\section{Finite-Sample Confidence Intervals}
\label{sec::cis}  
\subsection{Finite-Sample Confidence Intervals in  Dimension One}
\label{sec::ci_one_dim}
Throughout this subsection, let $r \geq 1$ and $\delta \in [0,1/2)$ be given, let $P, Q \in \calP(\bbR)$ be probability distributions
with respective CDFs
$F, G$, and let
$X_1, \dots, X_n \sim P$ and $Y_1, \dots, Y_m \sim Q$ be i.i.d. samples. Let 
$F_n(x) = \frac 1 n \sum_{i=1}^n I(X_i \leq x)$ and $G_m(x) = \frac 1 m \sum_{j=1}^m I(Y_j \leq x)$
denote their corresponding empirical CDFs, for all $x \in \bbR$. We derive confidence intervals 
$C_{nm} \subseteq \bbR$ for the $\delta$-trimmed Wasserstein distance, with the following non-asymptotic coverage guarantee
\begin{equation}
\label{eq::one_dim_coverage}
\inf_{P,Q \in \calP(\bbR)} \bbP\big(W_{r,\delta}(P,Q) \in C_{nm}\big) \geq 1-\alpha,
\end{equation}
for some pre-specified level $\alpha \in (0,1)$.
Our approach hinges on the fact that the one-dimensional Wasserstein distance may be expressed as the $L^r$ norm
of the quantile functions of $P$ and $Q$ (cf. equation \eqref{eq::wass-onedim}), suggesting that a confidence interval may be derived via uniform control
of the empirical quantile process. Specifically, the starting point for our confidence intervals is a confidence band of the form
\begin{equation}
\label{eq::uniform_quantile}
\inf_{P\in \calP(\bbR)} \bbP\Big(F_n^{-1}\big(\gamma_{\alpha,n}(u)\big) \leq F^{-1}(u) \leq F_n^{-1}\big(\eta_{\alpha,n}(u)\big), \ \forall u \in (0,1)\Big) \geq 1-\alpha/2,
\end{equation}
for some sequences of functions $\gamma_{\alpha,n},\eta_{\alpha,n}: (0,1) \to \bbR$. 
The study of uniform quantile bounds of the form~\eqref{eq::uniform_quantile} is a classical topic (see for instance, the book of \citet{shorack2009}). We
discuss two prominent examples
that will form the basis of our development.
 
\begin{example}
\label{ex::DKW} 
By the Dvoretzky-Kiefer-Wolfowitz (DKW)
inequality \citep{dvoretzky1956, massart1990}, we have
\begin{equation}
\label{eq::dkw}
\bbP\Big( | F_n(x)-F(x)| \leq \beta_n, \ \forall x \in \bbR \Big) \geq 1-\frac{\alpha}{2}, 
\quad \beta_n = \sqrt{\frac 1 {2n} \log(4/\alpha)}.
\end{equation}
Inverting this inequality leads to the choice 
\begin{equation}
\label{eq::dkw_fns}
\gamma_{\alpha,n}(u) = u - \beta_n,\quad 
\eta_{\alpha,n}(u) = u + \beta_n, \quad u \in (0,1).
\end{equation} 
\end{example}

\begin{example}
\label{ex::rel_VC}
Scale-dependent choices of $\gamma_{\alpha,n}$ and $\eta_{\alpha,n}$ may be obtained
via the relative Vapnik-Chervonenkis (VC) inequality \citep{vapnik2013}. The latter
implies the inequality
\begin{equation}
\label{eq::relVC}
\bbP\left(|F_n(x) - F(x)| \leq \nu_{\alpha,n} \sqrt{F_n(x)(1-F_n(x))} , \ \forall x \in \bbR\right) \geq 1-\frac\alpha 2,
\end{equation}
where $\nu_{\alpha,n}:=\sqrt{\frac {16}{n} \left[ \log(16/\alpha) + \log(2n+1)\right]}$.
As shown in Appendix \ref{subsec::relVC_proof}, inverting inequality \eqref{eq::relVC} leads to the following choice,
for all $u \in (0,1)$,
\begin{equation}
\label{eq::relVC_fns}
\begin{aligned}
\gamma_{\alpha,n}(u) &= \frac{2u + \nu_{\alpha,n}^2-\nu_{\alpha,n} \sqrt{\nu_{\alpha,n}^2 + 4u(1-u)}}{2(1+\nu_{\alpha,n}^2)}, \\
\eta_{\alpha,n}(u) &= \frac{2u + \nu_{\alpha,n}^2 + \nu_{\alpha,n} \sqrt{\nu_{\alpha,n}^2 + 4u(1-u)}}{2(1+\nu_{\alpha,n}^2)}.
\end{aligned}
\end{equation}
\end{example}

Given sequences of functions $\gamma_{\alpha,n},\eta_{\alpha,n}$ satisfying equation \eqref{eq::uniform_quantile}, 
one has with probability at least $1-\alpha$, 
$A_{nm}(u) \leq |F^{-1}(u) - G^{-1}(u) | \leq B_{nm}(u)$ uniformly in $u \in [\delta,1-\delta],$
where,% for all $u \in [\delta,1-\delta]$,
\begin{align*}
A_{nm}(u) &= \big[F_n^{-1}\big(\gamma_{\alpha,n}(u) \big) - G_m^{-1}\big(\eta_{\alpha,m}(u)\big)\big] \vee 
			    \big[G_m^{-1}\big(\gamma_{\alpha,m}(u) \big) - F_n^{-1}\big(\eta_{\alpha,n}(u)\big)\big] \vee 0,\\
B_{nm}(u) &= \big[F_n^{-1}\big(\eta_{\alpha,n}(u)\big) - G_m^{-1}\big(\gamma_{\alpha,m}(u)\big)\big] \vee 
			         \big[G_m^{-1}\big(\eta_{\alpha,m}(u)\big) - F_n^{-1}\big(\gamma_{\alpha,n}(u)\big)\big].
\end{align*}
This observation readily leads to the following Proposition.
\begin{proposition}
\label{prop::one_dim_coverage}
Let $\delta \in [0,1/2)$ and $r \geq 1$. Then, the interval
\begin{equation}
\label{eq::one_dim_CI}
C_{nm} = \left[ \left(\frac 1 {1-2\delta}\int_{\delta}^{1-\delta} A_{nm}^r(u)du \right)^{\frac 1 r} ,
                \left(\frac 1 {1-2\delta}\int_{\delta}^{1-\delta} B_{nm}^r(u)du \right)^{\frac 1 r}\right],
\end{equation}                   
satisfies
$\displaystyle\inf_{P,Q \in \calP(\bbR)} \bbP\big( W_{r,\delta}(P,Q) \in C_{nm} \big) \geq 1-\alpha.$
\end{proposition}
Proposition \ref{prop::one_dim_coverage} establishes the finite-sample coverage of the confidence interval
$C_{nm}$, under no assumptions on the distributions $P, Q$. We emphasize, however, 
that for distributions $P, Q$ with unbounded support, the interval $C_{nm}$ only has finite length under the following condition.
\begin{enumerate}[label=\textbf{A1($\delta; \alpha$)},leftmargin=1.65cm,listparindent=-\leftmargin]   %[label=\textbf{(A1($\delta$)}]      %[\textbf{A1($\delta$)}]
\item \label{assm::delta}%The trimming constant $\delta$ satisfies 
We have  
$\gamma_{\alpha,n\wedge m}(\delta)>0$ and $\eta_{\alpha,n\wedge m}(1-\delta) < 1.$
\end{enumerate}
If $\gamma_{\alpha,n}, \eta_{\alpha,n}$ are chosen via the  DKW inequality  \eqref{eq::dkw_fns}, 
these inequalities imply the choice $\delta \gtrsim (n\wedge m)^{-1/2}$, while
if they are chosen via the relative VC inequality \eqref{eq::relVC_fns},
one must take $\delta \gtrsim \log(n\wedge m) (n\wedge m)^{-1}$. These choices exclude the untrimmed case $\delta=0$, 
for which statistical inference for the Wasserstein distance is not possible without any assumptions on the tail behaviour of $P$ and~$Q$.
If explicit bounds on the quantile functions of $P$ and $Q$ are known near the boundary of the unit interval---which
is for instance the case when an upper bound on the moments of $P$ and $Q$ is known---these may be used
to replace the confidence band $[F_n^{-1}(\gamma_{\alpha,n}(u)), F_n^{-1}(\eta_{\alpha,n}(u))]$ 
by one of finite length for values of $u \in [0,1]$ satisfying $\gamma_{\alpha,n}(u) < 0$
and $\eta_{\epsilon,n}(u) > 1$. Doing so would lead to a confidence interval $C_{nm}$ of finite length
for the untrimmed Wasserstein distance. Since our goal is assumption-free inference, however, we do not
pursue this avenue here and we therefore assume \hyperref[assm::delta]{\textbf{A1($\delta;\alpha$)}} holds throughout the sequel. 

\subsection{Finite-Sample Confidence Intervals in General Dimension}
We now use Proposition \ref{prop::one_dim_coverage} to derive a confidence interval for $\sw_{r,\delta}(P,Q)$, where $P,Q \in \calP(\bbR^d)$. In analogy to Section \ref{sec::ci_one_dim}, 
a natural  approach is to choose functions $\widebar\gamma_{\alpha,n}$
and $\widebar\eta_{\alpha,n}$ such that
\begin{equation}
\label{eq::vc_half_spaces}
\bbP\left(F_{\theta,n}^{-1}\big(\widebar\gamma_{\alpha,n}(u)\big) 
		\leq F_\theta^{-1}(u) \leq F_{\theta,n}^{-1}\big(\widebar \eta_{\alpha,n}(u)\big), \ \forall u \in (0,1), \theta \in \bbS^{d-1}\right) \geq 1-\frac \alpha 2,
\end{equation}
uniformly in $P\in \calP(\bbR^d)$, where $F_{\theta,n}(x)= (1/n)\sum_{i=1}^n I(X_i^\top \theta \leq x)$
for all $x \in \bbR$ and $\theta \in\bbS^{d-1}$, and $F_\theta^{-1}$ denotes the quantile
function of $P_\theta = {\pi_\theta}_\# P$.
Such a bound can be obtained, for instance, by applying the VC inequality \citep{vapnik2013}
to the empirical process indexed by the set of half-spaces in $\bbR^d$.
An assumption-free confidence interval for $\sw_{r,\delta}(P,Q)$ with finite-sample coverage
may then be constructed by following the same lines as in the previous section.
Due to the uniformity of equation \eqref{eq::vc_half_spaces} over the unit sphere, however, it can be seen that the length of such
an interval is necessarily dimension-dependent. In what follows, we instead show that it is possible to obtain a 
confidence interval with dimension-independent length by exploiting the fact that the Sliced Wasserstein distance is a mean with respect
to  $\mu$.

Let $\theta_1, \dots, \theta_N$ be an i.i.d.\ sample from the distribution $\mu$, for
some integer $N \geq 1$, and let $\mu_N = (1/N) \sum_{i=1}^N \delta_{\theta_i}$ denote the corresponding empirical measure.
Consider the following Monte Carlo approximation of the Sliced Wasserstein distance
between the distributions $P$ and $Q$,
$$\sw^{(N)}_{r,\delta}(P,Q) = \left(\int_{\bbS^{d-1}} W_{r,\delta}^r(P_{\theta}, Q_{\theta})d\mu_N(\theta)\right)^{\frac 1 r}
    = \left(\frac 1 N \sum_{j=1}^N W_{r,\delta}^r(P_{\theta_j}, Q_{\theta_j})\right)^{\frac 1 r}.$$
For any $\theta \in \bbS^{d-1}$, let $[\ell_{N,nm}(\theta), u_{N,nm}(\theta)]$ be the 
confidence interval in equation 
\eqref{eq::one_dim_CI} for $W_{r,\delta}(P_{\theta}, Q_{\theta})$, at level $1-\alpha/N$.
Let
$$L_{N,nm} = \int_{\bbS^{d-1}} \ell_{N,nm}^r(\theta) d\mu_N(\theta),\quad
  U_{N,nm} = \int_{\bbS^{d-1}} u_{N,nm}^r(\theta) d\mu_N(\theta),$$
  and set
\begin{equation}
\label{eq::MC_CI}
C_{nm}^{(N)} = \left[ L_{N,nm}^{\frac 1 r}, U_{N,  nm}^{\frac 1 r}\right].
%\left[\left(\int_{\bbS^{d-1}} \ell_{N,nm}^r(\theta) d\mu_N(\theta)\right)^{\frac 1 r}, 
%                 \left(\int_{\bbS^{d-1}} u_{N,nm}^r(\theta )d\mu_N(\theta)\right)^{\frac 1 r} \right].
\end{equation}
By a  Bonferroni correction, we obtain conditional coverage of $\sw^{(N)}_{r,\delta}(P,Q)$, i.e. almost surely,
$$\inf_{P,Q \in \calP(\bbR^d)} \bbP\Big( \sw^{(N)}_{r,\delta}(P,Q) \in C_{nm}^{(N)} ~\Big|~ \theta_1, \ldots, \theta_N \Big) \geq 1-\alpha.$$
We further obtain finite-sample coverage of the Sliced Wasserstein distance itself
by the following small enlargement of $C_{nm}^{(N)}$.
%The following result, ensuring coverage of the (trimmed) Sliced Wasserstein distance, 
%is now a straightforward consequence of the Central Limit Theorem, when the dimension $d$ is fixed.  
\begin{proposition}
\label{prop::sw_coverage}
Let $b,r\geq 1$ and $\delta \in (0,1/2)$.
Let $(M_N)_{N=1}^\infty$ be a nonnegative sequence such that $M_N \to \infty$ as $N \to \infty$.
Define 
\begin{equation}
\label{eq::MC_CI_expanded}
\widebar C_{nm}^{(N)} = \left[\left(L_{N,nm}-M_N/\sqrt N\right)^{\frac 1 r}, \left(U_{N,nm}+M_N/\sqrt N\right)^{\frac 1 r}\right].
\end{equation}
 Then, there is a constant $c > 0$ depending only on $r$ such that
$$\inf_{P,Q \in \calK_{2r}(b)} \bbP\Big(\sw_{r,\delta}(P,Q) \in \widebar C_{nm}^{(N)} \Big) \geq 1-\alpha - \frac{bc}{M_N^2\delta^r}.$$
\end{proposition}
Proposition \ref{prop::sw_coverage} ensures that an enlargement of
the interval  $C_{nm}^{(N)}$, of size less than $(M_N/\sqrt N)^{\frac 1 r}$,
will cover $\sw_{r,\delta}(P,Q)$ at level $1-\alpha-O(M_N^{-2})$, for any fixed sample sizes $n$ and $m$. Notice that 
$N$ is chosen by the practitioner, so that this enlargement can be made to be of lower order than the length of $C_{nm}^{(N)}$.
We shall therefore focus our analysis and numerical studies on the interval $C_{nm}^{(N)}$ rather than $\widebar C_{nm}^{(N)}$. 

Although the coverage of the above intervals requires no assumptions on $P$ and $Q$, apart from the mild
moment condition $P,Q \in \calK_{2r}(b)$, 
we now show that their length achieves the minimax rates established in Theorem \ref{thm::minimax_distance}, 
and is adaptive to the magnitude of 
$\sj_{r,\delta}(P), \sj_{r,\delta}(Q)$.

\subsubsection{Bounds on the Confidence Interval Length} 
In this section, we state a general upper bound (Theorem~\ref{thm::length}) on the length
of $C_{nm}^{(N)}$, 
depending on $\gamma_{\alpha,n}, \eta_{\alpha,n}$. We subsequently specialize this 
result through Corollaries~\ref{cor::length_examples}, \ref{cor::null} 
to illustrate the different rates of convergence which can be obtained under various choices of these functions, and under various 
conditions on the underlying distributions.

In what follows, we assume $\gamma_{\alpha,n}$ and $\eta_{\alpha,n}$ are both differentiable, invertible with differentiable inverses
over $(0,1)$,
and are respectively increasing and decreasing as functions of $\alpha$. 
Given $\epsilon \in (0, 1)$, for notational convenience we write
$\varepsilon := (\epsilon \wedge \alpha)/N$ and $a = \alpha/N$. In the sequel, we also omit explicitly indexing various quantities by the 
number of Monte Carlo samples $N$. 
Our upper bounds on the length of $C_{nm}^{(N)}$
will depend on the function
$$\tilde\kappa_{\varepsilon,n}(u) = \max\Big\{ |f^{-1}(u) - g^{-1}(u)| : f,g \in \{\gamma_{a,n}, \gamma_{\varepsilon,n}, \eta_{a,n},\eta_{\varepsilon,n}\}\Big\}, \quad u \in (0,1),$$
as measured by the following two sequences,
$$\kappa_{\varepsilon,n} {=} \sup_{\frac \delta 2\leq u \leq 1-\frac \delta 2} \tilde\kappa_{\varepsilon,n}(u), ~~ 
  V_{\varepsilon,n}(P) {=} \frac 1 {1-2\delta}\int_{\bbS^{d-1}}\int_{\delta/2}^{1-\delta/2} \left[\frac{\tilde\kappa_{\varepsilon,n}(u)}{p_\theta(F_\theta^{-1}(u))}\right]^r{ dud\mu(\theta)}.$$
Here, recall $p_\theta$ denotes the density of the absolutely continuous component of $P_\theta$.
Additional technical assumptions \ref{assm::B1}-\ref{assm::B4} regarding  $\gamma_{\alpha,n},\eta_{\alpha,n}, \kappa_{\varepsilon,n}$, 
appear in Appendix \ref{app::length}. For appropriate
choices of $\delta$ and $\epsilon$, these assumptions are satisfied by the choices of $\gamma_{\alpha,n}$ and $\eta_{\alpha,m}$ described
in Examples \ref{ex::DKW} and \ref{ex::rel_VC}, for which the corresponding values of $\kappa_{\varepsilon,n}$ and $V_{\varepsilon,n}(P)$ 
are derived in the following simple Lemma.
\begin{lemma} 
\label{lem::kappa_examples}
Let $\varepsilon \in (0,1)$. 
\begin{enumerate}
\item If $\gamma_{\varepsilon,n}$ and $\eta_{\varepsilon,m}$ are chosen as in equation \eqref{eq::dkw_fns}, then 
there exist constants $c_1,c_2 > 0$ depending only on $r$ such that
$$\kappa_{\varepsilon,n} \leq c_1 \sqrt{\frac {\log\left(4/\varepsilon\right)} {n} }, ~~ \text{and}~~  V_{\varepsilon,n}(P) \leq 
		c_2\left(\frac{\kappa_{\varepsilon,n}}{\sqrt \delta}\right)^r \sj_{r,\frac \delta 2}(P).$$
\item If $\gamma_{\varepsilon,n}$ and $\eta_{\varepsilon,m}$ are chosen as in equation \eqref{eq::relVC_fns}, 
then there exist constants $k_1,k_2 > 0$ depending only on $r$ such that
%whenever $\delta \geq \nu_{\epsilon,n}^2$,  
$$\kappa_{\varepsilon,n} \leq k_1\nu_{\varepsilon, n}, \quad \text{and} \quad 
V_{\varepsilon,n}(P) \leq k_2 \nu_{\varepsilon,n}^r \sj_{r,\frac  \delta 2}(P).$$
\end{enumerate}
\end{lemma}
The proof is a straightforward consequence of Examples \ref{ex::DKW}, \ref{ex::rel_VC}, together with the 
derivations in Appendices~\ref{app::pf_prop_empirical_sw}  and \ref{subsec::relVC_proof},  and is therefore omitted.
We now define the functional
$$U_{\varepsilon,n}(P) = \frac 1 {1-2\delta}\int_{\bbS^{d-1}} \bigg(\sup_{\substack{\delta\leq u \leq 1-\delta\\ |h| \leq \kappa_{\varepsilon,n}}}  
		\big|F_\theta^{-1}(u+h) - F_\theta^{-1}(u)\big|^{r-1}\bigg)d\mu(\theta).$$	
$U_{\varepsilon,n}(P)$ is an upper bound on the magnitude of the largest jump discontinuity of the quantile function $F_\theta^{-1}$,
averaged over directions $\theta \in \bbS^{d-1}$. 
When $\sj_{r,\delta}(P) < \infty$, the quantile function $F_\theta^{-1}$
is absolutely continuous for almost all $\theta \in \bbS^{d-1}$ (see Lemma \ref{lem::abs_cont} in Appendix \ref{app::preliminary}), 
implying that  $U_{\varepsilon,n}(P)$
decays to zero as $n \to \infty$. The lengths of our confidence intervals will now  depend on the quantities
\begin{equation*} 
\psi_{\varepsilon,nm}
  {=} \begin{cases}
\big( {\sw_{\infty,\delta}^{(r{-}1)}(P,Q)} {+} U_{\varepsilon,n}(P) {+} U_{\varepsilon,m}(Q)\big) \frac{\sqrt b\kappa_{\varepsilon,n}}{\sqrt \delta}, 
%\text{if }
~ \sj_{r,\frac \delta 2}(P){\vee} \sj_{r,\frac \delta 2}(Q){=}\infty&
\\[0.1in]
 \left( \sw_{r,\delta}^r(P,Q) {+} V_{\varepsilon,n}(P) {+} 
 V_{\varepsilon,m}(Q) \right)^{\frac{r-1}{r}}\left[V_{\varepsilon,n}(P)\right]^{\frac 1 r}, ~\text{otherwise},&
 \end{cases}
 \end{equation*}
and,
 \begin{equation*}
\varphi_{\varepsilon,nm}
 {=} \begin{cases}
 \big( \sw_{\infty,\delta}^{(r-1)}(P,Q) {+} U_{\varepsilon,n}(P) {+} U_{\varepsilon,m}(Q)\big)\frac{\sqrt b\kappa_{\varepsilon,m}}{\sqrt \delta}, 
%  \text{if } 
~ \sj_{r,\frac \delta 2}(P){\vee} \sj_{r,\frac \delta 2}(Q) {=}\infty
 \\[0.1in]
 \left( \sw_{r,\delta}^r(P,Q) {+} V_{\varepsilon,n}(P) {+} 
 V_{\varepsilon,m}(Q) \right)^{\frac{r-1}{r}}\left[V_{\varepsilon,m}(Q)\right]^{\frac 1 r}, ~\text{otherwise}.
\end{cases}
\end{equation*}
With this notation in place, we arrive at the following upper bound on the length of $C_{nm}^{(N)}$. 
Recall that $\lambda$ denotes the Lebesgue measure on $\bbR$.
To simplify our statement, we shall only consider the case where $\delta$ is bounded away from~$1/2$. 
\begin{theorem}%[Length of $C_{nm}^{(N)}$]
\label{thm::length}
Let $r,b \geq 1$ and $\alpha,\epsilon \in (0,1)$. Let $P,Q \in \widebar \calK_2(b)$, and define $\delta \in (0,\delta_0)$
for some $\delta_0\in (0,1/2)$.  
Recall that $\varepsilon = (\epsilon\wedge \alpha)/N$, 
and assume 
$\kappa_{\varepsilon,n\wedge m} \leq \frac \delta 2 \wedge(1-2\delta)$. Assume further that 
conditions \ref{assm::B1}-\ref{assm::B4} hold for some constants $K_1,K_2 > 0$.
Then, 
there exists  $c  > 0$ depending only on $K_1,K_2,\delta_0,r$
such that with probability at least $1-\epsilon$, 
$$\lambda(C_{nm}^{(N)}) \leq \Big\{\sw_{r,\delta}^r(P,Q) + c\big(\psi_{\varepsilon, nm} + \varphi_{\varepsilon,nm} + \varkappa_N\big) \Big\}^{1/r}
 -  \sw_{r,\delta}(P,Q).$$
Here,  $\varkappa_N$ denotes a random variable
depending only on $\mu_N$, such that $\bbE |\varkappa_N| \leq c_1N^{-1/2r} I(d\geq 2)$,
for a constant $c_1 > 0$ depending on $b,r,\delta$ and $K_1$--$K_2$.  
\end{theorem}
The proof of Theorem~\ref{thm::length} appears in Appendix~\ref{app::length}.
As we shall see, the presence of  $\sw_{\infty,\delta}^{(r-1)}(P,Q)$ 
or $\sw_{r,\delta}(P,Q)$ in the definition of $\varphi_{\varepsilon, nm}$ and $\psi_{\varepsilon, nm}$ implies 
distinct rates of decay for the confidence interval length, depending on whether $P,Q$
approach each other under the  Sliced Wasserstein distance. The fact that $\sw_{\infty,\delta}$ is a stronger
metric than $\sw_{r,\delta}$, and the presence of the functional $U_{\varepsilon,n}$,
will imply a second dichotomy in the rate of decay of the confidence interval length, based
on whether or not $\sj_{r,\delta/2}(P) \vee \sj_{r,\delta/2}(Q) <\infty$.%, as we shall explore in the sequel. 

The following result
specializes Theorem \ref{thm::length} to Examples \ref{ex::DKW} and \ref{ex::rel_VC}.
\begin{corollary}
\label{cor::length_examples}
%Assume the same conditions on $P$ and $Q$  as in Theorem \ref{thm::length}.
%Let $\epsilon \in (0,1)$.
Let $r,b \geq 1$ and $\alpha,\epsilon \in (0,1)$. Let $P,Q \in \widebar \calK_2(b)$, and define $\delta \in (0,\delta_0)$
for some $\delta_0\in (0,1/2)$.  
\begin{enumerate}
\item[(i)] Suppose $\gamma_{\alpha,n}, \eta_{\alpha,n}$ are chosen as in Example \ref{ex::DKW}. 
Then, there is a constant $c_\alpha > 0$ such that whenever $\delta \wedge(1-2\delta) \geq \sqrt{c_\alpha\log(N/\epsilon)/(n\wedge m)}$, 
we have with probability at least $1-\epsilon$,
\end{enumerate}
\begin{align*}
\lambda&(C_{nm}^{(N)}) \underset{\alpha,b,\delta_0,r}\lesssim  \varkappa_N^{\frac 1 r} {+} \begin{cases}
\delta^{-\frac 1 2}\log(N/\epsilon)^{\frac 1 {2r}} \Big(n^{-\frac 1 {2r}} + m^{-\frac 1 {2r}}\Big) , ~ \sj_{r,\frac \delta 2}(P){\vee} \sj_{r,\frac \delta 2}(Q) {=} \infty \\[0.08in]
\delta^{-\frac 1 {2}} \log(N/\epsilon)^{\frac 1 {2}} \Big(\frac{\sj_{r,\delta/2}^{\frac 1 r}(P)}{\sqrt n} + \frac{\sj_{r,\delta/2}^{\frac 1 r}(Q)}{\sqrt m}\Big), 
~ \mathrm{otherwise}.
\end{cases}
\end{align*}
\begin{enumerate}
\item[(ii)] Suppose $\gamma_{\alpha,n}, \eta_{\alpha,n}$ are chosen as in Example \ref{ex::rel_VC}.
Let $\beta_{\epsilon,nm} = \log(n\wedge m) + \log(N/\epsilon)$. 
Then, there is a   constant $k_\alpha > 0$ such that
whenever $\delta \wedge (1-2\delta) \geq \sqrt{k_\alpha\beta_{\epsilon,nm}/(n\wedge m)}$, we have with probability at least $1-\epsilon$,
\end{enumerate}
\begin{align*}
\lambda(C_{nm}^{(N)}) \underset{\alpha,b,\delta_0,r}\lesssim 
   \varkappa_N^{\frac 1 r}  {+} \begin{cases}
\delta^{-\frac 1 2}\beta_{\epsilon,nm}^{\frac 1 {2r}} \Big(n^{-\frac 1 {2r}} + m^{-\frac 1 {2r}}\Big), ~~ \sj_{r,\frac \delta 2}(P){\vee}\sj_{r,\frac\delta 2}(Q) {=}\infty  \\[0.08in]
\beta_{\epsilon,nm}^{\frac 1 {2}} \bigg(\frac{\sj_{r,\delta/2}^{\frac 1 r}(P)}{\sqrt n} + \frac{\sj_{r,\delta/2}^{\frac 1 r}(Q)}{\sqrt m}\bigg),
~ \mathrm{otherwise}.\\
\end{cases}
\end{align*}
\end{corollary}
Whenever the trimming sequence is chosen as $\delta \asymp (n\wedge m)^{-a}$ for some $a \in (0,1/2)$, 
notice that one may allow $\epsilon$ to vanish at an exponentially fast rate with respect to $n\wedge m$,
in both cases of Corollary~\ref{cor::length_examples}. 
The high-probability bounds in this result may then be turned into bounds on the expected confidence
interval lengths,
 similarly as in Proposition~\ref{prop::empirical_sw}, though we avoid doing so here for brevity.

Corollary \ref{cor::length_examples}(i) shows that the length of the DKW-based interval 
achieves the minimax lower bound of Theorem \ref{thm::minimax_distance}(i), up to a polylogarithmic factor
in $N$ and the approximation error $\varkappa_N$. It does not, however, achieve the optimal dependence
on $\delta$. 
This is a consequence of the DKW inequality not adapting to the variance of the distributions therein. 
Corollary \ref{cor::length_examples}(ii) instead shows that the  relative VC-based
interval has length depending on $\delta$ solely through the magnitude of the $\sj_{r,\delta/2}$ functional, at the expense of a polylogarithmic
term in $n, m$. In both cases, the confidence interval length scales polynomially with $\delta^{-1}$
when $\sj_{r,\delta/2}(P)\vee \sj_{r,\delta/2}(Q)= \infty$, suggesting that in this case, the practitioner 
should not let $\delta$ vanish with $n\wedge m$ at a rate faster  than logarithmic, to guarantee consistent inference.

When the distributions $P$ and $Q$ are assumed to be bounded away from each other in $\sw_{r,\delta}$, 
Theorem~\ref{thm::minimax_distance}(ii) suggests that the nonparametric rate $n^{-\frac 1 {2r}} + m^{-\frac 1 {2r}}$ in 
Corollary \ref{cor::length_examples} is improvable. This is indeed the case, as shown below.
\begin{corollary} 
\label{cor::null}
Suppose  $P,Q \in \widebar \calK_2(b)$ satisfy $\sw_{r,\delta}(P,Q) \geq \Gamma $ for some constant $\Gamma > 0$. Then, under the assumptions
of Theorem \ref{thm::length}, we have with probability at least $1 - \epsilon$,
$$\lambda(C_{nm}^{(N)}) \underset{\Gamma,  \delta_0,b,r}\lesssim  ~\varkappa_N + \begin{cases}
\delta^{-\frac r 2} \big(\kappa_{\varepsilon,n} + \kappa_{\varepsilon,m}\big), & \sj_{r,\delta/2}(P) \vee \sj_{r,\delta/2}(Q) =\infty \\[0.08in]
V_{\varepsilon,n}^{1 /r}(P) + V_{\varepsilon,m}^{1 /r}(Q), & \mathrm{otherwise}.
\end{cases}$$
\end{corollary}
For example, when $\gamma_{\epsilon,n},\eta_{\epsilon,n}$ are based on the DKW inequality
(Example \ref{ex::DKW}), Corollary \ref{cor::null} implies that the length of $C_{nm}^{(N)}$ achieves the parametric
rate $n^{-\frac 1 2} + m^{-\frac 1 2}$ with high probability
(ignoring factors depending only on $N$ and $\delta$), 
under the mere condition that $P$ and $Q$
are bounded away from each other. Theorem \ref{thm::minimax_distance}(ii) implies that this rate is minimax optimal. 
As before, adaptivity to the magnitudes of 
$\sj_{r,\delta/2}(P), \sj_{r,\delta/2}(Q)$, without further dependence
on $\delta$, is available using the relative VC-based interval in Example \ref{ex::rel_VC}.
Further implications of Theorem \ref{thm::length} are discussed in Appendix \ref{app::mixed_case}.

\section{Asymptotic Confidence Intervals}% and a Hybrid Bootstrap Approach}
\label{sec::asymptotic} 
We now discuss several existing asymptotic confidence intervals
for the one-dimensional Wasserstein distance, their extensions to the Sliced Wasserstein distance, 
and we compare them to our finite-sample confidence intervals in Section \ref{sec::cis}.
%We then discuss how the strengths of these methods can be combined.

In the context
of goodness-of-fit testing, \cite{munk1998} prove central limit theorems
of the form
$
\sqrt{\frac{nm}{n+m}}\big\{W_{2,\delta}^2(P_n, Q_m) - W_{2,\delta}^2(P,Q)\big\} \rightsquigarrow N(0, \sigma^2),
$
where $P$ and $Q$ are one-dimensional distributions, and $\sigma > 0$. They also construct
a consistent estimator of $\sigma$. 
These results assume that $P\neq Q$, and that each of %the distributions
$P$ and $Q$  satisfy the following condition,
\begin{itemize}
\item[(C)] $F$ is twice continuously differentiable, with density $p$, 
which is strictly positive over the real line.
Moreover, 
$$\sup_{x \in \bbR } F(x)(1-F(x)) \left|\frac{p'(x)}{p^2(x)} \right| < \infty.$$
\end{itemize}
Assumption (C) originates from strong approximation theorems for the
empirical quantile process \citep{csorgo1978}, and entails that $P$ and $Q$ have differentiable 
densities, whose supports are intervals. Under the weaker assumption that $P$ and
$Q$ merely admit continuous and positive densities on the real line,
and still retaining the assumption that $P \neq Q$, 
\cite{freitag2003, freitag2005, freitag2007} prove the consistency of the bootstrap
in estimating the distribution of $W_{2,\delta}^2(P_n, Q_m)$ in the one-dimensional case. 

The Wasserstein distance is well-defined between
any pairs of (possibly mutually singular) distributions with sufficient moments, unlike other classical metrics between probability distributions such as the 
Hellinger and $L^r$ metrics. 
Indeed, this feature of the Wasserstein distance 
is a primary motivation for its use in statistical applications. 
Smoothness assumptions such as (C) are therefore prohibitive in inferential problems for
the Wasserstein distance, and motivated our development of assumption-light confidence intervals in the previous section. 
Nevertheless, when a smoothness assumption such as (C) happens to hold, asymptotic confidence intervals 
based on limit laws such as those of \cite{munk1998} above, or those based on the bootstrap, 
may have shorter length than those described in Section \ref{sec::cis}.

We show in Section~\ref{sec::bootstrap} that under some regularity conditions,
the bootstrap is valid in estimating the distribution of $\sw_{r,\delta}(P_n,Q_m)$ for all $r> 1$, thereby generalizing
the results of \cite{freitag2003, freitag2005, freitag2007} from the case $d=1$ and $r=2$. 
We then illustrate in Section~\ref{sec::hybrid},  
how the strengths of bootstrap can be combined with those of the finite-sample confidence intervals of Section~\ref{sec::cis}.
 
\subsection{Bootstrapping the Sliced Wasserstein Distance}
\label{sec::bootstrap}

Let $P,Q \in \calP(\bbR^d)$, and let 
$X_1, \dots, X_n \sim P$, $Y_1, \dots, Y_m \sim Q$ be i.i.d. samples
which are independent of each other. Furthermore, let $P_n$ and $Q_m$ 
denote their corresponding empirical measures, and let $P_n^*$ and $Q_m^*$ 
denote their bootstrap counterparts (that is, $P_n^*$ is the sampling distribution of a sample of size $n$ drawn from $P_n$). 
Lemma \ref{lem::hadamard_sw} in Appendix \ref{app::boot} establishes the Hadamard
differentiability of the Sliced Wasserstein distance at
pairs of distributions $(P,Q)$ satisfying certain regularity conditions. 
Limit laws for the empirical Sliced Wasserstein distance, together with consistency of the
bootstrap, then follow from the functional delta method \citep{vandervaart1996},
as outlined in Theorem~\ref{thm::bootstrap} below. We first introduce some notation.
Denote by $\text{BL}_1$ 
the set of 1-Lipschitz functions $f: \bbR \to \bbR$,
such that $\norm f_\infty \leq 1$. 
We also write for all $u \in [\delta,1-\delta]$ and all $\theta\in \bbS^{d-1}$, 
$$w(u,\theta)  =  \frac r {1-2\delta}  \text{sgn}\left(F_\theta^{-1}(u) - G_\theta^{-1}(u)\right)\big|F_\theta^{-1}(u) - G_\theta^{-1}(u)\big|^{r-1},$$
as well as, 
\begin{equation*}
\begin{multlined}[0.85\linewidth]
\sigma_P^2 =  \int_0^{1-\delta} \left( \int_{\bbS^{d-1}}  \int_{F_\theta^{-1}(\delta\vee t)}^{F_\theta^{-1}(1-\delta)} w(F_\theta(x),\theta)dx d\mu(\theta)\right)^2dt \\ - \left( \int_0^{1-\delta} \int_{\bbS^{d-1}} \int_{F_\theta^{-1}(\delta\vee t)}^{F_\theta^{-1}(1-\delta )} w(F_\theta(x),\theta)dx d\mu(\theta)dt\right)^2,
\end{multlined}
\end{equation*}
and,  
\begin{equation*}
\begin{multlined}[0.85\linewidth]
\sigma_Q^2 =  \int_0^{1-\delta} \left( \int_{\bbS^{d-1}}  \int_{G_\theta^{-1}(\delta\vee t)}^{G_\theta^{-1}(1-\delta)} w(G_\theta(x),\theta)dx d\mu(\theta)\right)^2dt \\ - \left( \int_0^{1-\delta} \int_{\bbS^{d-1}}  \int_{G_\theta^{-1}(\delta\vee t)}^{G_\theta^{-1}(1-\delta)} w(G_\theta(x),\theta)dx d\mu(\theta)dt\right)^2,
\end{multlined}
\end{equation*}
Finally, for any absolutely continuous distribution $P \in \calP(\bbR)$ with density $p$, we shall make use of the 
following trimmed version of the $J_\infty$ functional introduced by~\cite{bobkov2019},
$$J_{\infty,\delta}(P) % = \sup_{r \geq 1} J_{r,\delta}^{1/r}(P) 
= \esssup_{\delta \leq u \leq 1-\delta} \frac 1 {p(F^{-1}(u))}.$$
With this notation in place, the main result of this section is stated as follows. 
\begin{theorem}
\label{thm::bootstrap}
%Let $b,r\geq 1$ and $\delta \in (0,1/2)$ be fixed, and  $P,Q \in \widebar\calK_2(b)$. 
Let $\delta\in [0,1/2)$ and $r > 1$. Let $P,Q \in \widebar\calK_2$ be distributions such that for all $\theta \in \bbS^{d-1}$, $P_\theta, Q_\theta$ are absolutely continuous with respect to the Lebesgue measure, 
with respective families of densities $\{p_\theta\}_{\theta\in \bbS^{d-1}},\{q_\theta\}_{\theta\in \bbS^{d-1}}$ 
which are uniformly integrable 
over $\bbR$.
Assume that, 
\begin{equation}
\label{eq:assm_Jinfty_clt}
\sup_{\theta \in \bbS^{d-1}} J_{\infty,\delta/2}(P_\theta) \vee J_{\infty,\delta/2}(Q_\theta) < \infty.
\end{equation}
Then, the following statements hold as $n,m \to \infty$ such that
$\frac{n}{n+m} \to a \in (0,1)$.
\begin{enumerate}
\item[(i)] (Central Limit Theorem) We have, 
$$\sqrt {\frac{nm}{n+m}} \Big( \sw_{r,\delta}^r(P_n, Q_m) - \sw_{r,\delta}^r(P, Q)\Big) \rightsquigarrow N\big(0, a \sigma_P^2 + (1-a)\sigma^2_Q\big).$$ 
\item[(ii)] (Bootstrap Consistency) For $\bX=(X_1, \dots, X_n), \bY=(Y_1, \dots, Y_m)$, we~have 
\begin{equation*}
\begin{multlined}[\linewidth]
\sup_{h\in \mathrm{BL}_1} \Bigg| \bbE\left[h\left(\sqrt{\frac{nm}{n+m}} \left\{ \sw_{r,\delta}^r(P_n^*,Q_m^*) - \sw_{r,\delta}^r(P_n,Q_m)\right\}\right) \bigg| \bX,\bY\right] \\
 - 
	\bbE\left[h\left(\sqrt{\frac{nm}{n+m}} \left\{ \sw_{r,\delta}^r(P_n,Q_m) - \sw_{r,\delta}^r(P,Q)\right\}\right) \right]\Bigg|  \rightarrow  0,
	\end{multlined}
	\end{equation*}
in outer probability.
\end{enumerate}
\end{theorem}

Theorem~\ref{thm::bootstrap}(i) provides a central limit theorem for the empirical trimmed Sliced Wasserstein distance, centered
at its population counterpart. The primary assumptions required for this result are (a) the existence and uniform
integrability of the densities of $P_\theta$ and $Q_\theta$ along 
directions $\theta \in \bbS^{d-1}$, and (b) a uniform lower bound on these 
densities, over the compact sets $[F_\theta^{-1}(\delta/2), F_\theta^{-1}(1-\delta/2)]$, as measured
by the $J_{\infty,\delta/2}$ functional. Note that assumption (a) holds
if $P, Q$ admit upper bounded densities with respect to the Lebesgue measure over $\bbR^d$, but is strictly weaker; indeed, it can
be satisfied by non-atomic measures which are singular with respect to the Lebesgue measure on $\bbR^d$. 
Furthermore, we note that the assumption of uniform integrability
is vacuous in the special case $d=1$. Assumption (b) requires the bulk of the supports of $P_\theta$ and $Q_\theta$ 
to be connected. Such a condition is necessary for the limit in Theorem~\ref{thm::bootstrap}(i) to 
be a mean-zero Gaussian distribution,  as can be anticipated, for instance, from the lack of Hadamard differentiability
of the Wasserstein distance over finite spaces \citep{sommerfeld2018a}.
Nevertheless, our condition is stronger, 
since we have assumed $J_{\infty,\delta/2}(P_\theta)$ and 
$J_{\infty,\delta/2}(Q_\theta)$ are uniformly bounded in $\theta$. 
Inspired by our results on estimation
and finite sample inference for the Sliced Wasserstein distance, it is natural to ask
whether this condition can be replaced by, say, $\sj_{r,\delta}(P), \sj_{r,\delta}(Q) < \infty$.
Such a condition would
allow the densities $p_\theta$ and $q_\theta$ to approach zero at a sufficiently slow rate, which is currently precluded by 
our theorem. We leave this question open for future work. 

In the special case $d=1$ and $r=2$, the limiting variance obtained in Theorem~\ref{thm::bootstrap}(i) is equal
to the one obtained by~\cite{munk1998}, up to renormalizing their definition of the trimmed
Wasserstein distance, though our assumptions are significantly weaker since we do not require
the aforementioned condition~(C).
Nevertheless, their result allows the trimming constant $\delta$ to vanish, while we require  $\delta$ to be held fixed
and, in fact, positive, when $P$ and $Q$ have unbounded support.
In this regard, our result is closer to those of~\cite{freitag2003,freitag2005,freitag2007}, 
who prove, in particular, the Hadamard differentiability of the functional $W_{2,\delta}^2$ for a nonvanishing trimming constant $\delta$. 
Their results require $P$ and $Q$ to admit positive and continuously differentiable densities over the real line, which
is a strictly stronger assumption than those of Theorem~\ref{thm::bootstrap}. In particular, we require no smoothness conditions on the
various densities.

We next compare Theorem~\ref{thm::bootstrap}(i) to existing central limit theorems  
for untrimmed Wasserstein distances. 
Let $d=1$, and assume $P$ and $Q$ are compactly-supported, so that
one may take $\delta=0$ in Theorem~\ref{thm::bootstrap}(i).
In this case, the limiting variance  may be reformulated in terms of the expressions
$$\sigma_P^2 = \Var[\phi_0(X)], \quad \sigma_Q^2 = \Var[\psi_0(Y)],$$
where $X \sim P\in \calP(\bbR)$, $Y \sim Q\in \calP(\bbR)$, and for all $x,y \in \bbR$, 
$$ \phi_0(x) =  \int_{-\infty}^x w(F(t))dt,\quad \psi_0(y) =  \int_{-\infty}^y w(G(t))dt.$$
% - \left( \int_\delta^{1-\delta} \int_{\bbS^{d-1}}  \int_{F_\theta^{-1}(\delta)}^{F^{-1}(t)} w(F_\theta(x),\theta)dx d\mu(\theta)dt%\right)^2, $$
Here, we abbreviate $w(\cdot) \equiv w(\cdot,\theta)$ in the one-dimensional case. 
%Since $G^{-1}\circ F$ is the optimal transport map from $P$ to $Q$, 
It can be deduced from~\cite{gangbo1996} that $(\phi_0,\psi_0)$ forms an optimal pair of Kantorovich potentials 
in the dual $|\cdot|^r$-optimal transport problem~\citep{villani2003} from $P$ to $Q$. In particular, 
for $r = 2$, the limiting variance in Theorem~\ref{thm::bootstrap}(i)
reduces to the one obtained by~\cite{delbarrio2019a}, who derive central limit theorems for $W_2^2(P_n, Q_m)$, in general dimension $d \geq 1$. 
Their results are not centered at $W_2^2(P,Q)$ due to the large bias of empirical Wasserstein 
distances in general dimension, however,  it was shown by~\cite{delbarrio2019b} that when $d=1$, 
these limit theorems may be centered at the population Wasserstein distance under assumptions akin to condition (C). 
We also refer to~\cite{berthet2020a} and the recent work of \cite{hundrieser2022} 
for distinct assumptions under which such a result can be obtained.  

Theorem~\ref{thm::bootstrap}(ii) proves the consistency of the bootstrap in estimating the distribution
of $\sw_{r,\delta}^r(P_n, Q_m)$. 
Letting $F_{nm}^*$ denote  the CDF of $\sw_{r,\delta}^r(P_n^*, Q_m^*)-\sw_{r,\delta}^r(P_n, Q_m)$, it follows that 
an asymptotic $(1-\alpha)$-confidence interval for $\sw_{r,\delta}(P,Q)$ is given~by
\begin{align*}
%\label{eq::bootstrap_CI}
C^*_{nm}  =  
 & \left[\left(\sw_{r,\delta}^r(P_n,Q_m) - F_{nm}^*(1- \alpha /2)\right)^{\frac 1 r},        \left(\sw_{r,\delta}^r(P_n,Q_m) + F_{nm}^*( \alpha /2)\right)^{\frac 1 r}\right].
\end{align*}        
The CDF $F^*_{nm}$ is typically estimated via Monte Carlo simulation \citep{efron1994}.
The assumptions for the validity of $C^*_{nm}$ are those of Theorem~\ref{thm::bootstrap}, and in addition, 
the condition that $\sw_{r,\delta}(P,Q) > 0$,
which is necessary and sufficient for the limiting variance
$a \sigma_P^2 + (1-a)\sigma_Q^2$ in Theorem~\ref{thm::bootstrap}
to be positive. 
Failure of the bootstrap at the null $\sw_{r,\delta}(P,Q)=0$ is due to the Sliced Wasserstein distance being a functional
with first-order degeneracy~\citep{munk1998}, for which corrections such as those of \cite{chen2019},
or the  $m$-out-of-$n$ bootstrap \citep{sommerfeld2018a}, yield consistent procedures, but are practically less attractive as they introduce further tuning parameters. We illustrate in the sequel how our finite sample confidence intervals
can be combined with the bootstrap to relax this assumption. 
 
\subsection{A Hybrid Bootstrap Approach}
\label{sec::hybrid}
  
Let $C_{nm}^*$ denote the 
preceding bootstrap confidence interval at level $1-\alpha/2$, and let $C_{nm}^{\dagger}$ %and $C_{nm}^{\dagger}(\alpha/2)$ 
denote the assumption-light confidence interval for $\sw_{r,\delta}(P,Q)$ 
in equation \eqref{eq::MC_CI_expanded} at level $1-\alpha/2$. Assume that
the number of Monte Carlo replications~$N$ 
therein is taken to diverge as $n,m\to\infty$.
%Let  $(\epsilon_n)_{n=1}^\infty$ be any sequence of decreasing nonnegative real numbers. 
We define the $(1-\alpha)$-hybrid confidence interval as: 
\begin{align}
\label{eq::ci_hybrid}
C_{nm} = \begin{cases}
C_{nm}^{\dagger},~~~& \text{if } 0 \in C^\dagger_{nm} , \\
C_{nm}^*, & \text{otherwise}.
\end{cases}
\end{align} 
Roughly, we use the bootstrap interval if we are reasonably certain that $P$ and $Q$ are bounded away from
each other in Sliced Wasserstein distance,
and fall back on the  finite-sample interval otherwise. The following simple 
result characterizes the asymptotic coverage and length 
of the hybrid interval. In order to simplify our discussion, 
we write $C_{nm} = [a_{nm}^{1/r}, b_{nm}^{1/r}]$ and we focus on bounding
the length of the confidence interval $[a_{nm}, b_{nm}]$ for the $r$-th power of the $r$-Sliced Wasserstein
distance. We also assume that the finite-sample interval $C_{nm}^\dagger$ is defined in terms
of the DKW confidence band in Example~\ref{ex::DKW}.
\begin{proposition}
\label{prop::pretest}
Let $a,\alpha\in (0,1)$, $\delta \in (0,1/2)$, and assume the same conditions as Theorem~\ref{thm::bootstrap}. Then, 
the following holds assuming $\frac{n}{n+m} \to a$ when~$n,m\to~\infty$.
\begin{enumerate}
\item[(i)] (Coverage) We have, 
\begin{equation}
\label{eq::hybrid_coverage} 
\liminf_{n,m \to \infty} \bbP\big(\sw_{r,\delta}(P,Q) \in C_{nm} \big) \geq 1-\alpha.
\end{equation}
\item[(ii)] (Length) Let $N\asymp n^{r^2}$, and choose $M_N \asymp \log N$
in the definition of $C_{nm}^\dagger$.
Then,  we have with probability at least $1-\alpha$, 
$$(b_{nm} - a_{nm}) = O\left( \left(\frac{\log n}{n}\right)^{\frac r 2} + \frac{\sw_{r,\delta}(P,Q)}{\sqrt n}\right).$$ 
%The implicit constants in the above display depend, in particular, on $P$ and $Q$.
\end{enumerate}
\end{proposition}
Proposition~\ref{prop::pretest} establishes the asymptotic coverage of $C_{nm}$
under the same conditions as Theorem~\ref{thm::bootstrap}. In particular, it removes
the assumption $\sw_{r,\delta}(P,Q) > 0$, which is needed for the asymptotic coverage
of the bootstrap interval $C_{nm}^*$. We note that many other existing corrections 
of the bootstrap for functionals with first-order degeneracy, such as the $m$-out-of-$n$ bootstrap or the procedures
outlined in Section 2.1 of~\cite{verdinelli2021}, involve expanding the asymptotic length of the interval, leading to a loss of efficiency. 
In contrast, Proposition~\ref{prop::pretest} shows that with high (albeit, fixed) probability,
the hybrid interval achieves the rate-optimal asymptotic length both
at the null $\sw_{r,\delta}(P,Q)=0$ and away from the null, up to a polylogarithmic factor in $N$ (which can be removed when $d=1$). 
We emphasize that this adapativity is obtained without  tuning parameters, apart from 
the sequences $M_N,N$ whose precise choice does not greatly alter the properties of $C_{nm}$. Once again, 
the choice of these sequences is vacuous in the special case $d=1$. 

Though this methodology  inherits  benefits from both the bootstrap and finite-sample
confidence intervals, it is not assumption-free. In principle, it is possible to extend 
this procedure by 
empirically testing whether the conditions of Theorem~\ref{thm::bootstrap} are met, and to use the outcome of such
a test in the conditions of equation~\eqref{eq::ci_hybrid}. While doing so may allow for certain assumptions
to be relaxed, it could become impractical: for instance, we do not know of a test for the finiteness of the $J_{\infty,\delta/2}$ functional
which is free of tuning parameters.

\section{Simulation Study}
\label{sec::simulations}
We perform a simulation study to illustrate the coverage and length of the confidence intervals described in Sections
\ref{sec::cis} and \ref{sec::asymptotic}. 
All simulations were performed in Python 3.5 on a typical Linux machine
with twelve cores.
 Implementations for all
 confidence intervals described in this paper, along with code for reproducing the following simulations, are publicly
available\footnote{\url{https://github.com/tmanole/SW-inference}.}.

\paragraph{Comparison of Asymptotic and Finite-Sample Confidence Intervals.}
We compare 
the following confidence intervals: (i) The finite-sample interval in equation \eqref{eq::MC_CI} (or 
\eqref{eq::one_dim_CI} in the one-dimensional case), based on the DKW inequality from Example \ref{ex::DKW}, 
(ii)  The standard bootstrap confidence interval $C_{nm}^*$ in Section \ref{sec::bootstrap}, and
(iii) The hybrid interval in equation \eqref{eq::ci_hybrid}.
We also implemented the finite-sample interval \eqref{eq::MC_CI} with respect to the relative VC
inequality (Example \ref{ex::rel_VC}), however we rarely noticed an improvement over the 
DKW finite-sample interval in practice. This is likely due to the sub-optimal constants in the relative VC inequality
\eqref{eq::relVC}, unlike those in the DKW inequality of \cite{massart1990}, and consequently we do not consider this method in the present simulation study.

\begin{figure}[t]
\centering

\begin{subfigure}{0.33\textwidth}
  \centering
  \includegraphics[width=\linewidth]{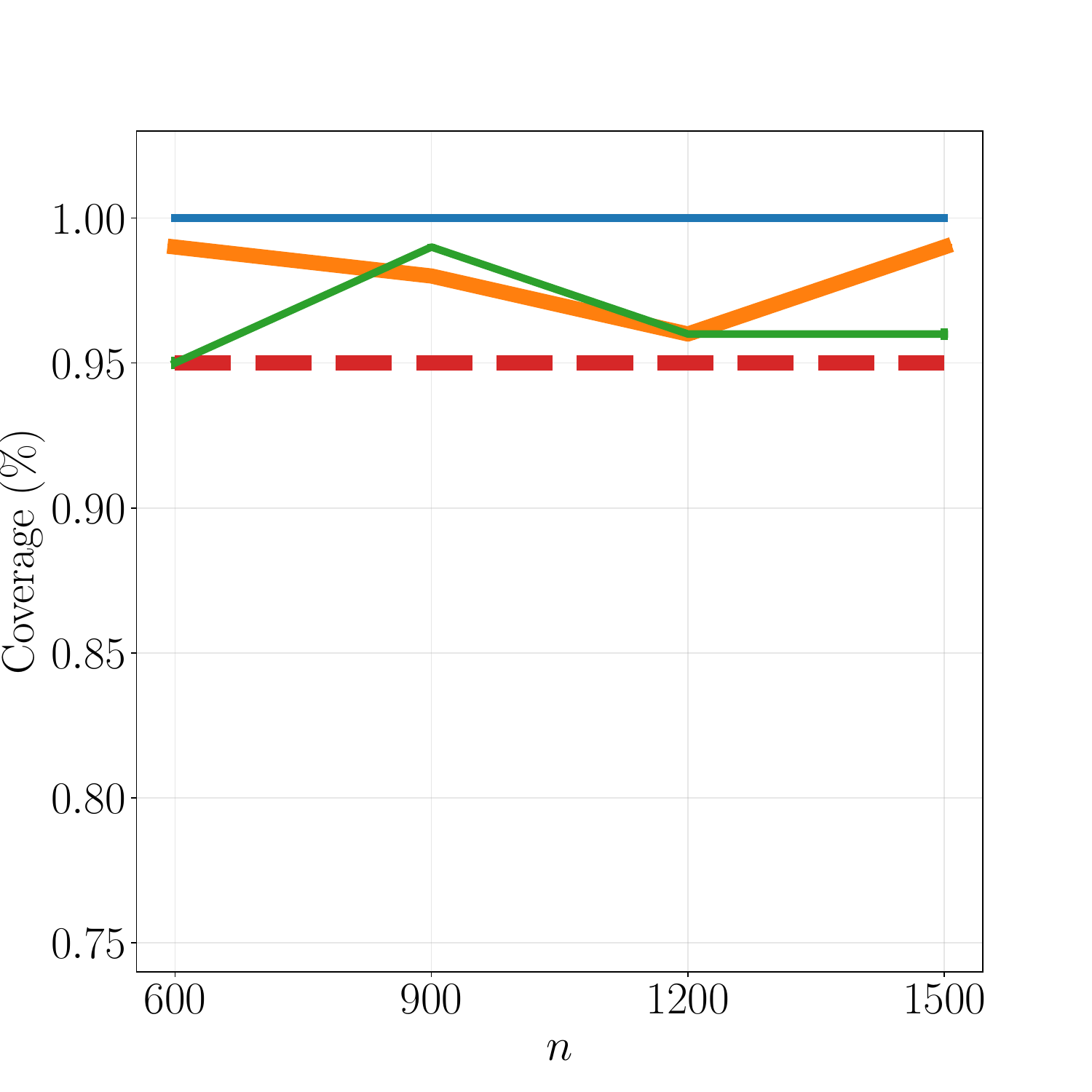}
  \caption{Model 1, Coverage.}
\end{subfigure}%
\begin{subfigure}{0.33\textwidth}
  \centering
  \includegraphics[width=\linewidth]{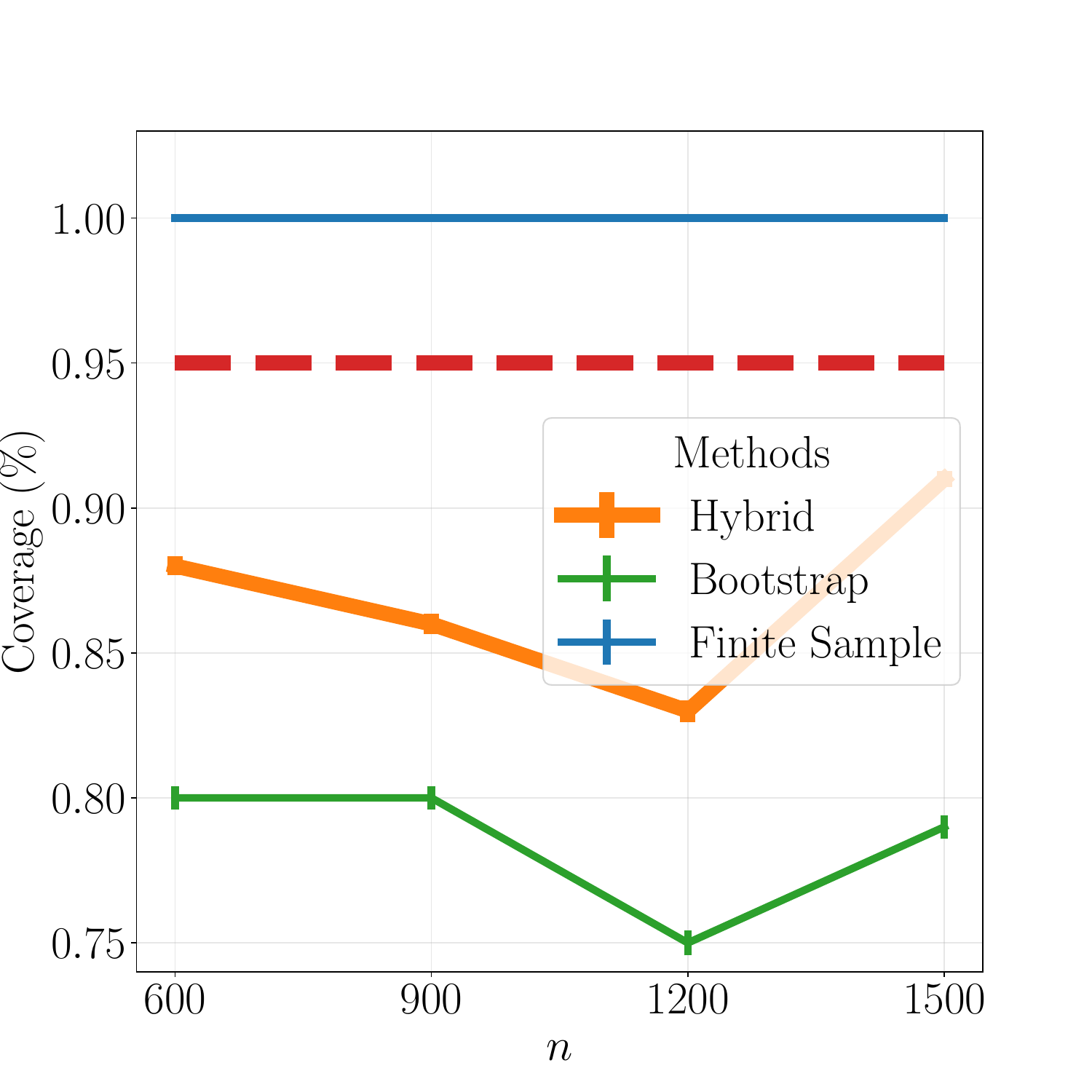}
  \caption{Model 2, Coverage.}
\end{subfigure}%
\begin{subfigure}{0.33\textwidth}
  \centering
  \includegraphics[width=\linewidth]{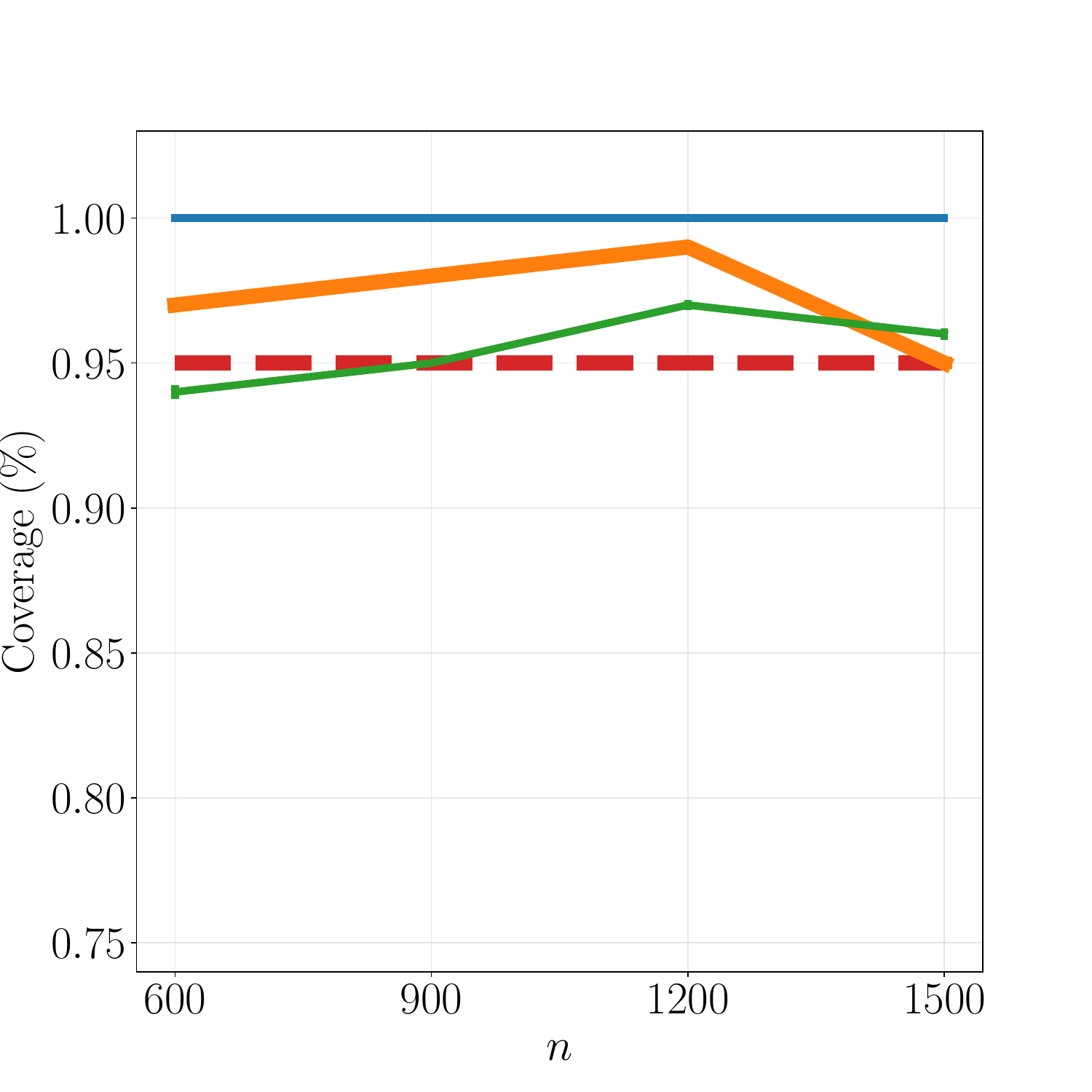}
  \caption{Model 3, Coverage.}
\end{subfigure}  \par\medskip

\begin{subfigure}{0.33\textwidth}
  \centering
  \includegraphics[width=\linewidth]{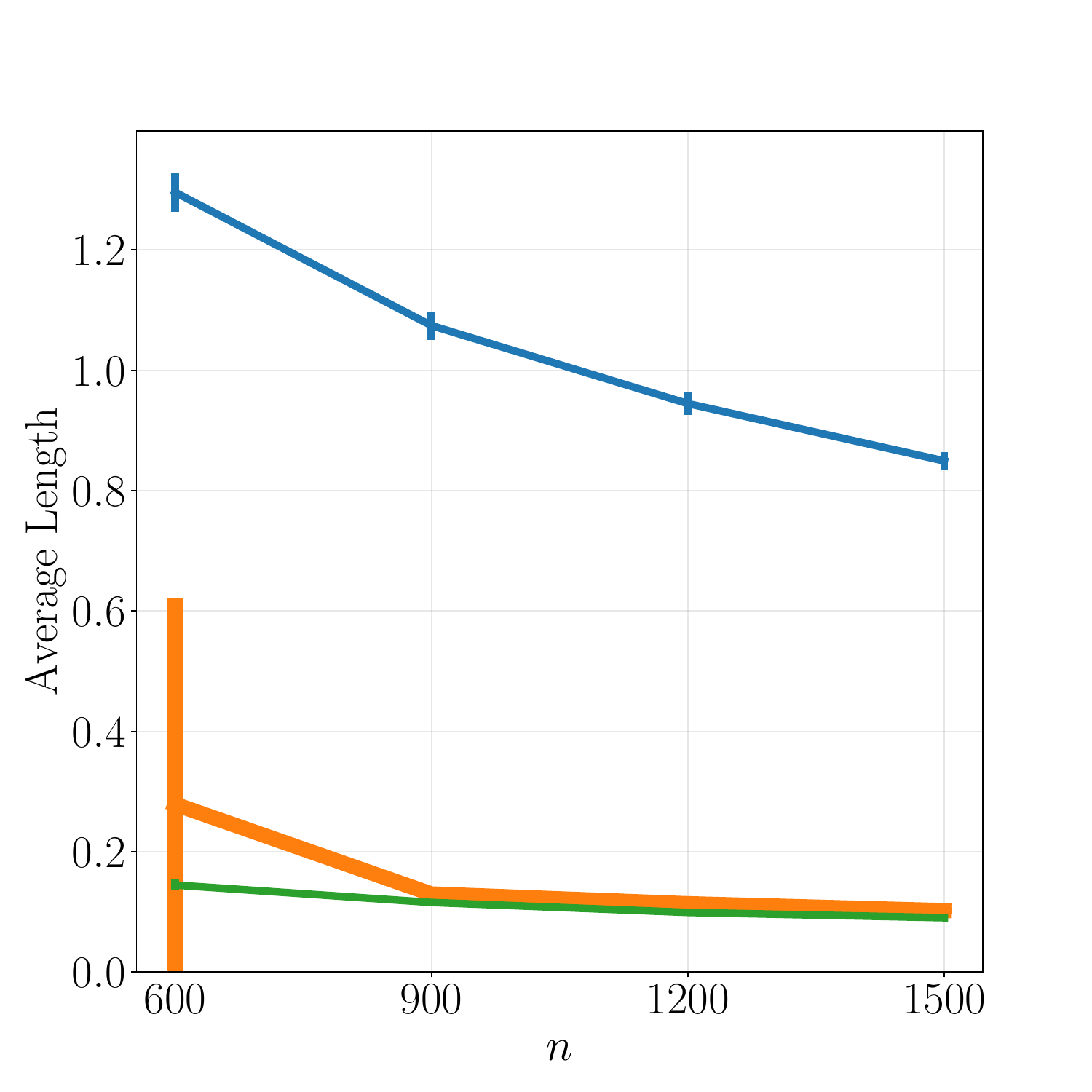}
  \caption{Model 1, Avg. Length.}
\end{subfigure}%
\begin{subfigure}{0.33\textwidth}
  \centering
  \includegraphics[width=\linewidth]{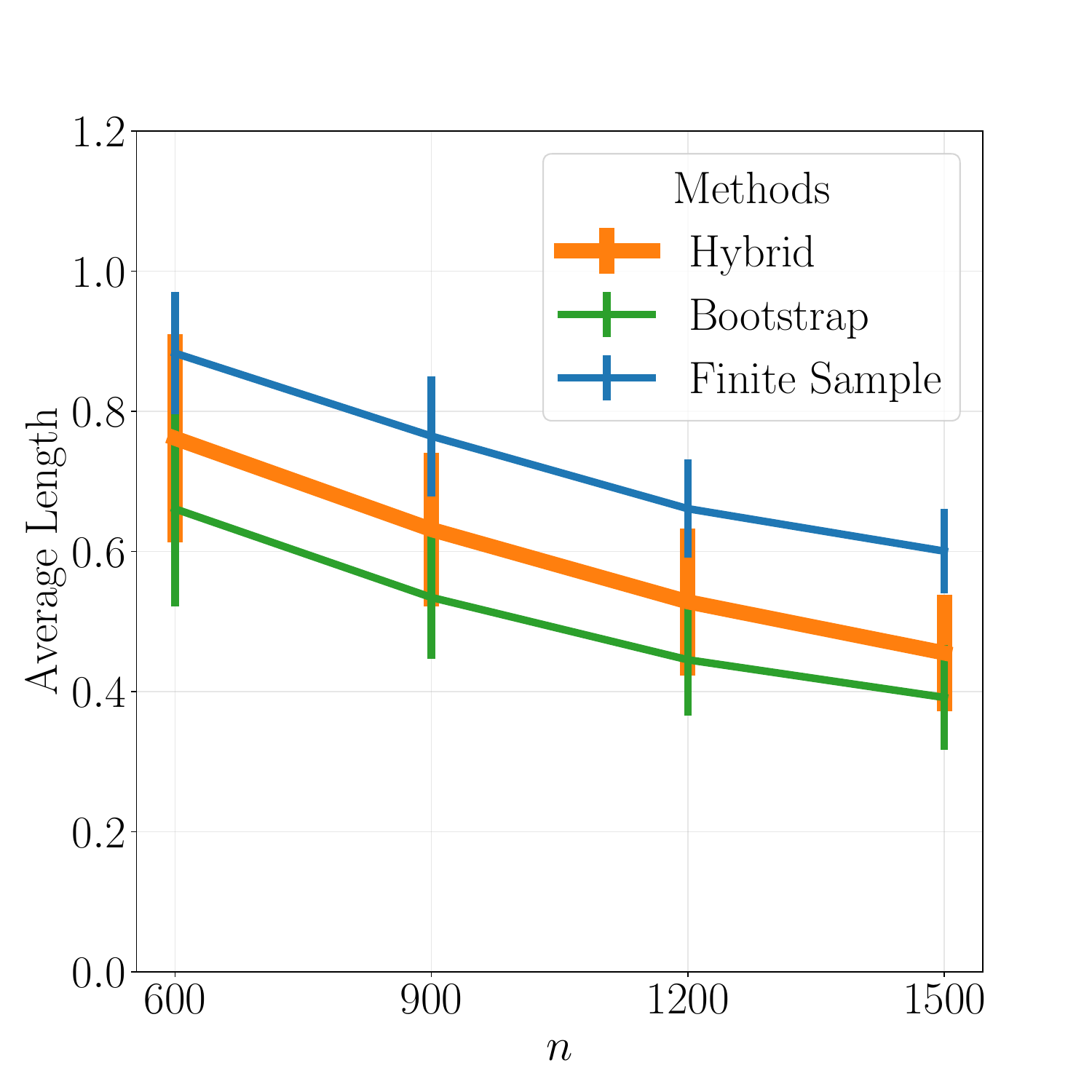}
  \caption{Model 2, Avg. Length.}
\end{subfigure}%
\begin{subfigure}{0.33\textwidth}
  \centering
  \includegraphics[width=\linewidth]{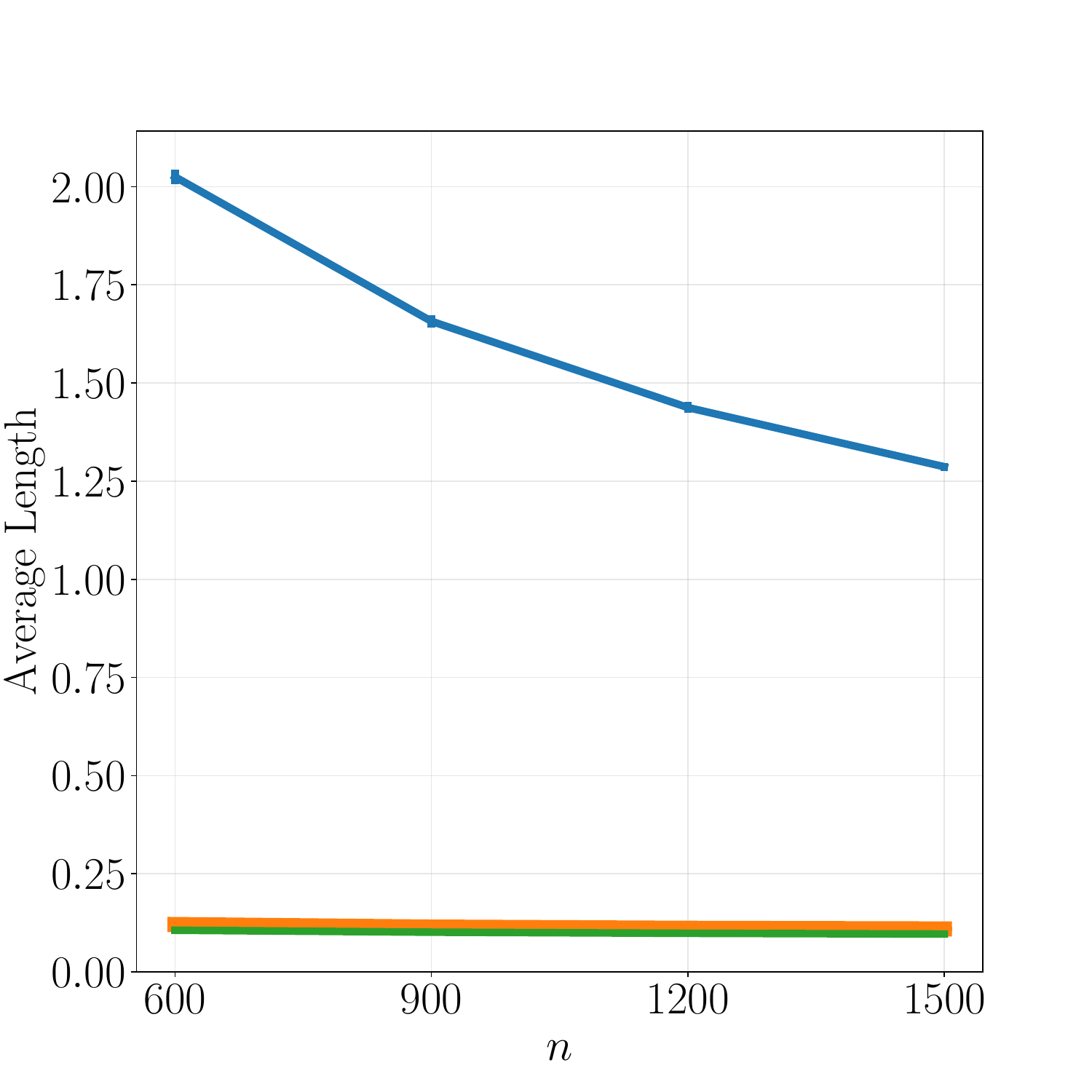}
  \caption{Model 3, Avg. Length.}
\end{subfigure}
\caption{\label{fig::mods1_3}  
Average  interval length and coverage for Models 1-3. 
Error bars represent one standard deviation.
Dotted line ({\bf \color{red}- - -}) indicate the nominal level $0.95$.
}
\end{figure}

\begin{figure}[t]
\centering

\begin{subfigure}{0.33\textwidth}
  \centering
  \includegraphics[width=\linewidth]{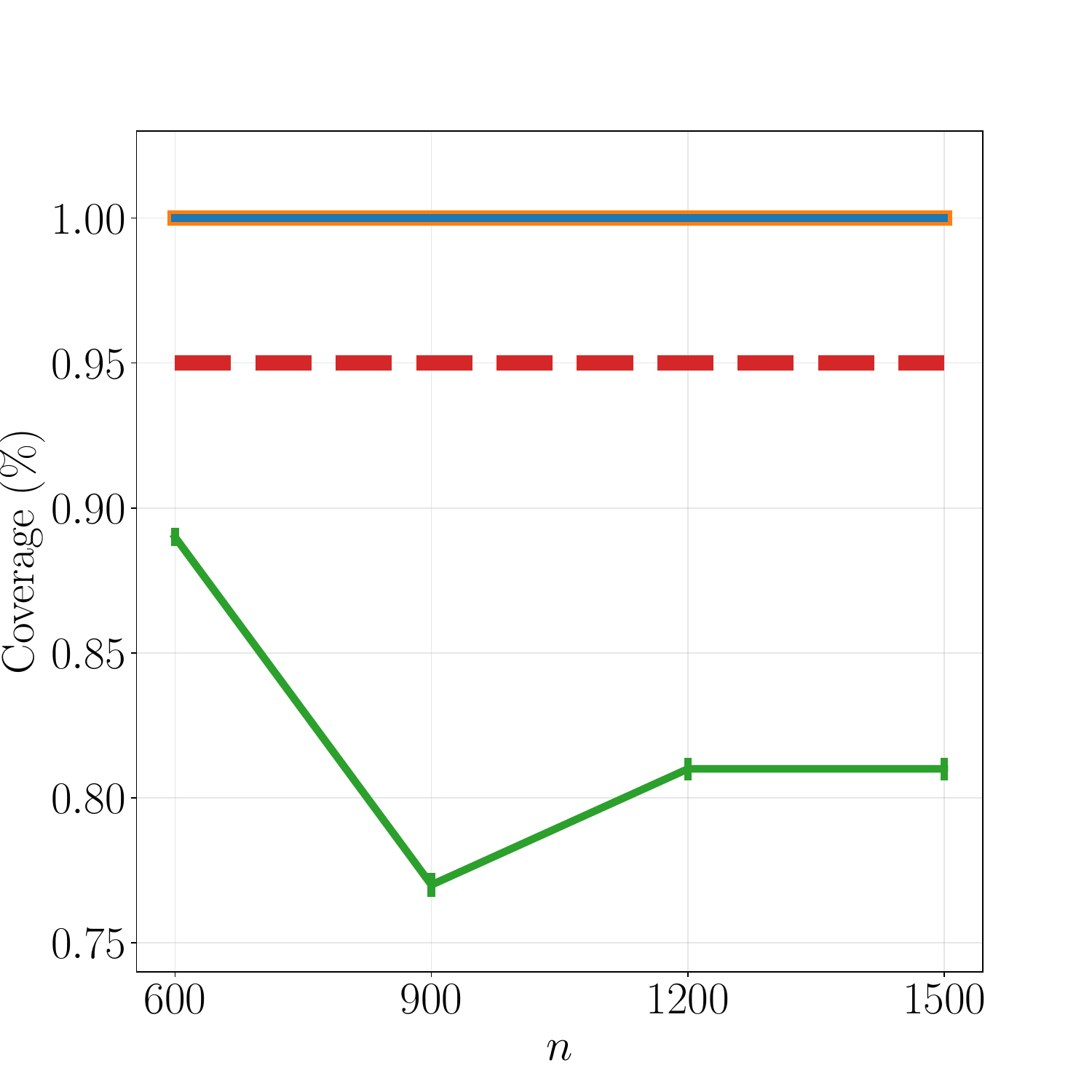}
  \caption{Model 4, Coverage.}
\end{subfigure}%
\begin{subfigure}{0.33\textwidth}
  \centering
  \includegraphics[width=\linewidth]{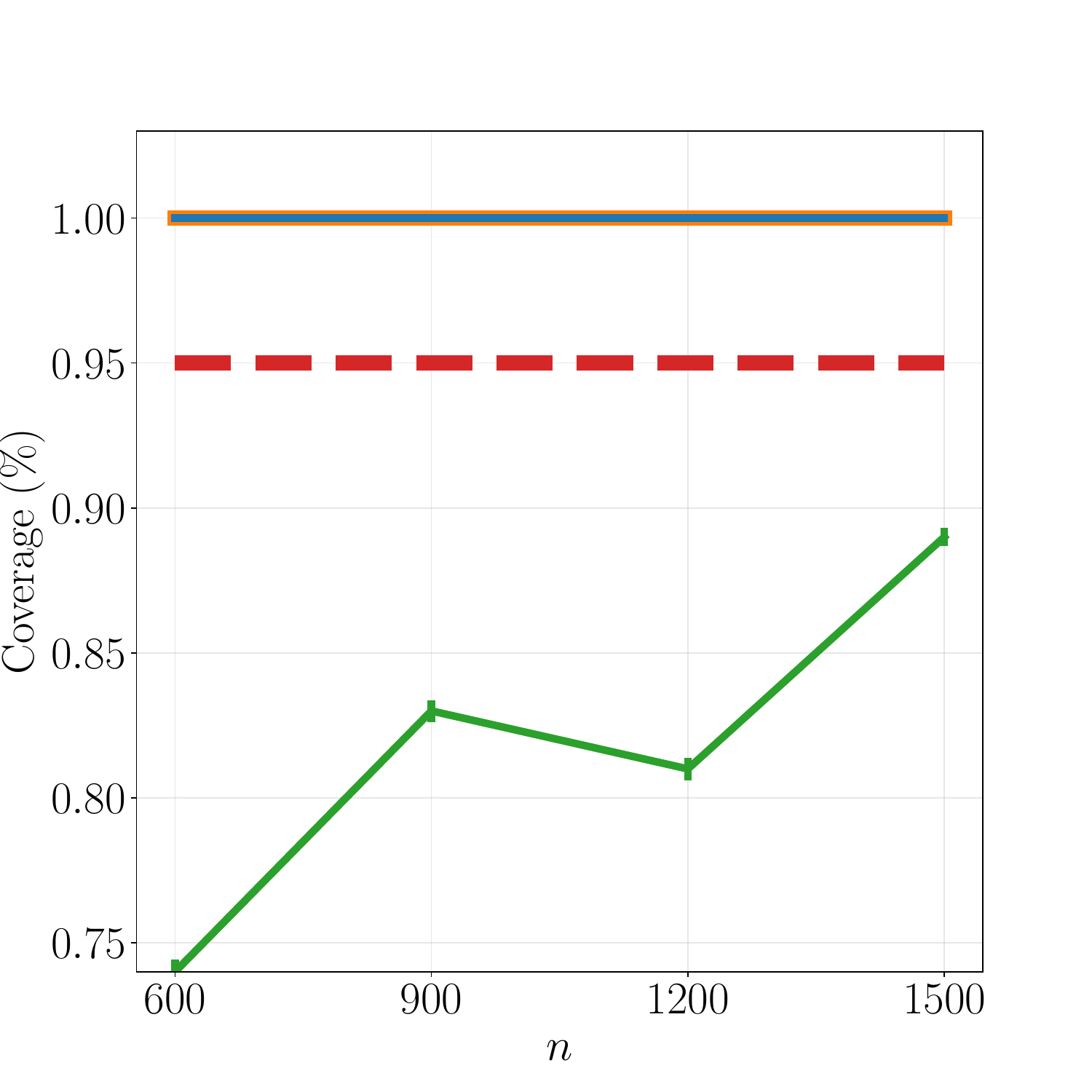}
  \caption{Model 5, Coverage.}
\end{subfigure}%
\begin{subfigure}{0.33\textwidth}
  \centering
  \includegraphics[width=\linewidth]{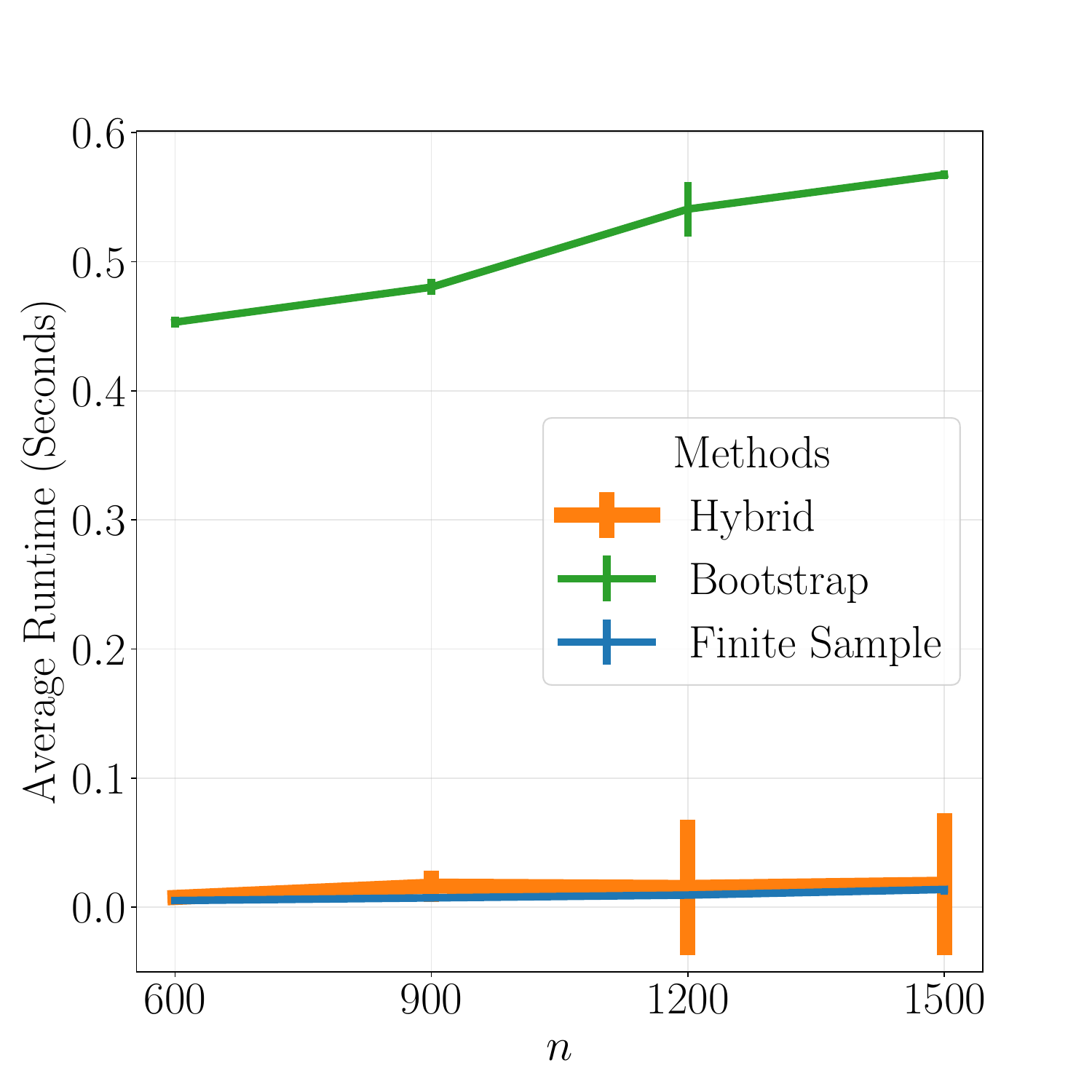}
  \caption{Model 4, Avg. Time.}
\end{subfigure}  \par\medskip

\begin{subfigure}{0.33\textwidth}
  \centering
  \includegraphics[width=\linewidth]{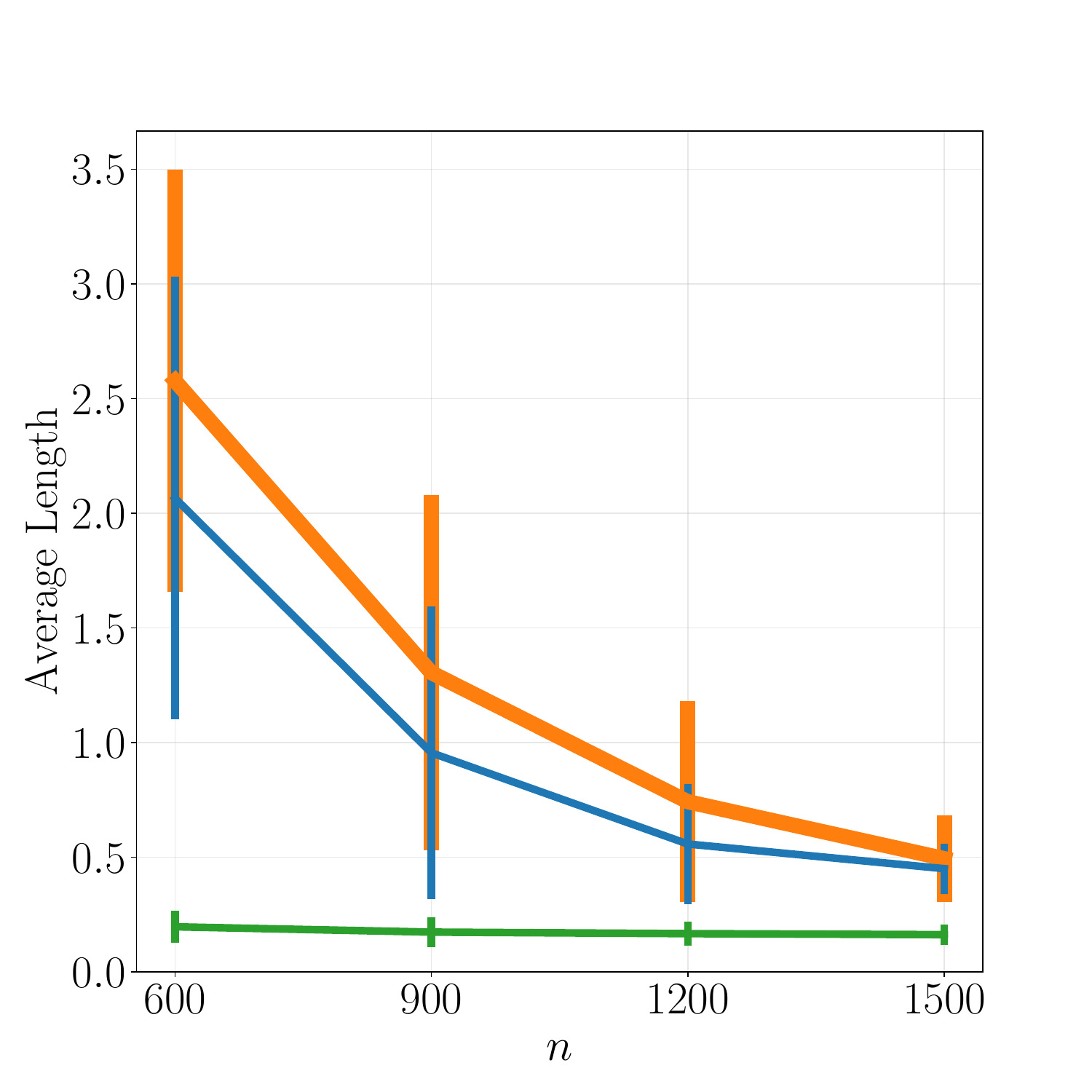}
  \caption{Model 4, Avg. Length.}
\end{subfigure}%
\begin{subfigure}{0.33\textwidth}
  \centering
  \includegraphics[width=\linewidth]{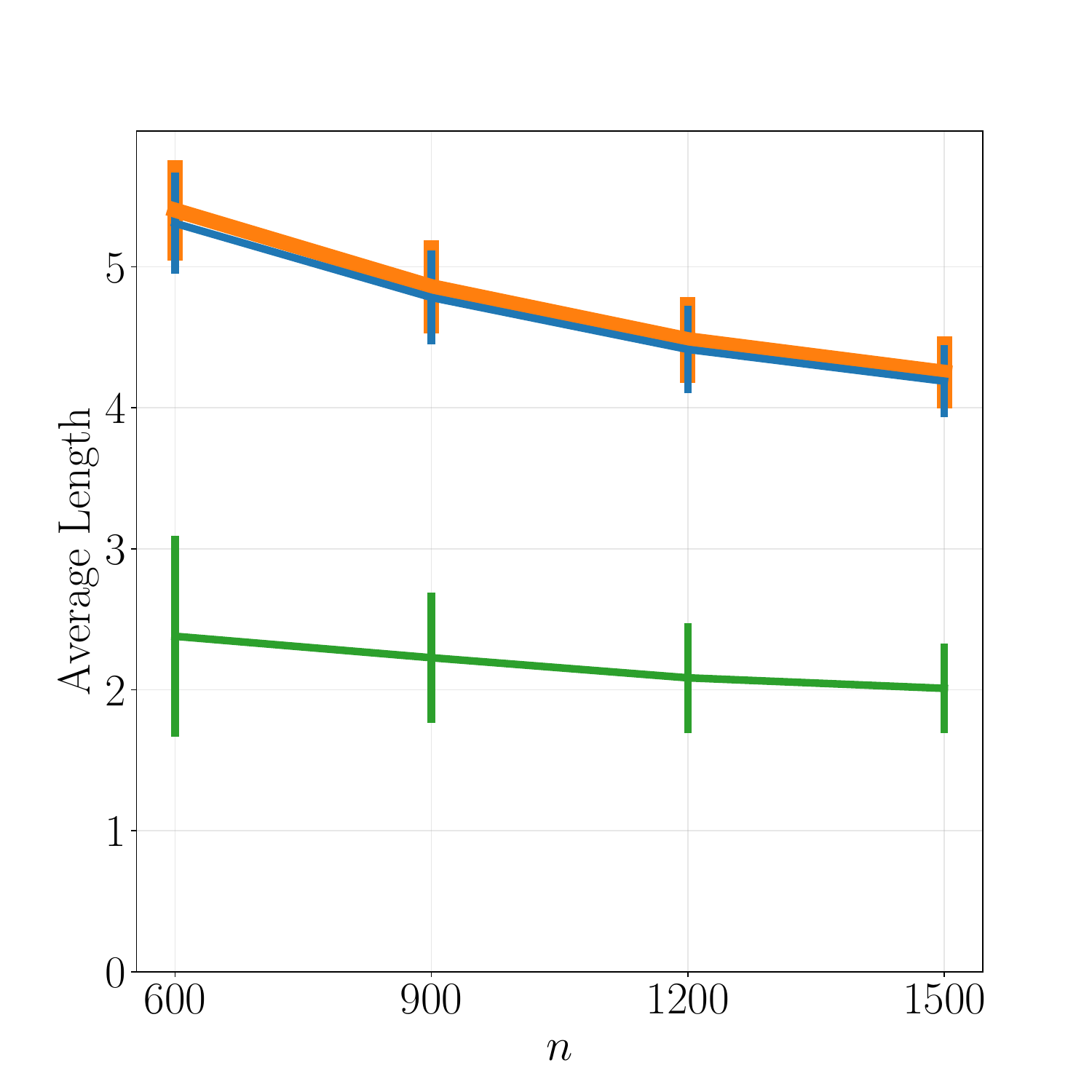}
  \caption{Model 5, Avg. Length.}
\end{subfigure}%
\begin{subfigure}{0.33\textwidth}
  \centering
  \includegraphics[width=\linewidth]{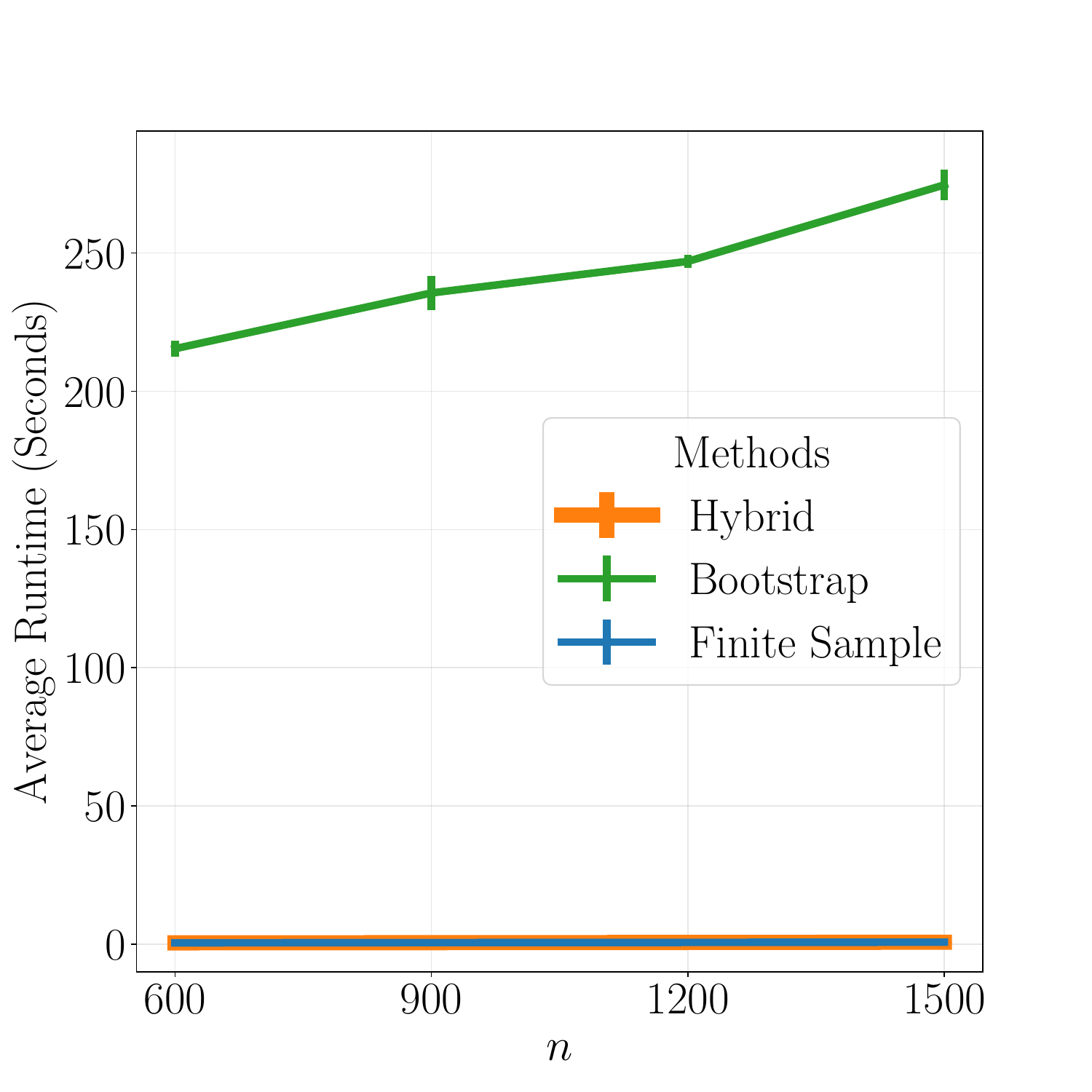}
  \caption{Model 5, Avg. Time.}
\end{subfigure}
\caption{\label{fig::mods4_5}  
Average   interval length, coverage, and runtime for Models 4 and 5. 
Error bars represent one standard deviation.
Dotted lines ({\bf \color{red}- - -}) indicate the nominal level $0.95$.}
\end{figure}

\begin{table}[bp]
\centering
\begin{tabular}{c|cc}
\hline 
	Mod.  & $P$ & $Q$ \\
	\hline                                  % 13
	$1$ & $\frac 1 2 N\Big((-1, -1)^\top, I_2\Big) + \frac 1 2 N\Big( (1, 1)^\top, I_2\Big)$
	    & $N(0, I_2) $        \\[0.05in]    
	$2$ & $\frac{1 + n^{-1/2}}{2}\delta_2 + \frac{1 - n^{-1/2}}{2}\delta_4$ 
		& $\frac 1 2 \delta_2 + \frac 1 2 \delta_5$ \\[0.05in]  
	$3$ & $\bbT(\frac 1 2, 1)$    & $\bbT(\frac 1 2, 5)$ \\[0.05in]                                  
	$4$ & $.95 N(0,1) + .05 N(0.1)$ 
	    & $N(0,1)$  \\[0.05in]
	$5$ & $.55  N\Big((-5, -5)^\top{,} I_2\Big) {+} .45 N\Big( (5, 5)^\top{,} I_2\Big)$
	    & $\frac 1 2 N\Big((-5, -5)^\top{,} I_2\Big) {+} \frac 1 2 N\Big( (5, 5)^\top{,} I_2\Big)$ \\
\end{tabular}
%\end{minipage}
%}
\caption{\label{tab::models} Parameter settings for Models 1-5.
%For any $r > 0$, $\mu_r$ denotes the uniform distribution over the sphere
%$\{\theta \in \bbR^d: \norm \theta = r\}$. 
For any $R > r > 0$, $\bbT(r,R)$ denotes the uniform distribution over the torus
$\big\{\big( (R + r \cos\theta) \cos\psi,(R + r \cos\theta) \sin\psi, r \sin\theta\big)^\top:$ $0 \leq \theta,\psi \leq 2\pi\big\}
\subseteq \bbR^3$.
We sample from $\bbT(r,R)$  using the Algorithm 1 of \cite{diaconis2013}.
%$\Sigma_a$ denotes the $d \times d$ Toeplitz matrix with $(i,j)$-th entry equal
%to $a^{|i-j|}$.
}\end{table}
We generate 100 samples of size $n=m= 600, 900, 1200$ and $1500$, from each of the pairs of distributions $(P,Q)$ described in Table~\ref{tab::models}. 
We choose the level $\alpha=.05$, the trimming constant $\delta=.1$ and the Monte Carlo sample size $N=500$, 
for which Assumption~\hyperref[assm::delta]{\textbf{A1($\delta;\alpha/N$)}}
 is met for the sample sizes under consideration.  
The number of bootstrap replications is set to $B=1,000$, and we set $r=2$
except where otherwise specified. %100 replications are performed for each model and sample size.
The empirical coverage and average lengths of the three confidence intervals are reported
in Figure \ref{fig::mods1_3} for Models 1-3, and in Figure \ref{fig::mods4_5} for Models 4-5.

Model 1 satisfies the regularity conditions required for the validity of the bootstrap, and 
we indeed observe its valid coverage for all sample sizes considered. When $n=600$, 
the finite-sample interval does not distinguish $\sw_{r,\delta}(P,Q)$ 
from zero on all replications, hence the hybrid confidence interval exhibits length and coverage between
those of the finite-sample and bootstrap intervals. For larger sample sizes, the length of the hybrid interval
essentially coincides with that of the the bootstrap interval. The distributions in Model 2 are not absolutely continuous, 
thus the bootstrap and hybrid intervals are seen to  undercover the true Wasserstein distance. Those of Model 3 are also not
absolutely continuous with respect to the Lebesgue measure on $\bbR^3$, since the supports
of $P$ and $Q$ are two-dimensional manifolds. Nevertheless, the linear projections of $P$ and $Q$
are absolutely continuous, with positive and distinct densities, thus ensuring that the conditions for all three
methods are met.
Models 4 and 5 consist of pairs of measures admitting Sliced Wasserstein distance near zero, 
causing the bootstrap to undercover. 
We also report the average runtime under these models in Figures \ref{fig::mods4_5}(c) and \ref{fig::mods4_5}(f),
showing a clear computational advantage of our finite-sample intervals over the bootstrap.

\begin{figure}[t]
\centering

\begin{subfigure}{0.25\textwidth}
  \centering
  \includegraphics[width=\linewidth]{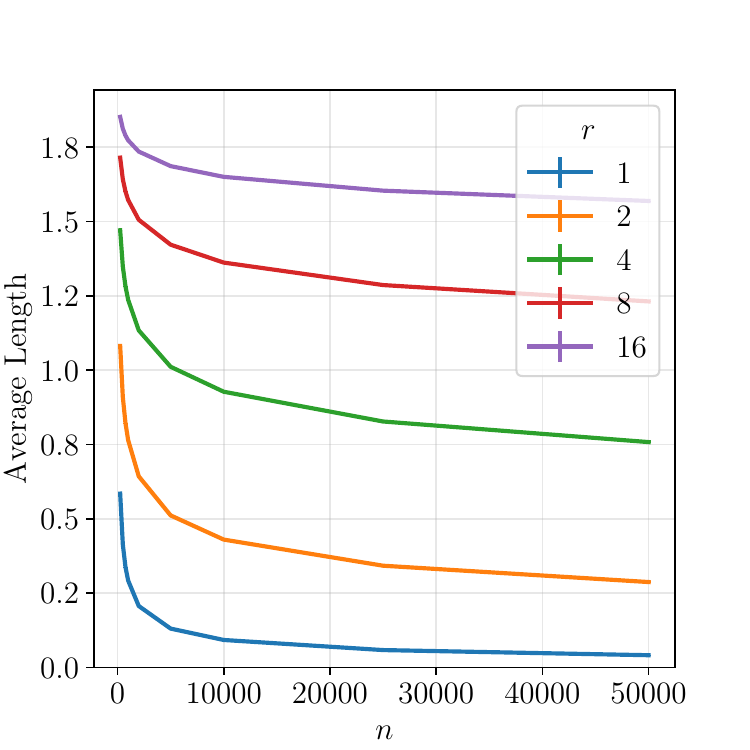}
  \caption{\centering Model 6(i), $\Delta = 0$.}
\end{subfigure}%
\begin{subfigure}{0.25\textwidth}
  \centering
  \includegraphics[width=\linewidth]{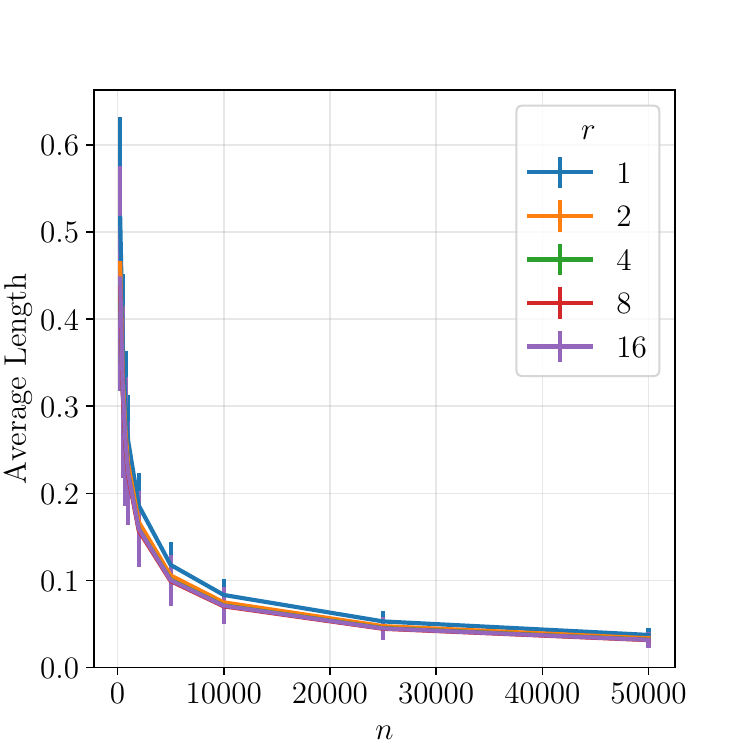}
  \caption{\centering  Model 6(ii), $\Delta=0$.}
\end{subfigure}%  
\begin{subfigure}{0.25\textwidth}
  \centering
  \includegraphics[width=\linewidth]{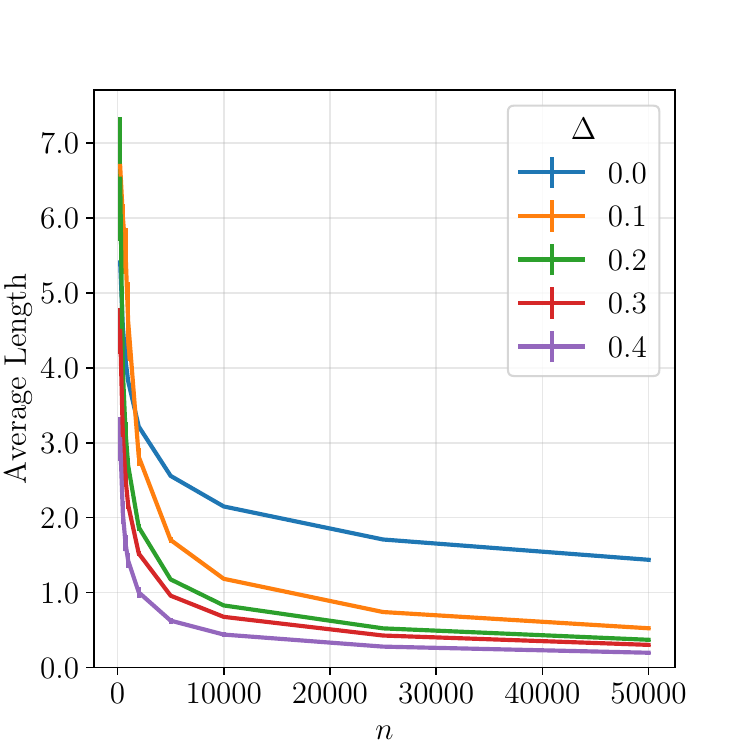}
  \caption{\centering  Model 6(i),{\color{white}d} $r=2$.}
\end{subfigure}%
\begin{subfigure}{0.25\textwidth}
  \centering
  \includegraphics[width=\linewidth]{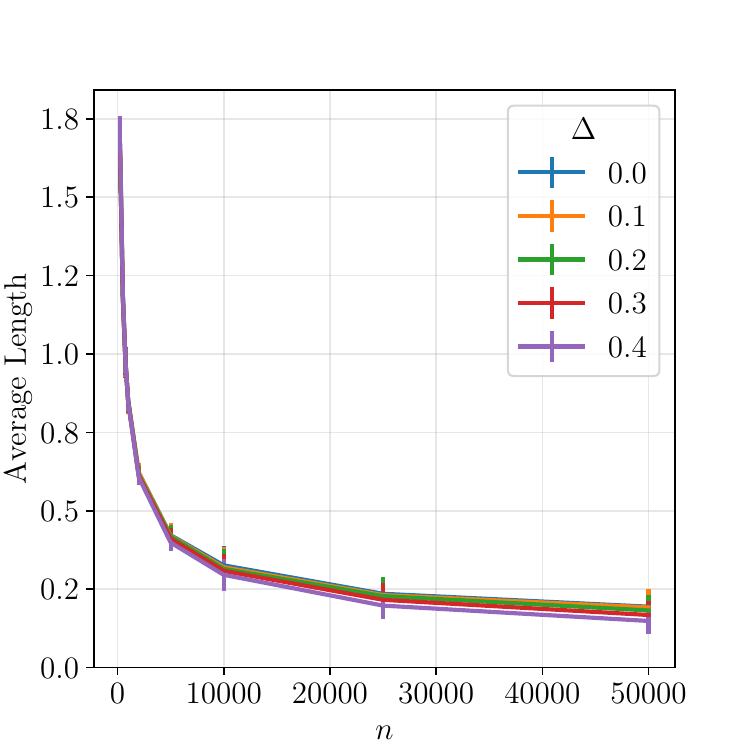}
  \caption{\centering Model 6(ii), $r=2$.}
\end{subfigure}%
\caption{\label{fig:sims} 
%\begin{subfigure}{0.333\textwidth}
%  \centering
%  \includegraphics[width=.92\linewidth]{figures/plot_model1_pi4}
%  \caption{Model 1, $\pi=0.4$.}
%\end{subfigure}\par\medskipLog-log scale plots for the simulation results under Model A ($\pi \in \{0.1, 0.25, 0.4\}$) and Model S.
\label{fig::asymp}
Average length of our finite-sample confidence interval for Models 6(i) and 6(ii), for varying values of $r$ and $\Delta$, based on 100 replications
for each sample size considered. Error bars represent one standard deviation of the confidence interval length. }
\end{figure}

\paragraph{Adaptivity and Asymptotic Confidence Interval Length.} 
We now illustrate the behaviour  predicted by Theorem \ref{thm::length}, and Corollaries thereafter, 
regarding the asymptotic length of our finite-sample confidence intervals. Consider the following two pairs of distributions,
\begin{alignat*}{4}
&\text{\bf Model 6(i).} ~~ &&P_1 &&= \frac {\delta_{-5} + \delta_5} 2, 
    ~~ &&Q_{1,\Delta} = \left(\frac 1 2 + \Delta\right) \delta_{-5} + \left(\frac 1 2 - \Delta\right) \delta_{5},
    % \label{eq::model_infinite}
     \\
&\text{\bf Model 6(ii).} ~~&&P_2 &&= U(-5, 5), 
      &&Q_{2,\Delta} =  \left(\frac 1 2 + \Delta\right)U(-5, 0) + 
                                \left(\frac 1 2 - \Delta\right)U(0, 5),
    % \label{eq::model_finite}
\end{alignat*}
for any $\Delta \in [0,1/2]$.
Notice that $J_{r,\delta}(P_1) = J_{r,\delta}(Q_{1,\Delta}) = \infty$, while
$J_{r,\delta}(P_2)$, $J_{r,\delta}(Q_{2,\Delta}) < \infty$. 
We report the average length of the finite-sample confidence interval \eqref{eq::one_dim_CI}, for each of Models~6(i) and 6(ii), 
based on 100 samples of
sizes $n=m\in[250,\  50,000]$. In Figures \ref{fig::asymp}(a) and \ref{fig::asymp}(b), 
we do so for $\Delta=0$ under varying
orders $r \in [1, 16]$ of the Wasserstein distance, while in Figures \ref{fig::asymp}(c) and \ref{fig::asymp}(d), we do so for a range of values  $\Delta \in [0, .4]$
under the fixed order $r=2$. 
In the former case, the average confidence interval length for the pair $(P_1, Q_{1,0})$
decays at an increasingly slow rate as $r$ increases, while that of the pair $(P_2, Q_{2,0})$
remains nearly unchanged. This behaviour was predicted by Corollary \ref{cor::length_examples}, indicating 
that the finite-sample interval length scales at the $n^{-1/2r}$ rate in general, but does so at the 
faster rate $n^{-1/2}$ for distributions admitting finite $\sj_{r,\delta}$ values. 
When $r=2$ is fixed, Figures \ref{fig::asymp}(c) and \ref{fig::asymp}(d) 
similarly exhibit increasing interval lengths for the pair $(P_1, Q_{1,\Delta})$ as $\Delta$ decreases to zero, 
yet nearly identical lengths for the pair $(P_2, Q_{2,\Delta})$. This behaviour is in line with Corollary \ref{cor::null}.

\section{Application to Likelihood-Free Inference}
\label{sec::applications} 
In a wide range of statistical applications, the likelihood function for a parametric model of interest may be intractable, 
though samples from the model can be easily generated. Examples 
include proton-proton collisions in particle physics \citep{brehmer2020}%,brehmer2020a}
, predator-prey dynamics in ecology \citep{lotka1920, lotka1920a}, 
inference for cosmological parameters in astronomy \citep{dalmasso2020a, dalmasso2020}, 
and network dynamics in queuing theory \citep{ebert2021}. In such applications, the practitioner
typically has access to a parametrized stochastic simulator for the data generating process, which produces samples 
from a distribution $P^\eta \in \calP(\bbR^d)$
with unknown closed form, and 
which depends on some physically meaningful parameters $\eta\in\Theta \subseteq \bbR^D$. 
The goal of \textit{likelihood-free inference} is
to characterize the values of $\eta$ for which an observed sample $X_1, \dots, X_n$ is likely to have been generated by the simulator.

Approximate Bayesian computation (ABC; \cite{sisson2018}) is arguably the most popular family of methodologies for likelihood-free inference. 
ABC methods repeatedly simulate parameter values in $\Theta$, and accept those for which the simulator $P^\eta$ produces a similar synthetic sample 
to the observed sample. The similarity of two samples is typically measured on the basis of summary statistics of the datasets. 
These summary statistics are often application-specific, and can be difficult to specify. Furthermore, 
due to the intractability of the likelihood, summary statistics can rarely be chosen as sufficient statistics for
$\eta$, making information loss inevitable. These considerations have motivated the development of methods which 
replace tailored summary statistics by distances between empirical measures of the synthetic and observed samples
\citep{park2016, gutmann2018, jiang2018}. In particular, 
\cite{bernton2019} and \cite{nadjahi2020} suggest the use of the Wasserstein and Sliced Wasserstein distances for this purpose. 

In what follows, we propose a simple alternative to such ABC methods, which provides
frequentist guarantees
for likelihood-free inference. Using the method developed in Section \ref{sec::cis}, 
we build confidence sets for the simulator's parameters on the basis
of minimizing the (Sliced) Wasserstein distance between the empirical measures
of an observed sample and synthetic samples from the simulator.
We focus on the case where $X_1, \dots, X_n$ is an i.i.d. sample from a distribution $P \in \calP(\bbR^d)$.
Fix
%$$\Theta_0 = \argmin_{\eta \in \Theta} \sw_{r,\delta}(P_\eta, P).$$
$\eta_0 \in  \argmin_{\eta \in \Theta} \sw_{r,\delta}(P,P^\eta),$
and let $\epsilon_0 = \sw_{r,\delta}(P,P^{\eta_0})$. 
Here $\eta_0$ denotes an $\sw_{r,\delta}$-projection of the distribution $P$ onto the family $\{P^\eta\}_{\eta \in \Theta}$.
If the simulator is correctly-specified, we have $P=P^{\eta_0}$ and $\epsilon_0 = 0$, whereas
if the simulator is misspecified, the set 
%$\Theta_\epsilon = 
$\{\eta \in \Theta: \sw_{r,\delta}(P,P^\eta) \leq \epsilon\},$
is empty for sufficiently small values of $\epsilon \geq 0$. 

We propose to construct confidence sets for $\eta_0$. 
For any $\eta \in \Theta$, and for any synthetic sample $Y_1^\eta, \dots, Y_m^\eta \sim P^\eta$,  
let $\big[\bar \ell_{nm}^{(N)}(\eta), \bar u_{nm}^{(N)}(\eta)\big]$ 
be a $(1-\alpha)$-confidence interval for $\sw_{r,\delta}(P, P^\eta)$,  
obtained via equation \eqref{eq::MC_CI_expanded} on the basis of $Y^\eta_1, \dots, Y^\eta_m$
and the observed sample $X_1, \dots, X_n$.
A confidence set for $\eta_0$ is then easily given by the following Proposition.
\begin{proposition}
\label{prop::likelihood_free}
Let $b,r\geq 1$ and $\delta \in (0,1/2)$ be fixed. Given any fixed
real number $\epsilon \geq \epsilon_0$, define 
$\widebar C_{nm}^{(N)} = \{\eta \in \Theta: \bar \ell_{nm}^{(N)}(\eta) \leq \epsilon\}.$
%\end{equation}
Then, 
$$\displaystyle\inf_{\substack{P \in \calK_{2r}(b) \\ \{P^\eta\}_{\eta\in\Theta} \subseteq \calK_{2r}(b)}}
\bbP\Big(\eta_0 \in \widebar C_{nm}^{(N)}\Big) \geq 1-\alpha-O(M_N^{-2}).$$
\end{proposition} 
Proposition \ref{prop::likelihood_free} provides a $(1-\alpha)$-confidence set
for the projection parameter $\eta_0$. In the well-specified setting $\epsilon_0=0$, 
$\widebar C_{nm}^{(N)}$ is simply
a confidence set for the parameter corresponding to the data-generating distribution $P=P^{\eta_0}$. 
We emphasize that no assumptions were made in the statement of Proposition \ref{prop::likelihood_free}
beyond the mild moment assumption $P,P^\eta \in \calK_{2r}(b)$ 
(which can be removed when $d=1$). 
The intractability 
of the likelihood function makes such assumption-lean inference particularly attractive.

In practice, computation of the lower confidence bounds $\bar \ell_{nm}^{(N)}(\eta)$ in  
%equation \eqref{eq::likelihood_free_CI} 
Proposition \ref{prop::likelihood_free} may be carried out over a finite grid $\{\eta_1, \dots, \eta_M\}\subseteq \Theta$ of candidate parameter values. 
While such a search may be computationally expensive, particularly for parameter spaces of high dimension $D$, 
it is akin to repeated sampling of parameters in ABC, or similar operations in other likelihood-free methods.
Nevertheless, efficient computation of the individual intervals $[\bar \ell_{nm}^{(N)}(\eta), \bar u_{nm}^{(N)}(\eta)]$ 
can dramatically reduce the computational burden of $\widebar C_{nm}^{(N)}$.
The simulation study in Section \ref{sec::simulations} suggests that the runtime of our finite-sample   intervals is considerably lower than
that of bootstrap-based methods (cf. Figures \ref{fig::mods4_5}(c) and \ref{fig::mods4_5}(f)).
  
\paragraph{Example: The Toggle Switch Model.} We illustrate our methodology in a systems biology model
used by \cite{bonassi2011,bonassi2015}. This model was analyzed by \cite{bernton2019} using an ABC method based on the Wasserstein distance,
and serves as a realistic example of likelihood-free inference with independent data. 
The toggle switch model describes the expression level of two genes across $n$ cells
over $T \in \bbN$ time points. Specifically, we let $(U_{i,t}, V_{i,t})$ denote their expression level in cell $i \in \{1, \dots, n\}$, and at
time $t \in \{1, \dots, T\}$. Given a starting value $(U_{i,0}, V_{i,0})$ for every $i=1, \dots, n$, the model is given by 
\begin{align}
\label{eq::toggle}
\begin{cases}
U_{i,t+1} = U_{i,t} + \frac{\alpha_1}{1 + V_{i,t}^{\beta_1}} - (1 + 0.03 U_{i,t}) + \frac 1 2 \xi_{i,t}  \\
V_{i,t+1} = V_{i,t} + \frac{\alpha_2}{1 + U_{i,t}^{\beta_2}} - (1 + 0.03 V_{i,t}) + \frac 1 2 \zeta_{i,t}
\end{cases},\quad t =1, \dots, T,
\end{align}
where $\alpha_1, \alpha_2 \in \bbR$ and $\beta_1, \beta_2 \geq 0$ are parameters, and $\xi_{i,t}, \zeta_{i,t}$ are independent standard Gaussian 
random variables. Following \cite{bernton2019}, $\xi_{i,t}$ and $\zeta_{i,t}$ are truncated such that $U_{i,t}, V_{i,t}$ remain nonnegative for all $i, t$. 
In applications, the full evolution 
\eqref{eq::toggle} is not observed, except for the %following
noisy measurement $X_i = U_{i,T} + \epsilon_i$ at time $T$,
%$$X_i = U_{i,T} + \epsilon_i, \quad i=1, \dots, n,$$
where $\epsilon_i \sim N(\mu, \mu\sigma/U_{i,T}^\gamma)$
are drawn conditionally on a realization
$U_{i,T}$ from the model~\eqref{eq::toggle}, for $i=1, \dots, n$,
where $\mu \in \bbR$ and $\sigma,\gamma\geq 0$. 
$X_1, \dots, X_n$ thus forms an i.i.d. sample from a distribution $P^\eta$ on $\bbR$ with respect to the parameter
$\eta=(\alpha_1, \alpha_2, \beta_1, \beta_2, \mu, \sigma, \gamma) \in \bbR^7$. A closed form for $P^\eta$ is unclear, but
the evolution \eqref{eq::toggle} makes simulation from $P^\eta$ simple, making this model a good candidate for likelihood-free inference.

We illustrate our methodology on simulated observations from this model in both a well-specified and a misspecified case. 
We treat the exponent parameters $\beta_1, \beta_2, \gamma$ as known, but possibly misspecified, 
and perform inference on $(\alpha_1, \alpha_2, \mu, \sigma)$.
In what follows, we set $U_{i,0} = V_{i,0} = 10$, $i=1, \dots, n$,  and we generate $n=2,000$ observations from
$P^{\eta_0}$ with $\eta_0 = (22, 12, 4,$ $4.5, 325, .25, .15)$, matching the parameter setting of \cite{bernton2019}.

\begin{figure}[h!]
\centering
\begin{minipage}{0.48\textwidth}
\includegraphics[width=\textwidth]{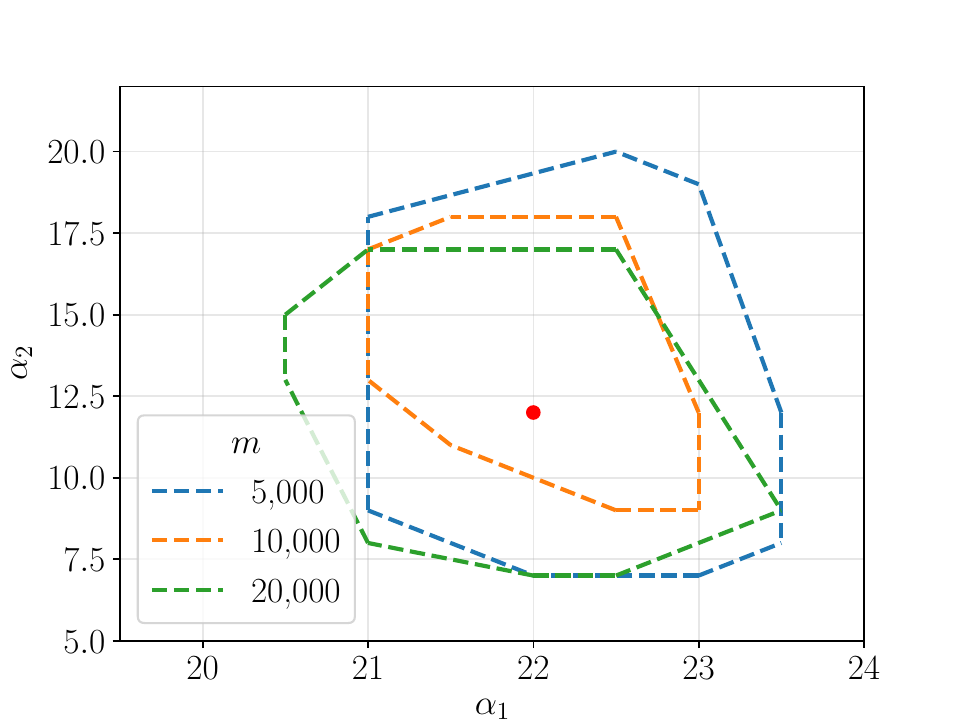}
\caption{Confidence sets for $\eta_0$ (in red ({\color{red} $\bullet$})) under the well-specified setting.%, for varying values of $m$.
         \label{fig::well_specified}\vspace{0.34in}}
\end{minipage}
\hspace{0.02\textwidth}
\begin{minipage}{0.48\textwidth}
\includegraphics[width=\textwidth]{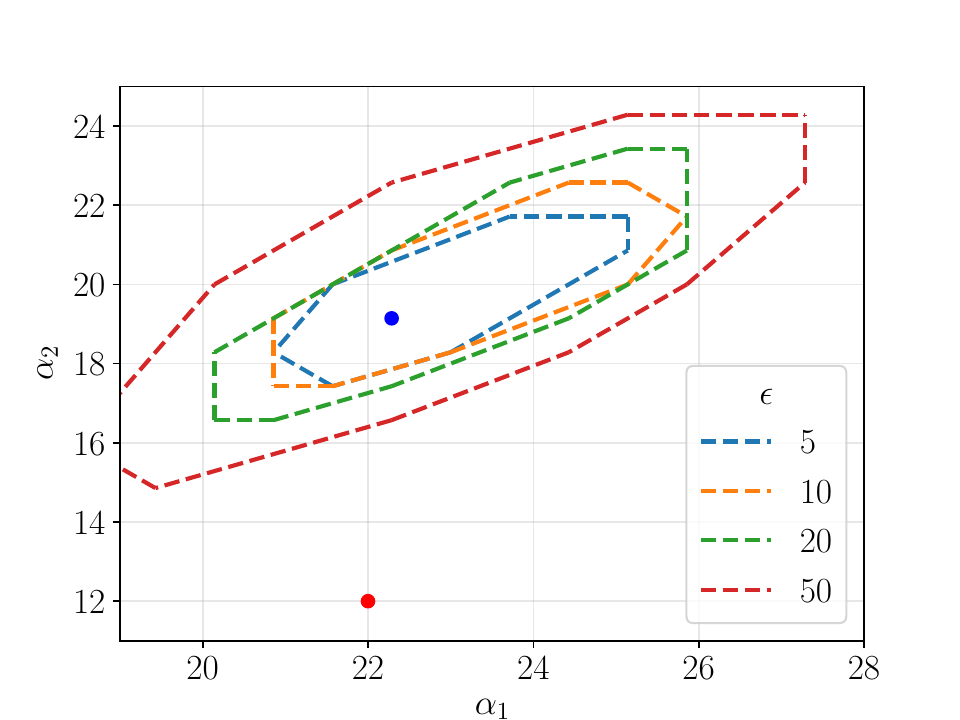}
\caption{Confidence sets for $\eta_0$ (in blue ({\color{blue}$\bullet$})) under  the misspecified setting. %, for varying values of $\epsilon$ 
		 The data-generating parameter is shown in red ({\color{red} $\bullet$}).
         \label{fig::misspecified}}
\end{minipage}
\end{figure}

\begin{itemize}
\item \textbf{Well-specified Setting.} Treat $\beta_1=4, \beta_2=4.5, \gamma=0.15$
as known and correctly specified. We compute the confidence set $C_{nm}$ of Proposition \ref{prop::likelihood_free}
with $r=1$ and $\epsilon=0$,
by repeatedly simulating $m$ observations from a grid of candidate values of $(\alpha_1, \alpha_2, \mu, \sigma) \in \bbR^4$,
for $m \in \{5 \cdot 10^3, 10^4, 2 \cdot 10^4\}$. The resulting two-dimensional confidence sets for the parameters
$(\alpha_1, \alpha_2)$, which are of primary interest, are reported in Figure \ref{fig::well_specified}. 
These confidence sets can be seen to cover the true parameter value, and naturally have decreasing 
area as $m$ increases.
\item \textbf{Misspecified Setting.} Using the same observed sample, 
we now misspecify the simulator with the values $\beta_1=2$ and $\beta_2=2$.
The resulting confidence set $C_{nm}$ is shown in Figure \ref{fig::misspecified}
for several choices of $\epsilon$, and can be seen to cover the projection parameter
$\eta_0$. The latter was approximated by $\heta_0 = \argmin_\eta W_{1,\delta}(P_M, P^\eta)$
over a grid of parameters $\eta$, where $P_M$ denotes a simulated  empirical measure based on $M=10,000$ observations.
\end{itemize}

\section{Conclusion and Discussion} 
\label{sec:conclusion}
Our aim in this paper has been to develop assumption-light 
finite-sample confidence intervals for the Sliced Wasserstein distance. 
After deriving minimax rates for estimating the Sliced Wasserstein distance, which are of independent interest, 
we bounded the length of our confidence intervals, showing that they achieve near minimax optimal length. Their length
is also shown to be adaptive to whether or not the underlying distributions are near the classical null, as well
as to their regularity, as measured by the magnitude of the functional $\sj_{r,\delta}$. 
These findings contrast asymptotic methods such as the bootstrap, whose validity we show is subject to certain prohibitive
assumptions on the underlying distributions, and whose asymptotic length does not enjoy the same adaptivity 
as that of our finite-sample intervals. 

Our work leaves open the problem of statistical inference for Wasserstein distances in dimension greater than one, for
which new techniques would have to be developed. Indeed, our work has hinged upon the representation of the one-dimensional
Wasserstein distance as the $L^r$ distance between quantile functions, which is unavailable in general dimension. 
For the same reason, our work does not shed light on statistical inference
for other modifications of the Wasserstein distance based on projections of distributions 
to low-dimensions greater than one, such as those summarized in Section \ref{sec::related_work}. 
We have shown that the Sliced Wasserstein distance can be estimated at dimension-independent rates, and it is of
interest to understand how this finding changes for other  low-dimensional modifications of the Wasserstein distance.

\section*{Acknowledgements}
The authors would like to thank an anonymous referee for  
thoughtful comments
which  significantly helped improve the content and presentation of this paper.  
Tudor Manole would like to thank Niccol\`o Dalmasso for conversations which inspired the likelihood-free inference application in Section~\ref{sec::applications}, and Arun Kumar Kuchibhotla for suggesting references on concentration inequalities for self-normalized empirical processes. 
This work was supported in part by the National Science Foundation grants DMS-1713003, DMS-2113684 and CIF-1763734. 
TM acknowledges the support of the Natural Sciences and Engineering Research Council of Canada (NSERC), 
 through a PGS D scholarship.  
 SB was supported by a Google Research Scholar Award and an Amazon Research Award.

%
%%%%%%%%%%%%%%%%%%%%%%%%%%%%%%%%%%%%%%%%%%%%%%%
%%% Supplementary Material, if any, should   %%
%%% be provided in {supplement} environment  %%
%%% with title and short description.        %%
%%%%%%%%%%%%%%%%%%%%%%%%%%%%%%%%%%%%%%%%%%%%%%%
%\begin{supplement}
%\stitle{ }
%\vspace{-0.3in}
%\sdescription{Proofs of results.} 
%\end{supplement}

\appendix
\noindent 

%\paragraph{Outline of Appendix}
%The remainder of this paper is organized as follows. 
%Appendix \ref{app::preliminary} collects several preliminary results which will be frequently used in the sequel. 
%Appendices \ref{app::pf_prop_empirical_sw}, \ref{app::main_lowerbound}, \ref{app:untrimmed}, \ref{app::length}, \ref{app::boot}, and \ref{app::pf_retest},
%respectively contain the proofs of our main results: Proposition \ref{prop::empirical_sw}, Theorem \ref{thm::minimax_distance}, 
%Theorem~\ref{thm:untrimmed_rate}, 
%Theorem~\ref{thm::length}, Theorem~\ref{thm::bootstrap}, and Proposition~\ref{prop::pretest}.
%All remaining results are proven in 
%Appendix \ref{app::secondary}. 
%Further discussion of Theorem \ref{thm::length} appears in Appendix \ref{app::mixed_case}.

\section{Preliminary Technical Results}
\label{app::preliminary}
In this section, we collect several preliminary results which will frequently be used in the sequel.
We begin with the following straightforward Lemma, which follows from Appendix A of \cite{bobkov2019}.
\begin{lemma}
\label{lem::abs_cont}
Let $P \in \calP(\bbR^d)$, $r \geq 1$, and $\delta \in (0,1/2)$. Let $F_\theta^{-1}$ denote the quantile
function of $P_\theta={\pi_\theta}_{\#} P$ for all $\theta \in \bbS^{d-1}$. If $\sj_{r,\delta}(P) < \infty$, 
then $F_\theta^{-1}|_{[\delta,1-\delta]}$ is absolutely continuous for $\mu$-almost all $\theta \in \bbS^{d-1}$. 
%In particular, for all such $\theta$, $F_\theta^{-1}$ is Lebesgue almost-everywhere differentiable, with strictly
%positive derivative wherever it exists.
\end{lemma}  
Furthermore, we describe the following characterization of distributions falling in the 
collections  $\calK_{r,\rho}(b)$ and $\widebar\calK_r(b)$.
 \begin{lemma}
\label{lem::calK}
Let $\delta \in (0,1/2)$ and $r,\rho,b \geq 1$.
Then,  
for all distributions $P \in \calK_{r,\rho}(b)$,   
$$\int_{\bbS^{d-1}} \big|F_\theta^{-1}(a)\big|^{r} d\mu(\theta) \leq  b(2/\delta)^{r/\rho}, \quad a \in\{ \delta, 1-\delta\}.$$ 
Furthermore, 
for all distributions $P \in \widebar\calK_r(b)$, we have
$$\sup_{\theta \in \bbS^{d-1}} \big|F_\theta^{-1}(a)\big| \leq  \left(\frac {2b} \delta\right)^{\frac 1 r}, \quad a \in \{\delta,1-\delta\}.$$
\end{lemma}
\paragraph{Proof of Lemma \ref{lem::calK}.}   
The claim is a simple consequence of Markov's inequality. 
Given $\theta \in \bbS^{d-1}$, it must hold that (a) $F_\theta^{-1}(\delta) < 0$ 
or (b)
$F_\theta^{-1}(1-\delta) > 0$. In the former case, since $F_\theta(F_\theta^{-1}(\delta)) \geq \delta$, we have
%$F_\theta(2F_\theta^{-1}(\delta)) \leq \delta$ and 
\begin{align} 
\label{eq::pf_calK_low_quantile}
\nonumber
\delta \leq \bbP\big\{X^\top\theta \leq F_\theta^{-1}(\delta)\big\}
&= \bbP\big\{-X^\top\theta \geq |F_\theta^{-1}(\delta)|\big\} \\
&\leq \bbP\big\{|X^\top\theta| \geq |F_\theta^{-1}(\delta)|\big\}
\leq \frac{\bbE_X[|X^\top\theta|^\rho]}{|F_\theta^{-1}(\delta)|^\rho}.
\end{align} 
while in the latter case, $F_\theta(F^{-1}_\theta(1-\delta)/2^{1/\rho})\leq 1-\delta$, therefore 
\begin{equation}
\label{eq::pf_calK_high_quantile}
\delta \leq \bbP\left\{X^\top\theta \geq \frac{F_\theta^{-1}(1-\delta)}{2^{1/\rho}}\right\} 
  \leq \bbP\left\{|X^\top\theta| \geq \frac{|F_\theta^{-1}(1-\delta)|} {2^{1/\rho}}\right\}
\leq \frac{2\bbE[|X^\top\theta|^\rho]}{|F_\theta^{-1}(1-\delta)|^\rho}
\end{equation}
We deduce, 
\begin{equation}
\label{eq::pf_calK_step}
\max_{a \in \{\delta,1-\delta\}} \big|F_\theta^{-1}(a)\big| \leq \left(\frac{2\bbE[|X^\top\theta|^\rho]}{\delta}\right)^{\frac 1 \rho}.
\end{equation}
Indeed, when both (a) and (b) hold, the above display follows from equations~\eqref{eq::pf_calK_low_quantile}
and~\eqref{eq::pf_calK_high_quantile}. When only (a) holds, it is clear that $|F_\theta^{-1}(\delta)| \geq |F_\theta^{-1}(1-\delta)|$,
thus the above display follows from equation~\eqref{eq::pf_calK_high_quantile}, and similarly when only case (b) holds.
Thus, since the above display holds for any $\theta \in \bbS^{d-1}$, we deduce that for both $a \in \{\delta,1-\delta\}$,
$$\int_{\bbS^{d-1}} \big|F_\theta^{-1}(a)\big|^r d\mu(\theta) 
 \leq \left(\frac{2 }{\delta}\right)^{\frac r \rho}\int_{\bbS^{d-1}}\bbE[|X^\top\theta|^\rho]^{\frac r \rho}d\mu(\theta)
 \leq \left(\frac{2 }{\delta}\right)^{\frac r \rho}b,$$
where we used the assumption $P \in \calK_{r,\rho}(b)$.    
 
To prove the second claim, one similarly has the following bound when
$P \in \widebar\calK_\rho(b)$, for $a \in \{\delta,1-\delta\}$,
\begin{align*}
\sup_{\theta \in \bbS^{d-1}} \big|F_\theta^{-1}(a)\big|
 &\leq \sup_{\theta \in \bbS^{d-1}}\left(\frac{2\bbE[|X^\top\theta|^\rho]}{\delta}\right)^{\frac 1 \rho} \\
 &\leq  \left(\frac 2 \delta \bbE\left[\sup_{\theta \in \bbS^{d-1}}|X^\top\theta|^\rho\right]\right)^{\frac 1 \rho} 
 =   \left(\frac 2 \delta \bbE\left[\norm X^\rho\right] \right)^{\frac 1 \rho} \leq \left(\frac{2b}{\delta}\right)^{\frac 1 \rho}.
\end{align*}  
\qed
 
\section{Proof of Propositions \ref{prop::ub_empirical} and \ref{prop::empirical_sw}}
\label{app::pf_prop_empirical_sw}

We shall begin by proving Proposition~\ref{prop::ub_empirical}(i) and Proposition~\ref{prop::empirical_sw}(ii). 
As we shall explain, Proposition~\ref{prop::ub_empirical}(ii) will then follow as a special case of Proposition~\ref{prop::empirical_sw}(ii),
and Proposition~\ref{prop::empirical_sw}(i) will follow from Proposition~\ref{prop::ub_empirical}(i). 
Our proofs will make use of Examples~\ref{ex::DKW} and \ref{ex::rel_VC} which appear
in Section~\ref{sec::cis} of the main text, and of their corresponding proofs in Appendix~\ref{app::secondary}.
We shall also make use of the following Lemma, which is proven in Appendix~\ref{app::lem_pf_empirical_sw_bound}.
\begin{lemma}
\label{lem::empirical_sw_bound}
Let $\delta \in (0,1/2)$. Then, for any $\rbar \in [r,r+1]$, there exists
a constant $B_r > 0$ depending only on $r$ such that
for all $\theta \in \bbS^{d-1}$ and all $\delta \geq 2(r+2)/n$,
$$ \max_{a \in \{\delta ,1-\delta\}} \bbE \left[\big|F_{\theta,n}^{-1}(a)\big|^{\rbar}\right]
 \leq B_r \left( 1+\max_{a \in \left\{\frac \delta 2,1-\frac \delta 2\right\}} \big|F_\theta ^{-1}(a)\big|^{\rbar} + 
  \left(\frac{\bbE\big|X^\top\theta\big|^2}{\delta}\right)^{\frac{\rbar}{2}}\right).$$
\end{lemma}
In the above result and throughout this section, recall
that $F_{\theta,n}^{-1}$ (resp. $G_{\theta,n}^{-1}$, $F_\theta^{-1}$, $G^{-1}_\theta$) denotes the quantile function of
the distribution $P_{\theta,n} = {\pi_\theta}_\# P_n$ (resp. 
$Q_{\theta,m} = {\pi_\theta}_\# Q_m, P_\theta = {\pi_\theta}_{\#} P, Q_\theta = {\pi_\theta}_\# Q$), for any $\theta \in \bbS^{d-1}$. 
Finally, throughout the remainder of this section, the symbol ``$\lesssim$'' is used
to hide universal constants depending only on $r$.

\paragraph{Proof of Proposition~\ref{prop::ub_empirical}(i).}
The claim is trivial if $\sj_{r,\delta/2}(P) = \infty$, thus assume otherwise.
We shall begin bounding the expectation of the following quantity
$$Z_n(\theta) = W_{r,\delta}^r(P_{\theta,n}, P_\theta),$$
for any fixed $\theta \in \bbS^{d-1}$. The result will then follow by integration
over $\bbS^{d-1}$. We begin with the following key high probability bound.
\begin{lemma}
\label{lem::high_prob_sw_null}
Let $y_0 = \sqrt \delta/4$ and $C_r > 0$ a constant depending only on $r$. Then, 
for all
$y \in (0,y_0]$, we have
$$\esssup_{\theta \in \bbS^{d-1}} \bbP \Big(Z_n(\theta) \geq C_r y^r J_{r,\delta/2}(P_\theta)\Big)\leq 
\frac{2n+1}{16}e^{-\frac{ny^{2}}{16}},$$
where the essential supremum is taken with respect to the measure $\mu$.
 \end{lemma}
Lemma~\ref{lem::high_prob_sw_null}
is proven in Appendix~\ref{app::high_prob_sw_null}. Now, let 
$T_\theta= C_r y_0^r J_{r,\delta/2}(P_\theta)$, so that
\begin{align}
\label{eq:decomp_Zn_theta}
\bbE[Z_n(\theta)]
 &= \bbE[Z_n(\theta)\cdot I(Z_n(\theta) >  T_\theta)] + \bbE[Z_n(\theta)\cdot I(Z_n(\theta) \leq T_\theta)].
\end{align}
To bound the first term, notice first that
$$Z_n(\theta) = \frac 1 {1-2\delta} \int_\delta^{1-\delta} \big|F_{\theta,n}^{-1}(u) - F_\theta^{-1}(u)\big|^r du 
\lesssim \max_{a \in \{\delta,1-\delta\}} \big| F_{\theta,n}^{-1}(a)\big|^r + \big|F_\theta^{-1}(a)\big|^r.$$
Now, let $s = 1/(1-1/\eta)$, where $\eta=\rbar/r$ and $\rbar \in (r,r+1]$. 
Then, by H\"older's inequality, we have uniformly in $\theta\in \bbS^{d-1}$,
\begin{align}
		 \label{eq::step_lemmaB1B2_pf_prop1}
\nonumber  \bbE\big[Z_{n}&(\theta)\cdot  I(Z_{n}(\theta) > T_\theta)\big]  \\
\nonumber
&\leq \norm{Z_{n}(\theta)}_{L^\eta(\bbP)} \norm{I(Z_{n}(\theta) > T_\theta)}_{L^s(\bbP)} \\
\nonumber
&\lesssim \max_{a \in \left\{\delta ,1- \delta \right\}}\Big(\big\|F_{\theta,n}^{-1}(a)\big\|_{L^{\rbar}(\bbP)}^{\rbar} +
 				 \big|F_{\theta}^{-1}(a)\big|^{\rbar}  \Big)^{\frac 1 \eta}
\norm{I(Z_{n}(\theta) > T_\theta)}_{L^s(\bbP)} \\
 &\lesssim  \left(1+ \max_{a \in \left\{\frac \delta 2,1-\frac \delta 2\right\}}\big|F_{\theta}^{-1}(a)\big|^{\rbar}  + 
		 \left(\frac{\bbE|X^\top\theta|^2}{\delta}\right)^{\frac{\rbar}{2}} \right)^{\frac 1 \eta}
		 \left(\frac{2n+1}{16}\right)^{\frac 1 s}e^{-\frac{ny_0^{2}}{16s}}\\
\nonumber  &\leq  \left(1+ \max_{a \in \left\{\frac \delta 2,1-\frac \delta 2\right\}}\big|F_{\theta}^{-1}(a)\big|^r  + 
		 \left(\frac{\bbE|X^\top\theta|^2}{\delta}\right)^{\frac{r}{2}} \right)
		 \left(\frac{2n+1}{16}\right)^{\frac 1 s}e^{-\frac{ny_0^{2}}{16s}},
 \end{align}
 where we invoked Lemmas~\ref{lem::calK}--\ref{lem::high_prob_sw_null}
 in equation~\eqref{eq::step_lemmaB1B2_pf_prop1}. 
 %{\color{blue}(TODO I think can do $r\rightarrow\rbar$ throughout).}
 Therefore, using the fact that $s \geq 1$, $P \in \calK_{r}(b)$, and invoking
 Lemma~\ref{lem::calK}, we obtain
\begin{align}
\label{eq:tail_prop1_pf_finite}
\int_{\bbS^{d-1}} \bbE\big[Z_{n}(\theta)\cdot  I(Z_{n}(\theta) > T_\theta)\big]d\mu(\theta)
 \lesssim ( b/\delta^{r/2}) ne^{-\frac{ny_0^{2}}{16s}} 
\lesssim   b n  e^{-c_0 n\delta}/\delta^{r/2},
\end{align}
for a universal constant $c_0 > 0$ depending only on $r$.
We further bound the second term in equation~\eqref{eq:decomp_Zn_theta}. 
Setting $t_\theta =T_\theta  \wedge C_r J_{r,\delta/2}(P_\theta) \left(\frac{16}{n}\log\left(\frac{2n+1}{16}\right)\right)^{r/2}$, we arrive at
\begin{align*}
\bbE\big[Z_{n}(\theta)&\cdot I(Z_{n}(\theta) < T_\theta)\big]  \\
 &=  \int_0^\infty \bbP\big(Z_{n}(\theta) \cdot I( Z_{n}(\theta) < T_\theta) \geq x\big)dx \\
 &\leq  \int_0^{T_\theta} \bbP\big(Z_{n}(\theta)  \geq x\big)dx \\
 &\leq t_\theta 
  +   \int_{t_\theta}^{T_\theta} \bbP\big(Z_{n}(\theta)  \geq x\big)dx  \\
 &\leq t_\theta + \frac{2n+1}{16}  \int_{t_\theta}^{T_\theta} 
 \exp\left\{-\frac n {16} \left(\frac{x}{J_{r,\delta/2}(P_\theta)C_r}\right)^{2/r}\right\}dx \\
 &= t_\theta + \frac{r(2n+1)C_r J_{r,\delta/2}(P_\theta)}{16} \left(\frac{4 }{\sqrt n}\right)^r 
 \int_{\sqrt{\log\left(\frac{2n+1}{16}\right)}}^\infty
 e^{-y^2}y^{r-1} dy \\
 &\lesssim t_\theta + \frac{r(2n+1)C_r J_{r,\delta/2}(P_\theta)}{16} \left(\frac{4 }{\sqrt n}\right)^r
 \int_{\sqrt{\log\left(\frac{2n+1}{16}\right)}}^\infty
 e^{-y^2/2}  dy,
\end{align*}
where we used the change of variable $y =  \frac{\sqrt n} {4} \left(\frac{x}{J_{r,\delta}(P_\theta) C_r}\right)^{1/r}$, and  
where the final inequality holds for all $n$ larger than a universal constant depending only on $r$. 
It follows that
\begin{align*}
\bbE\big[Z_{n}(\theta)\cdot I(Z_{n}(\theta) < T_\theta)\big] 
 \lesssim t_\theta + \frac{J_{r,\delta/2}(P_\theta)}{n^{r/2}} 
 \lesssim J_{r,\delta/2}(P_\theta)\left[\frac{1}{n}\log\left(\frac{2n+1}{16}\right)\right]^{r/2}.
\end{align*}
Putting this fact together with equation~\eqref{eq:tail_prop1_pf_finite}, we have by the Fubini-Tonelli Theorem,
\begin{align*}
\bbE\big[\sw_{r,\delta}^r(P_n,P)\big] 
 &= \int_{\bbS^{d-1}} \bbE[Z_{n}(\theta)] d\mu(\theta) \\
 &\lesssim 
\sj_{r,\delta/2}(P)\left(\frac{\log n}{n}\right)^{\frac r 2} + b n  e^{-c_0 n\delta}/\delta^{r/2}.
\end{align*}
The claim now follows since
$\bbE\big[\sw_{r,\delta}(P_n,P)\big] {\leq} \Big(\bbE\big[\sw_{r,\delta}^r(P_n,P)\big] \Big)^{\frac 1 r}$.\qed

\paragraph{Proof of Proposition~\ref{prop::empirical_sw}(ii).} 
As we shall show, the assumption $ \sw_{r,\delta}(P,Q) \geq \Gamma$ implies
that the deviation
$$\Delta_{nm} = \sw_{r,\delta}(P_n, Q_m) - \sw_{r,\delta}(P,Q)$$
is of same order
as
$$D_{nm} =  \sw_{r,\delta}^r(P_n, Q_m) - \sw_{r,\delta}^r(P,Q).$$
To bound this quantity, define for all $\theta \in \bbS^{d-1}$,
\begin{align*}
a_{\theta} &= \min\Big\{ F_\theta^{-1}(\delta), F_{\theta,n}^{-1}(\delta), G_\theta^{-1}(\delta),
G_{\theta,m}^{-1}(\delta)\Big\},\quad \\
b_{\theta} &= \max\Big\{ F_\theta^{-1}(1-\delta), F_{\theta,n}^{-1}(1-\delta), G_\theta^{-1}(1-\delta),
G_{\theta,m}^{-1}(1-\delta)\Big\},\\
M(\theta) &= \max 
    \Big\{
           \big|F_\theta^{-1}(\delta/2)\big|, 
           \big|G_\theta^{-1}(\delta/2)\big|,
           \big|F_\theta^{-1}(1-\delta/2)\big|,
           \big|G_\theta^{-1}(1-\delta/2)\big|,1
    \Big\},
\end{align*}
and
\begin{align*}
Z_{nm}(\theta) = (b_\theta - a_\theta)^{r-1}  \int_\delta^{1-\delta} 
    \Big[|F_{\theta,n}^{-1}(u) - F_\theta^{-1}(u)| + |G_{\theta,m}^{-1}(u) - G_\theta^{-1}(u)|\Big]du.
\end{align*}
The bulk of our proof is now contained in the following result.
\begin{lemma}
\label{lem::rth-deviation}
We have,
$$|D_{nm}| \leq \frac{r}{1-2\delta} \int_{\bbS^{d-1}} Z_{nm}(\theta)d\mu(\theta).$$
Assume further that $P,Q \in \calK_{r}(b)$. Let $y_0 =  \delta/2$. Then,
there exists a universal constant $A > 0$ such that for all $y \in (0,y_0]$
%and all $\theta \in \bbS^{d-1}$,
$$\sup_{\theta \in \bbS^{d-1}}\bbP(|Z_{nm}(\theta)| \geq A M^r(\theta) y) \leq 4\exp(-2(n\wedge m) y^2).$$
\end{lemma} 
Let $A,y_0 > 0$ be as in 
Lemma~\ref{lem::rth-deviation}. Similarly as in the proof of 
Proposition~\ref{prop::ub_empirical}, let $T_\theta = AM^r(\theta)y_0$. Then, we have,
\begin{align}
\label{eq::pf_deviation_Delta_decomp}
\bbE[Z_{nm}(\theta)] = \bbE\big[Z_{nm}(\theta)   I(Z_{nm}(\theta) \geq T_\theta)\big] + 
\bbE\big[Z_{nm}(\theta)  I(Z_{nm}(\theta) < T_\theta)\big].
\end{align}
Set $s = 1/(1-1/\eta)$, where again $\eta=\rbar/r$ and $\rbar \in (r,r+1]$, so that by H\"older's inequality, we have uniformly in $\theta\in \bbS^{d-1}$,
\begin{align}
\label{eq::pf_deviation_step_finite}
\nonumber  &\bbE\big[Z_{nm}(\theta)\cdot  I(Z_{nm}(\theta) \geq T_\theta)\big]  \\
\nonumber&\leq \norm{Z_{nm}(\theta)}_{L^\eta(\bbP)} \norm{I(Z_{nm}(\theta) \geq T_\theta)}_{L^s(\bbP)} \\
%\nonumber  &\leq \left\{B M^{\rbar}(\theta) + 
%		\frac 1 {n \wedge m} \Big( \left[\bbE|X^\top\theta|^{\rbar}\right] + \left[\bbE|Y^\top\theta|^{\rbar}\right]\Big)\right\}^\eta \bbP(Z_{nm}(\theta) \geq T_\theta)^{1/s} \\
\nonumber &\lesssim \Big(\norm{a_\theta}_{L^{\rbar}(\bbP)}^{\rbar}+\norm{b_\theta}_{L^{\rbar}(\bbP)}^{\rbar}\Big)^{\frac 1 \eta}\exp(-2(n\wedge m) y_0^2/s) \\
%\nonumber    & \qquad\quad \times\max_{a \in \left\{\frac \delta 2,1-\frac \delta 2\right\}}\left( 
% \big\|F_{\theta,n}^{-1}(a)\big\|_{L^{\rbar}(\bbP)}^{\rbar} +
% \big\|G_{\theta,m}^{-1}(a)\big\|_{L^{\rbar}(\bbP)}^{\rbar} + \left(\frac{\bbE|X^\top\theta|^2}{\delta}\right)^{\frac{\rbar}{2}} 
% + \left(\frac{\bbE|Y^\top\theta|^2}{\delta}\right)^{\frac{\rbar}{2}}
% 				   \right)^{\frac 1 \eta}\\
\nonumber 				    &\lesssim 
 \left( 
 M^{\rbar}(\theta) + \left(\frac{\bbE|X^\top\theta|^2}{\delta}\right)^{\frac{\rbar}{2}} 
 + \left(\frac{\bbE|Y^\top\theta|^2}{\delta}\right)^{\frac{\rbar}{2}}
 				   \right)^{\frac 1 \eta}\exp(-2(n\wedge m) y_0^2/s) \\
 &\lesssim 
 \left( 
 M^{r}(\theta) + \left(\frac{\bbE|X^\top\theta|^2}{\delta}\right)^{\frac{r}{2}} 
 + \left(\frac{\bbE|Y^\top\theta|^2}{\delta}\right)^{\frac{r}{2}}
 				   \right)\exp(-2(n\wedge m) y_0^2/s),
\end{align}
where we invoked Lemmas~\ref{lem::empirical_sw_bound} and~\ref{lem::rth-deviation}.
Now, notice that $\int_{\bbS^{d-1}} M^r(\theta)d\mu(\theta) \lesssim b/\delta^{r/2}$ by Lemma~\ref{lem::calK},
since $P,Q \in \calK_r(b)$. 
Therefore, integrating both sides of the above display with respect to $\mu$ leads to 
\begin{align*}
\int_{\bbS^{d-1}} &\bbE\big[Z_{nm}(\theta) I(Z_{nm}(\theta) \geq T_\theta)\big] d\mu(\theta) \lesssim 
 \exp(-c_1(n\wedge m)\delta^2) b/\delta^{r/2},
\end{align*}
for a constant $c_1 > 0$. 
We now bound the second term in equation~\eqref{eq::pf_deviation_Delta_decomp}. 
We again use Lemma~\ref{lem::rth-deviation} to obtain
\begin{align*}
\int_{\bbS^{d-1}} \bbE\big[Z_{nm}&(\theta)\cdot I(Z_{nm}(\theta) < T_\theta)\big] d\mu(\theta) \\
 &= \int_{\bbS^{d-1}} \int_0^\infty \bbP\big(Z_{nm}(\theta) \cdot I( Z_{nm}(\theta) < T_\theta) \geq x\big)dx d\mu(\theta)\\
 &\leq \int_{\bbS^{d-1}} \int_0^{T_\theta} \bbP\big(Z_{nm}(\theta)  \geq x\big)dx d\mu(\theta)\\
 &\leq 4\int_{\bbS^{d-1}}  \int_0^{T_\theta} \exp\left\{-2(n\wedge m) \left(\frac{x}{AM^r(\theta)}\right)^2\right\}dx d\mu(\theta) \\
 &\lesssim \frac 1 {\sqrt{n\wedge m}}\int_{\bbS^{d-1}} M^r(\theta)d\mu(\theta) \lesssim \frac{b}{\delta^{r/2}\sqrt{n\wedge m}}.
\end{align*}
Putting this fact together with equation~\eqref{eq::pf_deviation_step_finite}, 
and applying the Fubini-Tonelli Theorem and Lemma~\ref{lem::rth-deviation}, we arrive at
\begin{align}
\label{eq::D_nm_bound}
\nonumber 
\bbE|D_{nm}| &\lesssim \frac 1 {1-2\delta}
   \bbE\left[\int_{\bbS^{d-1}} Z_{nm}(\theta)d\mu(\theta)\right] \\
  &=  \frac 1 {1-2\delta} \int_{\bbS^{d-1}} \bbE\big[Z_{nm}(\theta)\big] d\mu(\theta)  
 \lesssim \frac{b}{\delta^{r/2}(1-2\delta)\sqrt{n\wedge m}}.
 \end{align}
Finally, the numerical inequality $|x^r - y^r| \geq y^{r-1} |x-y|$ for all $x,y > 0$ implies
\begin{equation}
\label{eq::rpow_to_sw}
\bbE |D_{nm}| \geq \sw_{r,\delta}^{r-1}(P,Q) \bbE |\Delta_{nm}|\geq \Gamma^{r-1} \bbE |\Delta_{nm}|,
\end{equation}
so that $\bbE |\Delta_{nm}| \lesssim \Gamma^{1-r} b(n\wedge m)^{-1/2}/(\delta^{r/2}(1-2\delta))$, as claimed.\qed
 \newline
 
\paragraph{Proof of Proposition~\ref{prop::ub_empirical}(ii).}
Introduce an i.i.d. sample $X_1',\dots,X_n' \sim P$ independent of $X_1, \dots, X_n$,
and let $P_n = \frac 1 n \sum_{i=1}^n \delta_{X_i'}$. 
It follows from convexity of the mapping $(P,Q) \mapsto \sw_{r,\delta}(P,Q)$,
similarly as in the proof of Theorem 4.3 of \cite{bobkov2019}, that
$$\bbE\big[\sw_{r,\delta}^r(P_n,P)\big] \leq \bbE\big[\sw_{r,\delta}^r(P_n, P_n')\big].$$
The claim now follows from equation~\eqref{eq::D_nm_bound} with $P=Q$, which implies the bound
$$\bbE\big[\sw_{r,\delta}^r(P_n,P_n')\big] \lesssim \frac{bn^{-1/2}}{\delta^{r/2}(1-2\delta)},$$
so that, $\bbE\big[\sw_{r,\delta}(P_n,P_n')\big] \lesssim  b^{1/r}n^{-1/2r}/(\sqrt \delta(1-2\delta)^{1/r}).$
\qed\\

\paragraph{Proof of Proposition~\ref{prop::empirical_sw}(i).} 
The claim is immediate from Proposition~\ref{prop::ub_empirical}(i), by the triangle inequality.\qed
 
 \paragraph{}  It remains to prove Lemmas~\ref{lem::empirical_sw_bound}--\ref{lem::rth-deviation}.

\subsection{Proof of Lemma~\ref{lem::empirical_sw_bound}} 
\label{app::lem_pf_empirical_sw_bound}
We shall make use of Bennett's inequality, which we recall as follows
using the notation of \cite{pollard2002}.
\begin{lemma}[Bennett's Inequality]
Let $Z_1, \dots, Z_n$ be i.i.d. random variables bounded above
by 1, and such that $\bbE[Z_1] = \mu \in \bbR$, $\Var[Z_1] = v > 0$. 
Then, 
$$\bbP\left\{ \sum_{i=1}^n (Z_i-\mu) \geq x\right\} \leq \exp\left\{ -\frac{x^2}{2W} \psi\left(\frac{x}{W}\right)\right\},$$
for all $x \geq 0$, where $W \geq nv$ is arbitrary, and where for all $t \geq -1$,
$$\psi(t) = \frac{(1+t)\log(1+t) - t}{t^2/2},\quad \text{for} \ t \neq 0, \ \text{and} \ \psi(0)=1.$$
\end{lemma}
Turning back to the proof, %  similarly as \cite{bobkov2019}, Theorem 2.14, 
fix $\theta \in \bbS^{d-1}$. Set $m_\theta = \bbE(|X^\top\theta|^2)$, and 
$$x_\theta = 2\big[|F_\theta^{-1}(\delta/2)| \vee |F_\theta^{-1}(1-\delta/2)| \vee \ell\sqrt{m_\theta/\delta}\vee 1\big],$$
for a constant $\ell > 0$ to be determined below. In particular, for all $x \geq x_\theta$,
$$F_\theta(x) \geq 1-\delta/2, \quad \text{and} \quad F_\theta(-x) \leq \delta/2.$$
Then, 
for all $x \geq x_\theta$, 
\begin{align*}
\bbP(-F_{\theta,n}^{-1}(\delta) > x)
 &\leq \bbP(F_{\theta,n}(-x) > \delta) \\
 &\leq \bbP(F_{\theta,n}(-x) - F_\theta(-x) > \delta/2) \\ 
 &\leq \exp\left\{ -\frac{(\delta n)^2}{8W_\theta} \psi\left(\frac{n\delta}{2W_\theta}\right)\right\},
\end{align*}
for any given $W_\theta \geq  n\Var[I(X^\top\theta \leq -x)] = nF_\theta(-x)(1-F_\theta(-x))$, by Bennett's inequality. 
Now, from Markov's inequality, one also has 
$$F_\theta(-x) = \bbP(X^\top\theta \leq -x) = \bbP(-X^\top\theta \geq x) \leq \bbP(|X^\top\theta| \geq x) \leq m_\theta/x^2,$$ 
thus, we may take $W_\theta = nm_\theta/x^2$. Furthermore,
\begin{align*}
\frac{(\delta n)^2}{8W_\theta}  \psi\left(\frac{n\delta}{2W_\theta}\right)  
 &= \frac{(\delta n)^2}{8W_\theta} \frac{(1+\frac{n\delta}{2W_\theta})\log(1+\frac{n\delta}{2W_\theta}) - \frac{n\delta}{2W_\theta}}{(\frac{n\delta}{2W_\theta})^2/2} \\
 &= W_\theta \left[\left(1+\frac{n\delta}{2W_\theta}\right)\log\left(1+\frac{n\delta}{2W_\theta}\right) - \frac{n\delta}{2W_\theta} \right] \\ 
 &\geq \frac{nm_\theta}{x^2} \left[\frac{x^2\delta}{2m_\theta}\log\left(1+\frac{x^2\delta}{2m_\theta}\right) - \frac{x^2\delta}{2m_\theta} \right] \\ 
 &\geq n \left[\frac{\delta}{2}\log\left(1+\frac{x^2\delta}{2m_\theta}\right) - \frac{\delta}{2} \right] \\ 
 &\geq \frac{n\delta}{4}\log\left(1+\frac{x^2\delta}{2m_\theta}\right),
  \end{align*}
where the last inequality holds for all $x \geq  x_\theta$ upon choosing the constant $\ell = \sqrt{2(e^2-1)}$
in the definition of $x_\theta$.
%therefore,
%$$\frac{(\delta n)^2}{8W} \psi\left(\frac{n\delta}{2W}\right)
%\geq %\frac{(\delta n)^2}{2W}  \frac{8W^2}{(n\delta)^2} \log\left(1 + \frac {n\delta} {2W}\right)
%W \log \left(1+\frac{n\delta}{2W}\right)
%= \frac{nb}{x} \log\left(1 +\frac{x\delta}{2b}\right),$$
Therefore, for all such $x$,
\begin{align*}
\bbP(-F_{\theta,n}^{-1}(\delta) > x)
\leq \left(1 +\frac{x^2\delta}{2m_\theta}\right)^{-\frac{n\delta}{4}}
\end{align*}
Applying a similar argument, we obtain that for all $x \geq x_\theta$, 
\begin{align*}
\bbP(F_{\theta,n}^{-1}(1-\delta) > x)
 &\leq \bbP(1-F_{\theta,n}(x) > \delta) \\
 &\leq \bbP(F_\theta(x) - F_{\theta,n}(x) > \delta/2) 
 \leq \left(1 +\frac{x^2\delta}{2m_\theta}\right)^{-\frac{n\delta}{4}}.
\end{align*}
Now notice that
$$|F_{\theta,n}^{-1}(1-\delta)| \leq F_{\theta,n}^{-1}(1-\delta) \vee (-F_{\theta,n}^{-1}(\delta)),$$
thus we arrive at
$$\bbP(|F_{\theta,n}^{-1}(1-\delta)| \geq x) \leq 2\left(1 +\frac{x^2\delta}{2m_\theta}\right)^{-\frac{n\delta}{4}},\quad x \geq x_\theta.$$ 
It follows that
\begin{align*}
\bbE|F_{\theta,n}^{-1}(1-\delta)|^{\bar r }
 &=\bar r \int_0^\infty x^{{\bar r -1}} \bbP(|F_{\theta,n}^{-1}(1-\delta)| \geq x)dx \\
 &\lesssim x_\theta^{\bar r } +\int_{x_\theta}^\infty x^{\bar r-1} \bbP(|F_{\theta,n}^{-1}(1-\delta)| \geq x)dx \\
 &\leq x_\theta^{\bar r} +2 \int_{x_\theta}^\infty x^{\bar r-1} \left(1 +\frac{x^2\delta}{2m_\theta}\right)^{-\frac{n\delta}{4}} dx\\
 &\leq x_\theta^{\bar r} +2 \int_{\sqrt{m_\theta/\delta}}^\infty x^{\bar r-1} \left(1 +\frac{x^2\delta}{2m_\theta}\right)^{-\frac{n\delta}{4}} dx,
 & (\text{Since } \ell > 1)\\
 &= x_\theta^{\bar r} +2\left(\frac{m_\theta}{\delta}\right)^{\frac{\rbar}{2}} \int_{1}^\infty y^{\bar r-1} \left( 1+\frac{y^2}{2}\right)^{-\frac{n\delta}{4}} dy\\ 
 &\lesssim x_\theta^{\bar r} +\left(\frac{m_\theta}{\delta}\right)^{\frac{\rbar}{2}} \int_{1}^\infty y^{\bar r-1-\frac{n\delta}{2}} dy\\  
 &= x_\theta^{\bar r} +\left(\frac{m_\theta}{\delta}\right)^{\frac{\rbar}{2}} \frac 1 {\frac{n\delta}{2} - \rbar}
\leq  x_\theta^{\bar r} +\left(\frac{m_\theta}{\delta}\right)^{\frac{\rbar}{2}},
 \end{align*}
 where we used the assumptions $\delta \geq 2(r+2)/n$ and $\rbar \leq r+1$ on the final line of the above display. 
Therefore,   $\bbE|F_{\theta,n}^{-1}(1-\delta)|^{\bar r } \lesssim x_\theta^{\bar r} + (m_\theta/\delta)^{\frac {\rbar}{2}}.$
 Upon repeating the same argument for the $\delta$-quantile, we obtain
$$\max_{a \in \{\delta,1-\delta\}} \bbE|F_{\theta,n}^{-1}(a)|^{\bar r}  \lesssim  x_\theta^{\bar r} +\left(\frac{m_\theta}{\delta}\right)^{\frac{\rbar}{2}}
\lesssim 1+\max_{a \in \{\frac\delta 2, 1-\frac \delta 2\}} |F_\theta^{-1}(a)|^{\rbar} + \left(\frac{m_\theta}{\delta}\right)^{\frac {\rbar}{2}}
.$$
This proves the claim. 
\qed

\subsection{Proof of Lemma~\ref{lem::high_prob_sw_null}}
\label{app::high_prob_sw_null}
Since $\sj_{r,\delta/2}(P)  < \infty$, it follows from  
Lemma \ref{lem::abs_cont}
that 
$F_\theta^{-1}$ is absolutely continuous over $[\delta/2,1-\delta/2]$ for $\mu$-almost every $\theta \in \bbS^{d-1}$. 
We fix any such  $\theta$ throughout the sequel.

To prove the claim, we shall make use of the following analogue of the relative VC inequality described in Example~\ref{ex::rel_VC}
\citep{bousquet2003}. 
For any given $\theta \in \bbS^{d-1}$ and $\epsilon\in(0,1)$, we have
$$\bbP\left(|F_{\theta,n}(x) - F_\theta(x)| \leq \nu_{\epsilon,n} \sqrt{F_\theta(x)(1-F_\theta (x))} , \ \forall x \in \bbR\right) \geq 1-\epsilon.$$
Notice here that the right-hand side of the inequality within the above probability involves the population
CDF $F_\theta$ rather than $F_{\theta,n}$. Similarly as in Example~\ref{ex::rel_VC} and Appendix~\ref{app::secondary}, the above
bound implies that
$$\bbP\Big(F_\theta^{-1}(\gamma_{\epsilon,n}(u)) \leq F_{\theta,n}^{-1}(u) \leq F_\theta^{-1}(\eta_{\epsilon,n}(u)), \ \forall u \in (0,1)\Big) \geq 1-\epsilon,$$
where, for all $ u \in (0,1)$,
\begin{equation}
\label{eq::relVC_fns_pf}
\begin{aligned}
\gamma_{\epsilon,n}(u) &= \frac{2u + \nu_{\epsilon,n}^2-\nu_{\epsilon,n} \sqrt{\nu_{\epsilon,n}^2 + 4u(1-u)}}{2(1+\nu_{\epsilon,n}^2)}, \\
\eta_{\epsilon,n}(u) &= \frac{2u + \nu_{\epsilon,n}^2 + \nu_{\epsilon,n} \sqrt{\nu_{\epsilon,n}^2 + 4u(1-u)}}{2(1+\nu_{\epsilon,n}^2)}.
\end{aligned}
\end{equation}
These functions are invertible over $[0,1]$, with inverses given by
\begin{align}
\begin{aligned}
\label{eq::relVC_inverse_fns_pf}
\eta_{\epsilon,n}^{-1}(t) &= t - \nu_{\epsilon,n}\sqrt{t(1-t)},\quad t \in
\left[\eta_{\epsilon,n}(0), \eta_{\epsilon,n}(1)\right] =  \left[\frac{\nu_{\epsilon,n}^2}{1+\nu_{\epsilon,n}^2},1\right],\\
\gamma_{\epsilon,n}^{-1}(t) &= t + \nu_{\epsilon,n}\sqrt{t(1-t)}, \quad t \in [\gamma_{\epsilon,n}(0), \gamma_{\epsilon,n}(1)]=\left[0,\frac{1}{1+\nu_{\epsilon,n}^2}\right].
\end{aligned}
\end{align} 
We shall further make use of the following elementary Lemma, proven in Appendix~\ref{app::pf_lem_technical_prop1}.
\begin{lemma}
\label{lem::technical_prop1}
Assume $\epsilon \in (0,1)$ is chosen such that $\nu_{\epsilon,n} \leq y_0:=\sqrt{\delta}/4$. Then, for all $ u \in [\delta/2,1-\delta/2]$,
$$\frac{\gamma_{\epsilon,n}(u)}{u} \geq \frac 1 2, \quad \frac{1-\eta_{\epsilon,n}(u)}{1-u} \geq \frac 1 2.$$
In particular, for all $x \in [\delta/2,1-\delta/2]$, and all $y \in [\gamma_{\epsilon,n}(x), \eta_{\epsilon,n}(x)]$,
$\frac{x(1-x)}{y(1-y)} \leq \frac 1 4.$
\end{lemma}
By Lemma~\ref{lem::technical_prop1}, the inequalities $\gamma_{\epsilon,n}(\delta) \geq \delta/2$ 
and $\eta_{\epsilon,n}(1-\delta) \leq 1-\delta/2$ hold whenever $\nu_{\epsilon,n} \leq y_0$,
and this last inequality is satisfied whenever $\epsilon \geq \epsilon_0 := \frac{2n+1}{16}\exp(-ny_0^{2}/16)$.
For all such $\epsilon$, define the event
$$E_\epsilon\equiv E_\epsilon(\theta) = 
\Big\{ F_\theta^{-1}(\gamma_{\epsilon,n}(u)) \leq F_{\theta,n}^{-1}(u) \leq F_\theta^{-1}(\eta_{\epsilon,n}(u)), \ \forall u \in [\delta, 1-\delta]\Big\}.$$
%defined for any $\epsilon$ such that 
%\begin{equation} 
%\label{eq::gamma_eta_range_pf_prop1}
%\eta_{\epsilon,n}(1-\delta) \leq 1-\delta/2,\quad \gamma_{\epsilon,n}(u)\geq \delta/2.
%\end{equation}
%Notice that
Over the event $E_\epsilon$, we have
\begin{align*}
Z_n(\theta)
 &=\frac 1 {1-2\delta}\int_\delta^{1-\delta} \big| F_{\theta,n}^{-1}(u) - F_{\theta}^{-1}(u)\big|^r du \\
 &\leq \frac 1 {1-2\delta} \int_\delta^{1-\delta} \big[ F_{\theta}^{-1}(\eta_{\epsilon,n}(u)) - F_{\theta}^{-1}(\gamma_{\epsilon,n}(u))\big]^r du.
\end{align*}
Now, recall that $\theta$ was chosen such that $F_\theta^{-1}$ is absolutely continuous over $[\delta/2,1-\delta/2]$, whence
%there exists $\zeta_{\epsilon,n}(u) \in [\gamma_{\epsilon,n}(u),\eta_{\epsilon,n}(u)]$,
$$F_{\theta}^{-1}(\eta_{\epsilon,n}(u)) - F_{\theta}^{-1}(\gamma_{\epsilon,n}(u)) = 
\int_{\gamma_{\epsilon,n}(u)}^{\eta_{\epsilon,n}(u)} \frac {dt} {p_\theta(F_\theta^{-1}(t))}.
$$  
Now, since $\nu_{\epsilon,n} \leq y_0$, we have for all $u\in[\delta,1-\delta]$,
\begin{align}
\label{eq::pf_sqrt_bound}
\nonumber \eta_{\epsilon,n}(u)-\gamma_{\epsilon,n}(u)
 &=   \frac{ \nu_{\epsilon,n} \sqrt{\nu_{\epsilon,n}^2 + 4u(1-u)}}{(1+\nu_{\epsilon,n}^2)} \\
 &\leq 2\nu_{\epsilon,n}\sqrt{u(1-u)} + \nu_{\epsilon,n}^2
 \leq 3\nu_{\epsilon,n}\sqrt{u(1-u)}.
\end{align}
Using Jensen's inequality, we deduce that, over the event $E_\epsilon$,
\begin{align*}
(1&-2\delta)Z_n(\theta) \\
 \hspace{-0.03in}&\leq \int_\delta^{1-\delta} 
 \left(\int_{\gamma_{\epsilon,n}(u)}^{\eta_{\epsilon,n}(u)}\frac{1}{p_\theta(F_\theta^{-1}(t))}dt\right)^rdu  \\
 &\leq  \int_\delta^{1-\delta}  {(\eta_{\epsilon,n}(u) - \gamma_{\epsilon,n}(u))^{r-1}}
 \int_{\gamma_{\epsilon,n}(u)}^{\eta_{\epsilon,n}(u)}\left(\frac{1}{p_\theta(F_\theta^{-1}(t))}\right)^rdtdu \\
 &=  \int_\delta^{1-\delta} \frac {(\eta_{\epsilon,n}(u) - \gamma_{\epsilon,n}(u))^{r-1}}{(u(1-u))^{r/2}} 
 \int_{\gamma_{\epsilon,n}(u)}^{\eta_{\epsilon,n}(u)}\left(\frac{\sqrt{t(1-t)}}{p_\theta(F_\theta^{-1}(t))}\right)^r \frac {(u(1-u))^{r/2}} {(t(1-t))^{r/2}} dtdu \\
 &\lesssim \nu_{\epsilon,n}^{r-1}  \int_\delta^{1-\delta}  
 \frac 1 {\sqrt{u(1-u)}} 
 \int_{\gamma_{\epsilon,n}(u)}^{\eta_{\epsilon,n}(u)}\left(\frac{\sqrt{t(1-t)}}{p_\theta(F_\theta^{-1}(t))}\right)^r  dtdu
   ~~\text{(By \eqref{eq::pf_sqrt_bound}, Lem.~\ref{lem::technical_prop1}) }\\
 &\leq \nu_{\epsilon,n}^{r-1} \int_{\delta/2}^{1-\delta/2}\left(\int_{\eta^{-1}_{\epsilon,n}(t)}^{\gamma_{\epsilon,n}^{-1}(t)} 
 \frac 1 {\sqrt{u(1-u)}} du\right)
  \left(\frac{\sqrt{t(1-t)}}{p_\theta(F_\theta^{-1}(t))}\right)^r   dt,
 \end{align*}
where we interchanged the order of integration, and used the fact that
$\delta/2 \leq \gamma_{\epsilon,n}(u) \leq \eta_{\epsilon,n}(u)\leq 1-\delta/2$ for all $u \in [\delta,1-\delta]$.
We further have
$$\int_{\eta^{-1}_{\epsilon,n}(t)}^{\gamma_{\epsilon,n}^{-1}(t)} 
 \frac 1 {\sqrt{u(1-u)}} du \leq \frac {\gamma_{\epsilon,n}^{-1}(t) - \eta^{-1}_{\epsilon,n}(t) } {\sqrt{u_t^*(1-u_t^*)}},\quad
 \text{for }  u_t^* \in \argmin_{\eta_{\epsilon,n}^{-1}(t) \leq u \leq \gamma_{\epsilon,n}^{-1}(t)} \sqrt{u(1-u)}.$$
By again applying Lemma~\ref{lem::technical_prop1}, and equation~\eqref{eq::relVC_inverse_fns_pf}, we obtain
 $$\frac {\gamma_{\epsilon,n}^{-1}(t) - \eta^{-1}_{\epsilon,n}(t) } {\sqrt{u_t^*(1-u_t^*)}} 
\leq \frac {\nu_{\epsilon,n} \sqrt{t(1-t)}} {\sqrt{u_t^*(1-u_t^*)}}  \leq \frac 1 2  \nu_{\epsilon,n}.$$
We have thus shown
\begin{align*}
Z_n(\theta)
 &\lesssim \frac{\nu_{\epsilon,n}^{r}}{1-2\delta}  \int_{\delta/2}^{1-\delta/2}
  \left(\frac{\sqrt{t(1-t)}}{p_\theta(F_\theta^{-1}(t))}\right)^r  dt = 
  \nu_{\epsilon,n}^{r}  J_{r,\delta/2}(P_\theta).
\end{align*}
Thus, setting $\epsilon = \frac{2n+1}{16}e^{-\frac{ny^{2}}{16}}$ for any $y \in (0,y_0]$, 
we have
$$\bbP \Big(Z_n(\theta) \geq C_r y^r J_{r,\delta/2}(P_\theta) \Big)\leq \epsilon,$$
for a universal constant $C_r > 0$ depending only on $r$.\qed

\subsection{Proof of Lemma~\ref{lem::technical_prop1}} 
\label{app::pf_lem_technical_prop1}
For all $u\in[\delta/2,1-\delta/2]$, and all $\epsilon$ such that $\nu_{\epsilon,n}^2 \leq \delta/16$, we have
\begin{align*}
\frac{\gamma_{\epsilon,n}(u)}{u}
 &\geq  \frac 1 {1 + \nu_{\epsilon,n}^2} - \frac{\nu_{\epsilon,n} \sqrt{\nu_{\epsilon,n}^2 + 4u(1-u)} - \nu_{\epsilon,n}^2}{2u(1+\nu_{\epsilon,n}^2)} \\
 &\geq  \frac 1 {1 + \nu_{\epsilon,n}^2} - \frac{\nu_{\epsilon,n} \sqrt{1-u}}{\sqrt u(1+\nu_{\epsilon,n}^2)} \\
 &\geq  \frac 1 {1 + \nu_{\epsilon,n}^2} - \frac{\sqrt 2 \sqrt \delta \sqrt{1-u}}{4\sqrt \delta (1+\nu_{\epsilon,n}^2)} \\
 &\geq  \frac 1 {1 + \nu_{\epsilon,n}^2} - \frac{1 }{2\sqrt 2  (1+\nu_{\epsilon,n}^2)} \\
 &\geq  \frac{2\sqrt 2 -1}{2\sqrt 2 (1+\nu_{\epsilon,n}^2)} \geq \frac{2\sqrt 2 -1}{2\sqrt 2 (1+1/16)} \geq 1/2.
\end{align*}
Similarly,  
\begin{align*}
\frac{1-\eta_{\epsilon,n}(u)}{1-u}
 &\geq  \frac{1}{1+\nu_{\epsilon,n}^2} - \frac{\nu_{\epsilon,n}^2 + \nu_{\epsilon,n}\sqrt{\nu_{\epsilon,n}^2+4u(1-u)}}{2(1-u)(1+\nu_{\epsilon,n}^2)} \\
 &\geq  \frac{1}{1+\nu_{\epsilon,n}^2} - \frac{\nu_{\epsilon,n}^2+\nu_{\epsilon,n}\sqrt{u(1-u)}}{(1-u)(1+\nu_{\epsilon,n}^2)} \\
 &=  \frac{1 - \nu_{\epsilon,n}^2/(1-u)}{1+\nu_{\epsilon,n}^2} - \frac{\nu_{\epsilon,n}\sqrt{u }}{\sqrt{1-u}(1+\nu_{\epsilon,n}^2)} \\
 &\geq  \frac{1 - \nu_{\epsilon,n}^2/\delta}{1+\nu_{\epsilon,n}^2} - \frac{\sqrt 2 \nu_{\epsilon,n}}{\sqrt{\delta}(1+\nu_{\epsilon,n}^2)} \\
 &\geq  \frac{1 - 1/16}{1+\nu_{\epsilon,n}^2} - \frac{1}{2\sqrt 2 (1+\nu_{\epsilon,n}^2)} \\
 &\geq  \frac{2\sqrt 2(1 - 1/16)-1}{2\sqrt 2(1+1/16)} \geq 1/2.
\end{align*}
In particular, for all $x \in [\delta/2,1-\delta/2]$, and all $y \in [\gamma_{\epsilon,n}(x), \eta_{\epsilon,n}(x)]$,
$$\frac{x(1-x)}{y(1-y)} \leq \frac{x}{\gamma_{\epsilon,n}(x)} \frac{1-x}{1-\eta_{\epsilon,n}(x)}
 \leq \frac 1 2.$$ 
The claim follows.
\qed

\subsection{Proof of Lemma~\ref{lem::rth-deviation}}
When $r > 1$, we have,
\begin{align*}
(1-2\delta)&\sw_{r,\delta}^r(P_n,Q_m) \\
 &= \int_{\bbS^{d-1}}
 \int_\delta^{1-\delta} \big|F_{\theta,n}^{-1}(u) - G_{\theta,m}^{-1}(u)\big|^r du d\mu(\theta)\\
 &=  \int_{\bbS^{d-1}} 
 \int_\delta^{1-\delta} \bigg\{
 	|F_\theta^{-1}(u) - G_\theta^{-1}(u)|^r   + 
 	r\sgn(\tilde F_{\theta,n}^{-1}(u) - \tilde G_{\theta,m}^{-1}(u))\\ &  \hspace{-0.52in} \times  |\tilde F_{\theta,n}^{-1}(u) - \tilde G_{\theta,m}^{-1}(u)|^{r-1}\big\{
 		(F_{\theta,n}^{-1}(u) - F_\theta^{-1}(u)) {-} (G_{\theta,m}^{-1}(u) - G_\theta^{-1}(u))\big\} 
 \bigg\}dud\mu(\theta),
\end{align*}
by a Taylor expansion of the map $(x,y) \mapsto |x-y|^r$ about $\big(F_\theta^{-1}(u), G^{-1}_\theta(u)\big)$, where 
$\tilde F_{\theta,n}^{-1}(u)$ (resp. $\tilde G_{\theta,m}^{-1}(u)$) is a real number on the line joining
$F_\theta^{-1}(u)$ and $F_{\theta,n}^{-1}(u)$ (resp. $G_\theta^{-1}(u)$ and $G_{\theta,m}^{-1}(u)$). 
We then have
\begin{align*} 
\nonumber 
|D_{nm}| &\leq   
 \frac r {1-2\delta} \int_{\bbS^{d-1}}\int_\delta^{1-\delta} |\tilde F_{\theta,n}^{-1}(u) - \tilde G_{\theta,m}^{-1}(u)|^{r-1}
    \\ &\qquad\qquad\qquad\quad \times \Big[|F_{\theta,n}^{-1}(u) - F_\theta^{-1}(u)| + |G_{\theta,m}^{-1}(u) - G_\theta^{-1}(u)|\Big]dud\mu(\theta) \\
    \nonumber
 &\leq   
 \frac r {1-2\delta} \int_{\bbS^{d-1}}  (b_\theta - a_\theta)^{r-1} \\ &\qquad\qquad\quad \times \int_\delta^{1-\delta} 
    \Big[|F_{\theta,n}^{-1}(u) - F_\theta^{-1}(u)| + |G_{\theta,m}^{-1}(u) - G_\theta^{-1}(u)|\Big]dud\mu(\theta)\\
 &= \frac r {1-2\delta} \int_{\bbS^{d-1}} Z_{nm}(\theta)d\mu(\theta).
\end{align*}
This proves the first claim when $r > 1$, and the same conclusion
holds trivially when $r=1$ by the triangle inequality. To prove the second claim,
given $\theta \in \bbS^{d-1}$,  define the following event
for any $t \in (0,\delta/2]$,
$$\begin{multlined}[0.95\textwidth]
E_t \equiv E_t(\theta) =  \Big\{ F_{\theta,n}^{-1}(u-t) \leq F_\theta^{-1}(u)\leq  F_{\theta,n}^{-1}(u+t) , \ \ \forall u \in [\delta,1-\delta]\Big\}\\
  \cap \Big\{G_{\theta,m}^{-1}(u-t) \leq G_\theta^{-1}(u) \leq G_{\theta,m}^{-1}(u+t), \ \ \forall u \in [\delta,1-\delta]\Big\}.
\end{multlined}$$
A union bound together with the Dvoretzky-Kiefer-Wolfowitz 
inequality (Example \ref{ex::DKW}) implies that, for all $t \in (0,\delta/2]$,  
$\bbP(E_t)\geq  1-4\exp(-2(n\wedge m) t^2).$
%Therefore, with probability at least $1-4\exp(-2(n\wedge m) t^2)$, for all $t \in (0,\delta/2]$,
Now, for all such $t$, the following inequalities hold over $E_t$,
\begin{align*}
\int_\delta^{1-\delta}&|F_{\theta,n}^{-1}(u) - F_\theta^{-1}(u)|du \\
 &\leq \int_\delta^{1-\delta} \big[F_\theta^{-1}(u+t) - F_\theta^{-1}(u-t)\big]du \\
 &= \int_{\delta-t}^{1-\delta-t} F_\theta^{-1}(u)du - \int_{\delta+t}^{1-\delta+t} F_\theta^{-1}(u) du \\
  &= \int_{\delta-t}^{\delta+t} F_\theta^{-1}(u)du + \int_{1-\delta-t}^{1-\delta+t} F_\theta^{-1}(u)du\\
 &\leq 2 t \big[|F_\theta^{-1}(\delta-t)| + |F_\theta^{-1}(\delta+t)| + |F_\theta^{-1}(1-\delta + t)| + |F_\theta^{-1}(1-\delta-t)|\big]\\
 &\leq 2 t \big[|F_\theta^{-1}(\delta/2)| + |F_\theta^{-1}(1-\delta/2)|\big].
\end{align*}
Over $E_t$, we also have for all $t \in (0,\delta/2]$,
\begin{align*}
\int_\delta^{1-\delta}|G_{\theta,n}^{-1}(u) - G_\theta^{-1}(u)|du 
 &\leq 2 t \big[|G_\theta^{-1}(\delta/2)| + |G_\theta^{-1}(1-\delta/2)|\big],
\end{align*}
and,
\begin{align*}
a_\theta &\geq F_\theta^{-1}(\delta - t)\wedge G_\theta^{-1}(\delta - t)
    	   \geq F_\theta^{-1}(\delta/2)\wedge G_\theta^{-1}(\delta/2), \\[0.08in]
b_\theta &\leq F_\theta^{-1}\big(1-\delta+ t\big)\vee G_\theta^{-1}(1-\delta+t)
\leq F_\theta^{-1}(1-\delta/2)\vee G_\theta^{-1}(1-\delta/2).
\end{align*}
Combining these facts,
we deduce that for a universal constant $A > 0$, 
we have with probability at least $1-4\exp(-2(n\wedge m) t^2)$,
$$|Z_{nm}(\theta)| \leq   A t M^r(\theta),$$
as was to be shown. \qed 
 
\section{Proof of Theorem \ref{thm::minimax_distance}}
\label{app::main_lowerbound}  
Throughout the proof, $\KL$ denotes the Kullback-Leibler divergence, and $\chi^2$ denotes the
$\chi^2$-divergence. % \citep{tsybakov2008}. 
In view of the identity $W_{r,\delta}(P,Q) = W_r(P^\delta, Q^\delta)$ stated in Section \ref{sec::background_wasserstein},
and its natural analogue for the Sliced Wasserstein distance, 
together with the fact that all distributions considered below are
compactly supported, there will be no loss of generality in assuming $\delta=0$ in what follows. 
 
At a high-level, our general approach is to carefully construct two pairs of distributions $(P_0,Q_0),(P_1,Q_1) \in \calO(\Gamma;s_1,s_2)$
such
that the corresponding product measures $(P_0^{\otimes n} \otimes Q_0^{\otimes m})$ and $(P_1^{\otimes n} \otimes Q_1^{\otimes m})$
are close in the $\KL$ distance, but such that $\sw_r(P_0,Q_0)$ and $\sw_r(P_1,Q_1)$ are sufficiently different. In particular, if we can 
show that 
$\KL\big(P_0^{\otimes n} \otimes Q_0^{\otimes m}, P_1^{\otimes n} \otimes Q_1^{\otimes m}\big) \leq \zeta < \infty,$ then
via an application of Le Cam's inequality (see for instance, Theorem 2.2 of \citet{tsybakov2008}), we obtain the minimax lower bound that,
\begin{align}
\label{eqn::lecam}
\calR_{nm}(\calO(\Gamma;s_1,s_2);r) \geq c_{\zeta} |\sw_r(P_0,Q_0) - \sw_r(P_1,Q_1)|,
\end{align}
where $c_{\zeta} > 0$ is a constant depending only on $\zeta$. We will use four separate constructions to handle various cases of the Theorem.

Let $\epsilon_n = k_r n^{-1/2}$, for a constant $k_r \in (0,1)$, possibly depending on $r$, 
to be determined below. We use the following pairs of distributions.
\begin{itemize}
\item {\bf Construction 1.} For a vector $A = (a, 0, \dots, 0) \in \bbR^d$, and for 
$g > 0$, 
we define:
\begin{alignat*}{2}
P_{01} &= \frac 1 2 \delta_0 + \frac 1 2 \delta_A, \quad  && Q_{01} = \frac 1 2 \delta_{g A } + \frac 1 2 \delta_{(1+g)A}\\
P_{11} &= \left(\frac 1 2 + \epsilon_n\right) \delta_0 + \left(\frac 1 2 - \epsilon_n\right)\delta_A,\quad &&Q_{11} = \frac 1 2 \delta_{g A } + \frac 1 2 \delta_{(1+g)A}.
\end{alignat*}
  
\item {\bf Construction 2.} 
%For $0 < s_1 \leq s_2$ 
For $\gamma_2 , \Delta> 0$ to be chosen in the sequel we let $P_{02}$, $P_{12}$, $Q_{02}$, $Q_{12} \in \calP(\bbR^d)$ be the
probability distributions 
of random vectors of the form $(X, 0, \dots, 0) \in \bbR^d$, with $X$ respectively distributed
according to the distributions
\end{itemize}
\begin{alignat*}{3}
\quad P_{02}^{(1)} &= 
U\left(0, \gamma_2^{1/ r}\right), \quad && Q_{02}^{(1)} = U\left(\Delta\gamma_2^{\frac 1 r}, (1+\Delta)\gamma_2^{\frac 1 r}\right)\\
P_{12}^{(1)}  &= \frac{1{+}\epsilon_n}{2} U\Big(0, \frac{\gamma_2^{\frac  1 r}}{2}\Big) + 
                 \frac{1{-}\epsilon_n}{2} U\Big(\frac{\gamma_2^{\frac 1 r}}{2}, \gamma_2^{\frac  1 r}\Big)
 ,\quad 
&&Q_{12}^{(1)} = U\left(\Delta \gamma_2^{\frac 1 r},(1+\Delta)\gamma_2^{\frac 1 r}\right).
\end{alignat*}
\begin{itemize}
\item {\bf Construction 3.} For $0 < s_1 \leq s_2$ we let $P_{03}$,~$P_{13}$,  $Q_{03}$, $Q_{13} \in \calP(\bbR^d)$ be the
probability distributions 
of random vectors of the form $(X, 0, \dots, 0) \in \bbR^d$, with $X$ respectively distributed
according to the distributions
\begin{align*}
P_{03}^{(1)} &= U\left(0, s_1^{1/r}\right), \qquad Q_{03}^{(1)} = U\left(0, s_2^{1/r}\right),\\
P_{13}^{(1)} &= U\left(0, s_1^{1/r}\right), \qquad Q_{13}^{(1)} = (1-\epsilon_m) U\left(0, s_2^{1/r}\right) + \epsilon_m\delta_{s_2^{1/r}}.
\end{align*} 

\item {\bf Construction 4. } For $0 < s_2 \leq s_1$ we let $P_{04}, P_{14}, Q_{04}, Q_{14} \in \calP_r(\bbR^d)$ be the
probability distributions 
of random vectors of the form $(X, 0, \dots, 0) \in \bbR^d$, with $X$ respectively distributed
according to the distributions
\begin{alignat*}{2}
P_{04}^{(1)} &= U\left(0, s_1^{1/r}\right), \qquad &&Q_{04}^{(1)} = U\left(0, s_2^{1/r}\right),\\
P_{14}^{(1)} &= (1-\epsilon_n)U\left(0, s_1^{1/r}\right) + \epsilon_n \delta_{s_1^{1/r}}, \qquad &&Q_{14}^{(1)} = U\left(0, s_2^{1/r}\right).
\end{alignat*} 
\end{itemize}
Construction 1 uses pairs of distributions with infinite $\sj_r$, while Constructions 2-4 use pairs of distributions with finite $\sj_r$. 
To compactly state our next result we define several terms,
\begin{alignat*}{2}
%I_r &:= \left(\int_{\bbS^{d-1}} |A^\top\theta|^r d\mu(\theta)\right)^{\frac 1 r} \\
t_r &:= \left(\int_{\bbS^{d-1}} |\theta_1|^r d\mu(\theta)\right)^{\frac 1 r}, \\
\beta &:= (s_2/s_1)^{1/r},   &&\bar\beta:= 1/\beta,\\
\Delta_\beta &:= \beta-1  ,   && \Delta_{\bar \beta}  := \bar\beta - 1. 
\end{alignat*}
With these definitions in place the following technical lemma describes the main features of our constructions.
\begin{lemma}
\label{lem:main}
There exists a choice of constant $k_r \in (0,1)$ for which the following statements hold.
\begin{itemize}
\item {\bf Construction 1.} Let $g := \Gamma / \left(\int_{\bbS^{d-1}} |A^\top\theta|^r d\mu(\theta)\right)^{1 /r}$. Then, there exists a constant $c_1 > 0$, possibly
depending on $r$, such that
\begin{align*}
\sw_r^r(P_{01}, Q_{01}) &= \Gamma^r,\\ 
\sw_r^r(P_{11},Q_{11}) &\geq \Gamma^r + c_1\epsilon_n.%~~~\text{for some } c_1 > 0.
%\KL(P_{01}^{\otimes n} \otimes Q_{01}^{\otimes m}, (P_{11}^{\otimes n} \otimes Q_{11}^{\otimes m})) &\leq \zeta < \infty. 
\end{align*}
Furthermore, there exists a choice of the vector $A$ for which $P_{01}$, $Q_{01}$, $P_{11}$, $Q_{11} \in \mathcal{O}(\Gamma; \infty,\infty)$. 
\item {\bf Construction 2.} 
%We have that $\sj_r(P_{03}), \sj_r(P_{13}), \sj_r(Q_{03}), \sj_r(Q_{13}) \leq \gamma_2,$ and
There exists a constant $c_2 > 0$, possibly depending on $r$,  such that
\begin{align*}
\sw_r^r(P_{02}, Q_{02}) &= \gamma_2 ( t_r\Delta)^r,\\ 
\sw_r^r(P_{12},Q_{12}) &\geq \gamma_2 t_r^r \Big\{ \Delta^r + c_2 \Delta^{r-1} \epsilon_n\Big\}.
%\KL(P_{01}^{\otimes n} \otimes Q_{01}^{\otimes m}, (P_{11}^{\otimes n} \otimes Q_{11}^{\otimes m})) &\leq \zeta < \infty. 
\end{align*}
Furthermore, $P_{02}, Q_{02}, P_{12}, Q_{12} \in \mathcal{O}(0; \gamma_2, \gamma_2)$. 
 \item {\bf Construction 3.} %We choose $\epsilon = n^{-1/2}$. 
Assume that $\bar\beta \in (0,1]$. Then,% We have that $\sj_r(P_{03}), \sj_r(P_{13}) \leq s_1$ and $\sj_r(Q_{03}), \sj_r(Q_{13}) \leq s_2,$ and
\begin{align*}
 \sw_r^r(P_{03}, Q_{03}) &= \frac{s_2 t_r^r|\Delta_{\bar\beta}|^r}{r+1}, \\
 \sw_r^r(P_{13}, Q_{13}) &\geq \frac{t_r^r s_2}{r+1} \Big\{ |\Delta_{\bar\beta}|^r  + r |\Delta_{\bar\beta}|^{r-1}\epsilon_m\Big\}.
 \end{align*} 
  Furthermore, $P_{03}, Q_{03}, P_{13}, Q_{13} \in \mathcal{O}(0; s_1, s_2)$. 
\item {\bf Construction 4.} %We choose $\epsilon = m^{-1/2}$. 
Assume that $\beta \in (0,1]$. Then, %We have that $\sj_r(P_{04}), \sj_r(P_{14}) \leq s_1$ and $\sj_r(Q_{04}), \sj_r(Q_{14}) \leq s_2,$ and
\begin{align*}
 \sw_r^r(P_{04}, Q_{04}) &= \frac{s_1 t_r^r|\Delta_\beta|^r}{r+1}, \\
 \sw_r^r(P_{14},Q_{14}) &\geq \frac{t_r^r s_1}{r+1} \Big\{|\Delta_\beta|^r  + r |\Delta_\beta|^{r-1}\epsilon_n\Big\}.
 \end{align*}
  Furthermore, $P_{04}, Q_{04}, P_{14}, Q_{14} \in \mathcal{O}(0; s_1, s_2)$. 
% \item 
\end{itemize}
 In each case, for some fixed universal constant $\zeta > 0$ we have that,
\begin{align*}
\KL\big(P_{0i}^{\otimes n} \otimes Q_{0i}^{\otimes m}, P_{1i}^{\otimes n} \otimes Q_{1i}^{\otimes m}\big) &\leq \zeta < \infty, \quad i =1, 2, 3, 4.
\end{align*} 
\end{lemma}
Taking this result as given, we can now complete the proof of the Theorem. Using Construction~1 with $\Gamma = 0$, we obtain
from~equation \eqref{eqn::lecam} that
\begin{align*}
\calR_{nm}(\calO(0; \infty,\infty);r) \geq  c_{\zeta} |\sw_r(P_{01},Q_{01}) - \sw_r(P_{11},Q_{11})| \gtrsim \epsilon_n^{1/r} \asymp n^{-1/2r}.
\end{align*}
Reversing the roles of $n$ and $m$ we obtain the first claim of part (i) of the Theorem. Choosing $\Gamma$ to be a strictly positive constant, we 
instead obtain
\begin{align*}
\calR_{nm}(\calO(\Gamma; \infty, \infty);r) 
 \gtrsim \Gamma \left| 1 - \left(1 + \frac{c_1 \epsilon_n}{\Gamma^r}\right)^{\frac 1 r}\right| 
 = \frac{c_1k_r n^{-1/2}}{r\Gamma^{r-1}}(1+o(1)),
\end{align*}
by a first-order Taylor expansion of the map $x \mapsto (1+x)^{\frac 1 r}$. The fact that $\Gamma$ is bounded away from zero then implies
$\calR_{nm}(\calO(\Gamma; \infty,\infty); r) \gtrsim n^{-1/2}$
which proves part (ii) of the theorem, again upon reversing the roles of $n$ and $m$. 
It thus only remains to prove the second claim of part (i).

Without loss of generality we assume that $n \leq m$ in the remainder of the proof, noting that, as above, we may always reverse the roles of $n$ and $m$ and repeat our constructions.  
We consider four cases.

\paragraph{ Case 1:  $-1 \leq \Delta_\beta \leq - \epsilon_n$.} 
 In this case, the condition $\Delta_\beta \leq 0$ implies $s_1 \geq s_2$. 
Since $n \leq m$, it therefore suffices to prove
$\calR_{nm}(\calO(0, s_1, s_2);r)
 \gtrsim s_1^{1/r}\epsilon_n.$ Furthermore, since $\beta \leq 1$, we may invoke Construction~4 to obtain
$$|\sw_r(P_{04}, Q_{04}) - \sw_r(P_{14}, Q_{14})| 
  \geq \frac {s_1^{1/r}t_r|\Delta_\beta|}{(r+1)^{1/r}} \left[\left(1 + \frac{r\epsilon_n}{|\Delta_\beta|}\right)^{\frac 1 r} - 1\right]
  \asymp s_1^{1/r} \epsilon_n,
$$
where we have used the assumption $|\Delta_\beta| \geq \epsilon_n$ in the last order assesment of the above display.
This fact together with equation~\eqref{eqn::lecam} yields the desired lower bound for Case~1.

\paragraph{Case 2:  $-\epsilon_n < \Delta_\beta \leq 0$.} 
The inequality $s_1 \geq s_2$ continues to hold, thus it suffices
to prove $\calR_{nm}(\calO(0; s_1, s_2);r)
 \gtrsim s_1^{1/r}\epsilon_n.$ 
 Notice further that
 $$s_2^{1/r}\epsilon_n = s_1^{1/r} \beta \epsilon_n > s_1^{1/r}(1-\epsilon_n)\epsilon_n = s_1^{1/r}\epsilon_n(1+o(1)).$$
It will therefore suffice to prove $\calR_{nm}(\calO(0; s_1, s_2);r)
 \gtrsim s_2^{1/r}\epsilon_n.$ 
We use Construction 2, and choose $\gamma_2 = s_2$, and $\Delta \in (0,1]$ to be a constant
larger than $\epsilon_n$. We observe that all distributions have $\sj_r$ value at most $s_2 = \min\{s_1,s_2\}$.
Furthermore, 
\begin{align*}
\big|\sw_r(P_{02},Q_{02}) - \sw_r(P_{12},Q_{12})\big|
 &\geq s_2^{1/r}  t_r \Delta\left\{ \left( 1 + \frac {c_2\epsilon_n} \Delta \right)^{\frac 1 r} - 1  \right\}.
\end{align*}
Since $\Delta \geq \epsilon_n$, it is a straightforward observation that the factor in braces of the above display is of order $\epsilon_n$, 
thus we have shown
\begin{align*}
\big|\sw_r(P_{02},Q_{02}) - \sw_r(P_{12},Q_{12})\big| \gtrsim s_2^{1/r} \epsilon_n,
\end{align*}
and this together with equation~\eqref{eqn::lecam} yields the desired lower bound for Case 2.

\paragraph{Case 3: $-1\leq\Delta_{\bar\beta} \leq -\epsilon_m $ and $s_1^{1/r}\epsilon_n \leq s_2^{1/r} \epsilon_m$.}
In this case, it suffices to prove that $\calR_{nm}(\calO(0, s_1, s_2);r)
 \gtrsim s_2^{1/r}\epsilon_m.$ Notice that $\bar\beta \leq 1$, hence 
we may use Construction 3 to obtain
$$|\sw_r(P_{03}, Q_{03}) - \sw_r(P_{13}, Q_{13})| 
  \geq \frac {s_2^{1/r}t_r |\Delta_{\bar\beta}|}{(r+1)^{1/r}} \left[\left(1 + \frac{\epsilon_m}{|\Delta_{\bar\beta}|}\right)^{\frac 1 r} - 1\right]
  \asymp s_2^{1/r} \epsilon_m,
$$
where we have used the assumption $|\Delta_{\bar\beta}| \geq \epsilon_m$ in the last order assesment of the above display.
This fact together with equation~\eqref{eqn::lecam} yields the desired lower bound for Case 1.
 
\paragraph{Case 4: $-\epsilon_m < \Delta_{\bar\beta} < 0$ or $s_1^{1/r}\epsilon_n > s_2^{1/r} \epsilon_m$.}
Notice that if the condition $\Delta_{\bar\beta} > -\epsilon_m$ is satisfied, it implies
$$s_1^{1/r}\epsilon_n = s_2^{1/r} \bar\beta\epsilon_n > (1-\epsilon_m)\epsilon_ns_2^{1/r}  \geq (1-\epsilon_m)\epsilon_m s_2^{1/r}
\gtrsim \epsilon_m s_2^{1/r}.$$
For this case, it will thus suffice to prove $\calR_{nm}(\calO(0; s_1, s_2);r)
 \gtrsim s_1^{1/r}\epsilon_n.$ Since $\Delta_{\bar\beta} \leq  0$, we observe that all distributions have $\sj_r$ 
 value at most $s_1 = \min\{s_1,s_2\}$. 
Invoking Construction 2 with $\gamma_2=s_1$, the remainder of the argument follows similarly as in Case 2.

It remains to establish Lemma~\ref{lem:main} and we turn our attention to this now.

\subsection{Proof of Lemma~\ref{lem:main}}
Bounding the KL divergence in each case is straightforward. We observe that for each $1 \leq i \leq 4$,
\begin{align*}
\KL(P_{0i}^{\otimes n} \otimes Q_{0i}^{\otimes m}, (P_{1i}^{\otimes n} \otimes Q_{1i}^{\otimes m})) &= n \KL(P_{0i}, P_{1i}) + m \KL(Q_{0i}, Q_{1i}) \\
&\leq n \chi^2 (P_{0i}, P_{1i}) + m \chi^2(Q_{0i}, Q_{1i}).
\end{align*}
The $\chi^2$ divergences in each construction can be computed in closed form. Doing so yields the bounds:
\begin{align*}
\KL\big(P_{0i}^{\otimes n} \otimes Q_{0i}^{\otimes m}, P_{1i}^{\otimes n} \otimes Q_{1i}^{\otimes m}\big) &\lesssim  n \epsilon_n^2 , \quad i=1,2,4\\
\KL\big(P_{03}^{\otimes n} \otimes Q_{03}^{\otimes m}, P_{13}^{\otimes n} \otimes Q_{13}^{\otimes m}\big) &\lesssim m \epsilon_m^2.
\end{align*}
Together with the definition of $\epsilon_n$, we obtain the desired bounds on the KL divergence.

 As a second preliminary, let us verify that, for appropriate choice of various parameters, the distributions we construct have appropriately bounded moments,
and belong to the class $\mathcal{K}_r(b)$ defined in equation~\eqref{eq::calK}.  
Notice first that the distributions $P_{01},Q_{01}, P_{11},Q_{11}$ have support with diameter bounded above by $(1+G)a$. 
Choosing $a$ (possibly depending on $G$ and hence $\Gamma$) such that 
this expression is bounded above by $b^{1/r}$ ensures $P_{01},Q_{01},P_{11},Q_{11} \in \calK_r(b)$. 
We are guaranteed that such a choice exists by using the assumption $\Gamma^r \leq c_r b$,
which ensures that $\Gamma$ cannot be too large.%, which in turn ensures that $\Delta$ cannot be too large.

Furthermore, the distributions $P_{ij}, Q_{ij}$ for $i=2,3,4$ and $j=0,1$ have supports with diameter bounded
above by $s (1+\Delta) \leq 2s$. The assumption $b \geq (2s)^{1/r}$ therefore guarantees 
$P_{ij}, Q_{ij} \in \calK_r(b)$ for $i=2,3,4$ and $j=0,1$.

We now consider each construction in turn, establishing the remaining claims. As a preliminary technical result, it will be useful to study the Wasserstein
distance between several pairs of univariate distributions.
\begin{lemma}
\label{lem::main_univariate}
\begin{enumerate}
\item Let $\Delta \geq \epsilon > 0$, and define the distributions
$$\nu = \frac{1+\epsilon}{2} U(0,1/2) + \frac{1-\epsilon}{2} U(1/2,1),$$
and $\rho = U(\Delta, 1+ \Delta)$. Then, 
$$W_r^r(\nu,\rho)  \geq \Delta^r + \frac{r}{4}\epsilon \Delta^{r-1} .$$
\item %Define $\Delta_\beta := \left(\frac{s_2}{s_1}\right)^{\frac 1 r} -1$.
Given $\xi \in (0,1]$, $\Delta_\xi = \xi-1$, define for all $\epsilon \in (0,1]$,
%$$\nu = \frac{1+\epsilon}{2} U(0,1/2) + \frac{1-\epsilon}{2} U(1/2,1),\quad \rho = U(0, \beta).$$
%Then, there exists a constant $c_1 > 0$ such that if $\Delta_\beta > 1$ then, 
%$$W_r^r(\nu,\rho) \geq \frac{\Delta_\beta^r}{r+1} +
%	 c_1 \Delta_\beta^r\epsilon.$$
$$\nu = U(0, \xi), \quad \rho = (1-\epsilon) U(0, 1) +  \epsilon \delta_1.$$
Then, 
$$W_r(\nu, \rho) \geq \frac{1}{r+1}\Big[|\Delta_\xi|^{r} +  r\epsilon |\Delta_\xi|^{r-1}\Big]$$ 
\end{enumerate}
\end{lemma}
We prove this result in Appendix~\ref{app::proof_main_univariate}. Taking this result as given, we can now compute the various Sliced 
Wasserstein distances and $\sj_r$ evaluations.

\paragraph{Computing the Sliced Wasserstein distances.} 
\begin{itemize}
\item {\bf Construction 1.} For any $\theta \in \bbS^{d-1}$, let $F_{01,\theta}^{-1}$, $F_{11,\theta}^{-1}$,
$G_{01,\theta}^{-1}$ and $G_{11,\theta}^{-1}$ denote the respective
quantile functions of $\pi_\theta \# P_{01}$, $\pi_\theta \# P_{11}$, $\pi_\theta \# Q_{01}$, 
$\pi_\theta \# Q_{11}$. We have
\begin{align*}
F_{01,\theta}^{-1}(u) &= \begin{cases}
0 \wedge A^\top\theta, & u \in (0, 1/2) \\
0 \vee A^\top\theta,   & u \in [1/2,1),
\end{cases}\\
F_{11,\theta}^{-1}(u) &= \begin{cases}
0 \wedge A^\top\theta, & u \in (0, 1/2+\epsilon_n) \\
0 \vee A^\top\theta,   & u \in [1/2+\epsilon_n,1),
\end{cases}\\G_{01,\theta}^{-1}(u) = G_{11,\theta}^{-1}&=  \begin{cases} 
g A^\top\theta \wedge (1+g)A^\top\theta, & u \in (0, 1/2) \\
g A^\top \theta \vee (1+g)A^\top\theta,  & u \in [1/2,1).
\end{cases}
\end{align*}
Therefore, 
\begin{align*}
\sw_r^r&(P_{01},Q_{01}) \\
 &= \int_{\bbS^{d-1}} \int_0^1
 		\big| F_{01,\theta}^{-1}(u) - G_{01,\theta}^{-1}(u) \big|^r du d\mu(\theta) \\
 &= \frac 1 2 \int_{\{\theta \in \bbS^{d-1}: A^\top\theta \geq 0\}} |g A^\top\theta|^r d\mu(\theta) \\ &~~+ 
    \frac 1 2 \int_{\{\theta \in \bbS^{d-1}: A^\top\theta < 0\}} |A^\top\theta - (1+g) A^\top \theta|^r d\mu(\theta) \\
 &= g^r \int_{\bbS^{d-1}} |A^\top\theta|^r d\mu(\theta)= \Gamma^r.
\end{align*}
Furthermore, 
\begin{align*}
\sw_r^r&(P_{11},Q_{11}) \\
 &= \int_{\bbS^{d-1}} \int_0^1
 		\big| F_{1,\theta}^{-1}(u) - G_{11,\theta}^{-1}(u) \big|^r du d\mu(\theta) \\
 &= \int_{\bbS^{d-1}} \Bigg(
			\int_0^{1/2} |0 \wedge A^\top\theta - g A^\top \theta \wedge (1+g) A^\top\theta|^r du \\ &\qquad \quad \ + 
			\int_{1/2}^{1/2+\epsilon_n} |0 \wedge A^\top\theta - g A^\top \theta \vee (1+g) A^\top\theta|^r du \\ &\qquad \quad \ + 
			\int_{1/2+\epsilon_n}^1 |0 \vee A^\top\theta - 
						g A^\top\theta \vee (1+g) \theta^\top A |^r du\Bigg)d\mu(\theta) \\
 &= (1-\epsilon_n)g^r \int_{\bbS^{d-1}} |A^\top\theta|^rd\mu(\theta) \\ & \qquad \quad \ + 
 	\epsilon_n \int_{\bbS^{d-1}} |0 \wedge A^\top\theta - 
 							g A^\top\theta \vee (1+g)A^\top\theta|^r d\mu(\theta)		\\
 &= (1-\epsilon_n)g^r I_r^r + 
 	\epsilon_n \int_{\bbS^{d-1}} |0 \wedge A^\top\theta - 
 							g A^\top\theta \vee (1+g)A^\top\theta|^r d\mu(\theta) \\
 &= \Gamma^r + c_1\epsilon_n,			
\end{align*}
for a positive constant $c_1 > 0$. 
It follows that 
$\sw_r(P_{11},Q_{11}) \geq$ $\sw_r(P_{01},Q_{01}) \geq \Gamma$, thus $(P_{01},Q_{01}),(P_{11},Q_{11}) \in \calO(\Gamma; \infty,\infty)$, and
\begin{align*}
\big|\sw_r(P_{01},Q_{01}) - \sw_r(P_{11},Q_{11})\big|
 &= \left|\Gamma - \big(\Gamma^r + c_1\epsilon_n\big)^{\frac 1 r} \right|.
\end{align*}

\item {\bf Construction 2.} 
We use the first part of Lemma~\ref{lem::main_univariate}, and let $\nu =$ $ \frac{1+\epsilon_n}{2} U(0,1/2)$ $+ \frac{1-\epsilon_n}{2} U(1/2,1),$
and $\rho = U(\Delta, 1+ \Delta)$.
Notice that 
if $X \sim \nu$, then $\gamma_2^{1/r} X \sim P_{12}^{(1)}$, and 
if $Y \sim \rho$, then $\gamma_2^{1/r}Y \sim Q_{02}^{(1)}$. Therefore, by Proposition 7.16 of \cite{villani2003}, 
$W_r({\pi_\theta}_{\#}P_{12}, {\pi_\theta}_{\#} Q_{12}) = |\theta_1| \gamma_2^{1/r} W_r(\nu,\rho)$. Thus,
\begin{align*}
\sw_r^r(P_{12},Q_{12})
 &= \int_{\bbS^{d-1}} W_r^r({\pi_{\theta}}_{\#} P_{12}, {\pi_{\theta}}_{\#} Q_{12}) d\mu(\theta)\\ 
 &= \int_{\bbS^{d-1}} |\theta_1|^r \gamma_2 W_r^r(\nu, \rho) d\mu(\theta) 
 \geq \gamma_2 t_r^r \left[ \Delta^r + \frac{r}{4} \Delta^{r-1} \epsilon_n\right],
\end{align*}
by Lemma \ref{lem::main_univariate}. 
Furthermore, it is easy to show that
\begin{align*}
\sw_r^r(P_{02},Q_{02}) = \gamma_2\Delta^r \int_{\bbS^{d-1}} |\theta_1|^rd\mu(\theta) = \gamma_2 ( t_r\Delta)^r.
\end{align*}    

\item {\bf Construction 3.} We use the second part of Lemma~\ref{lem::main_univariate}. We set 
$$\nu = U(0, \bar \beta), \quad \rho = (1-\epsilon_m) U(0, 1) +  \epsilon_m \delta_1.$$
Then, for all $ \epsilon \in (0,1]$, 
$$W_r(\nu, \rho) \geq \frac{1}{r+1}\Big[|\Delta_{\bar\beta}|^{r} +  r\epsilon_m |\Delta_{\bar\beta}|^{r-1}\Big]$$ 
We then obtain
\begin{align*}
\sw_r^r(P_{13}, Q_{13})
 &= \int_{\bbS^{d-1}} |\theta_1|^r s_2 W_r^r(\nu, \rho)  d\mu(\theta) \\
 &\geq \frac{t_r^r s_2}{r+1}\Big[|\Delta_{\bar\beta}|^{r} +  r\epsilon_m |\Delta_{\bar\beta}|^{r-1}\Big].
\end{align*}
On the other hand, it is easy to see that  
$$\sw_r^r(P_{03}, Q_{03}) = \frac{s_2 t_r^r|\Delta_{\bar\beta}|^r}{r+1}.$$
%where in the last inequality of the above display, we used the fact that $s_2^{1/r} \geq 2s_1^{1/r}$, due to 
%\eqref{eq::delta_b_ineq}.
\item {\bf Construction 4.} We again use the second part of Lemma~\ref{lem::main_univariate}, setting
$$\nu = U(0, \beta), \quad \rho = (1-\epsilon) U(0, 1) +  \epsilon_n \delta_1.$$
The rest follows by the same argument as for Construction 3. 
\end{itemize}

\paragraph{Computing the $\sj_r$ evaluations.} 
Our next step will be to compute the $\sj_r$ functionals for the various distributions we have constructed. We note that for Construction 1
our distributions are allowed to have infinite $\sj_r$ so we only need to consider Constructions 2-4.
 The calculations for Construction 3  and 4
follow along very similar lines to those of Construction 2, which we detail below.
\begin{itemize}
\item {\bf Construction 2.}  
We have
\begin{align*}
\sj_r(Q_{02}) = \sj_r(Q_{12}) &= \sj_r(P_{02}) \\
 &\leq \int_{\bbS^{d-1}} \int_0^1 \left(\frac{\sqrt{u(1-u)}}{1/(|\theta_1| \gamma_2^{1/r})} \right)^r du d\mu(\theta)  \\
 &\leq \int(|\theta_1| \gamma_2^{1/r})^r d\mu(\theta)
 \leq \gamma_2.
\end{align*} 
 Furthermore, 
\begin{align*}
\sj_r(Q_{13})
&\leq\int_{\bbS^{d-1}} \int_0^1\left(\frac{\sqrt{u(1-u)}}{(1-\epsilon_n)/(|\theta_1|\gamma_2^{1/r})}\right)^rdu d\mu(\theta)
\\ &\leq \frac{\gamma_2}{(1-\epsilon_n)^r}\int_0^1 [u(1-u)]^{\frac r 2} du.
\end{align*}
Choosing the constant $k_r > 0$ to satisfy $k_r < 1-\left(\int_0^1 [u(1-u)]^{\frac r 2} du\right)^{\frac 1 r}$
guarantees that the above display is bounded above by $\gamma_2$. 
\end{itemize}
To complete the proof it remains to prove Lemma~\ref{lem::main_univariate}. 

\subsubsection{Proof of Lemma~\ref{lem::main_univariate}}
\label{app::proof_main_univariate}
We prove each of the two claims in turn.

\paragraph{Proof of Claim (1).} 
Notice that the quantile functions of $\nu$ and $\rho$ are respectively given by
\begin{align*}
F^{-1}(u) &= \begin{cases}
\frac u {1+\epsilon}, & 0 \leq u \leq (1+\epsilon)/2, \\
\frac 1 2 + \frac {u-(1+\epsilon)/2} {1-\epsilon}, & (1+\epsilon)/2 \leq u \leq 1
\end{cases}, \\
G^{-1}(u) &= (u+\Delta)I(0 \leq u \leq 1).
\end{align*}
Thus,
\begin{align*}
W_r^r(\nu,\rho) 
 &= \int_0^1 |F^{-1}(u) - G^{-1}(u)|^r du\\
 &= 
  \int_0^{\frac{1+\epsilon}{2}} \left|\Delta+u - \frac u {1+\epsilon}\right|^r du  + 
  \int_{\frac{1+\epsilon}{2}}^1 \left| \Delta+u - \frac 1 2 - \frac {u-(1+\epsilon)/2} {1-\epsilon}  \right|^r du
  \\
 &= \text{(I)} + \text{(II)},
\end{align*}
say. We have,
\begin{align*}
 \text{(I)}
 = \int_0^{(1+\epsilon)/2} \left[\Delta + \frac \epsilon {1+\epsilon}u\right]^r du  
 = \frac{1+\epsilon}{\epsilon(r+1)} \left\{ \left(\Delta+ \frac\epsilon 2\right)^{r+1} - \Delta^{r+1} \right\}.
\end{align*}
Also,
\begin{align*}
\text{(II)}
 &= \int_{(1+\epsilon)/2}^1 \left(\Delta +\frac \epsilon {1-\epsilon} - u\frac \epsilon{1-\epsilon}\right)^r du \\
 &= -\frac{1-\epsilon}{(r+1)\epsilon} \left\{\left( (\Delta + \frac\epsilon{1-\epsilon} - \frac{\epsilon}{1-\epsilon}\right)^{r+1}
						- \left(\Delta + \frac\epsilon{1-\epsilon} - \frac{(1+\epsilon)\epsilon}{2(1-\epsilon)}\right)^{r+1}\right\}\\
 &= -\frac{1-\epsilon}{\epsilon  (r+1)} \left\{\Delta^{r+1} - \left(\Delta + \frac \epsilon 2 \right)^{r+1}\right\}.
\end{align*}
Thus,
\begin{align*}
\text{(I)}+\text{(II)}
 &= \left\{ \left(\Delta + \frac \epsilon 2\right)^{r+1} - \Delta^{r+1}\right\} \left(\frac{1+\epsilon}{\epsilon (r+1)} + \frac{1-\epsilon}{\epsilon(r+1)}\right\}\\
 &= \frac 2 {\epsilon (r+1)} \left\{ \left(\Delta + \frac \epsilon 2\right)^{r+1} - \Delta^{r+1}\right\} \\
 &= \frac 2 {\epsilon (r+1)} \left\{ \Delta^{r+1} + \frac{r+1}{2}\Delta^r \epsilon + \frac{r(r+1)}{8}(\Delta+\tilde\epsilon)^{r-1} \epsilon^2 - \Delta^{r+1}\right\},
\end{align*}
for some $\tilde\epsilon \in (0,\epsilon/2)$, by a first-order Taylor expansion.
Therefore, 
$$W_r^r(\nu,\rho) = \Delta^r + \frac{r}{4}(\Delta + \tilde\epsilon)^{r-1}\epsilon  \geq \Delta^r + \frac{r}{4} \Delta^{r-1} \epsilon,$$
and the claim follows. \hfill $\square$

\paragraph{Proof of Claim (2).}
The respective quantile functions of $\nu$ and $\rho$ are given by
$$F^{-1}(u) = \begin{cases}
\frac u {1-\epsilon}, & 0 \leq u \leq 1 - \epsilon, \\
1, & 1-\epsilon < u \leq 1
\end{cases}, \qquad 
G^{-1}(u) = \xi uI(0 \leq u \leq 1).$$ 
Thus,
\begin{align*}
\sw_r^r(\nu,\rho)
 &=\int \big|F^{-1}(u) - G^{-1}(u)\big|^r du \\
 &=\int_0^{1-\epsilon} \left|\frac u {1-\epsilon} - \xi u \right|^r du + 
   \int_{1-\epsilon}^1 \left| 1 - \xi u\right|^r du \\
 &=\int_0^{1-\epsilon} \left[\frac u {1-\epsilon} - \xi u \right]^r du + 
   \int_{1-\epsilon}^1 \left[ 1 - \xi u\right]^r du,\quad  (\text{Since } \xi \in(0,1])\\
 &=\frac{1-\epsilon}{r+1}(-\Delta_\xi + \epsilon\xi)^r  + 
   \frac 1 {\xi(r+1)}\Big[(-\Delta_\xi + \epsilon\xi)^{r+1} - (-\Delta_\xi)^{r+1}\Big] \\
 &=\frac 1 {r+1}\left[(-\Delta_\xi + \epsilon\xi)^r \left( 1 - \epsilon + \frac{-\Delta_\xi + \epsilon\xi}{\xi}\right)
 			- \frac{|\Delta_\xi|^{r+1}}{\xi}\right]\\
 &= \frac 1 {\xi(r+1)}\Big[ (|\Delta_\xi| + \epsilon\xi)^r - |\Delta_\xi|^{r+1} \Big] \\
 &= \frac 1 {\xi(r+1)}\Big[ (|\Delta_\xi| + \epsilon\xi)^r - |\Delta_\xi|^{r+1} \Big] \\
 &\geq \frac 1 {\xi(r+1)}\Big[ |\Delta_\xi|^r + r\epsilon\xi|\Delta_\xi|^{r-1}  - |\Delta_\xi|^{r+1} \Big] \\ 
 &= \frac {|\Delta_\xi|^r} {r+1} + \frac {r\epsilon|\Delta_\xi|^{r-1}} {r+1}.
\end{align*} 			
The claim follows. \hfill $\square$
 
\section{Proof of Theorem \ref{thm:untrimmed_rate}}
\label{app:untrimmed}
By the same argument as in the proof of Proposition~\ref{prop::empirical_sw},
it will suffice to prove equation~\eqref{eq:untrimmed_rate_power_r}.
 Let $P,Q \in \calK_{r,\rho}(b)$. We prove the claim in five steps. 

\paragraph{Step 0: Preparation.}
By the same argument as in the proof of Lemma~\ref{lem::calK}, notice that there exists a constant $C_{\rho} > 0$ such that
for any $\theta \in \bbS^{d-1}$,  
$$\big|F_\theta^{-1}(u)\big|\vee \big|G_\theta^{-1}(u)\big| \leq C_{\rho} \left(\frac{b_\theta}{u(1-u)}\right)^{\frac 1 \rho},\quad u \in (0,1),$$
where for $X \sim P$, $Y \sim Q$ and each $\theta \in \bbS^{d-1}$, 
$$b_\theta = \bbE [|X^\top \theta|^\rho]+ \bbE [|Y^\top \theta|^\rho].$$ 
Notice that the assumption $P,Q \in \calK_{r,\rho}(b)$ for $\rho > 2r$ implies that
$\int b_\theta^{\frac r \rho} d\mu(\theta) \leq 2b$. In particular, $b_\theta$ is finite
for almost all $\theta \in \bbS^{d-1}$. 
These statements may be applied analogously to the empirical measures $P_n$ and $Q_m$. 
Specifically, we~have
$$\big|F_{\theta,n}^{-1}(u)\big| \vee \big|G_{\theta,m}^{-1}(u)\big|\leq C_\rho \left(\frac{b_{\theta,nm}}{u(1-u)}\right)^{\frac 1 \rho},\quad u \in (0,1),$$
where we set 
$$b_{\theta,nm} := \int |x^\top\theta|^\rho dP_n(x)+  \int |y^\top\theta|^\rho dQ_m(y).$$
Combining these facts, we have for all $u \in (0,1)$, 
$$\big|F_\theta^{-1}(u)\big| \vee \big| G_\theta^{-1}(u)\big| \vee \big|F_{\theta,n}^{-1}(u)\big|\vee \big| G_{\theta,m}^{-1}(u)\big| 
\leq \psi_\theta(u) := C_\rho \left(\frac{b_\theta + b_{\theta,nm}}{u(1-u)}\right)^{\frac 1 \rho}.$$
We suppress the dependence of $\psi_\theta$ on $n$ and $m$ for ease of notation,
but we emphasize that $\psi_\theta(u)$ is a random variable. 

Our proof makes use of a uniform 
self-normalized concentration inequality for the empirical quantile process, which was introduced
in Section~\ref{app::high_prob_sw_null}, as part of the proof of Lemma~\ref{lem::high_prob_sw_null}, and 
also in Example~\ref{ex::rel_VC} of the main manuscript. 
Specifically, for any $\epsilon \in (0,1)$, let $\gamma_{\epsilon,n}, \eta_{\epsilon,n}$ be the sequences
given in equation~\eqref{eq::relVC_fns_pf}, with inverses given in 
equation~\eqref{eq::relVC_inverse_fns_pf}, and defined in terms of the quantity 
$$\nu_{\epsilon,n} = \sqrt{\frac {16}{n} \left[ \log(16/\epsilon) +  \log(2n+1)\right]},$$ 
for any given $\epsilon \in (0,1)$. Recall that for any $\theta\in \bbS^{d-1}$,
the event 
$$A_\epsilon = \left\{ u \in (0,1): F_\theta^{-1}(\gamma_{\epsilon,n}(u)) \leq F_{\theta,n}^{-1}(u) \leq F_\theta^{-1}(\eta_{\epsilon,n}(u))\right\}$$
satisfies $\bbP(A_\epsilon) \geq 1-\epsilon$.

\paragraph{Step 1:  First Reduction.}
 Apply a similar first-order Taylor expansion  
as in the proof of Lemma~\ref{lem::rth-deviation}, to the cost function $|\cdot|^r$,  to deduce that
$$\begin{multlined}[0.82\textwidth]
\big|\sw_{r}^r(P_n,Q_m) - \sw_r^r(P,Q)\big|  \\ 
  \leq   
 r \int_{\bbS^{d-1}}\int_0^1 |\tilde F_{\theta,n}^{-1}(u) - \tilde G_{\theta,m}^{-1}(u)|^{r-1} \\ 
    \times \Big[|F_{\theta,n}^{-1}(u) - F_\theta^{-1}(u)| + |G_{\theta,m}^{-1}(u) - G_\theta^{-1}(u)|\Big]dud\mu(\theta),\hspace{-0.4in}
\end{multlined}$$
for some $\tilde F_{\theta,n}^{-1}(u)$ on the segment joining $F_\theta^{-1}(u)$ and $F_{\theta,n}^{-1}(u)$, 
and some $\tilde G_{\theta,m}^{-1}(u)$ on the segment joining $G_\theta^{-1}(u)$ and $G_{\theta,m}^{-1}(u)$. 
By definition of $\psi_{\theta}$, we deduce
\begin{align*}
\big|\sw_{r}^r&(P_n,Q_m) - \sw_r^r(P,Q) \big| \\
 &\lesssim
  \int_{\bbS^{d-1}}\int_0^1 \psi_\theta^{r-1}(u)
    \Big[|F_{\theta,n}^{-1}(u) - F_\theta^{-1}(u)| + |G_{\theta,m}^{-1}(u) - G_\theta^{-1}(u)|\Big]dud\mu(\theta)\\
 &=
  \int_{\bbS^{d-1}}\int_0^{1/2} \psi_\theta^{r-1}(u)
    \Big[|F_{\theta,n}^{-1}(u) - F_\theta^{-1}(u)| + |G_{\theta,m}^{-1}(u) - G_\theta^{-1}(u)|\Big]dud\mu(\theta) \\    
 & + 
  \int_{\bbS^{d-1}}\int_{1/2}^1 \psi_\theta^{r-1}(u)
    \Big[|F_{\theta,n}^{-1}(u) - F_\theta^{-1}(u)| + |G_{\theta,m}^{-1}(u) - G_\theta^{-1}(u)|\Big]dud\mu(\theta).
\end{align*}
We shall bound the quantity 
$$T = \int_{\bbS^{d-1}}\int_{1/2}^1 \psi_\theta^{r-1}(u)|F_{\theta,n}^{-1}(u) - F_\theta^{-1}(u)|dud\mu(\theta),$$
 and a similar argument can  be used to bound the remaining terms in the penultimate display.  

Throughout the sequel, let $\epsilon \in (0,1)$ be chosen such that $\nu_{\epsilon,n} \leq 1$. Let \begin{equation}
\label{eq:beta_def}
\beta > \left(\frac 1 2 - \frac r \rho\right)^{-1} > 0,
\end{equation}
where we recall that $\rho > 2r$. 
Define $\delta_k = k^{-\beta}/2$ for all $k \geq 1$. Notice that 
$\delta_1 = 1/2$. Furthermore, let 
$$K \equiv K_n = 1\vee \left\lfloor (c\nu_{\epsilon,n})^{-\frac{\rho}{\beta(\rho-r)}}\right\rfloor$$
for a  constant
$c \geq 8$ to be specified below. We summarize a few algebraic facts in relation to the sequences
$\delta_k,\eta_{\epsilon,n}, \gamma_{\epsilon,n}$, 
which we prove in Section~\ref{app:pf_technical_deltak}. 
%We have the following properties. 
\begin{lemma}
\label{lem::technical_deltak}
There exists a choice of the constant $c \geq 8$, as well as a constant
$c_1 > 0$, both possibly depending on $\beta,\rho,r$, such that for all $n \geq 1$, and all $\epsilon \in (0,1)$ for which $\nu_{\epsilon,n} \leq 1$, the following properties hold.
\begin{enumerate}
\item[(i)] $c_1\nu_{\epsilon,n}^{\frac{\rho}{\rho-r}} \geq  \delta_K$.
Furthermore, if $K \geq 2$, then
$\delta_K \geq (c\nu_{\epsilon,n})^{\frac{\rho}{\rho-r}}/2$. 
\\[-0.1in] 
\item[(ii)] $1 - \gamma_{\epsilon,n}(1-\delta_k) 
\leq \delta_k + \frac{\nu_{\epsilon,n}^2}{1+\nu_{\epsilon,n}^2} + \nu_{\epsilon,n}\sqrt{\delta_k}$
for all $k=1, \dots, K$. 
\item[(iii)] $\frac{1-\delta_k}{1+\nu_{\epsilon,n}^2}\leq \eta_{\epsilon,n}(1-\delta_k) \leq 1-\delta_k + \nu_{\epsilon,n}^2 + \nu_{\epsilon,n}\sqrt{\delta_k}$, 
for all $k=1, \dots, K$.\\[-0.1in]
\item[(iv)] $\eta_{\epsilon,n}(1-\delta_k) \leq \gamma_{\epsilon,n}(1-\delta_{k+1}), $ for all $k=1, \dots, K-1$, if $K \geq 2$. 
\end{enumerate}
\end{lemma}
With these facts in place, consider the decomposition
\begin{align}
\label{eq:peeling_II}
T = \int_{\bbS^{d-1}}\int_{1/2}^1 \psi_\theta^{r-1}(u)|F_{\theta,n}^{-1}(u) - F_\theta^{-1}(u)|dud\mu(\theta) = %\sum_{k=1}^{K-1} 
\int_{\bbS^{d-1}}  T_\theta  d\mu(\theta),
\end{align}
where we recall that $\delta_1=1/2$, and we set
\begin{align}
T_\theta &= \int_{1/2}^1|F_{\theta,n}^{-1}(u) - F_\theta^{-1}(u)| \psi_\theta^{r-1}(u)du= T_{\theta,K}^*  +\sum_{k=1}^{K-1} T_{\theta,k},\quad \theta \in \bbS^{d-1},\\ 
T_{\theta,k} &= \int_{1-\delta_k}^{1-\delta_{k+1}} |F_{\theta,n}^{-1}(u) - F_\theta^{-1}(u)| \psi_\theta^{r-1}(u)du , 
~~ k=1, \dots, K-1, \\
 T_{\theta,K}^* &= \  \int_{1-\delta_K}^1 |F_{\theta,n}^{-1}(u) - F_\theta^{-1}(u)| \psi_\theta^{r-1} (u)du.
\end{align}
In the following two steps, we bound the preceding two terms for any fixed $\theta \in \bbS^{d-1}$ and $k=1, \dots, K-1$. 
 The symbol ``$\lesssim$'' will always hide
constants which do not depend on $\theta$ and $k$.

\paragraph{Step 2: Bounding $T_{\theta,K}^*$.}  
We have, 
\begin{align*}
%\nonumber 
%\label{eq:peeling_II_truncation}
T_{\theta,K}^* &= \int_{1-\delta_K}^1 |F_{\theta,n}^{-1}(u) - F_\theta^{-1}(u)| \psi_\theta^{r-1}(u)du 
 \lesssim  \int_{1-\delta_K}^1 \psi_\theta^r(u)du  \\
 &\lesssim (b_\theta + b_{\theta,nm})^{\frac{r}{\rho}}  \int_{1-\delta_K}^1 (1-u)^{-\frac{r}{\rho}}du  \\
 &\lesssim (b_\theta + b_{\theta,nm})^{\frac{r}{\rho}} \delta_K^{1-\frac r \rho}
\lesssim(b_\theta + b_{\theta,nm})^{\frac{r}{\rho}}   \nu_{\epsilon,n},
\end{align*} 
where the final inequality follows by Lemma~\ref{lem::technical_deltak}(i).

\paragraph{Step 3: Bounding $T_{\theta,k}$ in Probability.}  
The bulk of our work will now go into
bounding $T_{\theta,k}$, for any given 
$\theta \in \bbS^{d-1}$ and $k=1, \dots, K-1$. Notice that
this case is vacuous when $ K =1$. 
%We shall assume throughout the sequel that $\epsilon$ is chosen sufficiently large to ensure $\nu_{\epsilon,n} \leq 1$. 
The following derivations are performed over the event $A_\epsilon$. By its definition, we have for any $k=1,\dots, K-1$, 
\begin{align}
\nonumber
T_{\theta,k}
 &\leq \int_{1-\delta_k}^{1-\delta_{k+1}} \big[F_\theta^{-1}(\eta_{\epsilon,n}(u)) - F_\theta^{-1}(\gamma_{\epsilon,n}(u))\big] \psi_\theta^{r-1}(u)du \\
\nonumber
 &\leq \psi_\theta^{r-1}(1-\delta_{k+1}) \int_{1-\delta_k}^{1-\delta_{k+1}} \big[F_\theta^{-1}(\eta_{\epsilon,n}(u)) - F_\theta^{-1}(\gamma_{\epsilon,n}(u))\big]  du
 \\ 
\nonumber
 &= \psi_\theta^{r-1}(1-\delta_{k+1}) \Bigg[
     \int_{\eta_{\epsilon,n}(1-\delta_k)}^{\eta_{\epsilon,n}(1-\delta_{k+1})}  
     		F_\theta^{-1}(u)\frac{\partial \eta_{\epsilon,n}^{-1}(u)}{\partial u}du \\ \nonumber &\hspace{1.5in} - 
     \int_{\gamma_{\epsilon,n}(1-\delta_k)}^{\gamma_{\epsilon,n}(1-\delta_{k+1})}  
     		F_\theta^{-1}(u)\frac{\partial \gamma_{\epsilon,n}^{-1}(u)}{\partial u}du 
   \Bigg]\\ 
   &= \psi_\theta^{r-1}(1-\delta_{k+1}) (A_{\theta,k} + B_{\theta,k}), 
   \label{eq::Ak_plus_Bk} 
\end{align}
where
\begin{align*}
A_{\theta,k} &= 
     \int_{\eta_{\epsilon,n}(1-\delta_k)}^{\eta_{\epsilon,n}(1-\delta_{k+1})}  F_\theta^{-1}(u)du- 
     \int_{\gamma_{\epsilon,n}(1-\delta_k)}^{\gamma_{\epsilon,n}(1-\delta_{k+1})}  F_\theta^{-1}(u)du, \\
B_{\theta,k} &= \int_{\eta_{\epsilon,n}(1-\delta_k)}^{\eta_{\epsilon,n}(1-\delta_{k+1})}  F_\theta^{-1}(u)\left(\frac{\partial \eta_{\epsilon,n}^{-1}(u)}{\partial u}-1\right) du  \\ &-
     \int_{\gamma_{\epsilon,n}(1-\delta_k)}^{\gamma_{\epsilon,n}(1-\delta_{k+1})}  F_\theta^{-1}(u)\left(\frac{\partial \gamma_{\epsilon,n}^{-1}(u)}{\partial u}-1\right)du.
\end{align*}

{\bf Step 3.1: Bounding $A_{\theta,k}$.} Consider the decomposition
$$A_{\theta,k} {=} \left(
     \int_{\eta_{\epsilon,n}(1-\delta_k)}^{\gamma_{\epsilon,n}(1-\delta_{k+1})}+ 
     \int_{\gamma_{\epsilon,n}(1-\delta_{k+1})}^{\eta_{\epsilon,n}(1-\delta_{k+1})}   - 
     \int_{\gamma_{\epsilon,n}(1-\delta_k)}^{\eta_{\epsilon,n}(1-\delta_{k})} - 
     \int_{\eta_{\epsilon,n}(1-\delta_k)}^{\gamma_{\epsilon,n}(1-\delta_{k+1})}   
   \right) F_\theta^{-1}(u)du.$$
Using Lemma~\ref{lem::technical_deltak}(iv) and the fact that $\eta_{\epsilon,n} \geq \gamma_{\epsilon,n}$, the four lower bounds of integration
in the above display  
are less than their respective upper bounds. Therefore, 
\begin{align*}
A_{\theta,k} 
 &= \left(
 		 \int_{\gamma_{\epsilon,n}(1-\delta_{k+1})}^{\eta_{\epsilon,n}(1-\delta_{k+1})}   - 
        \int_{\gamma_{\epsilon,n}(1-\delta_k)}^{\eta_{\epsilon,n}(1-\delta_{k})} \right) F_\theta^{-1}(u)du\\
 &\leq \left(
 		 \int_{\gamma_{\epsilon,n}(1-\delta_{k+1})}^{\eta_{\epsilon,n}(1-\delta_{k+1})}   + 
        \int_{\gamma_{\epsilon,n}(1-\delta_k)}^{\eta_{\epsilon,n}(1-\delta_{k})} \right) \psi_\theta(u) du \\
 &\leq \psi_\theta(\eta_{\epsilon,n}(1-\delta_{k+1})) \Big[
 		 \big(\eta_{\epsilon,n}(1-\delta_{k+1}) - \gamma_{\epsilon,n}(1-\delta_{k+1})\big) \\ &\hspace{1.6in}+ 
		 \big(\eta_{\epsilon,n}(1-\delta_{k  }) - \gamma_{\epsilon,n}(1-\delta_{k  })\big)  \Big] \\		
 &\leq  	\psi_\theta(\eta_{\epsilon,n}(1-\delta_{k+1}))\left[ \frac{\nu_{\epsilon,n}\sqrt{\nu_{\epsilon,n}^2 + 4\delta_{k+1}}}{1+\nu_{\epsilon,n}^2} + \frac{\nu_{\epsilon,n}\sqrt{\nu_{\epsilon,n}^2 + 4\delta_{k}}}{1+\nu_{\epsilon,n}^2}
 \right] \\ 
 &\lesssim  	\psi_\theta(\eta_{\epsilon,n}(1-\delta_{k+1}))\left[    \nu_{\epsilon,n}^2  + \nu_{\epsilon,n} \sqrt{\delta_{k}} 
 \right]	.
\end{align*}	
By Lemma~\ref{lem::technical_deltak}(i), we deduce that %for all $k=1, \dots, K-1$, 
$$A_{\theta,k} \lesssim 		\psi_\theta(\eta_{\epsilon,n}(1-\delta_{k+1})) \nu_{\epsilon,n} \left[\delta_K^{\frac{\rho-r}{\rho}} + \delta_k^{\frac 1 2}\right] \lesssim \psi_\theta(\eta_{\epsilon,n}(1-\delta_{k+1})) \nu_{\epsilon,n}   \sqrt{\delta_k},$$
where we used the assumption   $\rho > 2r$. 
We next turn to bounding $B_{\theta,k}$. 

\paragraph{Step 3.2: Bounding $B_{\theta,k}$.} We have, 
\begin{align*}
\frac{\partial}{\partial u} \eta_{\epsilon,n}^{-1}(u) &= 1 - \frac {\nu_{\epsilon,n}}{2}\left[  \sqrt{\frac{1-u}{u}} - 
\sqrt{\frac{u}{1-u}}\right],\quad u \in \left[\eta_{\epsilon,n}(1/2), \eta_{\epsilon,n}(1)\right],
\\
\frac{\partial}{\partial u} \gamma_{\epsilon,n}^{-1}(u) &= 1 + \frac {\nu_{\epsilon,n}}{2}\left[  \sqrt{\frac{1-u}{u}} - 
\sqrt{\frac{u}{1-u}}\right], \quad u \in \left[\gamma_{\epsilon,n}(1/2), \gamma_{\epsilon,n}(1)\right].
\end{align*}
Since $\nu_{\epsilon,n}\leq 1$, notice that $\gamma_{\epsilon,n}(1/2)$ and $\eta_{\epsilon,n}(1/2)$ 
are bounded below by a positive universal
constant.  
Therefore, in the above display, the first terms in brackets are bounded above by positive
universal constants, uniformly over the stated ranges, leading to
\begin{align*}
\left|\frac{\partial}{\partial u} \eta_{\epsilon,n}^{-1}(u) - 1\right|
\lesssim \frac {\nu_{\epsilon,n}}{2} \sqrt{\frac{1}{1-u}},\quad u \in \left[\eta_{\epsilon,n}(1/2), \eta_{\epsilon,n}(1)\right],
\\
\left|\frac{\partial}{\partial u} \gamma_{\epsilon,n}^{-1}(u) - 1\right|
\lesssim \frac {\nu_{\epsilon,n}}{2}  
\sqrt{\frac{1}{1-u}}, \quad u \in \left[\gamma_{\epsilon,n}(1/2), \gamma_{\epsilon,n}(1)\right].
\end{align*}
Deduce that, 
\begin{align*}
  B_{\theta,k} 
 &= \int_{\eta_{\epsilon,n}(1-\delta_k)}^{\eta_{\epsilon,n}(1-\delta_{k+1})}  F_\theta^{-1}(u)\left(\frac{\partial \eta_{\epsilon,n}^{-1}(u)}{\partial u}-1\right) du \\ &\hspace{0.5in} -
     \int_{\gamma_{\epsilon,n}(1-\delta_k)}^{\gamma_{\epsilon,n}(1-\delta_{k+1})}  F_\theta^{-1}(u)\left(\frac{\partial \gamma_{\epsilon,n}^{-1}(u)}{\partial u}-1\right)du \\
  &\lesssim \nu_{\epsilon,n} 
   \int_{\eta_{\epsilon,n}(1-\delta_k)}^{\eta_{\epsilon,n}(1-\delta_{k+1})}  |F_\theta^{-1}(u)| \sqrt{\frac{1}{1-u}} \frac{du}{2}\\
    &\hspace{0.5in}+
        \nu_{\epsilon,n}
   \int_{\gamma_{\epsilon,n}(1-\delta_k)}^{\gamma_{\epsilon,n}(1-\delta_{k+1})}  |F_\theta^{-1}(u)| \sqrt{\frac{1}{1-u}}\frac{du}{2}  \\
%  &\leq \frac{\nu_{\epsilon,n}}{2} |F_\theta^{-1}(\eta_{\epsilon,n}(1-\delta_{k+1}))| \left(\int_{\eta_{\epsilon,n}(1-\delta_k)}^{\eta_{\epsilon,n}(1-\delta_{k+1})} +
%   \int_{\gamma_{\epsilon,n}(1-\delta_k)}^{\gamma_{\epsilon,n}(1-\delta_{k+1})}\right) \left(\sqrt{\frac{1-u}{u}} + \sqrt{\frac{u}{1-u}}\right)   du   \\  
%Now, since $\nu_{\epsilon,n} \leq 1$, $\gamma_{\epsilon,n}(1/2)$ and $\eta_{\epsilon,n}(1/2)$
%are bounded from below by a universal positive constant, thus
%\begin{align*}
%B_{\theta,k}
% &\lesssim \frac{\nu_{\epsilon,n}}{2} |F_\theta^{-1}(\eta_{\epsilon,n}(1-\delta_{k+1}))| \left(\int_{\eta_{\epsilon,n}(1-\delta_k)}^{\eta_{\epsilon,n}(1-\delta_{k+1})} +
%   \int_{\gamma_{\epsilon,n}(1-\delta_k)}^{\gamma_{\epsilon,n}(1-\delta_{k+1})}\right) \left(\sqrt{\frac{1-u}{u}} + \sqrt{\frac{u}{1-u}}\right)   du   \\  
%  &\leq \nu_{\epsilon,n}  \psi_\theta(\eta_{\epsilon,n}(1-\delta_{k+1}))
%  \left(\sqrt{\eta_{\epsilon,n}(1-\delta_{k+1})\big(1 - \eta_{\epsilon,n}(1-\delta_{k+1})\big)}
%  + \sqrt{\eta_{\epsilon,n}(1-\delta_{k+1})\big(1 - \eta_{\epsilon,n}(1-\delta_{k+1})\big)}\right) \\ 
  &\leq 2 \nu_{\epsilon,n}  \psi_\theta(\eta_{\epsilon,n}(1-\delta_{k+1})) \sqrt{1-\gamma_{\epsilon,n}(1-\delta_k)}\\
  &\lesssim \nu_{\epsilon,n}  \psi_\theta(\eta_{\epsilon,n}(1-\delta_{k+1})) \sqrt{\delta_k + \nu_{\epsilon,n}^2 + \nu_{\epsilon,n}\sqrt{\delta_k}}
  \quad (\text{By Lemma } \ref{lem::technical_deltak}\text{(ii)})
  \\  
  &\lesssim \nu_{\epsilon,n}  \psi_\theta(\eta_{\epsilon,n}(1-\delta_{k+1})) 
  \left[\sqrt{\delta_k} + \delta_k^{\frac{\rho-r}{\rho}} + \delta_k^{\frac{\rho-r}{2\rho}+\frac 1 4}\right] \quad (\text{By Lemma } \ref{lem::technical_deltak}\text{(i)})\\
  &\lesssim \nu_{\epsilon,n}  \psi_\theta(\eta_{\epsilon,n}(1-\delta_{k+1}))  \sqrt{\delta_k} ,
  \end{align*}   
where we again used the assumption $\rho > 2r$. 

\paragraph{ Step 4: Bounding $ T_{\theta }$ in Probability.} Combine the conclusions of Steps 3.1 and 3.2
with equation~\eqref{eq::Ak_plus_Bk} to deduce
\begin{align*}
&\sum_{k=1}^{K-1} T_{\theta,k} \\
 &\lesssim  \nu_{\epsilon,n} \sum_{k=1}^{K-1} \sqrt{\delta_k} \psi_\theta^{r-1}(1-\delta_{k+1}) \psi_\theta(\eta_{\epsilon,n}(1-\delta_{k+1}))\\
 &\lesssim  \nu_{\epsilon,n} \sum_{k=1}^{K-1} \sqrt{\delta_k}  \psi_\theta^r(\eta_{\epsilon,n}(1-\delta_{k+1}))\\
 &\lesssim  \nu_{\epsilon,n} (b_\theta + b_{\theta,nm})^{\frac r \rho} \sum_{k=1}^{K-1} \sqrt{\delta_k}   \big[1-\eta_{\epsilon,n}(1-\delta_{k+1})\big]^{-\frac r \rho}\\ 
&\leq  \nu_{\epsilon,n} (b_\theta + b_{\theta,nm})^{\frac r \rho} \sum_{k=1}^{K-1} \sqrt{\delta_k} 
 \left[\delta_k - \nu_{\epsilon_n}^2 -\nu_{\epsilon,n}\sqrt{\delta_k}\right]^{-\frac r \rho}
\quad (\text{By Lem. } \ref{lem::technical_deltak}\text{(iii)})\\  
&\leq   \nu_{\epsilon,n} (b_\theta + b_{\theta,nm})^{\frac r \rho} \sum_{k=1}^{K-1} \sqrt{\delta_k}  
 \left[\delta_k -  \Big(\delta_K^{\frac{2(\rho-r)}{\rho}}  +\delta_K^{\frac{\rho-r}{\rho}}\delta_k^{\frac 1 2 }\Big)/4\right]^{-\frac r \rho}
~ (\text{By Lem. } \ref{lem::technical_deltak}\text{(i)})\\    
 &\lesssim  \nu_{\epsilon,n} (b_\theta + b_{\theta,nm})^{\frac r \rho} \sum_{k=1}^{K-1} \delta_k^{\frac 1 2 - \frac r \rho},
\end{align*}
where we again used the fact that $\rho > 2r$ on the final line. 
By definition of $\beta$ in equation~\eqref{eq:beta_def},
the sequence $\delta_k^{\frac 1 2 - \frac r \rho}$ is summable, thus the summation
in the final line of the above display is bounded above by a finite constant depending 
only on $r,\rho,\beta$. Combine this fact with the conclusion of Step 2 to deduce that, for a constant $C = C(\beta,\rho,r)$, we have
\begin{equation}
\label{eq:Ttheta_step4}
T_\theta \leq C (b_\theta + b_{\theta,nm})^{\frac r \rho}\nu_{\epsilon,n},
\end{equation}
over the high-probability event $A_\epsilon$, for all $\epsilon$ such that $\nu_{\epsilon,n} \leq 1$. 

\paragraph{Step 4: Bounding $T_\theta$ in Expectation.} We now wish to turn
the preceding display into a bound on the expectation $\bbE[T_\theta]$.  
Notice that $\bbE[b_{\theta,nm}] =b_\theta$, thus by Markov's inequality, we have  for all $y > 0$, 
$$\bbP(b_{\theta,nm}^{r/\rho} \geq b_\theta^{r/\rho}  y) \leq \bbE[b_{\theta,nm}] y^{-\rho/r} /b_\theta \leq   y^{-\rho/r}.$$
Furthermore, by inverting the definition of 
$\nu_{\epsilon,n}$ in terms of $\epsilon$, equation~\eqref{eq:Ttheta_step4} implies that 
for all $0 < u \leq 1$, 
$$\bbP \left( T_\theta  \geq C u(b_\theta + b_{\theta,nm})^{r/\rho} \right) \leq \frac{2n+1}{16}\exp\left(-\frac{nu^2}{16}\right).
$$
Combine the preceding two displays to deduce that for all $0<u \leq 1$, 
\begin{align*}
\bbP&\left\{T_\theta > Cu b_\theta^{r/\rho} + C u (b_\theta/n)^{r/\rho} \exp\left(\frac{(r/\rho)n u^2}{16}\right)\right\}  \\
 &= \bbP\left\{T_\theta > Cu b_\theta^{\frac  r\rho} + C u  \left(\frac{b_\theta}{n}\right)^{\frac r \rho}e^{\frac{(r/\rho)n u^2}{16}} , b_{\theta,nm}^{r/\rho} \leq  (b_\theta/n)^{r/\rho}  e^{\frac{(r/\rho)n u^2}{16}} \right\} 
 \\ &\quad +  \bbP\left\{T_\theta > Cu b_\theta^{r/\rho} + C u  (b_\theta/n)^{r/\rho}e^{\frac{(r/\rho)n u^2}{16}}, b_{\theta,nm}^{r/\rho} >  (b_\theta/n)^{r/\rho}  e^{\frac{(r/\rho)n u^2}{16}}  \right\} \\
 &\leq \bbP\left\{T_\theta > C  u (b_\theta^{\frac r \rho}+b_{\theta,nm}^{\frac r \rho})\right\} + \bbP\left\{b_{\theta,nm}^{r/\rho} > (b_\theta/n)^{r/\rho}    \exp\left(\frac{(r/\rho)n u^2}{16}\right)\right\} \\ 
 &\leq \bbP\left\{T_\theta > C  u (b_\theta + b_{\theta,nm})^{\frac r \rho}\right\} + \bbP\left\{b_{\theta,nm}^{\frac r \rho} > (b_\theta/n)^{r/\rho}    \exp\left(\frac{(r/\rho)n u^2}{16}\right)\right\}  \\
 &\lesssim n\exp\left(-\frac{nu^2}{16}\right).
\end{align*}
Let $f_n(u) = Cu b_\theta^{r/\rho} + C u (b_\theta/n)^{r/\rho} \exp\big((r/\rho)n u^2/16\big)$. $f_n$ is strictly increasing over $\bbR_+$, 
thus it is invertible with inverse $f_n^{-1}$. Notice that $f_n^{-1}(0) = 0$ and $f_n^{-1}(u) \to \infty$
as $u \to \infty$. Furthermore, 
$$f_n'(u) \lesssim b_\theta^{r/\rho} + (b_\theta/n)^{r/\rho}\Big[ \exp\big((r/\rho) n u^2/16\big) + nu^2\exp\big((r/\rho) n u^2/16\big)\Big].$$
Now, let $t_n = \sqrt{\frac{16 \log n}{n} }$ and let $n$ be sufficiently
large to ensure $t_n \leq 1$. 
We have, 
\begin{align*}
\bbE\Big[  T_\theta  \cdot I(T_\theta & \leq f_n(1))\Big]  \\
 &= \int_0^{f_n(1)}\bbP( T_\theta \geq x)dx \\
 &\leq f_n(t_n) + \int_{f_n(t_n)}^{f_n(1)}\bbP( T_\theta \geq x)dx \\
 &= f_n(t_n) + \int_{t_n}^{1} \bbP( T_\theta  \geq f(u)) f'(u) du \\
 &\lesssim f_n(t_n) + b_\theta^{r/\rho}n\int_{t_n}^1 \exp\big(-n  u^2/16\big)   du  \\
 &\qquad\qquad  + \left(\frac{b_\theta}{n}\right)^{r/\rho}n\int_{t_n}^1 \exp\big(-n (1-r/\rho)u^2/16\big)   du 
 \\ &\qquad\qquad+ \left(\frac{b_\theta}{n}\right)^{r/\rho}n^2\int_{t_n}^1 u^2\exp\big(-n (1-r/\rho)u^2/16\big)    du \\
 &\lesssim f_n(t_n) + b_\theta^{r/\rho} \sqrt n \exp(-nt_n^2/16) \\ 
 &\qquad\qquad +\left(\frac{b_\theta}{n}\right)^{r/\rho}\sqrt n (1 + \sqrt nt_n)  \exp\big(-n (1-r/\rho)t_n^2/16\big),
\end{align*}
where the  bound on the final term can be obtained by integration by parts. Thus, 
$$\begin{multlined}[0.9\textwidth]
\bbE\Big[  T_\theta  \cdot I(T_\theta  \leq f_n(1))\Big]  
\\ \lesssim  b_\theta^{\frac r \rho}\sqrt{\frac{\log n}{n}}  + 
 \frac{b_\theta^{r/\rho}}{\sqrt n}+\left(\frac{b_\theta}{n}\right)^{r/\rho}\sqrt {n\log n} \cdot  n^{\frac r \rho - 1}
 \lesssim b_\theta^{\frac r \rho} \sqrt{\frac{\log n}{n}}.
\end{multlined}$$
 Finally, in order to control $\bbE\left[T_\theta \cdot I(T_\theta > f_n(1))\right]$, we use  the naive bound 
\begin{align*}
T_\theta
 &= \int_{1/2}^{1} |F_{\theta,n}^{-1}(u) - F_\theta^{-1}(u)| \psi_\theta^{r-1}(u)du \\
 &\lesssim \int_{1/2}^{1}   \psi_\theta^{r}(u)du 
 \lesssim (b_{\theta,nm} + b_\theta)^{\frac r\rho } \int_{1/2}^{1} \frac{du}{(1-u)^{r/\rho}}  
 \lesssim (b_{\theta,nm} + b_\theta)^{\frac r\rho }.
\end{align*}
Using the inequality
$2r < \rho$ and Jensen's inequality, we deduce
$$\bbE\left[T_\theta^2\right] \lesssim \bbE[(b_\theta+b_{\theta,nm})^{2r/\rho}] = b_\theta^{\frac{2r}{\rho}}+ \bbE[b_{\theta,nm}]^{\frac{2r}{\rho}}= 2b_\theta^{\frac{2r}{\rho}}.$$
Thus, using the Cauchy-Schwarz inequality, we arrive at
\begin{align*}
\bbE\Big[  T_\theta  \cdot I( T_\theta  > f_n(1))\Big] 
 \leq \sqrt{\bbE\big[T_\theta^2\big] \bbP\big( T_\theta  > f_n(1)\big)} 
 \lesssim b_\theta^{\frac{r}{\rho}} \sqrt n\exp(-n/32).
 \end{align*}
We deduce that
$$\bbE[T_\theta]  = 
\bbE \Big[ T_\theta  \cdot I(  T_\theta > f_n(1))\Big] + 
\bbE \Big[ T_\theta \cdot I( T_\theta  \leq   f_n(1))\Big]  \lesssim b_\theta^{\frac{r}{\rho}} (\log n/n)^{1/2}.$$

\paragraph{Step 5: Conclusion.} By the Fubini-Tonelli Theorem and the nonnegativity of $T_\theta$, we deduce from the above display that, 
\begin{align*}
\bbE[T] 
 &=\bbE\left[ \int_{\bbS^{d-1}} T_\theta d\mu(\theta)\right] \\
 &=\int_{\bbS^{d-1}} \bbE\left[ T_\theta\right] d\mu(\theta) 
 \lesssim (\log n/n)^{1/2} \int_{\bbS^{d-1}} b_\theta^{r/\rho} d\mu(\theta) \leq b (\log n/n)^{1/2},
\end{align*}
where the final inequality follows from the assumption that $P\in \calK_{r,\rho}(b)$.
The claim follows.\qed 

\subsection{Proof of Lemma~\ref{lem::technical_deltak}}
\label{app:pf_technical_deltak}
Part (i) is trivial from the definition of $K$. To prove parts (ii)--(iv), note that for all $k \geq 1$, 
\begin{align}
\nonumber 
\gamma_{\epsilon,n}(1-\delta_{k}) 
 &=\frac{2(1-\delta_{k}) + \nu_{\epsilon,n}^2 - \nu_{\epsilon,n} \sqrt{\nu_{\epsilon,n}^2 + 4\delta_{k}(1-\delta_{k})}}{2(1+\nu_{\epsilon,n}^2)} \\
 &\geq \frac{1-\delta_{k}}{1+\nu_{\epsilon,n}^2}  - \nu_{\epsilon,n} \sqrt{\delta_{k}}.
 \label{eq:technical_lem_deltak_gamma_bound}
\end{align}
Therefore,
$$\gamma_{\epsilon,n}(1-\delta_k) 
\geq 1-\delta_k - \frac{\nu_{\epsilon,n}^2}{1+\nu_{\epsilon,n}^2} - \nu_{\epsilon,n}\sqrt{\delta_k},$$ 
which proves claim (ii). Furthermore,  
\begin{align}
\nonumber 
\frac{1-\delta_k}{1+\nu_{\epsilon,n}^2} \leq \eta_{\epsilon,n}(1-\delta_k) 
 &=\frac{2(1-\delta_k) + \nu_{\epsilon,n}^2 + \nu_{\epsilon,n} \sqrt{\nu_{\epsilon,n}^2 + 4\delta_k(1-\delta_k)}}{2(1+\nu_{\epsilon,n}^2)} \\
 \nonumber 
 &\leq (1-\delta_k) + \frac 1 2 \left[2\nu_{\epsilon,n}^2 + \nu_{\epsilon,n}\sqrt{4\delta_k}\right] \\
 &= (1-\delta_k) + \nu_{\epsilon,n}^2 + \nu_{\epsilon,n}\sqrt{\delta_k},
 \label{eq:technical_lem_deltak_eta_bound}
\end{align} 
which proves claim (iii). 
To prove claim (iv), it follows from equations~\eqref{eq:technical_lem_deltak_gamma_bound}--\eqref{eq:technical_lem_deltak_eta_bound}
that it suffices to show  
$$\frac{1-\delta_{k+1}}{1+\nu_{\epsilon,n}^2}  - \nu_{\epsilon,n} \sqrt{\delta_{k+1}} 
\geq (1-\delta_k) + \nu_{\epsilon_n}^2 + \nu_{\epsilon,n}\sqrt{\delta_k}.$$
This assertion is equivalent to 
$$ \delta_k -\delta_{k+1}  
\geq  \nu_{\epsilon,n} \Big[\sqrt{\delta_{k+1}}  +\sqrt{\delta_k} \Big] +  \nu_{\epsilon,n}^2\left[ 2-\delta_k + \nu_{\epsilon_n}^2 + \nu_{\epsilon,n}\sqrt{\delta_k} + \nu_{\epsilon,n}\sqrt{\delta_{k+1}}\right],$$
which, in turn, will be satisfied if the following inequality holds, 
$$ \delta_k -\delta_{k+1}  
\geq  2\nu_{\epsilon,n}  \sqrt{\delta_k}  +  5\nu_{\epsilon,n}^2.$$
Using a first-order Taylor expansion of the map $x \mapsto x^{-\beta}$, notice that
$\delta_k - \delta_{k+1} \geq \beta (k+1)^{-(1+\beta)}/2$. This fact together with property (i) implies that it is enough to show
$$ \frac \beta 2(k+1)^{-(1+\beta)}
\geq  \frac {2(2\delta_K)^{\frac{\rho-r}{\rho}}} c  \frac{k^{-\frac \beta 2}}{\sqrt 2}   +  \frac 5 {c^2} (2\delta_K)^{2\frac{\rho-r}{\rho}},$$
for which, in turn, it suffices to show
$$ \frac \beta 2(k+1)^{-(1+\beta)}
\geq  \frac {2k^{-\beta \left(\frac 3 2 - \frac r \rho\right)} }{ c\sqrt 2}   +  \frac 5 {c^2} k^{-2\beta \frac{\rho-r}{\rho}}.$$
By definition of $\beta$ in condition~\eqref{eq:beta_def}, we have
$$1+\beta < \beta \left(\frac 3 2 - \frac r \rho\right) \vee 2\beta \frac{\rho-r}{\rho},$$
thus for a sufficiently large choice of $c$, depending only on $\beta,\rho,r$, the penultimate display holds
for all $k \geq 1$. The claim follows. 
 \qed

\section{Proof of Theorem \ref{thm::length}}
\label{app::length}
We begin by formally stating assumptions \ref{assm::B1}-\ref{assm::B4}, referenced in the statement of 
Theorem \ref{thm::length}.
\begin{enumerate}[label=\textbf{B\arabic*}]
\item \label{assm::B1} $\gamma_{\epsilon,n}(u), \eta_{\epsilon,n}(u)$, viewed as functions of $u \in [0,1]$,
 are nondecreasing, differentiable, invertible with differentiable inverses, 
and are also respectively nondecreasing and nonincreasing functions of $\epsilon \in (0,1)$
and $n \geq 1$. 
%\item \label{assm::B2} There exists a constant $K_1 > 0$ such that for all  
%$ f \in \{\gamma_{\tau/N, n}^{-1}, \eta_{\tau/N, n}^{-1}: \tau \in \{\epsilon,\epsilon\wedge\alpha\}\}$,  
%$$\sup_{\delta \leq u \leq 1-\delta} \left|\frac{\partial f(u) }{\partial u}- 1\right| \leq K_1 \kappa_{\varepsilon,n}.$$
 \item \label{assm::B3}  There exists a constant $K_1 > 0$ such that
for all $f,g \in \{\gamma_{\tau/N, n}, \eta_{\tau/N, n},\iota: \tau \in \{\epsilon,\epsilon\wedge\alpha\}\}$, 
with $\iota$ the identity function on $[0,1]$, we have
$\delta/2 \leq f(\delta )$, $f(1-\delta) \leq 1-\delta/2$, and
%$f([\delta,1-\delta]) \subseteq \text{dom}(g^{-1})$, and
%$\delta/2 \leq g^{-1}(f(\delta))$, $g^{-1}(f(1-\delta)) \leq 1-\delta/2$, and 
$$\sup_{\delta/2 \leq u \leq 1-\delta/2} \left| \frac{\partial g^{-1}(f(u))}{\partial u} - 1\right| \leq K_1 \kappa_{\varepsilon,n}.$$
  \item \label{assm::B4} There exists a constant $K_2  > 1$ such that for all $t \in [\delta/2,1-\delta/2]$
  and all $\gamma_{\varepsilon,n}(t) \leq 
  u \leq \eta_{\varepsilon,n}(t)$ 
we have
$$ \frac 1 { K_2} \leq 
 \frac{\gamma_{\varepsilon,n}^{-1}(u)- \eta_{\varepsilon,n}^{-1}(u)} {\eta_{\varepsilon,n}(t) - \gamma_{\varepsilon,n}(t)}\leq  K_2.$$
 
 \end{enumerate} 
 
\paragraph{Proof of Theorem \ref{thm::length}.}
Throughout the proof, the symbol ``$\lesssim$'' is used to hide constants possibly depending on $K_1,K_2,\delta_0,r$. 
Furthermore, the symbol $\varkappa_N$ is used to denote
a random variable depending only on $\theta_1, \dots, \theta_N$, whose definition may
change from line to line, but which always satisfies $\bbE_{\mu^{\otimes N}}[\varkappa_N] \leq C N^{-1/2r}I(d \geq 2)$,
where $C > 0$ denotes a constant possibly depending
on $K_1,K_2,\delta,b,r$.
We prove the claim in five steps. 

\paragraph{Step 0: Setup.} With probability at least $1-\epsilon$, uniformly in $j=1, \dots, N$
and $u \in (0,1)$, we have both
\begin{equation}
\label{pf::thm1_dkw1}
F_{\theta_j,n}^{-1}\big(\gamma_{\epsilon/N,n}(u)\big) \leq F_{\theta_j}^{-1}(u) \leq F_{\theta_j,n}^{-1}\big(\eta_{\epsilon/N,n}(u)\big),
\end{equation}
and,
\begin{equation}
\label{pf::thm1_dkw2}
G_{\theta_j,m}^{-1}\big(\gamma_{\epsilon/N,m}(u)\big) \leq G_{\theta_j}^{-1}(u) \leq G_{\theta_j,m}^{-1}\big(\eta_{\epsilon/N,m}(u)\big).
\end{equation}
All derivations which follow will be carried out on the event that the above two inequalities are satisfied, 
which has probability at least $1-\epsilon$. For notational simplicty, we will write $a = \alpha/N$, $e = \epsilon/N$, 
and we recall that $\varepsilon = e \wedge a$. 

Recall that for all $\theta \in \{\theta_1, \dots, \theta_N\}$, $F_{\theta,n}$, $G_{\theta,m}$ denote the empirical CDFs of $P_\theta$ and $Q_\theta$ respectively, 
and $F_{\theta,n}^{-1}$, $G^{-1}_{\theta,m}$ their corresponding quantile functions. We may write
$$\begin{multlined}[0.85\textwidth]
C_{nm}^{(N)} = \Bigg[ \left(\frac 1 {1-2\delta} \int_{\bbS^{d-1}} \int_\delta^{1-\delta} A_{\theta,nm}^r(u) dud\mu_N(\theta)\right)^{\frac 1 r},  \\
                        \left(\frac 1 {1-2\delta} \int_{\bbS^{d-1}} \int_\delta^{1-\delta} B_{\theta,nm}^r (u)dud\mu_N(\theta)\right)^{\frac 1 r}\Bigg],
\end{multlined}$$
where 
$$
\begin{multlined}[0.75\textwidth]
A_{\theta,nm}(u)  = \big[F_{\theta,n}^{-1}\big(\gamma_{a,n}(u) \big) - G_{\theta,m}^{-1}\big(\eta_{a,m}(u)\big)\big] \\ \vee 
			    \big[G_{\theta,m}^{-1}\big(\gamma_{a,m}(u) \big) - F_{\theta,n}^{-1}\big(\eta_{a,n}(u)\big)\big] \vee 0,
\end{multlined} $$
and
\begin{align*}
B_{\theta,nm}(u) {=} \big[F_{\theta,n}^{-1}\big(\eta_{a,n}(u)\big) - G_{\theta,m}^{-1}\big(\gamma_{a,m}(u)\big)\big] {\vee }
			         \big[G_{\theta,m}^{-1}\big(\eta_{a,m}(u)\big) - F_{\theta,n}^{-1}\big(\gamma_{a,n}(u)\big)\big].
\end{align*} 

We will first show that %for some $\widebar\varkappa_N$ depending on $\mu_N$, 
$$\left|\frac 1 {1-2\delta}\int_{\bbS^{d-1}} \int_\delta^{1-\delta}A_{\theta,nm}^r(u)dud\mu_N(\theta) {-} \left[\sw_{r,\delta}^{(N)}(P,Q)\right]^r\right| 
		{\lesssim} \psi_{\varepsilon,nm} {+} \varphi_{\varepsilon,nm} {+} \varkappa_N.$$	
A similar argument can be used to bound this expression with $A_{\theta,nm}$ replaced by $B_{\theta,nm}$, 
and will lead to the claim. 

 We will assume without loss of generality that $r > 1$ in what follows. 
As will be clear from the proof, the arguments of Steps 2 and 4 alone may easily be used to prove the claim when $r=1$.

\paragraph{Step 1: Taylor Expansion Setup.}
Let $u\in [\delta,1-\delta]$, and $\theta \in \{\theta_1, \dots, \theta_N\}$. By Taylor expansion of the map
$(x,y)\in \bbR^2 \mapsto (x-y)^r$, 
there exists $\tilde F_{\theta,n}^{-1}(u)$ (resp. $\tilde G_{\theta,m}^{-1}(u)$)
on the segment joining $F_{\theta,n}^{-1}(\gamma_{a,n}(u))$ and $F_\theta^{-1}(u)$
(resp. $G_\theta^{-1}(u)$ and $G^{-1}_{\theta,m}(\eta_{a,m}(u))$)
such that
$$
\big[F_{\theta,n}^{-1} (\gamma_{a,n}(u)) - G_{\theta,m}^{-1}(\eta_{a,m}(u)) \big]^r  = 
\left[ F_\theta^{-1}(u) - G_\theta^{-1}(u) \right]^r + \xi_{nm}(u),$$
where
$$
\begin{multlined}[\textwidth]
\xi_{\theta,nm}(u)  
 = r \left(\tilde F_{\theta,n}^{-1}(u) - \tilde G_{\theta,m}^{-1}(u)\right)^{r-1}\times \\ 
 		 \Big\{\left(F_{\theta,n}^{-1}\big(\gamma_{a,n}(u)\big) - F_\theta^{-1}(u)\right)
			  -\left(G_{\theta,m}^{-1}\big(\eta_{a,m}(u)\big) - G_\theta^{-1}(u)\right)\Big\}.
\end{multlined}
$$
Likewise, there exists $\widebar F_{\theta,n}^{-1}(u)$ (resp. $\widebar G_{\theta,m}^{-1}(u)$)
on the segment joining
$F_\theta^{-1}(u)$ and $F_{\theta,n}^{-1}\big(\eta_{a,n}(u)\big)$ (resp.
$G_{\theta,m}^{-1}\big(\gamma_{a,m}(u)\big)$ and $G_\theta^{-1}(u)$) 
such that
\begin{align*}
\big[G_{\theta,m}^{-1}(\gamma_{a,m}(u)) - F_{\theta,n}^{-1} \big(\eta_{a,n}(u)\big) \big]^r  
 &= \left[G_\theta^{-1}(u) - F_\theta^{-1}(u) \right]^r + \zeta_{\theta,nm}(u),
 \end{align*}
where
$$\begin{multlined}[\textwidth]
\zeta_{\theta,nm}(u) 
 = r \Big(\widebar G_{\theta,m}^{-1}(u) - \widebar F_{\theta,n}^{-1}(u) \Big)^{r-1} \times \\
 	\Big\{\left(F_{\theta,n}^{-1}\big(\eta_{a,n}(u)\big) - F_\theta^{-1}(u)\right) 
-
 \left(G_{\theta,m}^{-1}\big(\gamma_{a,m}(u)\big) - G_\theta^{-1}(u)\right)\Big\}.
\end{multlined}$$
Now, consider the numerical inequality
$ \big| (a^r+b)\vee((-a)^r+d)\vee 0 - |a|^r \big| \leq 3 (|b| + |d|)$, for all $a \in \bbR$, $b,d \in \bbR$, 
%Now, consider the numerical inequality
%$ \big| (a+b)\vee(-a+d)\vee 0 - |a| \big| \leq c (|b| \vee |d|)$, for all $a,b,d \in \bbR$, 
%where the constant $c$ can be taken to be 1 if $|a| \geq |b| \vee |d|$ and 3 otherwise.
Taking $a = (F_\theta^{-1} - G_\theta^{-1})$, $b = \xi_{\theta,nm}$ and $d = \zeta_{\theta,nm}$,  this inequality implies
\begin{align}
\label{eq::xi-zeta-bound}
\nonumber \Bigg|&\frac 1 {1-2\delta}\int_{\bbS^{d-1}} \int_\delta^{1-\delta}A_{\theta,nm}^r(u)dud\mu_N(\theta) - \left[\sw_{r,\delta}^{(N)}(P,Q)\right]^r\Bigg| \\
\nonumber  &\leq \frac 1 {1-2\delta} \int_{\bbS^{d-1}} \int_\delta^{1-\delta} 
 		  \bigg|
 			\Big[\big(F_\theta^{-1}(u) - G_\theta^{-1}(u)\big)^r + \xi_{\theta,nm}(u)  \Big] \\ \nonumber & \vee 
			\Big[\big(G_\theta^{-1}(u) - F_\theta^{-1}(u)\big)^r + \zeta_{\theta,nm}(u) \Big] \vee 0 - 
     	    \big|F_\theta^{-1}(u) - G_\theta^{-1}(u)\big|^r \bigg| dud\mu_N(\theta)\\ 
 &\lesssim  \int_{\bbS^{d-1}} \int_\delta^{1-\delta} |\zeta_{\theta,nm}(u)| du d\mu_N(\theta) +  \int_{\bbS^{d-1}} \int_\delta^{1-\delta} |\xi_{\theta,nm}(u)| du d\mu_N(\theta).
\end{align}
It will now suffice to bound the second term of the above display, and a similar bound will
hold for the first. 
Note that
\begin{align*}
\int_{\bbS^{d-1}} \int_\delta^{1-\delta} |\xi_{\theta,nm}(u)|du d\mu_N(\theta)
 &\leq r(\calI + \calJ),
 \end{align*}
 where,
 \begin{align}
 \calI &= 
  \int_{\bbS^{d-1}} \int_\delta^{1-\delta} \big|\tilde F_{\theta,n}^{-1}(u) - \tilde G_{\theta,m}^{-1}(u)\big|^{r-1} \big|F_{\theta,n}^{-1}\big(\gamma_{a,n}(u)\big) - F_\theta^{-1}(u)\big|du d\mu_N(\theta) \\
 \calJ &= 
  \int_{\bbS^{d-1}} \int_\delta^{1-\delta} \big|\tilde F_{\theta,n}^{-1}(u) - \tilde G_{\theta,m}^{-1}(u)\big|^{r-1} \big|G_{\theta,m}^{-1}\big(\eta_{a,m}(u)\big) - G_\theta^{-1}(u)\big|du d\mu_N(\theta).
 \label{pf::calJ_def} 
\end{align}
It will suffice to prove
that $\calI \lesssim\psi_{\varepsilon,nm}$ and 
$\calJ \lesssim \varphi_{\varepsilon,nm}$, up to terms depending only on $N$. 
We consider the cases $\sj_{r,\delta/2}(P) \vee \sj_{r,\delta/2}(Q) = \infty$
and $\sj_{r,\delta/2}(P) \vee \sj_{r,\delta/2}(Q) < \infty$ seperately.  

\paragraph{Step 2: Bounding $\calI$ and $\calJ$ when $\sj_{r,\delta/2}(P) \vee \sj_{r,\delta/2}(Q) =\infty$.}
We have,
$$
\begin{multlined}[0.9\textwidth]
\calI 
   \lesssim \int_{\bbS^{d-1}} \left(\sup_{\delta \leq u \leq 1-\delta} \left|\tilde F_{\theta, n}^{-1}(u) - \tilde G_{\theta,m}^{-1}(u)\right|^{r-1} \right) \times  \\
    \left(\int_\delta^{1-\delta}\left|F_{\theta,n}^{-1}(\gamma_{a,n}(u)) - F_\theta^{-1}(u)\right|du\right) d\mu_N(\theta).
\end{multlined}$$
We will bound each factor in the above integral, beginning with the second. Using inequality \eqref{pf::thm1_dkw1}, since $e \geq \varepsilon$, we have
for all $u \in [\delta,1-\delta]$ and $\theta \in \{\theta_1, \dots, \theta_N\}$,
$$\begin{multlined}[0.95\textwidth]
|F^{-1}_{\theta,n}(\gamma_{a,n}(u)) - F_\theta^{-1}(u)| \\ \leq 
	\big[F_\theta^{-1}\big(\gamma_{\varepsilon,n}^{-1}(\gamma_{a,n}(u))\big) - F_\theta^{-1}(u)\big] +
 \big[F_\theta^{-1}(u) - F_\theta^{-1}\big(\eta_{\varepsilon,n}^{-1}(\gamma_{a,n}(u))\big)\big].
\end{multlined}$$ 
Now, write
$x_n = 	\gamma_{\varepsilon,n}^{-1}(\gamma_{a,n}(\delta))$ and 
$y_n = \eta_{\varepsilon,n}^{-1}(\gamma_{a,n}(1-\delta))$.  By condition~\ref{assm::B3} and
the  definition of $\kappa_{\varepsilon,n}$, we have 
\begin{equation}
\label{eq:x_n_y_n}
|x_n-\delta| \vee |y_n - 1 - \delta| \leq \kappa_{\varepsilon,n},
\end{equation}
which, combined with the assumption that $\kappa_{\varepsilon,n} \leq \frac \delta 2 \wedge (1-2\delta)$, also implies 
$$\delta \leq x_n \leq 1-\delta \leq y_n\leq 1-\delta/2.$$
Thus,  for all $\theta \in \{\theta_1, \dots, \theta_N\}$,
\begin{align*}
 &\int_\delta^{1-\delta}\big[F_\theta^{-1}\big(\gamma_{\varepsilon,n}^{-1}(\gamma_{a,n}(u))\big) - F_\theta^{-1}(u)\big] du \\
 &~~=  \int_{x_n}^{y_n}
 	F_\theta^{-1}(u) \left(\frac{\partial \gamma_{a,n}^{-1}(\gamma_{\varepsilon,n}(u))}{\partial u}\right) du 
 		- \int_\delta^{1-\delta} F_\theta^{-1}(u) du \\
 &~~\leq  \int_{x_n}^{y_n}
 		  F_\theta^{-1}(u)   du  +   K_1\kappa_{\varepsilon,n} \int_{x_n}^{y_n} | F_\theta^{-1}(u) | du 
 		-  \int_\delta^{1-\delta} F_\theta^{-1}(u) du, \qquad  (\text{By \ref{assm::B3}}) \\ 		
 &~~\leq \int_{1-\delta}^{y_n} F_\theta^{-1}(u)du -  \int_{\delta}^{x_n} F_\theta^{-1} (u)du 
 		+ K_1\kappa_{\varepsilon,n}  |F_\theta^{-1}(y_n)|  \\			
 &~~\lesssim (y_n -1  	 +\delta) \Big[| F_\theta^{-1}(y_n) | \vee |F_\theta^{-1}(1-\delta)|\Big] \\ &~~+
 	   (x_n-\delta) \Big[ | F_\theta^{-1}(\delta)| \vee |F_\theta^{-1}(x_n)|\Big]   + \kappa_{\varepsilon,n} |F_\theta^{-1}(y_n) | \\			    		 		
&~~\lesssim 
  (y_n - 1+\delta+\kappa_{\varepsilon,n})|F_\theta^{-1}(1-\delta/2)|  + (x_n - \delta) |F_\theta^{-1}(\delta)|.  
\end{align*} 
Since $P \in \widebar\calK_2(b)$, the quantiles in the above display are bounded above by 
a universal multiple of $(b/\delta)^{1/2}$, by
Lemma~\ref{lem::calK}.
Thus, together with equation~\eqref{eq:x_n_y_n}, we arrive at
\begin{align*}
\int_\delta^{1-\delta}\big[F_\theta^{-1}\big(\gamma_{\varepsilon,n}^{-1}(\gamma_{a,n}(u))\big) - F_\theta^{-1}(u)\big] du 
\lesssim \kappa_{\varepsilon,n}(b/\delta)^{1/2},
\end{align*}
We similarly have,
\begin{align*}
\int_\delta^{1-\delta} \big[F_\theta^{-1}(u) - F_\theta^{-1}\big(\eta_{\varepsilon,n}^{-1}(\gamma_{a,n}(u))\big)\big] du 
\lesssim \kappa_{\varepsilon,n}(b/\delta)^{1/2}.
\end{align*}
Combining these facts, we obtain
\begin{equation}
\label{eq::calI_bound}
\calI \lesssim
 \kappa_{\varepsilon,n} (b/\delta)^{1/2}\int_{\bbS^{d-1}} \sup_{\delta \leq u \leq 1-\delta} \left|\tilde F_{\theta, n}^{-1}(u) - \tilde G_{\theta,m}^{-1}(u)\right|^{r-1} d\mu_N(\theta).
\end{equation}
We now bound the second factor in the above display. Since we have
$\tilde F_{\theta,n}^{-1}(u) \in [F_{\theta,n}^{-1}(\gamma_{a,n}(u)),F_\theta^{-1}(u)]$
and
$\tilde G_{\theta,m}^{-1}(u) \in [G_\theta^{-1}(u), G^{-1}_{\theta,m}(\eta_{a,m}(u))]$, we deduce that for any $u \in [\delta,1-\delta]$,
\begin{align}
\label{pf::thm1_step}
\Big| & \tilde F_{\theta,n}^{-1}(u) - \tilde G_{\theta,m}^{-1}(u)\Big| \\
%\nonumber  &\leq \Big[ G_{\theta,m}^{-1}(\eta_{a,m}(u)) - F_{\theta,n}^{-1}(\gamma_{a,n}(u))\Big] \vee \Big[ F_\theta^{-1}(u) - G_\theta^{-1}(u)\Big]\\
%\nonumber  &\leq \Big[ G_\theta^{-1}(\gamma_{\varepsilon,m}^{-1}(\eta_{a,m}(u))) - 
% 					F_\theta^{-1}(\eta_{\varepsilon,n}^{-1}(\gamma_{a,n}(u)))\Big] \vee \Big[ F_\theta^{-1}(u) - G_\theta^{-1}(u)\Big]  \\  
\nonumber &\leq  \Big| G_\theta^{-1}(u) - F_\theta^{-1}(u) \Big| + 
                    \Big| F_\theta^{-1}(u) - \tilde F_{\theta,n}^{-1}(u) \Big| +
	                \Big| \tilde G_{\theta,m}^{-1}(u) - G_\theta^{-1}(u) \Big| \\
\nonumber &\leq  \Big| G_\theta^{-1}(u) - F_\theta^{-1}(u) \Big|  \\ &\quad + 
                    \Big| F_\theta^{-1}(u) - F_\theta^{-1}\big(\eta_{\varepsilon,n}^{-1}(\gamma_{a,n}(u))\big) \Big| +
	                \Big| G_\theta^{-1}\big(\gamma_{\varepsilon,m}^{-1}(\eta_{a,m}(u))\big) - G_\theta^{-1}(u) \Big|.
\end{align}
It follows that
\begin{align*}
\int_{\bbS^{d-1}} \sup_{\delta \leq u \leq 1-\delta}  &\left|\tilde F_{\theta,n}^{-1}(u) - \tilde G_{\theta,m}^{-1}(u)\right|^{r-1} d\mu_N(\theta) \\
    &\lesssim \left\{ \sw_{\infty,\delta}^{(r-1)}(P,Q) + U_{\varepsilon,n}(P) + U_{\varepsilon,m}(Q) \right\} + \varkappa_N.
\end{align*} 
We conclude this section of the proof by combining the above display with equation \eqref{eq::calI_bound}. We then have
$$ \calI 
     \lesssim \kappa_{\varepsilon,n}
     (b/\delta)^{1/2}  \left\{\sw_{\infty,\delta}^{(r-1)}(P,Q) + U_{\varepsilon,n}(P) + U_{\varepsilon,m}(Q) + \varkappa_N\right\}
      \lesssim  \psi_{\varepsilon,nm} + \varkappa_N ,$$
and by a symmetric argument, 
$$ \calJ 
%   \lesssim   \big( \kappa_{\varepsilon,m} + \varkappa_N^{\big)\left\{ \sw_{\infty,\delta}^{r-1}(P,Q) + U_{e,n}(P) + U_{e,m}(Q) + \varkappa_N\right\}
      \lesssim  \varphi_{\varepsilon,nm} + \varkappa_N.$$
     
\paragraph{Step 3: Bounding $\calI$ and $\calJ$ when $\sj_{r,\delta/2}(P) \vee \sj_{r,\delta/2}(Q)  < \infty$.}
By means of H\"older's inequality, we have
\begin{align*}
\nonumber
\calI
 &\leq  \int_{\bbS^{d-1}}\int_\delta^{1-\delta}
 				\big|\tilde F_{\theta,n}^{-1}(u) - \tilde G_{\theta,m}^{-1}(u)\big|^{r-1} \big|F_{\theta,n}^{-1}(\gamma_{a,n}(u)) - F_\theta^{-1}(u)\big|du d\mu_N(\theta) \\
 &\nonumber\leq \left(\int_{\bbS^{d-1}}\int_\delta^{1-\delta} \big|\tilde F_{\theta,n}^{-1}(u) - \tilde G_{\theta,m}^{-1}(u)\big|^r du d\mu_N(\theta)\right)^{\frac{r-1}{r}}\times \\ &
              \qquad~~\left(\int_{\bbS^{d-1}}\int_\delta^{1-\delta} \big|F_{\theta,n}^{-1}\big(\gamma_{a,n}(u)\big) - F_\theta^{-1}(u)\big|^r dud\mu_N(\theta) \right)^{\frac 1 r} \\
\nonumber              
 &\nonumber \lesssim \left(\int_{\bbS^{d-1}} \int_\delta^{1-\delta} \big|\tilde F_{\theta,n}^{-1}(u) - \tilde G_{\theta,m}^{-1}(u)\big|^{r} du d\mu_N(\theta)\right)^{\frac{r-1}{r}} 
 		 \times \\ 
 		  & \qquad~~ \left(\int_{\bbS^{d-1}}\int_\delta^{1-\delta} \big[F_\theta^{-1}\big(\gamma_{\varepsilon,n}^{-1}(\gamma_{a,n}(u))\big) - 
 		  F_\theta^{-1}\big(\eta_{\varepsilon,n}^{-1}(\gamma_{a,n}(u))\big)\big]^r du d\mu_N(\theta) \right)^{\frac 1 r}\\
 &\nonumber = \left(\int_{\bbS^{d-1}} \int_\delta^{1-\delta} \big|\tilde F_{\theta,n}^{-1}(u) - \tilde G_{\theta,m}^{-1}(u)\big|^{r} du d\mu_N(\theta)\right)^{\frac{r-1}{r}} 
 		 \times \\ 
 		  & \hspace{-0.05in}\left(\int_{\bbS^{d-1}}\int_{\gamma_{a,n}(\delta)}^{\gamma_{a,n}(1-\delta)} \big[F_\theta^{-1}\big(\gamma_{\varepsilon,n}^{-1}(u)\big) - F_\theta^{-1}(\eta_{\varepsilon,n}^{-1}(u))\big]^r \left(\frac{\partial \gamma_{a,n}^{-1}(u)}{\partial u}\right)  dud\mu_N(\theta)\right)^{\frac 1 r} \\
 &\nonumber \leq  \left(\int_{\bbS^{d-1}} \int_\delta^{1-\delta} \big|\tilde F_{\theta,n}^{-1}(u) - \tilde G_{\theta,m}^{-1}(u)\big|^{r} du d\mu_N(\theta)\right)^{\frac{r-1}{r}} 
 		 \times \\ 
 		  & \hspace{-0.05in}\left(\int_{\bbS^{d-1}}\int_{\gamma_{a,n}(\delta)}^{\gamma_{a,n}(1-\delta)} \big[F_\theta^{-1}\big(\gamma_{\varepsilon,n}^{-1}(u)\big) - F_\theta^{-1}(\eta_{\varepsilon,n}^{-1}(u))\big]^r (1+K_1\kappa_{\varepsilon,n}) dud\mu_N(\theta)\right)^{\frac 1 r},
\end{align*} 
where we used condition~\ref{assm::B3} on the final line. 
We now reason similarly as in Appendix~\ref{app::high_prob_sw_null}.
  Since $\sj_{r,\delta/2}(P)  < \infty$, it follows from  
Lemma \ref{lem::abs_cont}
that 
$F_\theta^{-1}$ is absolutely continuous over $[\delta/2,1-\delta/2]$ for $\mu$-almost every $\theta \in \bbS^{d-1}$.  
We thus have for have almost surely that for each $\theta \in \{\theta_1, \dots, \theta_N\}$, 
\begin{align*}
&\int_{\gamma_{a,n}(\delta)}^{\gamma_{a,n}(1-\delta)} \big[ F_\theta^{-1}\big(\gamma_{\varepsilon,n}^{-1}(u)\big) - F_\theta^{-1}(\eta_{\varepsilon,n}^{-1}(u))\big]^r  du \\
 &\quad= \int_{\gamma_{a,n}(\delta)}^{\gamma_{a,n}(1-\delta)}\left( \int_{\eta_{\varepsilon,n}^{-1}(u)}^{\gamma_{\varepsilon,n}^{-1}(u)} \frac {dt} {p_\theta(F_\theta^{-1}(t))}\right)^rdu \\
 &\quad\leq \int_{\gamma_{a,n}(\delta)}^{\gamma_{a,n}(1-\delta)} (\gamma_{\varepsilon,n}^{-1}(u)- \eta_{\varepsilon,n}^{-1}(u))^{r-1}
 \int_{\eta_{\varepsilon,n}^{-1}(u)}^{\gamma_{\varepsilon,n}^{-1}(u)}
 \left(  \frac {1} {p_\theta(F_\theta^{-1}(t))}\right)^{r}dtdu \\
 &\quad\leq % \int_{\eta_{\varepsilon,n}^{-1}(\gamma_{a,n}(\delta))}^{\gamma_{\varepsilon,n}^{-1}(\gamma_{a,n}(1-\delta))}
 \int_{\delta/2}^{1-\delta/2}\left( \int_{\gamma_{\varepsilon,n}(t)}^{\eta_{\varepsilon,n}(t)}
   (\gamma_{\varepsilon,n}^{-1}(u)- \eta_{\varepsilon,n}^{-1}(u))^{r-1} du\right)\left(\frac {1} {p_\theta(F_\theta^{-1}(t))}\right)^rdt \\
 % ~ (\text{By } \eqref{eq:x_n_y_n})\\ 
 &\quad\leq %\int_{\gamma_{a,n}(\delta)}^{\gamma_{a,n}(1-\delta)} 
 \int_{\delta/2}^{1-\delta/2}\left(\sup_{\gamma_{\varepsilon,n}(t) \leq u \leq \eta_{\varepsilon,n}(t)} (\gamma_{\varepsilon,n}^{-1}(u)- \eta_{\varepsilon,n}^{-1}(u) )^{r-1} \right) \frac {\eta_{\varepsilon,n}(t) - \gamma_{\varepsilon,n}(t)} {[p_\theta(F_\theta^{-1}(t))]^r}dt \\
 &\quad\lesssim %\int_{\gamma_{a,n}(\delta)}^{\gamma_{a,n}(1-\delta)}
 \int_{\delta/2}^{1-\delta/2}\left( \frac {\eta_{\varepsilon,n}(t) - \gamma_{\varepsilon,n}(t)} {p_\theta(F_\theta^{-1}(t))}\right)^rdt
\lesssim \int_{\delta/2}^{1-\delta/2}%\int_{\gamma_{a,n}(\delta)}^{\gamma_{a,n}(1-\delta)} 
 \left(  \frac {\tilde \kappa_{\varepsilon,n}(t)} {p_\theta(F_\theta^{-1}(t))}\right)^rdt,
 %\leq \int_{\delta/2}^{1-\delta/2} \left(  \frac {\tilde \kappa_{\varepsilon,n}(t)} {p_\theta(F_\theta^{-1}(t))}\right)^rdt,
\end{align*}
where we repeatedly used assumption~\ref{assm::B4} on the final line.
%and again the fact that $[\gamma_{a,n}(\delta), \gamma_{a,n}(1-\delta)] \subseteq [\delta/2,1-\delta/2]$. 
We thus arrive at,
\begin{align}
\label{eq::pf_holder}
 \calI  
 &\lesssim \left( \int_{\bbS^{d-1}} \int_\delta^{1-\delta} \big|\tilde F_{\theta,n}^{-1}(u) - \tilde G_{\theta,m}^{-1}(u)\big|^{r} du d\mu_N(\theta)\right)^{\frac{r-1}{r}} 
 		 \big[ V_{\varepsilon,n}(P)\big]^{\frac 1 r}+ \varkappa_N.
\end{align}
By using similar calculations as in equations \eqref{pf::thm1_step} and \eqref{eq::pf_holder} to bound the first factor in the above display, we have
\begin{alignat*}{1}
&\int_{\bbS^{d-1}} \int_\delta^{1-\delta}\Big(\tilde F_{\theta,n}^{-1}(u ) - \tilde G_{\theta,m}^{-1}(u)\Big)^rdu d\mu_N(\theta)  \\
 &\qquad \lesssim\sw_{r,\delta}^r(P,Q)+ 
 			\int_{\bbS^{d-1}} \int_\delta^{1-\delta} \big| F_\theta^{-1}(u) - F_\theta^{-1}\big(\eta_{\varepsilon,n}^{-1}(\gamma_{a,n}(u))\big)\big|^r du d\mu_N(\theta)
 			\\& \qquad \qquad\qquad\qquad ~~ +
 			 \int_{\bbS^{d-1}} \int_\delta^{1-\delta} \big| G_\theta^{-1}\big(\gamma_{\varepsilon,m}^{-1}(\gamma_{a,m}(u))\big)-G_\theta^{-1}(u)\big|^r du d\mu_N(\theta) \\
 &\qquad \lesssim  \sw_{r,\delta}^r(P,Q) + V_{\varepsilon,n}(P) + V_{\varepsilon,m}(Q) + \varkappa_N, 			
\end{alignat*} 
Putting these facts together with equation \eqref{eq::pf_holder}, we arrive at
$$  \calI \lesssim \left(\sw_{r,\delta}^r(P,Q) + V_{\varepsilon,n}(P) +  V_{\varepsilon,m}(Q) \right)^{\frac{r-1}{r}} 
 		 \big[ V_{\varepsilon,n}(P)\big]^{\frac 1 r} + \varkappa_N = \psi_{\varepsilon,nm} + \varkappa_N.$$
Finally, by a symmetric argument, we also have $\calJ \lesssim \varphi_{\varepsilon,nm} + \varkappa_N.$

\paragraph{Step 4: Conclusion.}
Returning to equation \eqref{eq::xi-zeta-bound}, we have shown
$$\int_{\bbS^{d-1}} \int_{\delta}^{1-\delta} |\xi_{nm}(u)|du d\mu_N(\theta) \lesssim \psi_{\varepsilon,nm} + \varphi_{\varepsilon,nm} + \varkappa_N.$$
By the same arguments, we may obtain the same upper bound, up to universal constant factors, 
on the integral $\int_{\bbS^{d-1}} \int_{\delta}^{1-\delta} |\zeta_{nm}(u)|du d\mu_N(\theta)$ in 
equation \eqref{eq::xi-zeta-bound}.
%A symmetric argument can be used to bound $\int |\zeta_{nm}|$ in \eqref{eq::xi-zeta-bound}.
We deduce that, for some $c_1 > 0$ (possibly depending on $r$), we have
\begin{align*}
 &\left(\frac 1 {1-2\delta}\int_{\bbS^{d-1}} \int_\delta^{1-\delta} A_{nm}^r(u) du d\mu_N(\theta) \right)^{\frac 1 r}  \\
 &\qquad \qquad~~ \geq \left\{ \big[\sw_{r,\delta}^{(N)}(P,Q)\big]^r - c_1 \big(\psi_{\varepsilon,nm} + \varphi_{\varepsilon,nm} + \varkappa_N\big) \right\}^{\frac 1 r}\\
% \left(\frac 1 {1-2\delta}\int_{\bbS^{d-1}} \int_\delta^{1-\delta} A_{nm}^r(u) du d\mu_N(\theta) \right)^{\frac 1 r} 
 &\qquad \qquad~~ \geq \Big\{ \sw_{r,\delta}^r(P,Q) - c_1 \big(\psi_{\varepsilon,nm} + \varphi_{\varepsilon,nm} + \varkappa_N\big) \Big\}^{\frac 1 r}.
\end{align*}
By the same arguments, there exists a constant $c_2 > 0$  such that
$$\begin{multlined}[0.9\textwidth]
 \left(\frac 1 {1-2\delta}\int_{\bbS^{d-1}} \int_\delta^{1-\delta} B_{nm}^r(u) du d\mu_N(\theta) \right)^{\frac 1 r} \\
 \leq \left\{ \sw_{r,\delta}^r(P,Q) + c_2 \big(\psi_{\varepsilon,nm} + \varphi_{\varepsilon,nm} + \varkappa_N \big) \right\}^{\frac 1 r}.
 \end{multlined}
$$
The claim follows by choosing $c = c_1 \vee c_2$. \qed

\section{Proof of Theorem \ref{thm::bootstrap}} 
\label{app::boot}

The proof of this result has two main components. In Lemma~\ref{lem::hadamard_sw} we show that 
the Sliced Wasserstein distance is Hadamard differentiable under certain conditions. Theorem~\ref{thm::bootstrap} then
follows %from this Hadamard differentiability 
via an application of the functional delta method.

\paragraph{Hadamard Differentiability of the Sliced Wasserstein Distance.} 
Throughout this subsection, for a metric space $(T,\rho)$, $C[T]$ denotes the
set of real-valued continuous functions defined on~$T$, endowed with the supremum norm, and
$\ell^\infty(T) = \left\{f:T \to \bbR: \sup_{t \in T} |f(t)| < \infty\right\}.$
Let $D[I]$ denote the space of c\`adl\`ag functions defined over an interval $I = [a_1,a_2]\subseteq \bbR$, %also
endowed with the supremum norm. 
We will make use of the following result from \cite{vandervaart1996} (Lemma 3.9.20), regarding
the Hadamard differentiability of the quantile function at a fixed point $u \in (a_1,a_2)$. Let $\bbD_\psi$ 
denote
the set of nondecreasing maps $A\in D[I]$ such that the set $\{x \in I: A(x) \geq u\}$ is nonempty
for any given $u \in (0,1)$, and define the map
\begin{equation}
\label{eq:map_psi}
\psi: \bbD_\psi \subseteq D[I] \to \bbR, \quad \psi: A \mapsto A^{-1}(u) = \inf\{x \in I: A(x) \geq u\}.
\end{equation}
\begin{lemma}[\cite{vandervaart1996}]
\label{lem:quantile_hadamard}
Let $A \in \bbD_\psi$ satisfy the following two properties.
\begin{enumerate}
\item[(i)] $A$ is differentiable at a point $\xi_u \in (a_1,a_2)$ such that $A(\xi_u) = u$.
\item[(ii)] $A$ has strictly positive derivative at $\xi_u$. 
\end{enumerate}
Then, $\psi$ is Hadamard-differentiable at $A$ tangentially to the set
of functions $H \in D[I]$ which are continuous at $\xi_u$, with Hadamard derivative given by
$$\psi'_A(H) = -\frac{H(\xi_u)}{A'(\xi_u)}.$$
\end{lemma}
Now, define $\calH = \bbR \times \bbS^{d-1}$, identified with the set of half-spaces in $\bbR^d$. 
Let $\bbD = \ell^\infty(\calH)$, and let $\bbD_0$ denote the subspace of $\bbD$ consisting of maps $F: \calH \to \bbR$
such that $F(\cdot,\theta) \in C[\bbR]$ for  all $\theta \in \bbS^{d-1}$. 
 %, namely
%$$\calH = \left\{H_{\theta,x}: \theta \in \bbS^{d-1}, x \in \bbR\right\}, \quad
%   H_{\theta,x} = \{y\in \bbR^d: \langle y, \theta \rangle \leq x\}.$$
Furthermore, define $\bbD_\phi$ as the subset of maps $F:\calH \to \bbR$
such that $F(\cdot,\theta) \in D[\bbR]$ is a CDF for all $\theta \in \bbS^{d-1}$.
Define the map
$$\phi: \bbD_\phi^2 % \subseteq \ell^\infty(\calH) \times \ell^\infty(\calH) 
\to \bbR_+, \quad
  \phi: (F,G) \mapsto \frac 1 {1-2\delta} \int_{\bbS^{d-1}} \int_\delta^{1-\delta} \big| F^{-1}(u,\theta) - G^{-1}(u,\theta) \big|^rdud\mu(\theta),$$
for a fixed constant $\delta \in [0,1/2)$, where we interchangeably employ the notation $F^{-1}(u, \theta) = F_\theta^{-1}(u) = \inf\{x \in \bbR: F_\theta(x) \geq u\}$ 
and $F(\cdot,\theta) = F_\theta(\cdot)$ 
in this section only.

The Hadamard differentiability of $\phi$, tangentially to $\bbD_0$,
is established below.  
\begin{lemma}
\label{lem::hadamard_sw} 
%Let $P,Q\in \widebar\calK(b)$ for some $b \geq 1$, and assume 
%the measures ${\pi_\theta}_\# P$ and ${\pi_\theta}_\#Q$ are asolutely continuous for $\mu$-almost every $\theta \in \bbS^{d-1}$, 
%with respective densities $p_\theta$ and $q_\theta$, and respective CDFs
% $F(\cdot,\theta)$ and $G(\cdot,\theta)$.
Assume the same conditions as in Theorem~\ref{thm::bootstrap}. 
Then, the map $\phi$ is Hadamard differentiable at $(F, G)$, tangentially to $\bbD_0$, with Hadamard derivative given by
\begin{align*}
\phi'&: \bbD_0^2 \to \bbR, \\ 
&\phi'(H_1,H_2) = \frac r {1-2\delta} \int_{\bbS^{d-1}} \int_\delta^{1-\delta} \mathrm{sgn}\left(F^{-1}(u,\theta) - G^{-1}(u,\theta)\right)  \times 
\\ & %\qquad~~~
\big|F^{-1}(u,\theta) - G^{-1}(u,\theta)\big|^{r-1}
 	 \left(\frac{H_2(G^{-1}(u,\theta),\theta)}{q_\theta(G^{-1}(u,\theta))} - \frac{H_1(F^{-1}(u,\theta),\theta)}{p_\theta(F^{-1}(u,\theta))}\right)dud\mu(\theta).
 	 \end{align*}\end{lemma} 

\paragraph{Proof of Lemma \ref{lem::hadamard_sw}.}
Let $(H_{1k})_{k=1}^\infty, (H_{2k})_{k=1}^\infty\subseteq \bbD$ be sequences satisfying $F + t_kH_{1k}, G + t_kH_{2k} \in\bbD_\phi$ 
for all $k\geq 1$, and such that $H_{jk}$ converges uniformly to $H_j \in \bbD_0$, $j=1, 2$. 
Let $t_k \downarrow 0$ as $k \to \infty$, and define for all $k \geq 1$,
$$
\begin{multlined}[0.9\textwidth]
\Delta_k
 =\frac 1 {t_k(1-2\delta)} \int_{\bbS^{d-1}} \int_\delta^{1-\delta}\Big\{
			\big| (F + t_k H_{1k})^{-1}(u,\theta) - (G+t_kH_{2k})^{-1}(u,\theta) \big|^r 
			\\- 
 		    \big| F^{-1}(u,\theta) - G^{-1}(u,\theta)|^r \Big\}du d\mu(\theta).
\end{multlined} 		   
$$
We will prove that the limit of $\Delta_k$ exists when taking $k \to \infty$.
For all $r > 1$, the map $(x,y)  \in \bbR^2 \mapsto |x-y|^r$
is continuously differentiable. Thus, for all $u \in [\delta,1-\delta]$
and all $\theta \in \bbS^{d-1}$, 
there exists $\tilde F_k^{-1}(u,\theta)$ (resp, $\tilde G_k^{-1}(u,\theta)$) on the line joining
$F^{-1}(u,\theta)$ (resp. $G^{-1}(u,\theta)$) and $(F+t_kH_{1k})^{-1}(u,\theta)$ 
(resp. $(G + t_kH_{2k})^{-1}(u,\theta)$) such that 
\begin{equation}
\label{eq::Deltak}
\begin{multlined}[0.9\textwidth]
\Delta_k = \frac r {t_k(1-2\delta)}\int_{\bbS^{d-1}} \int_\delta^{1-\delta}
 	\varphi\big(\tilde F_k^{-1}(u,\theta); \tilde G_k^{-1}(u,\theta) \big) \\[0.05in] \times \Big\{ 
 		\big[(F+t_k H_{1k})^{-1}(u,\theta) - F^{-1}(u,\theta)\big] \\
 		 - \big[(G+t_k H_{2k})^{-1}(u,\theta) - G^{-1}(u,\theta)\big] \Big\}
 		dud\mu(\theta),
 		\end{multlined}
\end{equation} 
where
$\varphi(x;y) = \text{sgn}(x-y)|x-y|^{r-1}$.
We will now argue that each of the limits
\begin{align*}
B_1(u,\theta) &= \lim_{k\to\infty} B_{1k}(u,\theta), \quad B_{1k}(u,\theta) = \varphi\big(\tilde F_k^{-1}(u,\theta);\tilde G_k^{-1}(u,\theta)\big),\\
B_2(u,\theta) &= \lim_{k\to \infty} B_{2k}(u,\theta), \quad B_{2k}(u,\theta) = 
 		\frac{(F+t_k H_{1k})^{-1}(u,\theta) - F^{-1}(u,\theta)}{t_k},\\
B_3(u,\theta) &= \lim_{k\to \infty} B_{3k}(u,\theta), \quad B_{3k}(u,\theta) = 
 		\frac{(G+t_k H_{2k})^{-1}(u,\theta) - G^{-1}(u,\theta)}{t_k},
\end{align*}
exist  for $(\lambda \otimes \mu)$-almost every $(u,\theta) \in [\delta, 1-\delta] \times \bbS^{d-1}$.
Throughout the sequel, we write $I_\theta = [F^{-1}(\delta/2,\theta), F^{-1}(1-\delta/2,\theta)]$
for all $\theta \in \bbS^{d-1}$. By Lemma~\ref{lem::calK}, there exist finite constants $a_1 < a_2$, depending
on $\delta$ and the finite second moments of $P$ and $Q$, such that $\bigcup_{\theta \in \bbS^{d-1}} I_\theta \subseteq [a_1,a_2]$.
We fix $I = [a_1,a_2]$ in what follows.

Regarding the limit $B_1$, we make use of the following observation, which we  prove below 
in Section~\ref{app:pf_lem_quantile_technical_uniform}.
The conclusion of this assertion is stronger than necessary, but will be needed again in the sequel.
\begin{lemma}
\label{lem:quantile_technical_uniform}
Assume the same conditions as Theorem~\ref{thm::bootstrap}. Then, for all $u \in [\delta,1-\delta]$, we have as $k \to \infty$,
\begin{align*}
\sup_{ \theta \in \bbS^{d-1}} |(F+t_kH_{1k})^{-1}(u,\theta) - F^{-1}(u,\theta)| &\to 0, \\
\sup_{ \theta \in \bbS^{d-1}} |(G+t_kH_{2k})^{-1}(u,\theta) - G^{-1}(u,\theta)| &\to 0.
\end{align*}
\end{lemma}   
Notice further that the map $\varphi$ is continuous in both of its arguments, therefore  we obtain
from Lemma~\ref{lem:quantile_technical_uniform} that
$$B_1(u,\theta) = \varphi\big(F^{-1}(u,\theta);G^{-1}(u,\theta)\big),$$% = \text{sgn}\big(F^{-1}(u,\theta) - G^{-1}(u,\theta)\big) |F^{-1}(u,\theta) - G^{-1}(u,\theta)|^{r-1},$$
for $(\lambda \otimes \mu)$-almost every $(u,\theta)$.
 
We now turn to the limit $B_2$. Recall that $A := F(\cdot,\theta)$ is absolutely continuous
for any given $\theta \in \bbS^{d-1}$. %Fixing such a choice of $\theta$,  let $A = F(\cdot,\theta)$. 
 For any fixed $u \in [\delta,1-\delta]$, the existence of $B_2(u,\theta)$
would be implied by the Hadamard differentiability of the map $\psi$ in equation \eqref{eq:map_psi} at $A$,
tangentially to $\bbD_0$, 
sufficient conditions for which are given by conditions (i) and (ii) of Lemma~\ref{lem:quantile_hadamard}.  
Condition (i)  is immediately satisfied for almost all $u \in [\delta,1-\delta]$
due to the absolute continuity of $A$. % for almost all $u \in [\delta,1-\delta]$. 
Furthermore, 
the assumption $J_{\infty,\delta/2}(P)<\infty$ implies 
 that $p_\theta$ is nonzero at $A^{-1}(u)$
for almost every $u \in [\delta,1-\delta]$, implying that condition (ii) of Lemma~\ref{lem:quantile_hadamard} is satisfied
for all such $u$. We deduce from Lemma~\ref{lem:quantile_hadamard} the limit
$$B_2(u,\theta) = - \frac{H_1(A^{-1}(u),\theta)}{p_\theta(A^{-1 }(u))} = - \frac{H_1(F^{-1}(u,\theta),\theta)}{p_\theta(F^{-1 }(u,\theta))},$$
for $(\lambda \otimes \mu)$-almost every $(u,\theta)$. We similarly obtain, almost everywhere, 
$$B_3(u,\theta) =   - \frac{H_2(G^{-1}(u,\theta),\theta)}{q_\theta(G^{-1}(u,\theta))}.$$
With these facts in place, the claim will 
follow from the Dominated Convergence Theorem if we are able to interchange the limit as $k \to \infty$
with the integrations in equation~\eqref{eq::Deltak}. To this end, it will suffice to show that there exists %$f \in L^1(\bbS^{d-1})$ and
$K \geq 1$ such that
\begin{align}
\label{eq:dct_hdmrd}
\esssup_{\theta \in \bbS^{d-1}} \esssup_{\delta \leq u \leq 1-\delta} \sup_{k \geq K}  |B_{2k}(u,\theta)| <\infty.
\end{align}
A similar argument may then be used to obtain the same conclusion with $B_{2k}$ replaced by $B_{3k}$.
These facts will imply, in particular, that the maps $(F+t_kH_{1k})^{-1}$ and $(G+t_kH_{2k})^{-1}$ 
are uniformly bounded over $[\delta,1-\delta] \times \bbS^{d-1}$, which, together with the continuity of $\varphi$, 
then also implies that the above display holds with $B_{2k}$ replaced by $B_{1k}$.

It thus remains to prove equation~\eqref{eq:dct_hdmrd}. We shall make use of the following properties.
\begin{lemma}
\label{lem:uni_abs_cont}
Under the assumptions of Theorem~\ref{thm::bootstrap}, the following assertions hold.
\begin{enumerate}
\item[(i)] The family $\{F(\cdot,\theta)\}_{\theta\in \bbS^{d-1}}$
is uniformly absolutely continuous over $I$, in the sense that for all $t > 0$, there exists $\epsilon(t) > 0$ such that
for any $[\alpha,\beta] \subseteq I$, we have
$$(\beta-\alpha) \leq \epsilon(t) ~ \Longrightarrow ~ \sup_{\theta \in \bbS^{d-1}} \big|F(\alpha,\theta) - F(\beta,\theta)\big| \leq t.$$
\item[(ii)] We have, 
$$\gamma := \inf_{\theta \in \bbS^{d-1}} \inf_{\delta/2 \leq u \leq 1-3\delta/4} \Big[ F^{-1}(u + \delta/4,\theta) - F^{-1}(u, \theta)\Big] > 0.$$
\item[(iii)] There exists a constant $K \geq 1$ such that for all $k \geq K$, $u \in [\delta,1-\delta]$,
and $\theta \in \bbS^{d-1}$,
$$F^{-1}(3\delta/4,\theta) \leq (F+t_kH_{1k})^{-1}(u,\theta) \leq F^{-1}(1-\delta/4,\theta).$$
\item[(iv)] In particular,  $(F+t_kH_{1k})^{-1}(u,\theta)- \gamma \in I_\theta$ for all $\theta \in \bbS^{d-1}$. 
\end{enumerate}
\end{lemma}
Now, for all $u \in [\delta,1-\delta]$ and $\theta \in \bbS^{d-1}$, write
\begin{equation}
\label{eq:xi_abbrv}
\xi_{\theta,u} = F^{-1}(u,\theta), \quad \xi_{\theta,uk} =  (F+t_kH_{1k})^{-1}(u,\theta).
\end{equation}
For all $k \geq 1$, let $\epsilon_k := \epsilon(t_k)$ be the constant corresponding
to the choice $t = t_k$ in the statement of Lemma~\ref{lem:uni_abs_cont}(i). This defines a sequence $(\epsilon_k)_{k \geq 1}$, 
which we may assume is nonincreasing, and satisfies $\epsilon_k < \gamma$ for all $k \geq K$, without loss of generality. 
By definition of the quantile function, we have
$$(F+t_kH_{1k})(\xi_{\theta,uk} - \epsilon_k,\theta) \leq u \leq (F+t_kH_{1k})(\xi_{\theta,uk},\theta).$$
We use these inequalities to bound $\xi_{\theta,uk} - \xi_{\theta,u}$. Assume first that
$\xi_{\theta,uk} < \xi_{\theta,u}$. Then, by absolute continuity of $F(\cdot,\theta)$ for all $\theta \in \bbS^{d-1}$, we have
\begin{align*}
0 \leq  (F+t_kH_{1k})(\xi_{\theta,uk},\theta) - F(\xi_{\theta,u},\theta) 
 &=  t_kH_{1k}(\xi_{\theta,uk},\theta) - \int_{\xi_{\theta,uk}}^{\xi_{\theta,u}} p_\theta(x)dx. 
\end{align*}
By Lemma~\ref{lem:uni_abs_cont}(iii), we have $[\xi_{\theta,uk}, \xi_{\theta,u}]\subseteq  I_\theta$ for all $k\geq K$, $u \in [\delta,1-\delta]$, and $\theta \in \bbS^{d-1}$.  Therefore, we obtain 
\begin{align*}
0 \leq     t_kH_{1k}(\xi_{\theta,uk},\theta) - (\xi_{\theta,u}-\xi_{\theta,uk}) J_{\infty,\delta/2}^{-1}(P_\theta). 
\end{align*}
On the other hand, if $\xi_{\theta,uk} \geq  \xi_{\theta,u}$, we have
\begin{align*}
0&\geq  (F+t_kH_{1k})(\xi_{\theta,uk}-\epsilon_k) - F(\xi_{\theta,u}) \\
 &=  \Big[F(\xi_{\theta,uk}-\epsilon_k) - F(\xi_{\theta,uk})\Big] + 
     \Big[F(\xi_{\theta,uk}) - F(\xi_{\theta,u})\Big] + 
t_kH_{1k}(\xi_{\theta,uk}-\epsilon_k)\\
 &= - \int_{\xi_{\theta,uk}-\epsilon_k}^{\xi_{\theta,uk}} p_\theta(x)dx  + 
     \int_{\xi_{\theta,u}}^{\xi_{\theta,uk}} p_\theta(x)dx + 
t_kH_{1k}(\xi_{\theta,uk}-\epsilon_k)\\
 &\geq -t_k + 
     (\xi_{\theta,uk} - \xi_{\theta,u})J_{\infty,\delta/2}^{-1}(P_\theta) + 
t_kH_{1k}(\xi_{\theta,uk}-\epsilon_k),
\end{align*}
where on the final line,  we lower bounded the first term as follows.
Since $\epsilon_k < \gamma$, we have $\xi_{\theta,uk}-\epsilon_k \in I_\theta \subseteq I$
by Lemma~\ref{lem:uni_abs_cont}(iv). Again, we also have $\xi_{\theta,uk} \in I$ by Lemma~\ref{lem:uni_abs_cont}(iii).
Therefore, by the definition of $\epsilon_k = \epsilon(t_k)$,  
Lemma~\ref{lem:uni_abs_cont}(i) can be applied to obtain the stated lower bound. 
Combine the preceding two displays  to deduce, 
\begin{align*}  
\left| B_{2k}(u,\theta) \right| =  \left|\frac{\xi_{\theta,uk} - \xi_{\theta,u}}{t_k}\right|
 &\leq \left(1+ \left| H_{1k}(\xi_{\theta,uk}) \right| + \left| H_{1k}(\xi_{\theta,uk}-\epsilon_k) \right|\right) J_{\infty,\delta/2}(P_\theta)\\
 &\leq \left(1 + \|H_{1k}\|_{\infty} +\|H_{2k}\|_{\infty}\right) \sup_{\theta \in \bbS^{d-1}} J_{\infty,\delta/2}(P_\theta).
 \end{align*} 
Now, recall that for $j=1,2$, $H_{jk}$ converges uniformly to $H_j \in \bbD_0\subseteq \ell^\infty(\calH)$.
It must also follow that, up to modifying the value of 
$K \geq 1$, the function $H_{jk}$ is bounded, uniformly in  $k \geq K$. This fact combined with our assumption on $P$
implies that the right-hand side of the above display is bounded above by a finite constant
not depending on $k, u,\theta$. Equation~\eqref{eq:dct_hdmrd} readily follows, leading to the claim. 
\qed

%We now turn to the proof of the main result of this subsection.

\paragraph{Proof of Theorem \ref{thm::bootstrap}.}
The claim consists of two statements, to be proven in parallel. Since the set of half-spaces $\calH$  
forms a separable Vapnik-Chervonenkis class, it is Donsker, implying that
the empirical process $\bbG_{nm} = \sqrt{\frac{nm}{n+m}} (P_n-P, Q_m - Q)$ 
%indexed by $\calH$ 
converges weakly in $\bbD^2 = \ell^\infty(\calH) \times \ell^\infty(\calH)$,
\begin{align}
\label{eq::process_donsker}
\sup_{h\in \text{BL}_1(\bbD^2)} \Big| \bbE[h(\bbG_{nm})] - 
	\bbE[h(\bbG_{(P,Q)}) ]\Big| \longrightarrow 0,
\end{align}	
to a process $\bbG_{(P,Q)} := (\sqrt a \bbG_P, \sqrt{1-a} \bbG_Q)$, where $\bbG_P$ and $\bbG_Q$ denote independent
$P$- and $Q$-Brownian bridges respectively, and where we identify the set $\calH$ with the set of indicator
functions over $\calH$. Under this abuse of notation, notice that the process $\bbG_P$ takes the form
$\bbG_P(x,\theta) = \bbG \circ F(x,\theta)$ 
for a standard Brownian Bridge $\bbG$, for all $(x,\theta) \in \calH$. By assumption, for  
all $\theta \in \bbS^{d-1}$, $F(\cdot,\theta)$ is continuous, and since almost all sample
paths of $\bbG$ are continuous, we deduce that almost every sample path of $\bbG_P(\cdot,\theta)$
is also continuous. We deduce that $\bbG_P$ takes values in $\bbD_0$ almost surely, and similarly for $\bbG_Q = \widebar\bbG \circ  G$, where $\widebar \bbG$ is a standard Brownian Bridge independent of $\bbG$.

Furthermore, Theorem 3.6.3 of \cite{vandervaart1996} implies the same conditional limiting distribution for 
the bootstrap empirical process $\bbG_{nm}^* = \sqrt{\frac{nm}{n+m}} (P_n^* - P_n,Q_m^*-Q_m)$,
\begin{equation}
\label{eq::boot_donsker}
\begin{aligned}
\sup_{h\in \text{BL}_1(\bbD^2)} \Big| \bbE\big[h(\bbG_{nm}^*)\big| X_1, \dots, X_n, Y_1, \dots, Y_m\big] - 
	\bbE\big[h(\bbG_{(P,Q)}) \big]\Big| \to 0, \\
\bbE\big[h(\bbG_{nm}^*) | X_1, \dots, X_n, Y_1, \dots, Y_m\big]^* {-} 
 \bbE\big[h(\bbG_{nm}^*) | X_1, \dots, X_n, Y_1, \dots, Y_m\big]_* \to 0,
\end{aligned}
\end{equation}
in outer probability, where $h$ ranges over $\text{BL}_1(\bbD^2)$.
Now,  the Hadamard differentiability of $\phi$ (Lemma \ref{lem::hadamard_sw}) together with
equation~\eqref{eq::process_donsker}
and the functional delta method (see for instance Theorems 3.9.4 of \cite{vandervaart1996}) 
implies 
\begin{align}
\label{eq::boot_step1}
\sup_{h\in \text{BL}_1(\bbR)} \left| \bbE\left[h\left(\sqrt {\frac{nm}{n+m}} \big(\phi(P_n, Q_m) - \phi(P, Q)\big)\right) \right] - 
	\bbE\Big[h\left(\phi'\left(\bbG_{(P,Q)}\right)\right) \Big]\right| \to 0.
	\end{align}
Likewise, the delta method for the bootstrap (Theorem 3.9.11 of \cite{vandervaart1996}) and equations \eqref{eq::process_donsker}
and \eqref{eq::boot_donsker} imply
\begin{align} 
\nonumber\sup_{h\in \text{BL}_1(\bbR)} \bigg| \bbE&\left[h\left(\sqrt {\frac{nm}{n+m}} \left(\phi(P_n^*, Q_m^*) - \phi(P_n, Q_m)\right)\right) \bigg| X_1, \dots, X_n, Y_1, \dots, Y_m\right] 
      \\ &- 
	\bbE\Big[h\left(\phi'\left(\bbG_{(P,Q)}\right)\right) \Big]\bigg| \longrightarrow 0,
\label{eq::boot_step2}
	\end{align}
in outer probability.	
A combination of equations \eqref{eq::boot_step1} and \eqref{eq::boot_step2} readily leads 
to part (ii) of the claim. In view of equation~\eqref{eq::boot_step1}, part (i) of the claim will
 follow upon showing that the random variable $\phi'(\bbG_{(P,Q)})$ is equal in distribution to
$N(0,a\sigma_P^2+(1-a)\sigma_Q^2)$. 
In the sequel, write for all $u\in[\delta,1-\delta]$,
\begin{align*}
%w(u,\theta) &=  \frac r {1-2\delta}  \text{sgn}\left(F^{-1}(u,\theta) - G^{-1}(u,\theta)\right)\big|F^{-1}(u,\theta) - G^{-1}(u,\theta)\big|^{r-1},\quad \theta \in \bbS^{d-1} \\
w_P(u) &= \int_{\bbS^{d-1}} \frac{w(u,\theta)}{p_\theta(F_\theta^{-1}(u))}d\mu(\theta),\quad 
w_Q(u) = \int_{\bbS^{d-1}} \frac{w(u,\theta)}{q_\theta(G_\theta^{-1}(u))}d\mu(\theta). %,\quad u \in [\delta,1-\delta]. 
\end{align*}
Notice that
\begin{align*}
\phi'&(\bbG_{(P,Q)})  \\
 &=   \int_{\bbS^{d-1}} \int_\delta^{1-\delta} w(u,\theta)
 	 \left(\frac{\sqrt{1-a} \bbG_Q(G_\theta^{-1}(u),\theta)}{q_\theta(G_\theta^{-1}(u))} - 
 	 \frac{\sqrt a \bbG_P(F_\theta^{-1}(u),\theta)}{p_\theta(F_\theta^{-1}(u))}\right)dud\mu(\theta)\\
  &=   \int_{\bbS^{d-1}} \int_\delta^{1-\delta} w(u,\theta)
 	 \left(\frac{\sqrt{1-a} \widebar\bbG(u)}{q_\theta(G_\theta^{-1}(u))} - \frac{\sqrt a  \bbG(u)}{p_\theta(F_\theta^{-1}(u))}\right)dud\mu(\theta)\\
  &= \sqrt{1-a} \int_\delta^{1-\delta} w_Q(u) \widebar\bbG(u)du - \sqrt a \int_\delta^{1-\delta} w_P(u)   \bbG(u)du,
\end{align*}
where, on the final line, the interchange of order of integration is valid $\bbP$-almost surely.
Indeed, the sample paths of  $\bbG$ and $\widebar \bbG$ are almost surely continuous, whence bounded over $[\delta,1-\delta]$,
which in turn implies that the functions $w(u,\theta)\bbG(u)/p_\theta(F_\theta^{-1}(u))$ and 
$w(u,\theta)\widebar\bbG(u)/q_\theta(G_\theta^{-1}(u))$ are almost surely bounded,
due to the assumptions placed on $P$ and $Q$. 
We now make use of the following fact, which is proven below for completeness.
\begin{lemma}
\label{lem:brownian_bridges}
Let $f:[\delta,1-\delta] \to \bbR$ be a Lebesgue-measurable and bounded function. Then, the random variable
$\int_\delta^{1-\delta} f(u) \bbG(u)du$ has Gaussian distribution with mean zero and variance
$$\Var\left[\int_\delta^{1-\delta} f(u) \bbG(u)du\right] {=} \int_0^{1-\delta}\left(\int_{\delta\vee t}^{1-\delta} f(u)du\right)^2 dt
 - \left(\int_0^{1-\delta} \int_{\delta\vee t}^{1-\delta} f(u)du dt\right)^2.$$
\end{lemma}
By Lemma~\ref{lem:brownian_bridges} and the independence of $\bbG$ and $\widebar \bbG$, we obtain that $\phi'(\bbG_{(P,Q)})$ has Gaussian distribution with mean zero and variance
\begin{align*}
\Var&\left[\phi'(\bbG_{(P,Q)})\right]  \\
 &= a \left[\int_0^{1-\delta} \left(\int_{\delta\vee t}^{1-\delta} w_P(u)du\right)^2 dt
 - \left(\int_0^{1-\delta} \int_{\delta\vee t}^{1-\delta} w_P(u)du dt\right)^2\right] \\
 &+ (1-a)\left[\int_0^{1-\delta} \left(\int_{\delta\vee t}^{1-\delta} w_Q(u)du\right)^2 dt
 - \left(\int_0^{1-\delta} \int_{\delta\vee t}^{1-\delta} w_Q(u)du dt\right)^2\right].
\end{align*}
Finally, notice that
\begin{align*}
\int_{\delta\vee t}^{1-\delta} w_P(u)du 
 &= \int_{\delta\vee t}^{1-\delta} \int_{\bbS^{d-1}} \frac{w(u,\theta)}{p_\theta(F_\theta^{-1}(u))} d\mu(\theta) du \\
 &=   \int_{\bbS^{d-1}}  \int_{\delta\vee t}^{1-\delta}\frac{w(u,\theta)}{p_\theta(F_\theta^{-1}(u))} du d\mu(\theta)  \\ 
 &=   \int_{\bbS^{d-1}}  \int_{F_\theta^{-1}(\delta\vee t)}^{F_\theta^{-1}(1-\delta)} w(F_\theta(x),\theta)dx d\mu(\theta),
\end{align*}
 where, once again, the interchange of the order of integration is  valid due to the uniform boundedness of the integrands
 almost everywhere.
 A similar computation holds with  $w_P$ replaced by $w_Q$,
 thus $\Var[\phi'(\bbG_{(P,Q)})] = a \sigma_P^2 + (1-a) \sigma_Q^2$, and the claim follows. \qed 
 
\paragraph{} It remains to prove Lemmas~\ref{lem:quantile_technical_uniform}--\ref{lem:brownian_bridges}. 
 
\subsection{Proof of Lemma~\ref{lem:quantile_technical_uniform}}
\label{app:pf_lem_quantile_technical_uniform} 
Let $u \in [\delta,1-\delta]$. We prove the claim for $F$, and an identical argument may then be used for $G$.
We use the abbreviations in equation~\eqref{eq:xi_abbrv}. Let $\epsilon > 0$ be an arbitrary real number satisfying
$$\epsilon < \inf_{\theta \in \bbS^{d-1}} \big[ F^{-1}(1-\delta/2,\theta) - F^{-1}(1-\delta,\theta)\big].$$
The infimum on the right-hand side is strictly positive by Lemma~\ref{lem:uni_abs_cont}(ii), whose proof below is
a consequence of the uniform integrability of $\{p_\theta\}_{\theta \in \bbS^{d-1}}$, and does not require
the present result. 
By definition of $\epsilon$, we have $\xi_{\theta,u}, \xi_{\theta,u} +\epsilon\in I_\theta$, thus, by absolute
continuity of $F(\cdot,\theta)$, 
\begin{align*}
\inf_{\theta \in \bbS^{d-1}} \Big[F(\xi_{\theta,u}+\epsilon,\theta) - u\Big] 
 &=\inf_{\theta \in \bbS^{d-1}} \int_{\xi_{\theta,u}}^{\xi_{\theta,u}+\epsilon} p_\theta(x)dx  \\
 &\geq \epsilon \inf_{\theta \in \bbS^{d-1}} \essinf_{\xi_{\theta,u} \leq x \leq \xi_{\theta,u} + \epsilon} p_\theta(x) \\
 &\geq \epsilon \inf_{\theta \in \bbS^{d-1}} \essinf_{x \in I_\theta} p_\theta(x) >0,
\end{align*}
where the strict inequality follows from the fact that $\sup_\theta J_{\infty,\delta/2}(P_\theta) < \infty$. 
After repeating a symmetric argument, we deduce
\begin{equation}
\label{eq:well_separated_quantile_under_lb_density}
\sup_{\theta \in \bbS^{d-1}} F(\xi_{\theta,u} - \epsilon,\theta) < u < 
\inf_{\theta \in \bbS^{d-1}} F(\xi_{\theta,u} + \epsilon,\theta).
\end{equation}
On the other hand, by definition of quantile, we have for all $\epsilon > 0$ that, 
$$ \sup_{\theta \in \bbS^{d-1}} (F+t_kH_{1k}) (\xi_{\theta,uk} - \epsilon, \theta)
 \leq u \leq \inf_{\theta \in \bbS^{d-1}} (F+t_kH_{1k}) (\xi_{\theta,uk}, \theta).$$
Recall that $H_{1k}$ converges in $\ell^\infty(\calH)$ to $H_1 \in \ell^\infty(\calH)$, thus there exists $C > 0$ such that
$$ \sup_{\theta \in \bbS^{d-1}} F(\xi_{\theta,uk} - \epsilon, \theta) - Ct_k
 \leq u \leq \inf_{\theta \in \bbS^{d-1}} F(\xi_{\theta,uk}, \theta) + Ct_k.$$
Since $t_k \downarrow 0$, the above display contradicts equation~\eqref{eq:well_separated_quantile_under_lb_density} for all large enough $k$
if $\xi_{\theta,uk} < \xi_{\theta,u} - \epsilon$,
 or if $\xi_{\theta,uk} > \xi_{\theta,u} + \epsilon$.
We have therefore shown that for all $\epsilon > 0$ small enough, there exists a sufficiently
large $K\geq 1$ such that for all $k \geq K$,  $\sup_\theta |\xi_{\theta,uk} - \xi_{\theta,u}| \leq \epsilon$,
thus leading to the claim.\qed 
 
\subsection{Proof of Lemma~\ref{lem:uni_abs_cont}}
\label{app:pf_lem_uni_abs_cont}
Consider the family  $\calF = \{p_\theta\}_{\theta \in \bbS^{d-1}}$. 
$\calF$ is assumed to be uniformly integrable with respect to the Lebesgue
measure over $\bbR$, and hence is also uniformly integrable
with respect to the finite measure $\nu = \lambda|_{I}$ (that is, the restriction
of the Lebesgue measure to the bounded interval $I$).
Since $\nu$ does not possess any atoms, uniform integrability of $\calF$ is equivalent to 
$\calF$ having uniformly absolutely continuous integrals, by Proposition 4.5.3 of~\cite{bogachev2007}. 
Thus, for any $t > 0$, there exists $\epsilon(t) > 0$ such that
for any interval $[\alpha,\beta] \subseteq I$ for which $\beta-\alpha \leq \epsilon(t)$, we have  
$$\sup_{\theta \in \bbS^{d-1}} \big|F(\alpha,\theta) - F(\beta,\theta)\big| \leq  t.$$
This proves property (i). For part (ii), 
choose $t < \delta/4$ and fix the corresponding value of $\epsilon(t)$. Let $u \in [\delta/2, 1-\delta/4]$, and choose a sequence
$\theta_j \in \bbS^{d-1}$ such that for all $j \geq 1$, 
\begin{equation}
\label{eq:approximate_inf}
 \big|F^{-1}(u+\delta/4,\theta_j) - F^{-1}(u,\theta_j)\big| \leq \frac 1 j + \inf_{\theta \in \bbS^{d-1}} \big|F^{-1}(u+\delta/4,\theta) - F^{-1}(u,\theta)\big|.
 \end{equation}
Let $\alpha_j = F^{-1}(u, \theta_j)$
and $\beta_j = F^{-1}(u+\delta/4,\theta_j)$. Clearly, $F(\beta_j,\theta_j) - F(\alpha_j,\theta_j) = \delta/4 > t$, 
thus from the uniform absolute continuity of $\{F_\theta\}$ in part (i),  it must hold that % for some $\kappa > 0$ which does not depend on $u$, 
$$\beta_j - \alpha_j = F^{-1}(u+\delta/4,\theta_j) - F^{-1}(u, \theta_j) > \epsilon(t).$$
Since this property holds for all $j \geq 1$, we deduce from equation~\eqref{eq:approximate_inf} that
$$\inf_{\theta\in\bbS^{d-1}}\big|F^{-1}(u+\delta/4,\theta) - F^{-1}(u, \theta)\big| \geq \epsilon(t)/2.$$
Since $\epsilon(t)$ did not depend on $u$, we finally arrive at 
$$\gamma = \inf_{\delta/2 \leq u \leq 1-\delta/4} \inf_{\theta\in\bbS^{d-1}}\big|F^{-1}(u+\delta/4,\theta) - F^{-1}(u, \theta)\big| \geq \epsilon(t)/2 > 0,$$
which proves the second claim. To prove the third and fourth, notice simply that
\begin{align*}
(F+t_kH_{1k})^{-1}&(\delta,\theta) \\
 &\hspace{-0.4in}\geq F^{-1}(\delta,\theta) - \sup_{\theta \in \bbS^{d-1}}\big|F^{-1}(\delta,\theta) - (F+t_kH_{1k})^{-1}(\delta,\theta)\big| \\
 &\hspace{-0.4in}\geq F^{-1}(3\delta/4,\theta) + \gamma - \sup_{\theta \in \bbS^{d-1}}\big|F^{-1}(\delta,\theta) - (F+t_kH_{1k})^{-1}(\delta,\theta)\big|.
\end{align*}
By Lemma~\ref{lem:quantile_technical_uniform}, recall that there exists $K \geq 1$ such that 
$\sup_{\theta \in \bbS^{d-1}}\big|F^{-1}(\delta,\theta) - (F+t_kH_{1k})^{-1}(\delta,\theta)\big| \leq \gamma$, 
thus for all such $k$, we have
\begin{align*}
(F+t_kH_{1k})^{-1}(\delta,\theta)
 &\geq F^{-1}(3\delta/4,\theta).
\end{align*}
Similarly, up to modifying the value of $K$, we have for all $k \geq K$, 
\begin{equation}
\label{eq:belongs_toIdelta}
(F+t_kH_{1k})^{-1}(1-\delta,\theta) \leq F^{-1}(1-\delta/4,\theta).
\end{equation} 
The claim follows from here.
\qed

\subsection{Proof of Lemma~\ref{lem:brownian_bridges}}
\label{app:pf_lem_brownian_bridges}
We first prove the claim for step functions $f$. Let $M \geq 1$ be an integer and 
define a partition of $[\delta,1-\delta]$ via $\delta=s_0 < \dots < s_{M+1} = 1-\delta$.
Let $\alpha_0,\dots, \alpha_{M} \in \bbR$ and set  
$f = \sum_{i=0}^{M} \alpha_i I_{[s_i, s_{i+1})}$. 
Clearly, we may always rewrite $f$ in terms of any 
refinement of the partition $s_0, \dots, s_{M+1}$. 
Indeed, for any $K \geq 1$ and any set of real numbers  $0 = t_0 < \dots < t_{K+1} = 1-\delta$ 
containing $\{s_0,\dots, s_{M+1}\}$, we may find 
real numbers $a_0,\dots, a_K$ contained in $\{\alpha_0, \dots, \alpha_M\}$ such that
$f = \sum_{k=0}^{K} a_k I_{[t_k,t_{k+1})}.$ We must have $a_0 = 0$ when $t_0 = 0 < \delta$. 
Since $\bbG$ is almost surely continuous over $[\delta,1-\delta]$, and the function $f$ is piecewise
continuous, $f\bbG$ is almost surely Riemann integrable. Therefore, 
for any choice   of the partition $\{t_0, \dots, t_{K+1}\}$
with vanishing mesh as $K \to \infty$, we have
\begin{align}
\label{eq:riemann_sums}
\int_\delta^{1-\delta}  f(u) \bbG(u) du  = \lim_{K\to\infty} \sum_{k=1}^{K+1}  \bbG(t_{k-1})(t_{k}-t_{k-1})a_{k-1}.
\end{align}
Notice that,
\begin{align*}
%\int & f(u) \bbG(u) du  \\
 %&= \sum_{k=0}^K (t_{k+1}-t_k) \bbG(t_k)a_k   \\
 \sum_{k=1}^{K+1} & \bbG(t_{k-1})(t_{k}-t_{k-1})a_{k-1}    \\
 &= \sum_{k=1}^{K+1} \bbG(t_{k-1})\left[\sum_{j=k-1}^K a_j(t_{j+1}-t_j) -  \sum_{j=k}^K a_j(t_{j+1}-t_j)\right]  \\
 &= \sum_{k=0}^{K} \bbG(t_{k}) \sum_{j=k}^K a_j(t_{j+1}-t_j) -   \sum_{k=1}^{K+1} \bbG(t_{k-1}) \sum_{j=k}^K a_j(t_{j+1}-t_j) \\
 &= \sum_{k=1}^{K} \bbG(t_{k}) \sum_{j=k}^K a_j(t_{j+1}-t_j) -   \sum_{k=1}^{K} \bbG(t_{k-1}) \sum_{j=k}^K a_j(t_{j+1}-t_j) &
 \\ %(\text{As } \bbG(t_0) = \bbG(0) = 0) \\
 &= \sum_{k=1}^{K} (\bbG(t_{k}) -\bbG(t_{k-1}))\sum_{j=k}^K a_j(t_{j+1}-t_j)  \\
 &= \sum_{k=1}^{K} (\bbG(t_{k}) -\bbG(t_{k-1}))\int_{t_k}^{1-\delta} f(x)dx.
% &= \int_0^{1-\delta}F_K(t)d \bbG(t),
\end{align*}
Now, since $f$ is bounded, the function $t\in [0,1-\delta] \mapsto   \int_{\delta \vee t}^{1-\delta} f(x)dx$
is continuous  and bounded, 
%is integrable with respect to the process $\bbG$. Therefore, 
thus for any partition as in equation~\eqref{eq:riemann_sums}, we obtain
\begin{align*}
\int_\delta^{1-\delta}& f(u)\bbG(u)du  \\
  &=  \lim_{K\to \infty} \sum_{k=1}^{K+1} (\bbG(t_{k}) - \bbG(t_{k-1})) \int_{t_k}^{1-\delta} f(x)dx
=\int_0^{1-\delta} \int_{\delta \vee t}^{1-\delta} f(x)dx d\bbG(t) 
,
\end{align*}
where  convergence in the final equality 
is understood to hold in probability (see, for instance, Proposition 2.13 of~\cite{revuz2013}).  
%where for all $t \in [0,1-\delta]$, $F$ is the step function defined by 
%$$F_K(t) = \sum_{k=1}^K I(t_{k-1} \leq t \leq t_k) \int_{t_k}^{1-\delta} f(x)dx.$$
It now follows from Proposition 2.2.1 of~\cite{denker1985} that
the random variable on the right-hand side of the above display
has mean-zero Gaussian distribution, with variance
\begin{align*}
\int_0^{1-\delta} \left(\int_{\delta\vee t}^{1-\delta} f(u)du\right)^2 dt
 - \left(\int_0^{1-\delta} \int_{\delta\vee t}^{1-\delta} f(u)du dt\right)^2,
\end{align*}
which leads to the claim when $f$ is a step function. 

Assume now that $f$ is a Lebesgue measurable bounded function. By Theorem 4.3 of~\cite{stein2009}, there exists a sequence of step functions
$f_n$ converging pointwise to $f$, Lebesgue-almost everywhere on $[\delta,1-\delta]$. In view of the preceding 
result, we have in probability, 
$$
\begin{multlined}[0.9\textwidth]
\Bigg|\int_\delta^{1-\delta} f(u) \bbG(u)du - \int_0^{1-\delta} \int_{t\vee \delta}^{1-\delta} f(u)dud\bbG(t)\Bigg| \\
  \leq \int_\delta^{1-\delta} |f_n-f||\bbG| + \int_0^{1-\delta} \left(\int_{\delta \vee t}^{1-\delta}|f_n - f|\right) d\bbG(t).
\end{multlined}$$
Since $f$ is bounded, we may clearly take the functions $f_n$ to be uniformly bounded. 
Furthermore, since $\bbG$ is $\bbP$-almost surely bounded on the comapct set $[\delta,1-\delta]$, the first term on the right-hand side of the above display vanishes by the Dominated
Convergence Theorem, while the second vanishes in probability by Theorem 2.12 of~\cite{revuz2013}. 
Deduce that the identity $\int f(u) \bbG(u)du = \int_0^{1-\delta} \int_{t\vee \delta}^{1-\delta} f(u)dud\bbG(t)$ holds, with convergence holding 
in probability. The claim then follows as before, by~\cite{denker1985}. \qed

\section{Proofs of Additional Results}
\label{app::secondary}
 
\subsection{Proof of Proposition \ref{prop::sw_coverage}}
Given $\theta \sim \mu$, we have 
\begin{align*}
W_r^r(P_\theta, Q_\theta)
 &= \frac 1 {1-2\delta}\int_\delta^{1-\delta} \big| F_\theta^{-1}(u) - G_\theta^{-1}(u)\big|^rdu \\
 &\lesssim \max_{a \in \{\delta,1-\delta\}} \big|F_\theta^{-1}(a)\big|^r + 
      \max_{a \in \{\delta,1-\delta\}} \big|G_\theta^{-1}(a)\big|^r.
\end{align*}
Therefore, it follows from Lemma \ref{lem::calK} that
$$\sup_{P,Q \in \calK_{2r}(b)} \Var_\mu\big[W_{r,\delta}^r(P_\theta,Q_\theta)\big] 
\leq \sup_{P,Q \in \calK_{2r}(b)} \int_{\bbS^{d-1}} W_r^{2r}(P_\theta, Q_\theta)d\mu(\theta) \leq b(4/\delta)^r.$$\
Thus, denoting $S_N =\sw_{r,\delta}^{(N)}(P,Q)$, $S = \sw_{r,\delta}(P,Q)$ and $\Delta_N = M_N/\sqrt N$, we obtain
\begin{align*}
\bbP&\Big(S  \not\in \widebar C_{nm}^{(N)}\Big)   \\
 &=\bbP\Big(S \not\in \widebar C_{nm}^{(N)}, |S^r - S_N^r| > \Delta_N\Big) + 
  \bbP\Big(S  \not\in \widebar C_{nm}^{(N)}, |S^r - S_N^r| \leq \Delta_N\Big) 
       \\
 &\leq  \bbP\Big( |S^r - S_N^r| > \Delta_N\Big)  
 \\ &+ \bbP\Big(\big(\{S^r \leq L_{N,nm} - \Delta_N\} \cup\{ U_{N,nm} + \Delta_N\leq S^r\}\big)\cap \{ |S^r - S_N^r| \leq \Delta_N \}\Big)
 \\
 &\leq  \bbP\Big(|S^r-S_N^r| > \Delta_N\Big)  + \bbP\Big(S_N^r \not\in C_{nm}^{(N)} \Big)\\
 &= \bbP\Big(|S^r-S_N^r| > \Delta_N\Big)+ \bbE\left[ \bbP\Big(S_N \not\in C_{nm}^{(N)} ~\Big|~ \theta_1,\dots,\theta_N \Big)\right]  \\
  &\leq 
 \frac{\Var_{\theta\sim\mu}\big[W_{r,\delta}^r(P_\theta,Q_\theta)\big]}{N\Delta_N^2}  +\alpha
\leq \frac {bc/\delta^r}{M_N^2} + \alpha,
\end{align*}
for a constant $c > 0$ depending only on $r$, as claimed.
\qed  

\subsection{Proof of Corollary \ref{cor::length_examples}}
\label{app::proof_cor1}
Under the stated conditions on $\epsilon,\delta$, it can be seen
by direct verification that conditions~\ref{assm::B1}--\ref{assm::B4}
and the remaining conditions of Theorem~\ref{thm::length} hold
for the confidence bands of Examples~\ref{ex::DKW}--\ref{ex::rel_VC}, for constants $K_1, K_2$
possibly depending on $\alpha$. 
In what follows, the symbol ``$\lesssim$'' is used to hide constants possibly depending
on~$\alpha,b,\delta_0$ and~$r$. 

 Suppose first that $\sj_{r,\frac \delta 2}(P)\vee \sj_{r,\frac \delta 2}(Q) =\infty$. 
 Then, we use the bound %since $\kappa_{\varepsilon,n\wedge m} \leq \delta/2$ by assumption on $\varepsilon$, we have
\begin{align*}
\psi_{\varepsilon,nm}
 &\lesssim \big( \sw_{\infty,\delta}^{(r{-}1)}(P,Q) {+} U_{\varepsilon,n}(P) {+} U_{\varepsilon,m}(Q)\big)\frac{\kappa_{\varepsilon,n}}{ \sqrt\delta}\\
 &\lesssim  \left(\max_{a \in \{\delta/2, 1-\delta/2\}} \sup_{\bbS^{d-1}} \big|F_\theta^{-1}(a)\big|^{r-1} \right) \frac{\kappa_{\varepsilon,n}}{\sqrt \delta}
 \lesssim  \frac {\kappa_{\varepsilon,n}} {\delta^{r/2}} ,
 \end{align*}
 by Lemma~\ref{lem::calK}, under the assumption $P,Q \in \widebar\calK_{2}(b)$.
A similar argument can be used to bound $\varphi_{\varepsilon,nm}$,
leading to 
\begin{align*}
\lambda(C_{nm}^{(N)})
 &\leq \Big\{\sw_{r,\delta}^r(P,Q) + c\big(\psi_{\varepsilon, nm} + \varphi_{\varepsilon,nm} + \varkappa_N\big) \Big\}^{1/r}
 -  \sw_{r,\delta}(P,Q) \\
 &\leq c^{1/r} \big(\psi_{\varepsilon, nm} + \varphi_{\varepsilon,nm} + \varkappa_N\big)^{1/r} \\
 &\lesssim \varkappa_N^{1/r} +  \frac 1 {\sqrt\delta}\Big(\kappa_{\varepsilon,n}^{1/r} + \kappa_{\varepsilon,m}^{1/r}\Big),
 \end{align*}
 with probability at least $1-\epsilon$.
Parts (i) and (ii) now  follow from Lemma~\ref{lem::kappa_examples}  in the case $\sj_{r,\frac \delta 2}(P)\vee \sj_{r,\frac \delta 2}(Q) =\infty$.
 
Suppose now that  $\sj_{r,\frac \delta 2}(P)\vee \sj_{r,\frac \delta 2}(Q) < \infty$.
Using the shorthand notations
$$\Delta = V_{\varepsilon,n}(P) + V_{\varepsilon,m}(Q), \quad S = \sw_{r,\delta}(P,Q), $$
we have the bounds $\psi_{\varepsilon,nm},\varphi_{\varepsilon,nm} \lesssim (S^r + \Delta)^{\frac{r-1}{r}} \Delta^{\frac 1 r},$
whence
\begin{align*}
\lambda(C_{nm}^{(N)}) 
 &\lesssim  \varkappa_N^{\frac 1 r} + \Big\{S^r +  (S^r + \Delta)^{\frac{r-1}{r}} \Delta^{\frac 1 r} \Big\}^{\frac 1 r} - S,
\end{align*}
with probability at least $1-\epsilon$.
If $S^r \lesssim \Delta$, the right-hand side of the above display is clearly of order $\varkappa_N^{1/r}+\Delta^{1/r}$. Likewise, if
$S^r\gtrsim \Delta$, we obtain similarly as in equation~\eqref{eq::rpow_to_sw},
\begin{align*}
\lambda(C_{nm}^{(N)}) 
 &\lesssim  \varkappa_N^{\frac 1 r} + S\bigg[\Big\{1 +  (1 + (\Delta/S^r))^{\frac{r-1}{r}} (\Delta/S^r)^{\frac 1 r} \Big\}^{\frac 1 r} - 1\bigg] \\
 &\leq  \varkappa_N^{\frac 1 r} + S(1 + (\Delta/S^r))^{\frac{r-1}{r}} (\Delta/S^r)^{\frac 1 r}  \\
 &\lesssim  \varkappa_N^{\frac 1 r} + \Delta^{\frac 1 r},
\end{align*}
with probability at least $1-\epsilon$.
The conclusion of the above display thus holds irrespective of $S$ and $\Delta$. The claim now follows
by substituting the expressions for $V_{\varepsilon,n}(P)$ stated in Lemma~\ref{lem::kappa_examples} for each of parts (i) and (ii). 
\qed

\subsection{Proof of Corollary \ref{cor::null}}
We reason similarly as in the proof of Corollary~\ref{cor::length_examples}.
By Theorem \ref{thm::length}, we have with probability at least $1-\epsilon$,
\begin{align*}
\lambda(C_{nm}) 
 &\leq \Big\{\sw_{r,\delta}^r(P,Q) + c\big(\psi_{\varepsilon, nm} + \varphi_{\varepsilon,nm} + \varkappa_N\big) \Big\}^{1/r}
 -  \sw_{r,\delta}(P,Q) \\
 &\lesssim   \sw_{r,\delta}^{1-r}(P,Q)\big(\psi_{\varepsilon, nm} + \varphi_{\varepsilon,nm} + \varkappa_N\big)  \\
 &\leq \Gamma^{1-r} \big(\psi_{\varepsilon, nm} + \varphi_{\varepsilon,nm} + \varkappa_N\big),
\end{align*}
since $\sw_{r,\delta}(P,Q) \geq \Gamma$. The claim now follows by invoking
similar bounds on $\psi_{\varepsilon,nm}$ and $\varphi_{\varepsilon,nm}$
as in the proof of Corollary~\ref{cor::length_examples}. \qed

\subsection{Proof of Proposition \ref{prop::pretest}}
\label{app::pf_retest} 
Notice that  $\sigma_P, \sigma_Q > 0$ whenever $\sw_{r,\delta}(P,Q) > 0$.
In view of Theorem~\ref{thm::bootstrap}, % and the fact that $\sigma_P, \sigma_Q > 0$ when $\sw_{r,\delta}(P,Q) > 0$, 
it is then a standard result that the percentile bootstrap interval $C_{nm}^*$ satisfies
\begin{equation}
\label{eq:boot_coverage}
\liminf_{n,m\to \infty} \bbP(\sw_{r,\delta}(P,Q) \in C_{nm}^*) \geq 1-\alpha/2,
\end{equation}
under the assumptions of Theorem~\ref{thm::bootstrap} and the assumption $\sw_{r,\delta}(P,Q) > 0$
(see, for instance, Lemma~23.3 of~\cite{vandervaart1998}). Therefore, when these assumptions hold, we have
\begin{align*}
\bbP&(\sw_{r,\delta}(P,Q) \not\in C_{nm}) \\
 &= \bbP (\sw_{r,\delta}(P,Q) \not\in C_{nm}^*, 0 \not\in C_{nm}^\dagger) +  \bbP(\sw_{r,\delta}(P,Q)\not\in C_{nm}^\dagger,0\in C_{nm}^\dagger) 
\\
 &\leq \bbP(\sw_{r,\delta}(P,Q) \not\in C_{nm}^*) +  \bbP(\sw_{r,\delta}(P,Q)\not\in C_{nm}^\dagger) 
 \leq \alpha + o(1),
\end{align*} 
where on the final line, we use equation~\eqref{eq:boot_coverage} and Proposition~\ref{prop::sw_coverage}.
Notice that the assumptions of Proposition~\ref{prop::sw_coverage}   are satisfied, since $\widebar \calK_2 \subseteq \calK_{2r}$. 
On the other hand, when  $\sw_{r,\delta}(P,Q) = 0$, we obtain
\begin{align*}
\bbP&(\sw_{r,\delta}(P,Q) \not\in C_{nm}) \\
 &= \bbP (\sw_{r,\delta}(P,Q) \not\in C_{nm}^*, 0 \not\in C_{nm}^\dagger) +  \bbP(\sw_{r,\delta}(P,Q)\not\in C_{nm}^\dagger,0\in C_{nm}^\dagger) 
\\
 &\leq 2\bbP(0 \not\in C_{nm}^\dagger) \leq \alpha+o(1).
\end{align*} 
This proves the stated asymptotic coverage property of the confidence interval $C_{nm}$.
In order to bound its length, note that it is a direct consequence of Theorem~\ref{thm::bootstrap} 
that the bootstrap quantiles $F^*_{nm}(\alpha/2)$ and $F^*_{nm}(1-\alpha/2)$ 
are of the order $O_p(n^{-1/2})$ as $n/(n+m) \to a\in (0,1)$.
Thus, 
$$b_{nm}^* - a_{nm}^* = O_p(n^{-1/2}),$$
where we write $C_{nm}^* = [(a_{nm}^*)^{1/r}, (b_{nm}^*)^{1/r}]$. 
Furthermore, under the conditions of Theorem~\ref{thm::bootstrap}, and for $N \asymp n^{r^2}$
and $M_N \asymp \log N$,  it can be deduced from Theorem~\ref{thm::length} 
and Corollary~\ref{cor::length_examples} that when $\sw_{r,\delta}(P,Q) = 0$, we have
$$b_{nm}^\dagger - a_{nm}^\dagger = O_p\left(\left(\frac{\log n}{n}\right)^{\frac r 2}\right)$$
where we write $C_{nm}^\dagger = [(a_{nm}^\dagger)^{1/r}, (b_{nm}^\dagger)^{1/r}]$.  
Finally, as in the above proof of coverage of $C_{nm}$, when
$\sw_{r,\delta}(P,Q) > 0$ we have $C_{nm} = C_{nm}^*$ with probability at least $1-\alpha/2 - o(1)$, 
while when $\sw_{r,\delta}(P,Q) = 0$ we have $C_{nm} = C_{nm}^\dagger$ with probability at least $1-\alpha/2-o(1)$. 
Combine these facts to deduce that for any $\epsilon > 0$, there exist constants $C,n_0 > 0$ such that
for all $n \geq n_0$,
$$b_{nm} - a_{nm} \leq C \left(  \left(\frac{\log n}{n}\right)^{\frac r 2} + \frac{\sw_{r,\delta}(P,Q)}{\sqrt n}\right)$$
with probability at least $1-\alpha/2 - \epsilon$. Choosing $\epsilon = \alpha/2$ leads to the claim.
 \qed  
 
\subsection{Proof of Proposition \ref{prop::likelihood_free}}
The proof is straightforward. We have,
\begin{align*}
\bbP\Big(\eta_0 \not\in \widebar C_{nm}^{(N)}\Big)
 &\leq \bbP\Big(\bar \ell_{nm}^{(N)}(\eta_0) > \epsilon \Big) \\
 &\leq \bbP\Big(\bar \ell_{nm}^{(N)}(\eta_0) > \epsilon_0 \Big)
 = \bbP\Big(\bar \ell_{nm}^{(N)}(\eta_0) > \sw_{r,\delta}(P,P_{\eta_0}) \Big).
 \end{align*}
The claim now follows from Proposition \ref{prop::sw_coverage}.
\qed 
 
\subsection{Proof of Example \ref{ex::rel_VC}}
\label{subsec::relVC_proof}
We begin by proving the validity of the inequality in equation~\eqref{eq::relVC}. 
Let $\calA$ be a collection of sets, and let $\calS_\calA(n)$ denote the shattering
number \citep{vapnik2013} of $\calA$. The relative VC inequality is then given by
$$\bbP\left(\sup_{A \in \calA} \frac{|P_n(A) - P(A)|}{\sqrt{P_n(A)}} \geq t\right) \leq 4\calS_\calA(2n) e^{-nt^2/4},\quad t > 0.$$
Letting $\calA = \{(-\infty,x] : x \in \bbR\}$ and 
$\calA = \{[x, \infty) : x \in \bbR\}$ respectively, we obtain
\begin{align*}
\bbP&\left(\sup_{x \in \bbR} \frac{|F_n(x) - F(x)|}{\sqrt{F_n(x)}} \geq t\right) \leq 4(2n+1) e^{-nt^2/4},\\
\bbP&\left(\sup_{x \in \bbR} \frac{|F_n(x) - F(x)|}{\sqrt{1-F_n(x)}} \geq t\right) \leq 4(2n+1) e^{-nt^2/4},
\end{align*}
for all $t > 0$. By a union bound
and the fact that $u(1-u)\geq \frac 1 2 (u\wedge (1-u))$ for all $u \in [0,1]$, we arrive at
$$\bbP\left(\sup_{x \in \bbR} \frac{|F_n(x) - F(x)|}{\sqrt{F_n(x)(1-F_n(x))}} \geq t \right) \leq 8(2n+1) e^{-\frac{nt^2}{16}}.$$ 
Setting $t=\nu_{\alpha,n}:=\sqrt{\frac {16}{n} \left[ \log(16/\alpha) + \log(2n+1)\right]}$, we deduce that
with probability at least $1-\alpha/2$,
\begin{equation}
\label{eq::relVC_appendix}
|F_n(x) - F(x)| \leq  \nu_{\alpha,n} \sqrt{F_n(x)(1-F_n(x))},\quad \forall x \in \bbR.
\end{equation}
This proves the validity of equation~\eqref{eq::relVC}. 

We now invert equation~\eqref{eq::relVC_appendix} to obtain the functions $\gamma_{\alpha,n}$ and $\eta_{\alpha,n}$ which lead to a
quantile confidence band. We will require the following definitions of lower CDF and upper quantile function, 
$$\widebar F(x) := \lim_{y \to x^-} F(y) = \bbP(X_1 < x), \quad \widebar F^{-1}(u) = \inf\big\{x \in \bbR: \widebar F(x) > u\big\},$$
with empirical analogues given by
$$\widebar F_n(x) := \lim_{y \to x^-} F_n(y) = \frac 1 n \sum_{i=1}^n I(X_i < x), \quad \widebar F_n^{-1}(u) = \inf\big\{x \in \bbR: \widebar F_n(x) > u\big\}.$$
Notice that $F$ and $\widebar F^{-1}$ are right continuous, whereas $\widebar F$ and $F^{-1}$ are left continuous.
Furthermore, we make use of the following elementary inequalities relating quantile functions and CDFs,
\begin{align}
\label{eq::qu_cdf_1} F_n(x) \geq u \Longrightarrow& x \geq F_n^{-1}(u), \\
\label{eq::qu_cdf_3} \widebar F_n(x) \leq u \Longrightarrow& x \leq \widebar F_n^{-1}(u),\\
\label{eq::qu_cdf_2} F(x) \geq u   \Longleftrightarrow& x \geq F^{-1}(u), \\
\label{eq::qu_cdf_4} \widebar F(x) \leq u  \Longleftrightarrow& x \leq \widebar F^{-1}(u).
\end{align}
We now turn to the proof. The calculations which follow are elementary, but tedious.
Let $v = F(x)$. By equation~\eqref{eq::relVC_appendix}, we have with probability at least $1-\alpha/2$
that for all $x \in \bbR$,
\begin{align*}
F_n(x)& + \nu_{\alpha,n} \sqrt{F_n(x)(1-F_n(x))}  \geq v 
\geq F_n(x) - \nu_{\alpha,n}\sqrt{F_n(x)(1-F_n(x))}\\
 &\Longrightarrow F_n(x)(1-F_n(x))\nu_{\alpha,n}^2 \geq (v-F_n(x))^2 \\
 &\Longrightarrow (F_n(x) - F_n(x)^2) \nu_{\alpha,n}^2 \geq v^2 - 2v F_n(x) + F_n(x)^2 \\
 &\Longrightarrow F_n(x)^2 (1+\nu_{\alpha,n}^2) - F_n(x)\big(2v + \nu_{\alpha,n}^2\big) + v^2 \leq 0\\
 &\Longrightarrow F_n(x) \geq  \frac{2v + \nu_{\alpha,n}^2}{2(1+\nu_{\alpha,n}^2)} - 
 	\frac{\sqrt{[2v+\nu_{\alpha,n}^2\big]^2 - 4(1+\nu_{\alpha,n}^2)v^2}}{2(1+\nu_{\alpha,n}^2)}\\
% &\Longrightarrow F_n(x) \in \left\{ \frac{ 2v + \nu_{\alpha,n}^2}{2(1+\nu_{\alpha,n}^2)} + \iota 
% 	\frac{\nu_{\alpha,n} \sqrt{\nu_{\alpha,n}^2 + 4v(1-v)}}{2(1+\nu_{\alpha,n}^2)}: \iota \in \{-1,1\}\right\} \\
 &\Longrightarrow F_n(x) \geq \frac{2v + \nu_{\alpha,n}^2}{2(1+\nu_{\alpha,n}^2)} -  \frac{\nu_{\alpha,n} \sqrt{\nu_{\alpha,n}^2 + 4v(1-v)}}{2(1+\nu_{\alpha,n}^2)} = \gamma_{\alpha,n}(v)\\
 &\Longrightarrow x \geq F_n^{-1}(\gamma_{\alpha,n}(v)) ,\qquad  (\text{By } \eqref{eq::qu_cdf_1}) \\
 &\Longrightarrow x \geq F_n^{-1}(\gamma_{\alpha,n}(F(x))).
\end{align*}
Now, let $u \in (0,1)$. Setting $x = F^{-1}(u)$ and using the fact that $F\circ F^{-1}(u) \geq u$ by 
equation~\eqref{eq::qu_cdf_2}, the above display implies
$$ F^{-1}(u) \geq F_n^{-1}(\gamma_{\alpha,n}(F\circ F^{-1}(u))) \geq F_n^{-1}(\gamma_{\alpha,n}(u)),$$
uniformly in $u \in (0,1)$, with probability at least $1-\alpha/2$. 

We now turn to an upper confidence bound on $F^{-1}(u)$. Upon taking limits from the left in equation~\eqref{eq::relVC_appendix}, 
we obtain 
$$\widebar F_n(x) - \nu_{\alpha,n} \sqrt{\widebar F_n(x)(1-\widebar F_n(x))} \leq \widebar F(x)
 \leq \widebar F_n(x) + \nu_{\alpha,n} \sqrt{\widebar F_n(x)(1-\widebar F_n(x))} $$
uniformly in $x \in \bbR$, on the same
event of probability at least $1-\alpha/2$. Thus, letting $v = \widebar F(x)$, we have
\begin{align*}
%\widebar F_n&(x) - \nu_{\alpha,n} \sqrt{\widebar F_n(x)(1-\widebar F_n(x))}  \leq v \\
% &\Longrightarrow 
 \nu_{\alpha,n}^2&\widebar F_n(x)(1-\widebar F_n(x)) \geq (\widebar F_n(x) - v)^2 \\
 &\Longrightarrow \nu_{\alpha,n}^2 \widebar F_n(x) -\nu_{\alpha,n}^2 \widebar F_n(x)^2 \geq \widebar F_n(x)^2 - 2v \widebar F_n(x) + v^2 \\
 &\Longrightarrow \widebar F_n(x)^2(1+\nu_{\alpha,n}^2) - (\nu_{\alpha,n}^2 + 2v) \widebar F_n(x) + v^2 \leq 0 \\
 &\Longrightarrow \widebar F_n(x) \leq 
 	\frac{\nu_{\alpha,n}^2 + 2v}{2(1+\nu_{\alpha,n}^2)} + \frac{\sqrt{[2v+\nu_{\alpha,n}^2\big]^2 - 4(1+\nu_{\alpha,n}^2)v^2}}{2(1+\nu_{\alpha,n}^2)}\\
 &\Longrightarrow \widebar F_n(x) \leq \frac{\nu_{\alpha,n}^2 + 2v + \nu_{\alpha,n} \sqrt{\nu_{\alpha,n}^2 + 4v(1-v)}}{2(1+\nu_{\alpha,n}^2)} = \eta_{\alpha,n}(v)\\
 &\Longrightarrow x \leq \widebar F_n^{-1}(\eta_{\alpha,n}(v)) ,\qquad  (\text{By } \eqref{eq::qu_cdf_3}) \\
 &\Longrightarrow x \leq \widebar F_n^{-1}(\eta_{\alpha,n}(\widebar F(x))).
\end{align*}
Therefore, setting $x = \widebar F^{-1}(u)$ for $u \in (0,1)$, and using the fact that
$\widebar F \circ \widebar F^{-1}(u) \leq u$ by equation~\eqref{eq::qu_cdf_4}, we obtain 
$$\widebar F^{-1}(u) \leq \widebar F_n^{-1}\big(\eta_{\alpha,n}(\widebar F( \widebar F^{-1}(u))\big) 
\leq \widebar F_n^{-1}\big(\eta_{\alpha,n}(u)\big).$$
Upon taking limits from the right, this implies
$$F^{-1}(u) \leq F_n^{-1}\big(\eta_{\alpha,n}(u)\big),$$
uniformly in $u$ with probability at least $1-\alpha/2$. We conclude that
$$\bbP\Big(F_n^{-1}\big(\gamma_{\alpha,n}(u)\big) \leq F^{-1}(u) \leq F_n^{-1}\big(\eta_{\alpha,n}(u)\big), \ \forall u \in (0,1) \Big) \geq 1-\alpha/2.$$
The validity of equation~\eqref{eq::relVC_fns} follows. 
\qed 
\section{Further Interpretation of Theorem \ref{thm::length}}
\label{app::mixed_case}

We now briefly discuss Theorem \ref{thm::length} 
in the situation where $\sj_{r,\delta/2}(P) < \infty$ but $\sj_{r,\delta/2}(Q) =\infty$.  
In this case, Corollary \ref{cor::null} continues to imply that the parametric rate is achievable for distributions $P$ and
$Q$ which are bounded away from each other. When $P$ and $Q$ are permitted to approach each other, 
the conditions $\sj_{r,\delta/2}(P) <\infty,\sj_{r,\delta/ 2}(Q) =\infty$ need not be violated so long as $P \neq Q$. 
In this specialized situation, Corollary \ref{cor::length_examples} merely implies the nonparametric rate $n^{-\frac 1 {2r}} + m^{-\frac 1 {2r}}$,
though one may expect a better dependence on $n$ as a result of the condition $\sj_{r,\delta/2}(P) < \infty$, 
as well as a better dependence on $m$ by virtue of $Q$ being in proximity to $P$.
 
We focus on the one-dimensional case for simplicity, and we omit superscripts depending on $N=1$ below. In this setting, 
upper bounds on the $\infty$-Wasserstein distance in terms of $r$-Wasserstein distances have been
established by \cite{bouchitte2007}:   
If $P$ 
is absolutely continuous with respect to the Lebesgue measure, with strictly
positive density $p$ over a compact set $\calX\subseteq \bbR$, then for all $r \geq 1$, 
\begin{equation}
\label{eq::Linf_ineq}
W_\infty^{r+1}(P,Q) \leq C_{r}(\calX)  W_r^r(P,Q)\sup_{x \in \calX} \left|\frac 1 {p(x)}\right|,
\end{equation}
where $C_{r}(\calX) > 0$ is a constant depending only on $r,\calX$.  
Equipped with this result, we arrive at the following Corollary to Theorem~\ref{thm::length}.

\begin{corollary}
\label{cor::mixed}
Let $P,Q \in \calP(\bbR)$ be absolutely continuous with respect to the Lebesgue measure, and assume the same conditions 
and notation as Theorem \ref{thm::length}. 
Assume there exists a fixed constant $s > 0$ such that $\text{J}_{\infty,\delta/2}(P) \leq s$, and let $\delta$ be fixed.
Assume further that $W_{r,\delta}(P,Q) \lesssim \kappa_{\varepsilon,n\wedge m}$.
Then, there exists a constant $C > 0$ depending on $b,r,\delta ,s$ such that
with probability at least $1-\epsilon$, we have
$$\lambda(C_{nm}^{(N)}) \leq C \big[\kappa_{\epsilon,n\wedge m}\big]^{\frac{r^2+1}{r(r+1)}}.$$
\end{corollary}
For example, when $\gamma_{\epsilon,n},\eta_{\epsilon,n}$ are based on the DKW inequality
as in Example~\ref{ex::DKW}, 
and $r=2$, Corollary \ref{cor::mixed} implies the  rate $n^{-\frac 5 {12} } + m^{-\frac 5 {12}}$,
which is strictly improves the $n^{-\frac 1 {4}} + m^{-\frac 1 {4}}$ rate implied by Corollary \ref{cor::length_examples}.
Moreover, the rate implied by Corollary \ref{cor::mixed} approaches the parametric rate for large $r$, 
which we conjecture to be the minimax rate for this scenario.

\paragraph{Proof of Corollary \ref{cor::mixed}.}
Let $P^\delta,Q^\delta$ denote the $\delta$-trimmings of $P$ and $Q$, in the sense of equation~\eqref{eq::trimmed_W}.
The quantile functions of $P^\delta, Q^\delta$ are given by
$$F^{-1}_\delta(u) = F^{-1}(u(1-2\delta)+\delta), \quad G^{-1}_\delta(u) = G^{-1}(u(1-2\delta) + \delta).$$
Now, write $h=\kappa_{\epsilon,n\wedge m}$ for simplicity, and let $P_h^\delta, Q_h^\delta$ denote distributions
with quantile functions $u\mapsto F_\delta^{-1}(u+h)$ and $u \mapsto G_\delta^{-1}(u+h)$ respectively, which
are well defined since $h\leq \delta/2$ by assumption.
Then, by definition of $U_{\epsilon,m}$, 
$$U_{\epsilon,m}(Q) \leq W_\infty^{r-1}(Q^\delta, Q_h^\delta).$$
Furthermore, the finiteness of $\text{J}_{\infty,\delta/2}(P)$ implies
$W_\infty(P^\delta,P_h^\delta) \lesssim h$, and we have,
\begin{align*}
W_r(P_h^\delta, Q_h^\delta) 
 &= W_r (P^\delta,Q^\delta) + 
     \int_{\delta+h}^{1-\delta+h} \big|F^{-1} - G^{-1}\big|^r -  \int_{\delta}^{1-\delta} \big|F^{-1}  - G^{-1} \big|^r   \\
 &\lesssim W_r(P^\delta,Q^\delta) + h.
\end{align*}
Combining these facts, together with equation \eqref{eq::Linf_ineq} and the bound $\text{J}_{\infty,\delta/2}(P) \leq s$, we arrive at
\begin{align*}
U_{\epsilon,m}^{\frac 1 {r-1}}(Q) 
 &\leq  W_\infty(Q^\delta, Q_h^\delta)\\
 &\leq W_\infty(Q^\delta, P^\delta) + W_\infty(P^\delta, P_h^\delta) + W_\infty(P_h^\delta, Q_h^\delta) \\
 &\lesssim W_r^{\frac r {r+1}}(P^\delta,Q^\delta) +h + W_r^{\frac r {r+1}}(P_h^\delta, Q_h^\delta)\\
% &\lesssim 2\sw_r^{\frac r {r+1}}(P^\delta,Q^\delta) + h + h^{\frac r {r+1}} \\
% &\lesssim W_r^{\frac r {r+1}}(P^\delta,Q^\delta) +  h^{\frac r {r+1}} \\
&\lesssim W_{r,\delta}^{\frac r {r+1}}(P,Q) +  h^{\frac r {r+1}}. 
\end{align*}
Thus, by Theorem \ref{thm::length},   under the assumption
$W_{r,\delta}(P,Q) \lesssim h $,
we have with probability at least $1-\epsilon$,
\begin{align*}
\lambda&(C_{nm}) \\ 
 &\lesssim \Big\{W_{r,\delta}^r(P,Q) + h\Big(W_{\infty,\delta}^{r-1}(P,Q) + U_{\epsilon,n}(P) + U_{\epsilon,m}(Q) \Big) \Big\}^{1/r}
 -  W_{r,\delta}(P,Q) \\
 &\lesssim W_{r,\delta}(P,Q)\left\{ \left[1 {+} \frac{h}{W_{r,\delta}^r(P,Q)}\Big(W_{r,\delta}^{\frac{r(r-1)}{r+1}}(P,Q) + h^{r-1}% + W_{r,\delta}^{\frac {r(r-1)} {r+1}}(P,Q)
  +  h^{\frac {r(r-1)} {r+1}} \Big)
 \right]^{\frac 1 r} {-}  1 \right\}\\ 
%  &\lesssim W_{r,\delta}(P,Q) \left\{\left[1 +  h \left(W_{r,\delta}^{\frac{(r-1)}{r+1}}(P,Q) %+ W_{r,\delta}^{\frac {(r-1)} {r+1}}(P,Q) 
% 	+  \frac{h^{\frac {r(r-1)} {r+1}}}{W_{r,\delta}^r(P,Q)} \right)\right]^{1/r}
% -  1\right\} \\
 &\lesssim W_{r,\delta}(P,Q) \left( \frac{h^{\frac {r^2+1} {r+1}}}{W_{r,\delta}^r(P,Q)} \right)^{1/r}
\asymp  h^{\frac {r^2+1} {r(r+1)}}.
\end{align*}
The claim follows.\qed

%% if your bibliography is in bibtex format, uncomment commands:
\bibliographystyle{imsart-nameyear} % Style BST file (imsart-number.bst or imsart-nameyear.bst)
\bibliography{manuscript_EJS2001}       % Bibliography file (usually '*.bib')

\end{document}